%% file: 1study.tex
\documentclass[a4paper,12pt,reqno]{amsbook}

\usepackage{graphicx}
\usepackage{amsmath}
\usepackage{amssymb}
\usepackage{amsfonts}
\usepackage{epsfig}
\usepackage{epigraph}

\makeatletter

\@addtoreset{equation}{chapter}

\newtheorem{theorem}{Theorem}[chapter]

\newtheorem{corollary}[theorem]{Corollary}
\newtheorem{remark}[theorem]{Remark}
\newtheorem{proposition}[theorem]{Proposition}
\newtheorem{definition}[theorem]{Definition}

\newtheorem{example}[theorem]{Example}
\def\PerfProof{{\it Proof.\ }}

\def\PerfProof{{\it Proof.\ }}

\makeatother

\makeindex

\begin{document}

\title{Notes on Coxeter Transformations and the McKay correspondence}
       \author{Rafael Stekolshchik}
         \thanks{{\large email: rs2@biu.013.co.il}}

\date{}

\begin{abstract}
We consider the Coxeter transformation in the context of the McKay
correspondence, representations of quivers, and Poincar\'{e}
series.

We study in detail the Jordan forms of the Coxeter transformations
and prove shearing formulas due to Subbotin and Sumin for the
characteristic polynomials of the Coxeter transformations. Using
shearing formulas we calculate characteristic polynomials of the
Coxeter transformation for the diagrams $T_{2,3,r}, T_{3,3,r},
T_{2,4,r}$, prove J.~S.~Frame's formulas, and generalize
R.~Steinberg's theorem on the spectrum of the affine Coxeter
transformation for the multiply-laced diagrams. This theorem is
the key statement in R.~Steinberg's proof of the McKay
correspondence. For every extended Dynkin diagram, the spectrum of
the Coxeter transformation is easily obtained from R.~Steinberg's
theorem.

In the study of representations $\pi_n$ of $SU(2)$,
 we extend B.~Kostant's construction of a vector-valued generating
function
$$
  P_G(t)=([P_G(t)]_0, [P_G(t)]_1, \dots, [P_G(t)]_r).
$$
B.~Kostant's construction appears in the context of the McKay
correspondence and gives a way to obtain multiplicities of
indecomposable representations $\rho_i$ of the binary polyhedral
group $G$ in the decomposition of $\pi_n|G$. In the case of
multiply-laced graphs, instead of indecomposable representations
$\rho_i$ we use restricted representations and induced
representations of $G$ introduced by P.~Slodowy. Using
B.~Kostant's construction we generalize to the case of
multiply-laced graphs W.~Ebeling's theorem which connects the
Poincar\'{e} series $[P_G(t)]_0$ and the Coxeter transformations.
According to W.~Ebeling's theorem
\begin{equation*}
        [P_G(t)]_0 = \frac{\chi(t^2)}{\tilde{\chi}(t^2)},
\end{equation*}
where $\chi$ is the characteristic polynomial of the Coxeter
transformation and $\tilde{\chi}$ is the characteristic polynomial
of the corresponding affine Coxeter transformation.

Using the Jordan form of the Coxeter transformation we prove a
criterion of V.~Dlab and C.~M.~Ringel of regularity of quiver
representations, consider necessary and sufficient conditions of
this criterion for extended Dynkin diagrams and for diagrams with
indefinite Tits form.
\end{abstract}

\vspace{7mm} 


\maketitle


\newpage
~\\ ~\vspace{52mm}
\begin{center}
{\large In memory of V.~F.~Subbotin}
\end{center}

\tableofcontents
\newpage
\listoffigures
\listoftables
\newpage


\input 2prelim.tex
\input 3jordan.tex
\input 4shearing.tex
\input 5steinberg.tex
\input 6regular.tex

\appendix
\input 7mckay.tex
\input 8hurwitz.tex


\renewcommand{\appendixname}{}

  \printindex

\end{document}

%% file: 2prelim.tex

\chapter{\sc\bf Introduction}

\setlength{\epigraphwidth}{67mm}

\epigraph{
 ...A second empirical procedure for finding the exponents was
discovered by H.~S.~M.~Coxeter. He recognized that the exponents
can be obtained from a particular transformation $\gamma$ in the
Weyl group, which he had been studying, and which we take the
liberty of calling a Coxeter-Killing transformation...} {B.
~Kostant, \cite[p.974]{Kos59}, 1959.}

\section{The three historical aspects of the Coxeter transformation}
\label{history_1} The three areas, where the Coxeter
transformation plays a dramatic role, are:

$\bullet$ the theory of Lie algebras of the compact simple Lie groups;

$\bullet$ the representation theory of algebras and quivers;

$\bullet$ the McKay correspondence.

 \index{root system}
 \index{highest root}
 \index{Lie group}
 \index{Cartan subalgebra}
A {\it Coxeter transformation} or a {\it Coxeter element} is
defined as the product of all the reflections in the root system
of a compact Lie group. Neither the choice of simple roots nor the
ordering of reflections in the product affects its conjugacy
class, see \cite[Ch.5,\S6]{Bo}, see also Remark \ref{conj_class}.
H.~S.~M.~Coxeter studied these elements and their eigenvalues in
\cite{Cox51}.

Let $h$ be the order of the Coxeter transformation (called the
{\it  Coxeter number}), $|\varDelta|$ the number of roots in the
corresponding root system $\varDelta$, and $l$ the number of
eigenvalues of the Coxeter transformation, i.e., the rank of the
Cartan subalgebra. Then
\begin{equation}
\label{coxeter_result_1}
               hl = |\varDelta|.
\end{equation}
This fact was empirically observed by H.~S.~M.~Coxeter in
\cite{Cox51} and proved by A.~J.~Coleman in \cite{Col58}. One can
also find a proof of this fact
 in B.~Kostant's work \cite[p.1021]{Kos59}.

H.~S.~M.~Coxeter also observed that the order of the Weyl group is
equal to
\begin{equation}
\label{coxeter_result_2}
    (m_1 + 1)(m_2 + 1)...(m_l + 1),
\end{equation}
where the $m_i$ are the exponents of the eigenvalues of the
Coxeter transformation, and the factors $m_i + 1$ are the degrees
of $l$ basic polynomial invariants of the 
Weyl group. Proofs of these facts were obtained by C.~Chevalley
\cite{Ch55} and other authors; for historical notes, see
\cite{Bo}, \cite{Kos59}.

Let $\varDelta_{+}\subset \varDelta$ be the subset of positive
roots, let
\begin{equation*}
 \beta = n_1\alpha_1 + \dots + n_l\alpha_l,
\end{equation*}
where $\alpha_i \in \varDelta_{+}$, be the highest root in the
root system $\varDelta$. The Coxeter number $h$ from
(\ref{coxeter_result_1}) and coordinates $n_i$ of $\beta$ are
related as follows:
\begin{equation}
\label{coxeter_result_3}
    h = n_1 + n_2 + \dots + n_l + 1.
\end{equation}
Observation (\ref{coxeter_result_3}) is due to H.~S.~M.~Coxeter
\cite[p.234]{Cox49}, see also \cite[Th.1.4.]{Stb59}, \cite[Ch.6,
1, \S11,Prop.31]{Bo}.

The Coxeter transformation is important in the study of
representations of algebras, quivers, partially ordered sets
(posets) and lattices. Their distinguished role in this area is
related to the construction of the Coxeter functors given by
I.~N.~Bernstein, I.~M.~Gelfand, V.~A.~Ponomarev in \cite{BGP73}.
Further revelation of the role of Coxeter functors for
representations of algebras is due to V.~Dlab and C.~M.~Ringel
\cite{DR76}; for a construction of the functor $D{\rm Tr}$, see
M.~Auslender , M.~I.~Platzeck and I.~Reiten \cite{AuPR79},
\cite{AuRS95}. For an application of the Coxeter functors in the
representations of posets, see \cite{Drz74}; for their
applications in the representations of the modular lattices, see
\cite{GP74}, \cite{GP76}.

\index{orbit structure of the Coxeter transformation} Another area
where the  affine Coxeter transformations appeared is the McKay
correspondence --- a one-to-one correspondence between finite
subgroups of $SU(2)$ and simply-laced extended Dynkin diagrams.
Affine Coxeter transformations play the principal role in
R.~Steinberg's work \cite{Stb85} on the proof of the McKay
correspondence. B.~Kostant (\cite{Kos84}) obtains multiplicities
of the representations related to the concrete nodes of the
extended Dynkin diagram from the orbit structure of the affine
Coxeter transformation, see \S\ref{McKay}, \S\ref{orbit_str}.

In this work we only consider two areas of application of the
Coxeter transformation: representations of quivers and the McKay
correspondence. We do not consider other areas where the Coxeter
transformation plays an important role, such as the theory of
singularities of differentiable maps (the monodromy operator
coincides with a Coxeter transformation, see, e.g., \cite{A'C75},
\cite{Gu76}, \cite{ArGV86}, \cite{EbGu99}), Alexander polynomials,
pretzel knots, Lehmer's problem, growth series of Coxeter groups
(see, e.g., \cite{Lev66}, \cite{Hir02}, \cite{GH01},
\cite{McM02}).

\section{A brief review of this work}~

In Ch.\ref{chapter_prelim}, we recall some common definitions and
notions.

In Ch.\ref{jordan}, we establish general results about the Jordan form and the
spectrum of the
Coxeter transformation.

In Ch.\ref{chapter_shearing}, we give the eigenvalues of the
affine Coxeter transformation. After that we prove some shearing
formulas concerning the characteristic polynomials of the Coxeter
transformation. The main shearing formula is the Subbotin-Sumin
{\it splitting along the edge} formula which is extended in this
chapter to multiply-laced diagrams. One of applications of
shearing formulas is a construction of characteristic polynomials
of the Coxeter transformation for the hyperbolic Dynkin diagrams
$T_{2,3,r}$, $T_{3,3,r}$ and $T_{2,4,r}$. Two {\it Frame formulas}
from \cite{Fr51} (see Remark \ref{frame_formula_1} and Proposition
\ref{frame_formula_2}) are easily obtained from shearing formulas.

In Ch.\ref{chapter_steinberg}, we generalize a number of results
appearing in a context of the McKay correspondence to
multiply-laced diagrams. First, we consider R.~Steinberg's theorem
playing the key role in his proof of the McKay correspondence.
Essentially, R.~Steinberg observed that the orders of eigenvalues
of the affine Coxeter transformation corresponding to the {\it
extended} Dynkin diagram $\tilde\varGamma$ coincide with the
lengths of branches of the corresponding Dynkin diagram
$\varGamma$, \cite[p.591,$(*)$]{Stb85}.

Further, in Ch.\ref{chapter_steinberg}, we move on to B.~Kostant's
construction of a vector-valued {\it generating function} $P_G(t)$
\cite{Kos84}. Let $G$ be a binary polyhedral group, let $\rho_i$,
where $i = 0,\dots, r$, be irreducible representations of $G$
corresponding by the McKay correspondence to simple roots
$\alpha_i$ of the extended Dynkin diagram, let $\pi_n$, where $n =
0,1,\dots$, be indecomposable representations of $SU(2)$ in the
symmetric algebra ${\rm Sym}^{n}(\mathbb{C}^2)$, where ${\rm
Sym}^{n}(V)$ is the $n$th symmetric power of $V$. Let $m_i(n)$ be
{\it multiplicities} in the decomposition
\begin{equation}
  \label{def_pi_n}
   \pi_n|G = \sum\limits_{i=0}^r{m_i(n)}\rho_i;
\end{equation}
set
\begin{equation}
  \label{def_v_n}
v_n = \sum\limits_{i=0}^rm_i(n)\alpha_i.
\end{equation}
Then the Kostant generating function is defined as the following vector-valued function:
\begin{equation}
  \label{def_generating_funct}
   P_G(t)  =  ([P_G(t)]_0, [P_G(t)]_1, \ldots , [P_G(t)]_r)^t
           := \sum\limits_{n=0}^\infty{v_n}t^n.
\end{equation}
In particular, $[P_G(t)]_0$ is the Poincar\'{e} series of the
algebra of invariants ${\rm Sym}(\mathbb{C}^2)^G$, i.e.,
\begin{equation}
    [P_G(t)]_0 = P({\rm Sym}(\mathbb{C}^2)^G,t).
\end{equation}
B.~Kostant obtained explicit expression for the series $[P_G(t)]_0$,
and therefore found a way to calculate the multiplicities
$m_i(n)$. In Ch.\ref{chapter_steinberg}, we extend B.~Kostant's
construction to the case of  multiply-laced graphs. For this
purpose, we use P.~Slodowy's generalization  \cite[App.III]{Sl80}
of the McKay correspondence to the multiply-laced case. The main
idea of P.~Slodowy is to consider the pair of binary polyhedral
groups $H \triangleleft G$ and their {\it restricted
representations} $\rho\downarrow^G_H$ and {\it induced
representations} $\tau\uparrow^G_H$ instead of the representations
$\rho_i$.

 \index{induced representation $\tau\uparrow^G_H$}
 \index{restricted representation $\rho\downarrow^G_H$}
 \index{Slodowy correspondence}
 In Appendix \ref{chapter_slodowy}, we study in
detail P.~Slodowy's generalization for the case of the binary
octahedral group $\mathcal{O}$ and the binary tetrahedral group
$\mathcal{T}$, where $\mathcal{T} \triangleleft \mathcal{O}$. The
generalization of the McKay correspondence to the multiply-laced
case is said to be the {\it Slodowy correspondence}.

Generally, one can speak about the {\it McKay-Slodowy
correspondence}.

Finally, in Ch.\ref{chapter_steinberg} we
generalize to the multiply-laced case W.~Ebeling's theorem
\cite{Ebl02}: it relates the Poincar\'{e} series $[P_G(t)]_0$ and
the Coxeter transformations. According to W.~Ebeling's theorem,
\begin{equation}
        [P_G(t)]_0 = \frac{\chi(t^2)}{\tilde{\chi}(t^2)},
\end{equation}
where $\chi$ is the characteristic polynomial of the Coxeter
transformation ${\bf C}$ and $\tilde{\chi}$ is the characteristic
polynomial of the corresponding affine Coxeter transformation
${\bf C}_a$, see Theorem \ref{theorem_ebeling}.

The results of Ch.\ref{chapter_shearing} and Ch.\ref{chapter_steinberg}
are published for the first time.

In Ch.\ref{chap_regular}, using results on the Jordan form of the
Coxeter transformation we prove a criterion of V.~Dlab and
C.M.~Ringel of regularity of quiver representations, we consider
necessary and sufficient conditions of this criterion for extended
Dynkin diagrams and for diagrams with indefinite Tits form.

\section{The Spectrum and the Jordan form}
\subsection{The Jordan form and reduction to the golden pair of matrices}
In this work we review the research started 30 years ago in the
teamwork with V.F.~Subbotin and discuss the Jordan form, and the
spectrum of the Coxeter transformation. We show that the study of
eigenvalues of the Coxeter transformation reduces to the study of
the {\it golden pair} of matrices satisfying conditions of the
Perron-Frobenius theorem \cite{MM64}, \cite{Ga90} and having other
nice properties.

In \S3 we recall (see (\ref{symmetric_B}), (\ref{matrix_K}), (\ref{matrix_T}))
that Cartan matrix $K$ can be represented in the
form
\begin{equation}
 K =
   \begin{cases}
        2{\bf B} & \text{for $K$ symmetric,  see  (\ref{symmetric_B})} \\
        T{\bf B} & \text{for $K$ symmetrizable, see
        (\ref{matrix_K}),}
   \end{cases}
\end{equation}
where $T$ is a diagonal matrix, {\bf B} is a symmetric matrix.
 \index{multiply-laced case}
 \index{simply-laced case}
 For simply-laced (resp. multiply-laced) graphs, the {\it golden pair}
of matrices is
$$
  DD^t \text{ and } D^tD \text{ (resp. } DF \text{ and } FD),
$$
where the matrices $D$ and $F$ are found from the expressions
\begin{equation}
  \label{symmetric_B_1}
{\bf B} =
\begin{pmatrix}
    I_m & D   \\
    D^t & I_k
\end{pmatrix}
\end{equation}
for the simply-laced case and
\begin{equation}
  \label{matrix_K_1}
K =
\begin{pmatrix}
    2I_m & 2D   \\
    2F & 2I_k
\end{pmatrix}
\end{equation}
for the multiply-laced case.

The mentioned reduction method works only for trees (they have
symmetrizable Cartan matrices) and for graphs with even cycles.
 In \S\ref{biblio_notes} we give
some bibliographical remarks regarding generalized Cartan matrices.

According to Proposition \ref{gold_pair}, relation
(\ref{Coxeter_eigenvalues}), the corresponding eigenvalues
$\lambda^{\varphi_i}_{1,2}$ of Coxeter transformations are
\begin{equation*}
        \lambda^{\varphi_i}_{1,2} =
        2\varphi_i - 1 \pm 2\sqrt{\varphi_i(\varphi_i-1)},
\end{equation*}
where the $\varphi_i$ are the eigenvalues of $DD^t$ and $D^tD$.

One of the central results of this work is Theorem \ref{th_jordan}
(\S\ref{jordan}, Fig.~\ref{jordan_form}) on Jordan form
\cite{SuSt75, SuSt78, St85}:

\medskip

1) {\it The Jordan form of the Coxeter transformation is diagonal
if and only if the Tits form is non-degenerate.}

2) {\it If $\mathcal{B}$ is nonnegative definite, $\mathcal{B}
\geq 0$ ($\varGamma$ is an extended Dynkin diagram), then the
Jordan form of the Coxeter transformation contains only one
$2\times2$ Jordan block. All eigenvalues $\lambda_i$ lie on the
unit circle.}

3) {\it If $\mathcal{B}$ is indefinite and degenerate, then the
number of two-dimensional Jordan blocks coincides with
$\dim\ker{\bf B}$. All other Jordan blocks are $1\times1$. There
is simple maximal $\lambda^{\varphi_1}_1$ and simple minimal
$\lambda^{\varphi_1}_2$ eigenvalues such that}
$$
   \lambda^{\varphi_1}_1 > 1,   \hspace{5mm} \lambda^{\varphi_1}_2 < 1.
$$
\medskip

C.~M.~Ringel \cite{Rin94} generalized
the result of Theorem \ref{th_jordan}
for non-symmetrizable Cartan matrices, see \S\ref{biblio_notes}.

\subsection{An explicit construction of eigenvectors. The eigenvalues are roots of
unity} \label{expl_creation} An important point of this work is an
explicit construction of eigenvectors and adjoint vectors of the
Coxeter transformation --- the vectors that form a Jordan basis,
see Proposition \ref{eigenvectors_basis}, \cite{SuSt75, SuSt78}.
This construction is used to obtain the necessary and sufficient
condition for regularity of representations \cite{SuSt75, SuSt78}.
This condition was also found by V.~Dlab and C.M.~Ringel,
\cite{DR76}, see \S\ref{criteria_DR}.

The eigenvalues for all cases of extended Dynkin diagrams are
easily calculated using a generalized R.~Steinberg theorem
 (Theorem \ref{gen_Steinberg}) and Table
\ref{table_char_ext_polynom_Dynkin}.

\index{roots of unity} Theorem \ref{roots_unity} and Table
\ref{table_eigenvalues} summarize this calculation as follows
(\cite{SuSt79, St81, St82a, St85}):
\medskip

{\it the eigenvalues of the Coxeter transformation for any
extended Dynkin diagram are {\bf roots of unity.}}

\medskip

The case of $\widetilde{A}_n$ is considered in the
\S\ref{biblio_notes}. According to formula (\ref{char_polyn_An}),
the characteristic polynomial of the Coxeter transformation for
$\widetilde{A}_n$ is (\cite{MeSu82}, \cite{Men85}, \cite{Col89},
 \cite{Shi00}, \cite{BT97}):
$$
   (\lambda^{n-k+1}-1)(\lambda^{k}-1),
$$
where $k$ is the number characterizing the conjugacy class of the
Coxeter transformation. For $\widetilde{A}_n$, there are $[n/2]$
characteristic polynomials \cite{Col89}, (see
\S\ref{biblio_notes}), \index{roots of unity} but, for all these
cases, the eigenvalues of the Coxeter transformation are roots of
unity.

\subsection{Study of the Coxeter transformation and the Cartan matrix}
To study of the Coxeter transformation is almost the same
as to study of the Cartan matrix. The Cartan matrix and
the matrix of the Coxeter transformation
(more precisely, the bicolored representative of the conjugacy
class of the Coxeter transformation, \S\ref{bicolored_C})
are constructed from the same blocks,
see relations (\ref{symmetric_B}), (\ref{matrix_K}).

By Proposition \ref{fixed_points}, the eigenvalues $\lambda$ of
the Coxeter transformation and the eigenvalues $\gamma$ of the
matrix {\bf B} of the quadratic Tits form are related as follows:
$$
    \frac{(\lambda+1)^2}{4\lambda} = (\gamma-1)^2.
$$
By Corollary \ref{diag_form}, the Jordan form of the Coxeter
transformation is diagonal if and only if the Tits form is
nondegenerate.

It seems that the Coxeter transformation contains more information
than the Cartan matrix, namely, the Coxeter transformation
contains additional information about orientation, see
\S\ref{orientation} and considerations on the graphs containing
cycles in \S\ref{biblio_notes}.

\index{fixed points}
\index{anti-fixed points}
The Coxeter transformation and
the Cartan matrix are also related  to
the {\it fixed points} and {\it anti-fixed
points}  of the powers of the Coxeter transformation.  This connection
is given by means of the Chebyshev polynomials in
Theorem \ref{fixed_main_proposition} (\cite{SuSt82}).

\subsection{Monotonicity of the dominant eigenvalue of the golden pair}
According to Corollary \ref{corollary_dominant}, the matrices
$DD^t$ and $D^tD$ have a common simple positive eigenvalue
$\varphi_1$, the maximal eigenvalue. This eigenvalue is said to be
the {\it dominant} eigenvalue. \index{dominant eigenvalue} The
dominant eigenvalue $\varphi_1$ is a certain characteristic of the
graph $\varGamma$. In Proposition \ref{trend_fi}, we show that if
any edge is added to $\varGamma$, then this characteristic only
grows. The same is true for the maximal eigenvalue
$\lambda^{\varphi_1}_1$, see Proposition \ref{gold_pair}, relation
(\ref{Coxeter_eigenvalues}).

\medskip
{\bf Problem}. Is the dominant eigenvalue an invariant of the
graph in the class of trees, $\mathbb{T}$, i.e., is there
one-to-one correspondence between the dominant eigenvalue and the graph
if the Tits form of the graph is indefinite?

If there exist two graphs $\varGamma_1$ and $\varGamma_2$ with the
same dominant eigenvalue $\varphi_1$, and $\theta$ is the
assignment of the dominant eigenvalue to a graph,
then what class of graphs $\mathbb{T}/\theta$  we obtain modulo the relation
given by $\theta$?

\medskip

This problem is solved in this work for
the
diagrams $T_{2,3,r}$, $T_{3,3,r}$, $T_{2,4,r}$,
see Propositions \ref{polyn_T_1}, \ref{polyn_T_2}, \ref{polyn_T_3}
and Tables \ref{table_char_E_series}, \ref{table_char_E6_series},
\ref{table_char_E7_series} in \S\ref{section_T_pqr}.

\section{Shearing formulas and the diagrams $T_{p,q,r}$}
\subsection{Shearing formulas for the characteristic polynomial}

 \index{splitting along the edge}
 \index{Subbotin-Sumin formula! - simply-laced case}
 \index{formula! - Subbotin-Sumin, simply-laced case}
There is a number of recurrent formulas used to calculate the
characteristic polynomial of the Coxeter transformation of a given
graph in terms of characteristic polynomials of the Coxeter
transformation of the graph's components. Subbotin and Sumin
proved the formula of {\it splitting along the edge}
\cite{SuSum82}:

\begin{equation}
 \label{sub_sum_intro_1}
   \mathcal{X}(\varGamma, \lambda) =
       \mathcal{X}(\varGamma_1, \lambda)\mathcal{X}(\varGamma_2, \lambda) -
       \lambda\mathcal{X}(\varGamma_1\backslash\alpha, \lambda)
              \mathcal{X}(\varGamma_2\backslash\beta, \lambda).
\end{equation}
The proof of the Subbotin-Sumin formula is given
in Proposition
\ref{split_edge}.

Another formula (\cite{KMSS83}) is given in Proposition \ref{gluing}.

V.~Kolmykov kindly informed me that the following statement holds:

\begin{proposition}
If $\lambda$ is the eigenvalue of Coxeter transformations for
graphs $\varGamma_1$ and $\varGamma_2$, then $\lambda$ is also the
eigenvalue of the graph $\varGamma$ obtained by gluing
$\varGamma_1$ and $\varGamma_2$.
\end{proposition}

For details, see Proposition \ref{merge_2G}.

\index{splitting along the weighted edge}
In Proposition \ref{split_edge_2} we generalize
formula (\ref{sub_sum_intro_1}) to the multiply-laced case.
The formula of {\it splitting along the weighted edge} holds:
\index{multiply-laced case}
\begin{equation}
  \label{splitting_intro_2}
   \mathcal{X}(\varGamma, \lambda) =
       \mathcal{X}(\varGamma_1, \lambda)\mathcal{X}(\varGamma_2, \lambda) -
       \rho\lambda\mathcal{X}(\varGamma_1\backslash\alpha, \lambda)
              \mathcal{X}(\varGamma_2\backslash\beta, \lambda),
\end{equation}
where $\rho$ is the ratio of the lengths of the roots,
corresponding to the endpoints of the split edge. Corollary
\ref{split_edge_corol_2} deals with the case where $\varGamma_2$
contains only one point. In this case, we have
\begin{equation}
  \label{splitting_intro_3}
   \mathcal{X}(\varGamma, \lambda) =
       -(\lambda + 1)\mathcal{X}(\varGamma_1, \lambda) -
       \rho\lambda\mathcal{X}(\varGamma_1\backslash\alpha, \lambda).
\end{equation}

A formula similar to (\ref{splitting_intro_2}) can be proved not
only for the Coxeter transformation, but, for example, for the
Cartan matrix, see \`E.~B.~Vinberg's paper \cite[Lemma
5.1]{Vin85}. See also Remark \ref{frame_formula_1} concerning the
works of J.~S.~Frame and S.~M.~Gussein-Zade.

\subsection{An explicit calculation of characteristic polynomials}
We use recurrent formulas (\ref{sub_sum_intro_1}) --
(\ref{splitting_intro_3}) to calculate the characteristic
polynomials of the Coxeter transformation for the Dynkin diagrams
and extended Dynkin diagrams.

\begin{table} 
  \centering
  \vspace{2mm}
  \caption{\hspace{3mm}Characteristic polynomials, Dynkin diagrams}
   \renewcommand{\arraystretch}{1.5}
  \begin{tabular} {|c|c|c|}
  \hline \hline
        Dynkin     & Characteristic  & Form with \cr
        diagram &  polynomial     & denominator \\
   \hline \hline
     & & \cr
       ${A}_n$
     & $\lambda^n + \lambda^{n-1} + \dots + \lambda + 1$
     & $\displaystyle\frac{\lambda^{n+1} - 1}{\lambda - 1}$ \cr
     & & \\
  \hline
     & & \cr
       ${B}_n, {C}_n$
     & $\lambda^{n} +  1$ & \cr
     & &
     \\
  \hline
     & & \cr
       ${D}_n$
     & $\lambda^{n} +  \lambda^{n-1} + \lambda + 1$
     & $(\lambda + 1)(\lambda^{n-1} + 1)$ \cr
     & & \\
  \hline
    & & \cr
       ${E}_6$
     & $\lambda^{6} +  \lambda^{5} - \lambda^3 + \lambda + 1$
     & $\displaystyle\frac{(\lambda^{6} + 1)}{(\lambda^2 + 1)}
       \frac{(\lambda^{3} - 1)}{(\lambda - 1)}$ \cr
     & &  \\
  \hline
     & & \cr
       ${E}_7$
     & $\lambda^7 + \lambda^6 - \lambda^4 - \lambda^3 + \lambda + 1$
     & $\displaystyle\frac{(\lambda + 1)(\lambda^{9} + 1)}{(\lambda^3 + 1)}$ \cr
     & & \\
  \hline
     & & \cr
       ${E}_8$
     & $\lambda^8 + \lambda^7 - \lambda^5 - \lambda^4 - \lambda^3 + \lambda + 1$
     & $\displaystyle\frac{(\lambda^{15} + 1)(\lambda + 1)}{(\lambda^5 +
1)(\lambda^3  + 1)}$ \cr
     & & \\
  \hline
     & & \cr
       ${F}_4$
     & $\lambda^4  - \lambda^2  + 1$
     & $\displaystyle\frac{\lambda^{6} + 1}{\lambda^2 + 1}$ \cr
     & & \\
  \hline
     & & \cr
       ${G}_2$
     & $\lambda^2  - \lambda  + 1$
     & $\displaystyle\frac{\lambda^{3} + 1}{\lambda + 1}$ \cr
     & & \\
  \hline \hline
\end{tabular}
  \label{table_char_polynom_Dynkin}
\end{table}

\begin{table} 
  \centering
  \vspace{2mm}
  \caption{\hspace{3mm}Characteristic polynomials,
                        extended Dynkin diagrams}
  \renewcommand{\arraystretch}{1.3} 
  \begin{tabular} {|c|c|c|c|}
  \hline \hline
      Extended & Characteristic  & Form with & Class \cr
      Dynkin   & polynomial      &  $\chi_i$ &   $g$   \cr
      diagram  &                 &           &       \\
  \hline \hline
     & & & \cr
       $\tilde{D}_4$
     & $(\lambda - 1)^2(\lambda + 1)^3$
     & $(\lambda - 1)^2\chi^3_1$ & 0 \\
  \hline
     & & & \cr
       $\tilde{D}_n$
     & $(\lambda^{n-2} -  1)(\lambda - 1)(\lambda + 1)^2$
     & $(\lambda - 1)^2\chi_{n-3}\chi^2_1$ & 0 \\
  \hline
     & & & \cr
       $\tilde{E}_6$
     & $(\lambda^3 -  1)^2(\lambda + 1)$
     & $(\lambda - 1)^2\chi^2_2\chi_1$ & 0 \\
  \hline
     & & & \cr
       $\tilde{E}_7$
     & $(\lambda^4 -  1)(\lambda^3 -  1)(\lambda + 1)$
     & $(\lambda - 1)^2\chi_3\chi_2\chi_1$ & 0 \\
  \hline
     & & & \cr
       $\tilde{E}_8$
     & $(\lambda^5 -  1)(\lambda^3 -  1)(\lambda + 1)$
     & $(\lambda - 1)^2\chi_4\chi_2\chi_1$ & 0 \\
  \hline
     & & & \cr
       $\widetilde{CD}_n, \widetilde{DD}_n$
     & $(\lambda^{n-1} -  1)(\lambda^2 -1)$
     & $(\lambda-1)^2\chi_{n-2}\chi_1$  & 1 \\
  \hline
     & & & \cr
       $\tilde{F}_{41}, \tilde{F}_{42}$
     & $(\lambda^2 -  1)(\lambda^3 - 1)$
     & $(\lambda - 1)^2\chi_2\chi_1$ &  1 \\
  \hline
     & & & \cr
       $\tilde{B}_n, \tilde{C}_n, \widetilde{BC}_n$
     & $(\lambda^n -  1)(\lambda - 1)$
     & $(\lambda - 1)^2\chi_{n-1}$ &  2\\
  \hline
     & & & \cr
       $\tilde{G}_{21}, \tilde{G}_{22}$
     & $(\lambda-1)^2(\lambda  + 1)$
     & $(\lambda-1)^2\chi_1$ & 2 \\
  \hline
     & & & \cr
       ${A}_{11}, {A}_{12}$
     & $(\lambda - 1)^2$
     & $(\lambda - 1)^2$ &  3 \\
  \hline \hline
\end{tabular}
  \label{table_char_ext_polynom_Dynkin}
\end{table}

\index{characteristic polynomial of the Coxeter transformation! - Dynkin
diagrams}
The characteristic polynomials of the Coxeter transformations
for the Dynkin diagrams are presented in Table
\ref{table_char_polynom_Dynkin}. For calculations, see
\S\ref{case_Dynkin_diagr}.

\index{characteristic polynomial of the Coxeter transformation! -
extended Dynkin diagrams} The characteristic polynomials of the
Coxeter transformations for the extended Dynkin diagrams are presented
in Table \ref{table_char_ext_polynom_Dynkin}. The polynomials
$\chi_n$ from Table \ref{table_char_ext_polynom_Dynkin} are
characteristic polynomials of the Coxeter transformation for the
Dynkin diagram ${A}_n$:
$$
   \chi_n = (-1)^n\mathcal{X}({A}_n),
$$
defined (see \S\ref{shearing} and Remark \ref{up_to_sign}) to be
\begin{equation*}
  \begin{split}
    \chi_n =
       \frac{\lambda^{n+1} - 1}{\lambda - 1} =
       \lambda^{n} + \lambda^{n-1} + ... + \lambda^2 + \lambda + 1.
  \end{split}
\end{equation*}
For calculations of characteristic polynomials
for the extended Dynkin diagrams, see
\S\ref{Steinberg_simply_laced} and \S\ref{sect_gen_Steinberg}.

\subsection{Formulas for the diagrams $T_{2,3,r}, T_{3,3,r}, T_{2,4,r}$}
 For the three classes of diagrams ---
$T_{2,3,r}, T_{3,3,r}, T_{2,4,r}$ --- explicit formulas of
characteristic polynomials of the Coxeter transformations are
obtained, see \S\ref{hyperbolic}, Fig.~\ref{T_pqr_diagram}.

\index{${E}_n$-series} \index{characteristic polynomial of the
Coxeter transformation! - $T_{2,3,r}$} The case of $T_{2,3,r}$,
where $r \geq 2$, contains diagrams ${D}_5, {E}_6, {E}_7, {E}_8,
\tilde{E}_8, {E}_{10}$, and so we call these diagrams the
${E}_n$-{\it series}, where $n = r+3$. The diagram $T_{2,3,7}$ is
{\it hyperbolic}, (see \S\ref{hyperbolic}) and, for all $r \geq
3$, we have
\begin{equation}
  \chi(T_{2,3,r}) =
     \lambda^{r+3} + \lambda^{r+2} -
     \sum\limits_{i=3}^{r}\lambda^{i} + \lambda + 1,
\end{equation}
see (\ref{E_series}) and Table \ref{table_char_E_series}. The
spectral radius of $\chi(T_{2,3,r})$ converges to the {\it Zhang
number}\footnote{Hereafter we give all such numbers with six
decimal points.} \index{Zhang number} \index{spectral radius}
\begin{equation}
  \label{zhang_1}
   \sqrt[3]{\frac{1}{2} + \sqrt{\frac{23}{108}}} +
   \sqrt[3]{\frac{1}{2} - \sqrt{\frac{23}{108}}} \approx 1.324717...,
\end{equation}
as $r \rightarrow \infty$, see Proposition \ref{polyn_T_1} and
Remark \ref{rem_Zhang}. The number (\ref{zhang_1}) is also the
smallest {\it Pisot number},\index{Pisot number} see
\S\ref{sect_numbers}.

\index{${E}_{6,n}$-series} \index{characteristic polynomial of the
Coxeter transformation! - $T_{3,3,r}$} The case of $T_{3,3,r}$,
where $r \geq 2$, contains diagrams ${E}_6, \tilde{E}_6$, and so
we call these diagrams the ${E}_{6,n}$-{\it series}, where $n =
r-2$. The diagram $T_{3,3,4}$ is {\it hyperbolic}, (see
\S\ref{hyperbolic}) and, for all $r \geq 3$, we have
\begin{equation}
  \chi(T_{3,3,r}) =
     \lambda^{r+4} + \lambda^{r+3}
     -2\lambda^{r+1}
     -3\sum\limits_{i=4}^{r}\lambda^{i}
     -2\lambda^3 + \lambda + 1,
\end{equation}
see (\ref{E6_series}) and Table \ref{table_char_E6_series}. The
spectral radius of $\chi(T_{3,3,r})$ converges to
\begin{equation*}
   \frac{\sqrt{5} + 1}{2} \approx 1.618034...,
\end{equation*}
as $r \rightarrow \infty$, see Proposition \ref{polyn_T_2}.

\index{${E}_{7,n}$-series}
\index{characteristic polynomial of the
Coxeter transformation! - $T_{2,4,r}$}
\index{spectral radius}

The case of $T_{2,4,r}$, where $r \geq 2$, contains diagrams
${D}_6, {E}_7, \tilde{E}_7$, and so we call these diagrams the
${E}_{7,n}$-{\it series}, where $n = r-3$. \index{hyperbolic
diagrams} The diagram $T_{2,4,5}$ is {\it hyperbolic}, (see
\S\ref{hyperbolic}) and, for all $r \geq 4$, we have
\begin{equation}
  \chi(T_{2,4,r}) =
     \lambda^{r+4} + \lambda^{r+3}
     -\lambda^{r+1}
     -2\sum\limits_{i=4}^{r}\lambda^{i}
     -\lambda^3 + \lambda + 1,
\end{equation}
see (\ref{E7_series}) and Table \ref{table_char_E7_series}. The
spectral radius of $\chi(T_{2,4,r})$ converges to
\begin{equation*}
   \frac{1}{3} + \sqrt[3]{\frac{58}{108} + \sqrt{\frac{31}{108}}} +
   \sqrt[3]{\frac{58}{108} - \sqrt{\frac{31}{108}}} \approx 1.465571... \,
\end{equation*}
as $r \rightarrow \infty$, see Proposition \ref{polyn_T_3}.

\section{Coxeter transformations and the McKay correspondence}
\subsection{The generalized R.~Steinberg theorem}
Here we generalize R.~Steinberg's theorem concerning the
mysterious connection between lengths of branches of any Dynkin
diagram and orders of eigenvalues of the affine Coxeter
transformation. R.~Steinberg proved this theorem for the
simply-laced case in \cite[p.591,$(*)$]{Stb85}; it was a key
statement in his explanation of the phenomena of the McKay
correspondence.

 \index{multiply-laced case}
 We prove the R.~Steinberg theorem for the simply-laced case in
 \S\ref{Steinberg_simply_laced}, Theorem \ref{Steinberg}.
 The multiply-laced case (generalized
R.~Steinberg's theorem) is proved in \S\ref{sect_gen_Steinberg},
Theorem \ref{gen_Steinberg}. Essentially, the generalized
R.~Steinberg theorem immediately follows from Table
\ref{table_char_ext_polynom_Dynkin}.

\begin{theorem} [The generalized R.~Steinberg theorem]
The affine Coxeter transformation with the extended Dynkin diagram
$\tilde{\varGamma}$ has the same eigenvalues as the product of
some Coxeter transformations of types ${A}_i$, where $i \in \{p-1,
q-1, r-1\}$ and matches to
  $(3-g)$ branches of the corresponding Dynkin diagram $\varGamma$,
  where $g$ is the class number (\ref{class_g}) of
  $\tilde\varGamma$. In other words,
\begin{equation*}
 \begin{split}
  & \text{ For } g = 0, \text{ the product }
    \chi_{p-1}\chi_{q-1}\chi_{r-1}, \text{ is taken}. \\
  & \text{ For } g = 1, \text{ the product }
    \chi_{p-1}\chi_{q-1}, \text{ is taken}. \\
  & \text{ For } g = 2, \text{ the product consists of only from }
     \chi_{p-1}. \\
  & \text{ For } g = 3, \text{ the product is trivial (= 1)}. \\
 \end{split}
\end{equation*} \end{theorem}

For details, see Theorem \ref{gen_Steinberg}.

\subsection{The Kostant generating functions and W.~Ebeling's theorem}
 \label{multipl_kostant}
Now we consider B.~Kostant's construction  of the vector-valued
{\it generating function} $P_G(t)$ \cite{Kos84}. Let $G$ be a
binary polyhedral group, and $\rho_i$, where $i = 0,\dots, r$,
irreducible representations of $G$ corresponding due to the McKay
correspondence to simple roots $\alpha_i$ of the extended Dynkin
diagram; let $\pi_n$, where $n = 0,1,\dots$, be irreducible
representations of $SU(2)$ in ${\rm Sym}^n(\mathbb{C}^2)$. Let
$m_i(n)$ be {\it multiplicities} in the decomposition
\begin{equation}
   \pi_n|G = \sum\limits_{i=0}^r{m_i(n)}\rho_i,
\end{equation}
and so $m_i(n) = <\pi_n|G, \rho_i>$; set
\begin{equation}
   v_n = \sum\limits_{i=0}^rm_i(n)\alpha_i =
    \left (
      \begin{array}{c}
         m_0(n) \\
         \dots  \\
         m_r(n).
      \end{array}
    \right )
    \vspace{5mm}
\end{equation}
Then
\begin{equation}
  \label{def_P_G_2}
   P_G(t) = \sum\limits_{n=0}^\infty{v_n}t^n =
   \left (
     \begin{array}{c}
        \sum\limits_{n=0}^\infty{m_0(n)}t^n \\
        \dots \\
        \sum\limits_{n=0}^\infty{m_r(n)}t^n
     \end{array}
   \right ).
\end{equation}
Thus, $P_G(t)=([P_G(t)]_0, [P_G(t)]_1, \ldots , [P_G(t)]_r)^t$ is
a vector-valued series. In particular, $[P_G(t)]_0$ is the
Poincar\'{e} series of the algebra of invariants ${\rm
Sym}(\mathbb{C}^2)^G$, i.e.,
\begin{equation}
    [P_G(t)]_0 = P({\rm Sym}(\mathbb{C}^2)^G,t).
\end{equation}

B.~Kostant obtained explicit formulas for the series $[P_G(t)]_0$,
and therefore a way to calculate the multiplicities $m_i(n)$.
B.~Kostant's construction is generalized in
Ch.~\ref{chapter_steinberg} to the multiply-laced case. For this
purpose, we use the P.~Slodowy generalization
\cite[App.III]{Sl80} of the McKay correspondence to the
multiply-laced case.

The main idea of P.~Slodowy is to consider a pair of binary
polyhedral groups $H \triangleleft G$ and their {\it restricted
representations} $\rho\downarrow^G_H$ and {\it induced
representations} $\tau\uparrow^G_H$ instead of representations
$\rho_i$. In Appendix \ref{chapter_slodowy}, we study in detail
P.~Slodowy's generalization for the pair $\mathcal{T}
\triangleleft \mathcal{O}$, where $\mathcal{O}$ is the binary
octahedral group and $\mathcal{T}$ is the binary tetrahedral
group. We call the generalization of the McKay correspondence to
the multiply-laced case the {\it Slodowy correspondence}. Finally,
in Ch.~\ref{chapter_steinberg}, we generalize to the
multiply-laced case W.~Ebeling's theorem \cite{Ebl02} which
relates the Poincar\'{e} series $[P_G(t)]_0$ and the Coxeter
transformations.

First, we prove the following proposition due to B.~Kostant \cite{Kos84}. It
holds for the McKay operator and also for the Slodowy operator.

\begin{proposition} [B.~Kostant \cite{Kos84}]
  If $B$ is either the McKay operator $A$ or the Slodowy operator
  $\tilde{A}$ or $\tilde{A}^\vee$, then
 \begin{equation*}
      Bv_n = v_{n-1} + v_{n+1}.
 \end{equation*}
\end{proposition}
For details, see Proposition \ref{kostant_prop} in \S\ref{sect_char_McKay}.

\begin{theorem} [generalized W.~Ebeling theorem, \cite{Ebl02}]
  \label{gen_ebeling_1}
  Let $G$ be a binary polyhedral group and $[P_G(t)]_0$ the Poincar\'{e}
  series (\ref{poincare_alg_inv}) of the algebra of invariants
${\rm Sym}(\mathbb{C}^2)^G$. Then
\begin{equation*}
       [P_G(t)]_0 = \frac{\det{M_0}(t)}{\det{M}(t)},
\end{equation*}
where
\begin{equation*}
   \det{M}(t) = \det|t^{2}I - {\bf C}_a|, \hspace{5mm}
   \det{M_0}(t) = \det|t^{2}I - {\bf C}|,
\end{equation*}
and where {\bf C} is the Coxeter transformation and ${\bf C}_a$ is
the corresponding affine Coxeter transformation.
\end{theorem}
In Theorem \ref{gen_ebeling_1} the Coxeter transformation {\bf C}
and the affine Coxeter transformation ${\bf C}_a$ are related to
the binary polyhedral group $G$. For the multiply-laced case, we
consider a pair of binary polyhedral groups $H \triangleleft G$
and {\bf C}, ${\bf C}_a$ are related again to the group $G$. We
generalize W.~Ebeling's theorem for the multiply-laced case, see
Theorem \ref{theorem_ebeling}. For definition of the Poincar\'{e}
series for the multiply-laced case, see (\ref{poincare_alg_inv_2})
from \S\ref{generating_fun} and Remark \ref{poincare_m_case}.

\section{The regular representations of quivers}
\subsection{Finite-type, tame and wild quivers}
  \label{tame_quivers}
In \cite{BGP73}, Bernshtein, Gelfand and Ponomarev introduced {\it
regular representations} of quivers, i.e., representations that
never vanish under the Coxeter transformations. More precisely,
let $V$ be any representation in the category of representations
$\mathcal{L}(\varGamma, \varDelta)$ of the graph $\varGamma$ with
an orientation $\varDelta$, and let $\varPhi^+$, $\varPhi^-$ be
the Coxeter functors (\cite{BGP73}, \cite{DR76}). The
representation $V$ is said to be {\it regular} if
\begin{equation}
  \label{regular_repr}
   (\varPhi^+)^k(\varPhi^-)^k(\dim{V}) = (\varPhi^+)^k(\varPhi^-)^k(\dim{V}) = \dim{V}
   \text{ for all } k \in \mathbb{Z}.
\end{equation}

 \index{finite-type quivers}
 \index{theorem! - Gabriel}

The Dynkin graphs do not have regular representations; in the
category of graph representations, only a finite set of
irreducible representations is associated to any Dynkin diagram.
Graphs with such property are called {\it finite type quivers}.
According to P.~Gabriel's theorem \cite{Gab72}, {\it a quiver is
of finite type if and only if it is a (simply laced) Dynkin
diagram}.

In the category of all representations of a given quiver, the
regular representations are the most complicated ones; they have
been completely described only for the extended Dynkin diagrams
\index{tame quivers} which for this reason,  in the representation
theory of quivers, were called {\it tame quivers} (\cite{Naz73},
\cite{DR76}).

\index{wild quivers} Any quiver with indefinite Tits form
$\mathcal{B}$ is {\it wild}, i.e., the description of its
representations contains the problem of classifying pairs of
matrices up to simultaneous similarity; this classification is
hopeless in a certain sense, see \cite{GP69}, \cite{Drz80}.

By (\ref{regular_repr}), the representation $V$ is regular if and
only if
\begin{equation}
  \label{regular_repr_1}
    {\bf C}^k(\dim{V}) > 0 \text{ for all } k \in \mathbb{Z},
\end{equation}
see \cite{DR76}, \cite{St75}, \cite{SuSt75}, \cite{SuSt78}.

\subsection{The Dlab-Ringel defect $\delta$ and $\varDelta$-defect $\rho$}

In \cite{St75}, \cite{SuSt75}, \cite{SuSt78}, the linear form
$\rho_{\varDelta{'}}(z)$ was considered, see Definition
\ref{rho_defect}. For the simply-laced case, it is defined as
 \index{simply-laced case}
\begin{equation}
 \label{rho_simply_intro}
  \rho_{\varDelta{'}}(z) = <Tz, \tilde{z}^1>,
\end{equation}
   and for the multiply-laced case as
\index{multiply-laced case}
\index{defect of the vector}
\begin{equation}
 \label{rho_multi_intro}
  \rho_{\varDelta{'}}(z) = <Tz, \tilde{z}^{1\vee}>.
\end{equation}
 Here ${z}^{1}$ is the eigenvector of the Coxeter transformation
 corresponding to eigenvalue $1$, and $z^{1\vee}$ is the
 eigenvector corresponding to eigenvalue $1$
 of the Coxeter transformation for the dual diagram
 $\varGamma^\vee$;
\index{dual vector} $\tilde{v}$ denotes the {\it dual} vector to
$v$ obtained from $v$ by changing the sign of the
$\mathbb{Y}$-component (see Definition \ref{conjug_vector}).

The linear form $\rho_{\varDelta{'}}(z)$ is said to be the
$\varDelta{'}${\it-defect} of the vector $z$, or the
   {\it defect of the vector $z$ in the orientation} $\varDelta{'}$.
\index{Dlab-Ringel definition of defect}
\index{$\varDelta{'}$-defect}

In \cite{DR76}, V.~Dlab and C.~M.~Ringel introduced the {\it
defect} $\delta_{\varDelta{'}}$ obtained as the solution to the
following equation, see (\ref{DR_defect}):
\begin{equation*}
     {\bf C}^*_{\varDelta{'}}\delta_{\varDelta{'}} = \delta_{\varDelta{'}}.
\end{equation*}

In Proposition \ref{defect_coinciding}, \cite{St85} we show that
the Dlab-Ringel defect $\delta_{\varDelta{'}}$ coincides with the
$\varDelta{'}$-defect $\rho_{\varDelta{'}}$.

\subsection{Necessary and sufficient conditions of regularity}

In Theorem \ref{necess_codnd_reg} proved by Dlab-Ringel
\cite{DR76} and by Subbotin-Stekolshchik \cite{SuSt75},
\cite{SuSt78}, we show the necessary condition of regularity of
the representation $V$:

\index{regular vector}
\begin{center}
{\it If $\dim{V}$ is the regular vector for
the extended Dynkin diagram $\varGamma$ \\
in the orientation $\varDelta{'}$, then
\begin{equation}
  \label{nes_suf_cond_1}
      \rho_{\varDelta{'}}(\dim{V}) = 0.
\end{equation}
}
\end{center}

In Proposition \ref{prop_sufficient} (for the bicolored
orientation) and in Proposition \ref{arbitr_orient} (for an
arbitrary orientation $\varDelta{'}$) we show that (see
\cite{St82})
\begin{center}
{\it the condition (\ref{nes_suf_cond_1}) is also sufficient if
$\dim{V}$ is a positive root in the root system associated with
the extended Dynkin diagram $\varGamma$}.
\end{center}

To prove the sufficient condition of regularity, we study the
transforming elements interrelating Coxeter transformations for
different orientations. Here, Theorem \ref{th_transf},
\cite{St82}, plays a key role. Proposition \ref{prop_transf} and
Theorem \ref{th_transf} yield the following:

1) Let $\varDelta{'}, \varDelta{''}$ be two arbitrary orientations
of the graph $\varGamma$ that differ by the direction of $k$
edges. Consider a chain of orientations, in which every two
adjacent orientations differ by the direction of one edge:
\begin{equation*}
     \varDelta{'} = \varDelta_0, \varDelta_1, \varDelta_2, \dots,
     \varDelta_{k-1}, \varDelta_k = \varDelta{''}.
\end{equation*}
   Then, in the
   Weyl group, there exist elements $P_i$ and $S_i$, where $i = 1,2,...,k$,
such that
\begin{equation*}
  \begin{split}
    & {\bf C}_{\varDelta_0} = P_1{S}_1,  \\
    & {\bf C}_{\varDelta_1} = S_1{P}_1 = P_2{S}_2, \\
    & \dotfill
    \\
    & {\bf C}_{\varDelta_{k-1}} = S_{k-1}P_{k-1} = P_k{S}_k,   \\
    & {\bf C}_{\varDelta_k} = S_k{P}_k.
  \end{split}
\end{equation*}

2) ${T}^{-1}{\bf C}_{\varDelta{'}}T =
 {\bf C}_{\varDelta{''}}$ for the following $k+1$ transforming elements
$T:=T_i$:
\begin{equation*}
 \begin{split}
   & T_1 = P_1{P}_2{P}_3...P_{k-2}{P}_{k-1}P_k, \\
   & T_2 = P_1{P}_2{P}_3...P_{k-2}{P}_{k-1}S^{-1}_k, \\
   & T_3 = P_1{P}_2{P}_3...P_{k-2}{S}^{-1}_{k-1}S^{-1}_k, \\
   &    \dotfill
  \\
   & T_{k-1} =
            P_1{P}_2{S}^{-1}_3...{S}^{-1}_{k-2}S^{-1}_{k-1}S^{-1}_k,    \\
   & T_k =
            P_1{S}^{-1}_2{S}^{-1}_3...{S}^{-1}_{k-2}{S}^{-1}_{k-1}{S}^{-1}_k, \\
   & T_{k+1} =
      {S}^{-1}_1{S}^{-1}_2{S}^{-1}_3...{S}^{-1}_{k-2}{S}^{-1}_{k-1}{S}^{-1}_k.
 \end{split}
\end{equation*}
 In addition, for each reflection $\sigma_\alpha$, there exists a $T_i$ whose
 decomposition does not contain this reflection.

3) The following relation holds:
\begin{equation*}
    T_p{T}^{-1}_q = {\bf C}^{q-p}_{\varDelta{'}}.
\end{equation*}

\medskip

For the graph $\varGamma$ with indefinite Tits form $\mathcal{B}$,
we prove the following
 {\it necessary condition of regularity},
see Theorem \ref{regul_indef_Tits} \cite{SuSt75}, \cite{SuSt78}:

 \index{conjugate vector}
 If $z$ is the regular vector in the
orientation $\varDelta{'}$, then
\begin{equation*}
      <Tz, \tilde{z}^m_1> \hspace{3mm} \leq \hspace{3mm} 0,
      \hspace{7mm}
      <Tz, \tilde{z}^m_2> \hspace{3mm} \geq \hspace{3mm} 0,
\end{equation*}
where $\tilde{z}^m_1$ and $\tilde{z}^m_2$ are the {\it dual
vectors} (Definition \ref{conjug_vector}) to the vectors $z^m_1$
and $z^m_2$ corresponding to the maximal eigenvalue ${\varphi}^m =
{\varphi}^{max}$ of $DD^t$ and $D^tD$, respectively, see
\S\ref{necessary_cond}.

Similar results were obtained by Y.~Zhang in
\cite[Prop.1.5]{Zh89}, and by J.A.~de la Pe\~{n}a and M.~Takane in
\cite[Th.2.3]{PT90}.

For an application of this necessary condition to the star graph,
see \S\ref{star_example}.

\chapter{\sc\bf Preliminaries}
\label{chapter_prelim}

\epigraph{
 ...Having computed the $m$'s several years earlier ..., I recognized
them\footnotemark[1]  in the Poincar\'{e} polynomials  while
listening to Chevalley's address at the International Congress in
1950. I am grateful to A.~J.~Coleman for drawing my attention to
the relevant work of Racah, which helps to explain the
``coincidence''; also, to J.~S.~Frame for many helpful
suggestions...} {H.~S.~M.~Coxeter,  \cite[p.765]{Cox51}, 1951}

\footnotetext[1]{The numbers $m_1,\dots,m_n$ are said to be {\it
exponents} of the Weyl group. See \S\ref{reflections}.}

\section{The Cartan matrix and the Tits form}
\label{cartan}
 \hspace{1mm}  Let $K$ be a matrix with the following properties:
   \cite{Mo68}, \cite{Kac80}
\begin{align*}
 & (C1) \hspace{3mm}  {k}_{ii} = 2 \text{ for } i=1,..,n, \\
 & (C2)  \hspace{3mm}  -{k}_{ij} \in \mathbb{Z}_+ = \{0,1,2,...\}
      \text{ for } i \neq j, \\
 & (C3)  \hspace{3mm} {k}_{ij} = 0  \hspace{2mm} \text{implies} \hspace{2mm}
       {k}_{ij} = 0 \text{ for } i,j=1,...,n.
\end{align*}
\index{generalized Cartan matrix} \index{symmetrizable matrix}
Such a matrix is called a {\it generalized Cartan matrix}. A
matrix $M$ is said to be {\it symmetrizable} if there exists an
invertible diagonal matrix $T$ with positive integer coefficients
and a symmetric matrix {\bf B} such that $M = T${\bf B}. For
example, any generalized Cartan matrix $K$ whose diagram contains
no cycles is symmetrizable (\cite[\S3]{Mo68}). Any Cartan matrix
whose diagram is a simply-laced diagram (even with cycles) is
symmetrizable because it is symmetric. In particular,
$\widetilde{A}_n$ has a symmetrizable Cartan matrix. In this work
we consider only diagrams without cycles, so the diagrams we
consider  have symmetrizable Cartan matrices.

\index{valued graph} Let $\varGamma_0$ (resp. $\varGamma_1$) be
the set of vertices (resp. edges) of the graph $\varGamma$. A {\it
valued graph} $(\varGamma, d)$ (\cite{DR76}, \cite{AuRS95}) is a
finite set $\varGamma_1$ (of edges) rigged with non-negative
integers $d_{ij}$ for all pairs $i, j\in\partial \varGamma_1$ of
the endpoints of the edges in such a way that $d_{ii} = 0$ and
there exist numbers $f_i\in \mathbb{N}$ satisfying
\begin{equation}
   \label{valued_graph}
     d_{ij}f_j = d_{ji}f_i  \hspace{3mm}
           \text{ for all }\hspace{1mm} i,j \in \in\partial \varGamma_1.
\end{equation}
The rigging of the edges of $\varGamma_1$ is depicted by symbols
\begin{displaymath}
   i\stackrel{(d_{ij}, d_{ji})}{\line(1,0){40}}j
\end{displaymath}
If $d_{ij} = 1 = d_{ji}$, we simply write
\begin{displaymath}
   i \hspace{1mm} \line(1,0){40} \hspace{1mm} j
\end{displaymath}
There is, clearly, a one-to-one correspondence between valued
graphs and symmetrizable Cartan matrices, see \cite{AuRS95}.

A non-symmetric bilinear form $\mathcal{Q}$ defined on
vectors--dimensions of the quiver representation  was
independently introduced by Ringel \cite{Rin76} and Ovsienko and
Roiter \cite{OR77}. The form $\mathcal{Q}$ has the following
property important for the representations of quivers
\cite{Rin76}:
\begin{equation}
    \mathcal{Q}(\dim{X}, \dim{Y}) =
     \dim{\rm Hom}({X}, {Y}) -
      \dim{\rm Ext}^1({X}, {Y}).
\end{equation}
The {\it Tits form} is the quadratic form $\mathcal{B}$ given by the relation
\index{Tits form}
\begin{equation}
    \mathcal{B}(\alpha) = \mathcal{Q}(\alpha, \alpha).
\end{equation}
The corresponding symmetric bilinear form is
\index{Cartan matrix! - multiply-laced case}
\index{Cartan matrix! - simply-laced case}
\begin{equation}
 \label{bilinear}
    (\alpha, \beta) =
      \frac{1}{2}( \mathcal{Q}(\alpha, \beta) +
       \mathcal{Q}(\beta, \alpha) ).
\end{equation}
Let {\bf B} be the matrix of the quadratic form $\mathcal{B}$. In
the simply-laced case,
\begin{equation}
   K = 2{\bf B}
\end{equation}
with the symmetric Cartan matrix $K$. So, in the simply-laced
case, the Tits form is the Cartan-Tits form. In the multiply-laced
case, the symmetrizable matrix $K$ factorizes
\begin{equation}
 \label{symmetrizable_K}
\begin{array}{cc}
     & K = T{\bf B} \vspace{3mm}, \\
     & {k}_{ij} = \displaystyle\frac{2(a_i, b_j)}{(a_i, a_i)}, \hspace{5mm}
       {k}_{ji} = \displaystyle\frac{2(a_i, b_j)}{(b_j, b_j)},
\end{array}
\end{equation}
where $T$ is a diagonal matrix with positive integers
on the diagonal, the vectors $a_i$ and $b_j$ are defined in \S\ref{conjugacy},
see also \S\ref{bicolored_C}.

\section{The Schwartz inequality}
Let $\mathcal{B}$ be the quadratic Tits form associated with an arbitrary
tree graph (simply or multiply laced) and {\bf B} the matrix of $\mathcal{B}$.
Let
\begin{equation}
    \ker{\bf B} = \{x \mid {\bf B}x = 0\}.
\end{equation}
Since $(x, y)  = <x, {\bf B}y> =
<{\bf B}x, y>$, it follows that
\begin{equation}
  \label{inclusion_B}
    \ker{\bf B} = \{x \mid
          (x,y) = 0 \text{ for all } y \in \mathcal{E}_\varGamma \}
        \hspace{2mm}
        \subseteq \hspace{2mm} \{x \mid
           \mathcal{B}(x) = 0 \}.
\end{equation}
We will write $\mathcal{B} > 0$ (resp. $\mathcal{B} \geq 0$) if
$\mathcal{B}$ is positively (resp. nonnegative) definite. Let $x,
y \in \mathbb{R}^n$, where $n$ is the number of vertices in
$\varGamma_0$. \index{Schwartz inequality} Then the following {\it
Schwartz inequality} is true:
\begin{equation}
 \label{ineq_Sch}
  (x, y)^2 \leq \mathcal{B}(x)\mathcal{B}(y).
\end{equation}
To prove (\ref{ineq_Sch}), it suffices to consider the inequality
$(x + \alpha{y}, x + \alpha{y}) \geq 0$ which is true for all
$\alpha \in \mathbb{R}$. Then the discriminant of the polynomial
$$
 (x, x) + 2\alpha(x, y) +
\alpha^2(y, y)
$$
 should be non-positive, whence (\ref{ineq_Sch}).

Over $\mathbb{C}$, there exists $x \not\in \ker{\bf B}$ such that
$\mathcal{B}(x) = 0$. For example, the eigenvectors of the Coxeter
transformation with eigenvalues $\lambda \neq \pm{1}$
satisfy this condition because
\begin{equation*}
   \mathcal{B}(x) = \mathcal{B}(Cx) = \lambda^2\mathcal{B}(x).
\end{equation*}
If $\mathcal{B}(x) = 0$
and $x \in \mathbb{R}^n$, then from (\ref{ineq_Sch}) we
get $(x, y) \leq 0$ and
$(x, -y) \leq 0$, i.e., $(x, y) = 0$
for all $y \in \mathbb{R}^n$. In other words, $x \in \ker{\bf B}$.
Taking (\ref{inclusion_B}) into account we see that if $\mathcal{B} \geq 0$,
then
\begin{equation}
    \ker{\bf B} =  \{x \mid \mathcal{B}(x) = 0 \}.
\end{equation}

\section{The Dynkin diagrams and the extended Dynkin diagrams}
\begin{figure}[h]
\centering
\includegraphics{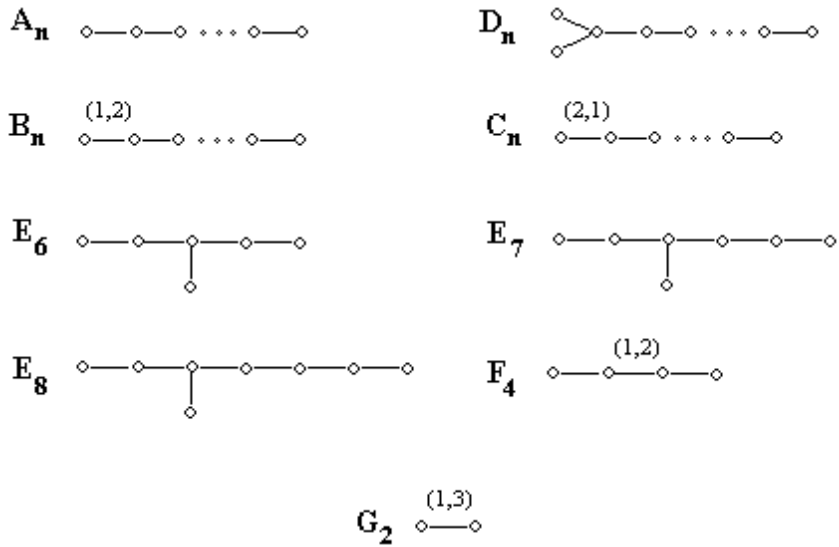}
\caption{\hspace{3mm} The Dynkin diagrams}
\label{dynkin_diag}
\end{figure}

For the proof of the following well-known characterization of the
Dynkin diagrams and the extended Dynkin diagrams, see, e.g.,
\cite{Bo}, \cite{DR76}:

1) $(\varGamma, d)$ is a Dynkin diagram if and only if its
quadratic form is positive definite, see Fig.~\ref{dynkin_diag}.

\begin{figure}[h]
\centering
\includegraphics{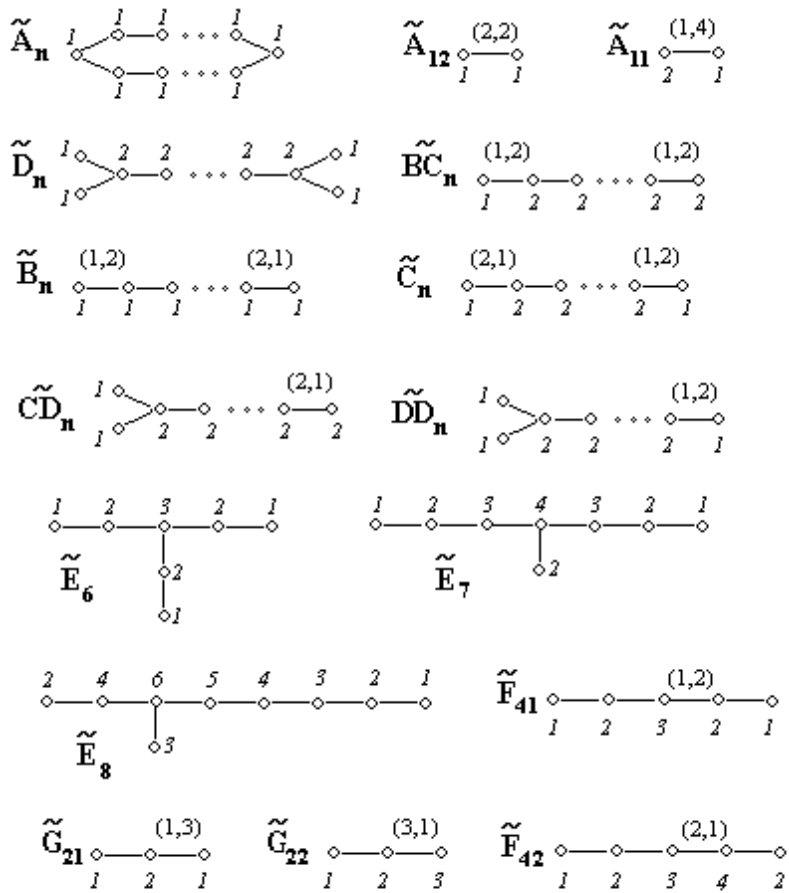}
\caption{\hspace{3mm} The extended Dynkin diagrams}
\label{euclidean_diag}
\end{figure}

2) $(\varGamma, d)$ is an extended Dynkin diagram if and only if
its quadratic form is nonnegative definite, see
Fig.~\ref{euclidean_diag}.


\begin{remark}
 \label{operations}
 {\rm \indent 1)  All quadratic forms $\mathcal{B}$ fall into 3
non-intersecting sets

   (a) $\mathcal{B}$ is positive definite, $\mathcal{B} > 0$
                (the Dynkin diagrams),

   (b) $\mathcal{B}$ is nonnegative definite, $\mathcal{B} \geq 0$
                (the extended Dynkin diagrams),

   (c) $\mathcal{B}$ is indefinite.
\medskip

2) Consider two operations:

\index{operation \lq\lq Add Edge"}
  2.1) \lq\lq Add Edge":  add a vertex and connect it with $\varGamma$
by only one edge.  We denote the new graph
$\stackrel{\wedge}{\varGamma}$. \index{operation \lq\lq Remove
Vertex"}

2.2) \lq\lq Remove Vertex": remove a vertex (the new graph may
contain more than one component).
 We denote the new graph $\stackrel{\vee}{\varGamma}$.

3) The set (b) ($\mathcal{B} \geq 0$) is not stable under these operations:
 \begin{displaymath}
  \begin{array}{cc}
   \wedge  :  \hspace{1mm} \mathcal{B} \geq 0 \hspace{3mm}
     \Longrightarrow  & \hspace{3mm} \mathcal{B} \text{ is indefinite}, \\
   \vee  :  \hspace{1mm} \mathcal{B} \geq 0 \hspace{3mm}
     \Longrightarrow  & \hspace{3mm} \mathcal{B} > 0.
  \end{array}
 \end{displaymath}

The set (a) is invariant under $\vee$ and sometimes is invariant
under $\wedge$. The set (c) is invariant under $\wedge$
and sometimes is invariant under $\vee$.

4) If the graph $\varGamma$ with indefinite form $\mathcal{B}$ is
obtained from any Dynkin diagram by adding an edge, then the same
graph $\varGamma$ can be obtained by adding an edge (or maybe
several edges) to an extended Dynkin diagram.}
\end{remark}

\section{The real and imaginary roots}
  \label{roots}
 \index{imaginary roots}
 \index{real roots}
 \index{root system}
 \index{Lie algebra}
 Recall now the definitions of {\it imaginary}
and {\it real} roots in the infinite root system associated with
infinite dimensional Kac-Moody Lie algebras. We mostly follow
V.~Kac's definitions \cite{Kac80}, \cite{Kac82}, \cite{Kac93}.

We consider the vector space $\mathcal{E}_{\varGamma}$ over
$\mathbb{Q}$; set
$$
     \dim{\mathcal{E}_\varGamma} = |\varGamma_0|.
$$
Let
$$
 \alpha_i = \{0,0,...0,\stackrel{i}{1},0,...0,0\} \in \mathcal{E}_{\varGamma}
$$
be the basis vector corresponding to the vertex ${v}_i \in
\varGamma_0$. The space $\mathcal{E}_{\varGamma}$ is spanned by
the vectors $\{\alpha_i \mid i \in \varGamma_0\}$; the vectors
$\alpha_1,...,\alpha_n$ form the basis in
$\mathcal{E}_{\varGamma}$. Let
$$
   \mathcal{E}_+ =
   \{ \alpha = \sum{k_i\alpha_i} \in \mathcal{E}_\varGamma  \mid
    k_i \in \mathbb{Z}_+, \; \sum{k_i} > 0 \}
$$
be the set of all non-zero elements in $\mathcal{E}_{\varGamma}$
with non-negative integer coordinates in the basis
$\{\alpha_1,...,\alpha_n\}$.

Define the linear functions $\phi_1,~\dots,~\phi_n$ on
$\mathcal{E}_{\varGamma}$ by means of the elements of the Cartan
matrix (\ref{symmetrizable_K}):
\begin{equation}
          \phi_i(\alpha_j) = k_{ij} \hspace{1mm}.
\end{equation}

\index{positive root system} The {\it positive root system}
$\varDelta_+$ associated to the Cartan matrix $K$ is a subset in
$\mathcal{E}_+$ defined by the following properties (R1)--(R3):
\vspace{3mm}

  (R1)  $\alpha_i \in \varDelta_+$ and
       $2\alpha_i \notin \varDelta_+$ for
         $i = 1,..., n$.

  (R2) If $\alpha \in \varDelta_+$ and
         $\alpha \neq  \alpha_i$,  then
       $\alpha + k\alpha_i \in \varDelta_+$ for
       $k \in \mathbb{Z}$ if and only if $-p \leq k \leq q$,
       where $p$ and $q$ are some non-negative integers such that
        $p - q = \phi_i(\alpha)$.

  (R3)  Any $\alpha \in \varDelta_+$ has a connected support.

\medskip

We define endomorphisms $\sigma_1,~...,~\sigma_n$ of
$\mathcal{E}_{\varGamma}$ by the formula
\begin{equation}
          \sigma_i(x) = x - \phi_i(x)\alpha_i.
\end{equation}

Each endomorphism $\sigma_i$ is the reflection in the hyperplane $\phi_i = 0$ such that \\
$\sigma_i(\alpha_i) = -\alpha_i$. These reflections satisfy the
following relations:
$$
    \sigma_i^2 = 1,  \hspace{3mm} (\sigma_i\sigma_j)^{n_{ij}} = 1,
$$
where $n_{ij}$ are given in the Table \ref{exponents}.
\begin{table}[h]
  \centering
  \renewcommand{\arraystretch}{1.9}
  \vspace{2mm}
  \caption{\hspace{3mm}The exponents $n_{ij}$}
  \begin{tabular} {||c|c|c|c|c|c||}
  \hline \hline
      \quad ${k}_{ij}{k}_{ji}$ \quad  & \quad 0 \quad & \quad 1 \quad &
      \quad 2 \quad & \quad 3 \quad & \quad $\geq 4$  \quad  \\
  \hline \hline
      $n_{ij}$  &  2 & 3 & 4 & 6 & $\infty$ \\
  \hline \hline
  \end{tabular}
  \label{exponents}
\end{table}

\index{Weyl group} \index{simple roots} The group $W$ generated by
the reflections $\sigma_1,~...,~\sigma_n$ is called the {\it Weyl
group}. The vectors $\alpha_1,~...,~\alpha_n$ are called {\it
simple roots}; we denote by $\varPi$ the set of all simple roots.
Let $W(\varPi)$ be the orbit of $\varPi$ under the $W$-action.

\index{root system} \index{fundamental set}
Set
$$
M= \{\alpha \in \mathcal{E}_+ \mid
    \phi_i(\alpha) \leq 0 \text{ for } i=1,...,n,
    \text{ and } \alpha \text{ has a connected support } \}.
$$
The set $M$ is called the {\it fundamental set}. Let $W(M)$ be the
orbit of $M$ under the $W$-action. We set
$$
  \varDelta^{re}_+ = \bigcup\limits_{w \in W}(w(\varPi)\cap\mathcal{E}_+),
  \hspace{15mm}
  \varDelta^{im}_+ = \bigcup\limits_{w \in W}(w(M)).
$$
The elements of the set $\varDelta^{re}_+$ are called {\it real
roots} and the elements of the set $\varDelta^{im}_+$ are called
{\it imaginary roots}. By \cite{Kac80}, the {\it system of
positive roots} $\varDelta_+$ is the disjoint union of the sets
$\varDelta^{re}_+$ and $\varDelta^{im}_+$:
$$
     \varDelta^{re}_+ =  \varDelta^{re}_+ \cup \varDelta^{im}_+. 
$$

If the Tits form $\mathcal{B}$  (or, which is the same, the Cartan
matrix $K$) is positive definite, then the root system is finite,
it corresponds to a simple finite dimensional Lie algebra. In this
case, the root system consists of real roots.

If the Tits form $\mathcal{B}$ is nonnegative definite, $\mathcal{B} \geq 0$,
we have an infinite root system whose
imaginary root system is one-dimensional:
$$
      \varDelta^{im}_+ = \{\delta, 2\delta, 3\delta,...\}, \hspace{4mm}
      {\rm where} \hspace{1mm} \delta = \sum\limits_{i}k_i\alpha_i,
$$
the coefficients $k_i$ being the labels of the vertices from
Fig.~\ref{euclidean_diag}.

The numerical labels at the vertices are the coefficients of the
imaginary root which coincides with the fixed point ${z}^1$ of the
Coxeter transformation, see \S\ref{explicitly}.

\index{fixed points} \index{nil-roots} The elements $k\delta$ ($k
\in \mathbb{N}$) are called {\it nil-roots}. Every nil-root is a
fixed point for the Weyl group.

\section{The hyperbolic Dynkin diagrams and hyperbolic Cartan matrices}
  \label{hyperbolic}
\index{hyperbolic Dynkin diagram} \index{hyperbolic Cartan matrix}
\index{hyperbolic Weyl group} \index{compact hyperbolic Weyl
group} \index{strictly hyperbolic Dynkin diagram} A connected
graph $\varGamma$ with indefinite Tits form is said to be a {\it
hyperbolic Dynkin diagram} (resp. {\it strictly hyperbolic Dynkin
diagram}) if every subgraph $\varGamma{'} \subset \varGamma$ is a
Dynkin diagram or an extended Dynkin diagram (resp. Dynkin
diagram). The corresponding Cartan matrix $K$ is said to be the
{\it hyperbolic Cartan matrix} (resp. {\it strictly hyperbolic
Cartan matrix}), see \cite[exs. of \S{4.10}]{Kac93}. The
corresponding Weyl group is said to be a {\it hyperbolic Weyl
group} (resp. {\it compact hyperbolic Weyl group}), see
\cite[Ch.5, exs. of \S{4}]{Bo}.

\begin{figure}[h]
\centering
\includegraphics{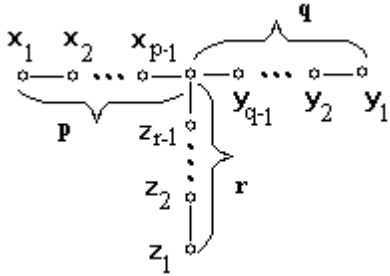}
\caption{\hspace{3mm} The diagram $T_{p,q,r}$}
\label{T_pqr_diagram}
\end{figure}

Consider the Dynkin diagram $T_{p,q,r}$ depicted on
Fig.~\ref{T_pqr_diagram}. By eq.(\ref{matrix_D}) the quadratic
Tits form $\mathcal{B}$ can be expressed in the following form
\begin{equation}
 \begin{split}
    \frac{1}{2}\mathcal{B}(z) =
    & x_1^2 + x_2^2  + \dots + x_{p-1}^2
      - x_1x_2  - \dots - x_{p-2}x_{p-1} - x_{p-1}u_0 + \\
    & y_1^2 + y_2^2 + \dots + y_{q-1}^2
      - y_1y_2  - \dots - y_{q-2}y_{q-1} - y_{q-1}u_0 + \\
    & z_1^2 + z_2^2 + \dots + z_{r-1}^2
      - z_1z_2  - \dots - z_{r-2}z_{r-1} - z_{r-1}u_0 +  u_0^2,
  \end{split}
\end{equation}
where $u_0$ is the coordinate of the vector $z$ corresponding to
the branch point of the diagram. Let
\begin{equation}
 \begin{split}
  V(x_1,...,& x_{p-1},u_0) = \\
  & x_1^2 + x_2^2  + \dots + x_{p-1}^2
      - x_1x_2  - \dots - x_{p-2}x_{p-1} - x_{p-1}u_0,  \\
  V(y_1,...,& y_{q-1},u_0) = \\
  & y_1^2 + y_2^2 + \dots + y_{q-1}^2
      - y_1y_2  - \dots - y_{q-2}y_{q-1} - y_{q-1}u_0,  \\
  V(z_1,..., & z_{r-1}, u_0) = \\
  & z_1^2 + z_2^2 + \dots + z_{r-1}^2
      - z_1z_2  - \dots - z_{r-2}z_{r-1} - z_{r-1}u_0.  \\
  \end{split}
\end{equation}
Then
\begin{equation}
    \frac{1}{2}\mathcal{B}(z) =
       V(x_1,...,x_{p-1}) + V(y_1,....y_{q-1}) + V(z_1,....z_{r-1}) +
       u_0^2.
\end{equation}
It is easy to check that
\begin{equation}
 \label{sums_i_i1}
 \begin{split}
    2V(x_1,...,x_{p-1}, u_0) = &
       \sum\limits_{i=1}^{p-1}\frac{i+1}{i}(x_i - \frac{i}{i+1}x_{i+1})^2
       - \frac{p-1}{p}u^2_0, \\
    2V(y_1,...,y_{p-1},  u_0) = &
       \sum\limits_{i=1}^{q-1}\frac{i+1}{i}(y_i - \frac{i}{i+1}y_{i+1})^2
       - \frac{q-1}{q}u^2_0, \\
    2V(z_1,...,z_{p-1},  u_0) = &
       \sum\limits_{i=1}^{r-1}\frac{i+1}{i}(z_i - \frac{i}{i+1}z_{i+1})^2
       - \frac{r-1}{r}u^2_0.
  \end{split}
\end{equation}
Denote by $U(x)$, $U(y)$, $U(z)$ the corresponding sums in (\ref{sums_i_i1}):
\begin{equation}
 \label{U_xyz}
 \begin{split}
    2V(x_1,...,x_{p-1},u_0) = U(x) - \frac{p-1}{p}u^2_0, \\
    2V(y_1,...,y_{p-1},u_0) = U(y) - \frac{q-1}{q}u^2_0, \\
    2V(z_1,...,z_{p-1},u_0) = U(z) - \frac{r-1}{r}u^2_0.
  \end{split}
\end{equation}
Thus,
\begin{equation}
    \mathcal{B}(z) = U(x) + U(y) + U(z)
      - \frac{p-1}{p}u^2_0 - \frac{q-1}{q}u^2_0 - \frac{r-1}{r}u^2_0
        + 2u_0^2,
\end{equation}
or
\begin{equation}
    \mathcal{B}(z) = U(x) + U(y) + U(z)
      + \left(\frac{1}{p} + \frac{1}{q} + \frac{1}{r} - 1\right)u^2_0.
\end{equation}
 Set
\begin{equation}
   \mu = \frac{1}{p} + \frac{1}{q} + \frac{1}{r};
\end{equation}
then we have
\begin{equation}
  \label{B_by_mu}
    \mathcal{B}(z) = U(x) + U(y) + U(z) + (\mu - 1)u^2_0.
\end{equation}

We see that $\mathcal{B} > 0$ for $\mu > 1$, i.e., the triples $p,
q, r$ are as follows:

1) $(1,q,r)$ for any  $q,r \in \mathbb{Z}$, i.e., the diagrams
$A_n$,

2) $(2,2,r)$, the diagrams ${D}_n$,

3) $(2,3,r)$, where $r = 3,4,5$, i.e., the diagrams ${E}_6, {E}_7,
{E}_8$.

\noindent Similarly, $\mathcal{B} \geq 0$ for $\mu = 1$, i.e., for
 the following triples $p, q, r$:

1) $(3,3,3)$, the diagram $\tilde{E}_6$,

2) $(2,4,4)$, the diagram $\tilde{E}_7$,

3) $(2,3,6)$, the diagram $\tilde{E}_8$.

\medskip

The Tits form $\mathcal{B}$ is indefinite for $\mu < 1$. Since
$U(x)$, $U(y)$, $U(z)$ in (\ref{B_by_mu}) are positive, we see
that, in this case, the signature of $\mathcal{B}$ is equal to
$(n-1, 1)$, see \cite[exs.4.2]{Kac93}. It is easy to check that
the triples $p, q, r$ corresponding to hyperbolic graphs
$T_{p,q,r}$ are only the following ones:

1) $(2,3,7)$, the diagram ${E}_{10}$, $1 - \mu =
\displaystyle\frac{1}{42}$,

2) $(2,4,5)$, $1 - \mu = \displaystyle\frac{1}{20}$, \vspace{2mm}

3) $(3,3,4)$, $1 - \mu =  \displaystyle\frac{1}{12}$.

\section{Orientations and the associated Coxeter transformations}
\label{orientation} Here, we follow the definitions of
Bernstein-Gelfand-Ponomarev \cite{BGP73}. Let us consider the
graph $\varGamma$ endowed with an orientation $\varDelta$.

\index{sink-admissible vertex} \index{source-admissible vertex}
The vertex $v_i$ is called {\it sink-admissible} (resp. {\it
source-admissible}) in the orientation $\varDelta$ if every arrow
containing $v_i$ ends in (resp. starts from) this vertex. The
reflection $\sigma_i$ can be applied only to vertices which are
either sink-admissible or source-admissible. The reflection
$\sigma_i$ acts on the orientation $\varDelta$ by reversing all
arrows containing the vertex $v_i$.

\index{sink-admissible sequence}
\index{source-admissible sequence}
Consider now a sequence of vertices and the corresponding reflections.
 A sequence of vertices
$$
  \{v_{i_n}, v_{i_{n-1}},~...,~v_{i_3}, v_{i_2}, v_{i_1}\}
$$
is called {\it sink-admissible}, if the vertex $v_{i_1}$ is
sink-admissible in the orientation $\varDelta$, the vertex
$v_{i_2}$ is sink-admissible in the orientation
$\sigma_{i_1}(\varDelta)$, the vertex $v_{i_3}$ is sink-admissible
in the orientation $\sigma_{i_2}\sigma_{i_1}(\varDelta)$, and so
on. {\it Source-admissible sequences} are similarly defined.

\index{fully sink-admissible sequence}
\index{fully source-admissible sequence}
A sink-admissible (resp. source-admissible) sequence
$$
  \mathcal{S} = \{v_{i_n}, v_{i_{n-1}},~...,~v_{i_3}, v_{i_2}, v_{i_1}\}
$$
 is called {\it fully sink-admissible} (resp. {\it fully
source-admissible}) if $\mathcal{S}$ contains every vertex $v \in
\varGamma_0$ exactly once. Evidently, the inverse sequence
$\mathcal{S}^{-1}$ of a sink-admissible sequence $\mathcal{S}$ is
source-admissible and vice versa.

Every tree has a fully sink-admissible sequence $\mathcal{S}$. To
every sink-admissible sequence $\mathcal{S}$, we assign the
Coxeter transformation depending on the order of vertices in
$\mathcal{S}$:
$$
{\bf C} = \sigma_{i_n}\sigma_{i_{n-1}}...\sigma_{i_2}\sigma_{i_1}.
$$

For every orientation $\varDelta$ of the tree, every fully
sink-admissible sequence gives rise to the same Coxeter
transformation ${\bf C}_\varDelta$, and every fully
source-admissible sequence gives rise to ${\bf C}^{-1}_\varDelta$.

\section{The Poincar\'{e} series}
  \label{sect_poincare}
\subsection{Graded algebras, symmetric algebras, algebras of invariants}
 \label{graded_algebras}
 \index{graded $k$-algebra}
Let $k$ be a field. We define a {\it graded $k$-algebra} to be a
finitely generated $k$-algebra $A$ (associative, commutative, and
with identity), together with a direct sum decomposition (as
vector space)
\begin{equation}
    A = A_0 \oplus A_1 \oplus A_2 \oplus \dots,
  \end{equation}
such that $A_0 = k$ and $A_iA_j \subset A_{i+j}$. The component $A_n$
is called the {\it $n$th homogeneous part of $A$} and the element $x \in A_n$
is said to be {\it homogeneous of degree $n$}, notation: $\deg{x} = n$.

 \index{graded $A$-module}
 \index{finitely generated module}

We define a {\it graded $A$-module} to be a finitely generated
$A$-module, together with a direct sum decomposition
\begin{equation}
    M = M_0 \oplus M_1 \oplus M_2 \oplus \dots,
  \end{equation}
such that $A_iM_j \subset M_{i+j}$.

\index{Poincar\'{e} series of a graded algebra}
\index{Poincar\'{e} series of a graded module}
The {\it Poincar\'{e} series} of a graded algebra
$A = \mathop{\oplus}\limits_{n=1}^{\infty}A_i$ is the formal series
\begin{equation}
 \label{poincare_alg}
   P(A,t) = \sum\limits_{n=1}^{\infty}(\dim{A_n})t^n.
\end{equation}

The {\it Poincar\'{e} series} of a graded $A$-module $M =
\mathop{\oplus}\limits_{n=1}^{\infty}M_i$ is the formal series
\begin{equation}
 \label{poincare_mod}
   P(M,t) = \sum\limits_{n=1}^{\infty}(\dim{M_n})t^n,
\end{equation}
see \cite{Sp77}, \cite{PV94}.

 \index{theorem! - Hilbert-Serre}
 \index{finitely generated module}
\begin{theorem} {\em (Hilbert, Serre, see \cite[p.117, Th.11.1]{AtMa69})}
The Poincar\'{e} series $P(M,t)$ of a finitely generated graded
$A$-module is a rational function in $t$ of the form
\begin{equation}
 \label{hilbert_serre}
   P(M,t) = \frac{f(t)}{\prod\limits_{i=1}^{s}(1 - t^{k_i})},
   \text{ where } f(t) \in \mathbb{Z}[t].
\end{equation}
\end{theorem}

In what follows in this section, any algebraically closed field
$k$ can be considered instead of $\mathbb{C}$. Set
\begin{equation}
  \label{R_algebra}
    R =  \mathbb{C}[x_1,\dots,x_n].
\end{equation}
The set $R_d$ of homogeneous polynomials of degree $d$ is a finite
dimensional subspace of $R$, and $R_0 = \mathbb{C}$. Moreover,
$R_dR_e \subset R_{d+e}$, and $R$ is a graded $\mathbb{C}$-algebra
with a direct sum decomposition (as a vector space)
  \begin{equation}
   \label{R_decomp}
    R = R_0 \oplus R_1 \oplus R_2 \oplus \dots
  \end{equation}

 \index{symmetric algebra}
Let $f_1,\dots,f_n \in V^*$ be the linear forms defined by
$f_i(x_j) = \delta_{ij}$.  Then the $f_i$, where $i=1,\dots,n$,
generate the {\it symmetric algebra} ${\rm Sym}(V^*)$ isomorphic
to $R$.

For any $a \in GL(V)$ and $f \in R$, define  $af \in R$ by the
rule
\begin{equation}
    (af)(v) =  f(a^{-1}v)\;\text{ for any $v\in V$};
\end{equation}
then, for any $a, b\in G$ and $f\in R$, we have
$$
  a(bf) = bf(a^{-1}v) = f(b^{-1}a^{-1}v) =
  f((ab)^{-1}v) = ((ab)f)(v),
$$
and $aR_d = R_d$.

Let $G$ be a subgroup in $GL(V)$. We say that $f \in R$ is a
$G$-{\it invariant} if $af = f$ for all $a \in G$.  The
$G$-invariant polynomial functions form a subalgebra $R^G$ of $R$,
which is a graded subalgebra, i.e., \index{$G$-invariant}
\index{algebra of invariants}
\begin{equation}
    R^G = \bigoplus R^G\cap{R_i} . 
\end{equation}
The algebra $R^G$ is said to be the {\it algebra of invariants} of the group
$G$.
If $G$ is finite, then
\index{theorem! - Moilen}
\begin{equation}
 \label{moilen}
   P(R^G,t) = \frac{1}{|G|}\frac{1}{\sum\limits_{g \in G}\det(1-gt)}
\hspace{1mm}.
\end{equation}
Eq. (\ref{moilen}) is a classical theorem of Moilen (1897), see
\cite[\S3.11]{PV94}, or \cite[Ch.5, \S5.3]{Bo}.

\subsection{Invariants of finite groups generated by reflections}
  \label{reflections}
  \index{pseudo-reflection}
  \index{reflection}

Let $G$ be a finite subgroup of $GL(V)$ and $g \in G$. Then $g$ is
called a {\it pseudo-reflection} if precisely one eigenvalue of
$g$ is not equal to 1. Any pseudo-reflection with determinant $-1$
is called a {\it reflection}.

 \index{theorem! - Sheppard-Todd-Chevalley-Serre}
\begin{theorem} {\em STCS (Sheppard and Todd, Chevalley, Serre)}
\label{STCS_theorem} Let $G$ be a finite subgroup of $GL(V)$.
There exist $n=\dim V$ algebraically independent homogeneous
invariants $\theta_1,\dots,\theta_n$ such that
\begin{equation*}
   R^G = \mathbb{C}[\theta_1,\dots,\theta_n]
\end{equation*}
if and only if $G$ is generated by pseudo-reflections.
\end{theorem}

For references, see \cite[Th.4.1]{Stn79}, \cite[Th.A,
p.778]{Ch55}.

Sheppard and Todd \cite{ShT54} explicitly determined all finite
subgroups of $GL(V)$ generated by pseudo-reflections and verified
the sufficient condition of Theorem STCS. Chevalley \cite{Ch55}
found the classification-free proof of this sufficient condition
for a particular case where $G$ is generated by reflections. Serre
observed that Chevalley's proof is also valid for  groups
generated by pseudo-reflections. Sheppard and Todd (\cite{ShT54})
proved the necessary condition of the theorem by a strong
combinatorial method; see also Stanley \cite[p.487]{Stn79}.

 \index{Betti numbers}
 \index{theorem! - Pontrjagin-Brauer-Chevalley}

The coefficients of the Poincar\'{e} series are called the {\it
Betti numbers}. For some bibliographical notes on Betti numbers,
see Remark \ref{betti_numbers}. The Poincar\'{e} polynomial of the
algebra $R^G$ is
\begin{equation}
  \label{betti_numb}
  (1 + t^{2p_1 + 1})(1 + t^{2p_2 + 1})\dots(1 + t^{2p_n + 1}),
\end{equation}
where $p_i + 1$ are the degrees of homogeneous basis elements of
$R^G$, see \cite{Cox51}, and \cite{Ch50}. See Table
\ref{table_Poincare_series_Dynkin} taken from
 \cite[p.781,Tab.4]{Cox51}.
\begin{table} 
  \centering
  \vspace{2mm}
  \caption{\hspace{3mm}The Poincar\'{e} polynomials for the simple compact Lie groups}
  \renewcommand{\arraystretch}{1.3} 
  \begin{tabular} {|c|c|}
  \hline \hline
       Dynkin   &   Poincar\'{e} polynomial \cr
       diagram  &                           \cr
  \hline \hline
     & \cr
       ${A}_n$
     & $(1 + t^3)(1 + t^5)\dots(1 + t^{2n+1})$ \\
  \hline
     & \cr
       ${B}_n$ or ${C}_n$
     & $(1 + t^3)(1 + t^7)\dots(1 + t^{4n-1})$ \\
  \hline
     & \cr
       ${D}_n$
     & $(1 + t^3)(1 + t^7)\dots(1 + t^{4n-5})(1 + t^{2n-1})$ \\
  \hline
     & \cr
       ${E}_6$
     & $(1 + t^3)(1 + t^9)(1 + t^{11})(1 + t^{15})(1 + t^{17})(1 + t^{23})$
\\
  \hline
     & \cr
       ${E}_7$
     & $(1 + t^3)(1 + t^{11})(1 + t^{15})(1 + t^{19})
        (1 + t^{23})(1 + t^{27})(1 + t^{35})$ \\
  \hline
     & \cr
       ${E}_8$
     & $(1 + t^3)(1 + t^{15})(1 + t^{23})(1 + t^{27})
        (1 + t^{35})(1 + t^{39})(1 + t^{47})(1 + t^{59})$ \\
  \hline
     & \cr
       ${F}_4$
     & $(1 + t^3)(1 + t^{11})(1 + t^{15})(1 + t^{23})$ \\
  \hline
     & \cr
       ${G}_2$
     & $(1 + t^3)(1 + t^{11})$ \\
  \hline \hline
\end{tabular}
  \label{table_Poincare_series_Dynkin}
\end{table}

Let $\lambda_1,...,\lambda_n$ be the eigenvalues of a Coxeter
transformation in a finite Weyl group. These eigenvalues can be
given in the form
\begin{equation*}
   \omega^{m_1},...,\omega^{m_n},
\end{equation*}
where $\omega = \exp^{2\pi{i}/h}$ is a primitive root of unity.
The numbers $m_1,...,m_n$ are called the {\it exponents} of the
Weyl group. H.~S.~M.~Coxeter observed that the {\it exponents}
$m_i$ and numbers $p_i$ in (\ref{betti_numb}) coincide (see the
epigraph to this chapter).

\index{theorem! - Coxeter-Chevalley-Coleman-Steinberg}
\index{exponents of the Weyl group}

\index{Lie group}
\begin{theorem} {\em CCKS (Coxeter, Chevalley, Coleman, Steinberg)}
Let
$$
 u_1,\dots,u_n
$$
  be homogeneous elements generating the algebra of invariants
$R^G$, where $G$ is the Weyl group corresponding to a simple
compact Lie group. Let $m_i + 1=\deg u_i$, where $i=1,2,...n$.
Then the exponents of the group $G$ are
$$
m_1,\dots,m_n.
$$
\end{theorem}

For more details, see \cite[Ch.5, \S6.2, Prop.3]{Bo} and
historical notes in \cite{Bo}, and \cite{Ch50}, \cite{Cox51},
\cite{Col58}, \cite{Stb85}.

\index{Betti numbers}
\begin{remark} {On the Betti numbers.}
\label{betti_numbers} {\rm
 Here, we mostly follow the works of C.~Chevalley
\cite{Ch50}, and A.~J.~Coleman \cite{Col58} and review of
C.~A.~Weibel \cite{Weib}.

Let $X$ be a topological space, and $H_k(X, \mathbb{Z})$ be the
$k$-th homology group of $X$. If $H_k(X, \mathbb{Z})$ is finitely
generated, then its rank is called the {\it $k$-th Betti number}
of $X$.

\index{Lie group}
\index{theorem! - Hopf}

L.~Pontrjagin \cite{Pon35} computed the homology of the four
classical Lie groups by means of combinatorial arguments;
R.~Brauer \cite{Br35} used  de Rham homology. H.~Hopf \cite{Ho41}
and H.~Samelson \cite{Sam41} showed that, for a compact Lie group
$\mathfrak{G}$, the Poincar\'{e} series $\sum{B_k}t^k$ of
$\mathfrak{G}$, (where $B_k$ is the $k$-th Betti number) is of the
form
\begin{equation}
 \label{hopf_result}
   P(t) = \prod\limits_{i=1}^{n}(1 + t^{k_i}).
\end{equation}
The relation (\ref{hopf_result}) is called Hopf's theorem
\cite{Ch50}. \footnote{As C.~A.~Weibel \cite[p.6]{Weib} writes:
\lq\lq ...today we would say that Hopf's result amounted to an
early classification of finite dimensional graded \lq\lq Hopf
algebras" over $\mathbb{Q}$" .}

 For the simple compact Lie groups, the
Poincar\'{e} series $P(t)$ are given in Table
\ref{table_Poincare_series_Dynkin}.

 J.~L.~Koszul \cite{Kosz50} proved Hopf's theorem
in a purely algebraic manner. C.~Chevalley \cite{Ch50} proved that
$k_i = 2p_i + 1$ in (\ref{hopf_result}), where $p_i + 1$ is the
degree of a minimal homogeneous invariant of the group
$\mathfrak{G}$. A.~Borel and C.~Chevalley \cite{BC55} have
simplified the calculations for the exceptional Lie groups. }
\end{remark}

%% file: 3jordan.tex

\chapter{\sc\bf The Jordan normal form of the Coxeter transformation}
 \label{jordan}

\setlength{\epigraphwidth}{75mm}

\epigraph{It turned out that most of the classical concepts of the
Killing-Cartan-Weyl theory can be carried over to the entire class
of Kac-Moody algebras, such as the Cartan subalgebra, the root
system, the Weyl group, etc. ... I shall only point out that
$\mathfrak{g}'(K)$ \footnotemark[1] does not always possess a
nonzero invariant bilinear form. This is the case if and only if
the matrix $K$ is {\it symmetrizable} ...} {V.~Kac,
\cite[p.XI]{Kac93}, 1993}

 \index{Cartan subalgebra}
\footnotetext[1]{Here, $\mathfrak{g}'(K)$ is the subalgebra
 $[\mathfrak{g}(K), \mathfrak{g}(K)]$ of
the Kac-Moody algebra $\mathfrak{g}(K)$ associated with the
generalized Cartan matrix $K$. One has
\begin{equation*}
  \mathfrak{g}(K) = \mathfrak{g}'(K) + \mathfrak{h},
\end{equation*}
where $\mathfrak{h}$ is the Cartan subalgebra, $\mathfrak{g}(K) =
\mathfrak{g}'(K)$ if and only if $\det{K} \neq 0$,
 \cite[\S1.3 and Th.2.2]{Kac93}.}

\section{The Cartan matrix and the Coxeter transformation}
In this subsection a graph $\varGamma$ and a partition $S = S_1
\coprod S_2$ of its vertices are fixed.

\subsection{A bicolored partition and a
bipartite graph} \label{bicolor_part}

 \index{bicolored partition}
 \index{bipartite graph}
 \index{bicolored orientation}

A partition $S = S_1 \coprod S_2$ of the vertices of the graph
$\varGamma$ is said to be {\it bicolored} if all edges of
$\varGamma$ lead from $S_1$ to $S_2$. A {\it bicolored partition}
exists if and only if all cycles in $\varGamma$ are of even order.
The graph $\varGamma$ admitting a bicolored partition is said to
be {\it bipartite} \cite{McM02}. An orientation $\varDelta$ is
said to be {\it bicolored}, if there is the corresponding
sink-admissible sequence
$$
   \{v_1,v_2,~...,~v_m,v_{m+1},v_{m+2},~...~v_{m+k}\}
$$
of vertices in this orientation $\varDelta$, such that
subsequences
$$
   S_1 = \{v_1, v_2,~...,~v_m\} \text{ and }
   S_2 = \{v_{m+1},v_{m+2},~...,~v_{m+k}\}
$$
form a bicolored partition. Any two generators $g'$ and $g''$ of
the Coxeter group $W(S)$ (contained in the same subpartition)
commute and the subgroups $W(S_1)$ and $W(S_2)$ are abelian. So,
the products $w_i\in W(S_i)$ for $i = 1,2$ of generators of
$W(S_i)$ are involutions, i.e.,
$$
    w_1^2 = 1,  \hspace{5mm}  w_2^2 = 1 \hspace{1mm}.
$$

For the first time (as far as I know), the technique of bipartite
graphs was used by R.~Steinberg in \cite{Stb59}, where he gave
classification-free proofs for some results of H.~S.~M.~Coxeter
\cite{Cox34} concerning properties of the order of the Coxeter
transformation; see also R.~Carter's paper \cite{Car70}.

\subsection{Conjugacy of Coxeter transformations}
 \label{conjugacy}
\indent All Coxeter transformations are conjugate for any tree or
forest $\varGamma$ \cite[Ch.5, \S6]{Bo}; see also Proposition
\ref{prop_transf} and Remark \ref{conj_class}. The Coxeter
transformations for the graphs with cycles are studied in
\cite{Col89}, \cite{Rin94}, \cite{Shi00}, \cite{BT97}, and in the
works by Menshikh and Subbotin of 1982--1985, see
\S\ref{biblio_notes}. Here, we consider only trees.
\index{bicolored Coxeter transformation}

Let us select one Coxeter transformation most simple to study. It
is one of the two {\it bicolored Coxeter transformations}
corresponding to a bipartite graph:
\begin{equation}
 \label{C_decomp}
  {\bf C} = w_1{w}_2  \hspace{5mm}
  {\rm or} \hspace{5mm}
  {\bf C}^{-1} = w_2{w}_1.
\end{equation}
From now on we assume that $S_1$ contains $m$ elements and $S_2$
contains $k$ elements, we denote by $a_1,~...,~a_m$ (resp.
$b_1,~...,~b_k$) basis vectors corresponding to vertices
$v_1,...,v_m$ of ${S}_1$ (resp. vertices $v_{m+1}...,v_{m+k}$ of
${S}_2$). We denote by $\mathcal{E}_{\varGamma_a}$ (resp.
$\mathcal{E}_{\varGamma_b}$) the vector space generated by the
$a_i, \text{ where } i=1,...,m$ (resp. by the $b_i, \text{ where }
i=1,...,k$). So,
$$
 \dim \mathcal{E}_{\varGamma_a} = m,\hspace{5mm}
 \dim \mathcal{E}_{\varGamma_b} = k \hspace{1mm}.
$$

\subsection{The Cartan matrix and the bicolored Coxeter transformation}
 \label{bicolored_C}
 The Cartan matrix and the bicolored Coxeter transformation
are constructed from the same blocks.
More exactly, the matrix {\bf B} and involutions
$w_i$,where $i = 1,2$,
are constructed from the same blocks.

In the simply-laced case
(i.e., for the symmetric Cartan matrix), we have
\begin{equation}
  \label{symmetric_B}
  \begin{split}
& K = 2{\bf B}, \text{\hspace{3mm}where\hspace{3mm}}
{\bf B} = \left (
\begin{array}{cc}
    I_m & D   \\
    D^t & I_k
\end{array}
\right ),  \\
& w_1 = \left (
\begin{array}{cc}
    -I_m & -2D   \\
    0 & I_k
\end{array}
\right ), \hspace{5mm}
w_2 = \left (
\begin{array}{cc}
    I_m & 0   \\
    -2D^t & -I_k
\end{array}
\right ),
\end{split}
\end{equation}
where the elements $d_{ij}$ that constitute matrix $D$
are given by the formula
\begin{equation}
  \label{matrix_D}
  d_{ij} = (a_i, b_j) =
  \left \{
   \begin{array}{cc}
     -\displaystyle\frac{1}{2} & \text{ if } |a_i - b_j| = 1 \hspace{1mm},
\vspace{2mm} \\
          0       & \text{ if } |a_i - b_j| > 1 \hspace{1mm}.
   \end{array}
  \right .
\end{equation}
\index{Cartan matrix! - multiply-laced case}
\index{Cartan matrix! - simply-laced case}

In the multiply-laced case (i.e., for the
symmetrizable and non-symmetric Cartan matrix $K$), we have
\begin{equation}
 \begin{split}
  \label{matrix_K}
& K = T{\bf B}, \hspace{3mm}\text{ where }\hspace{3mm}
K = \left (
\begin{array}{cc}
    2I_m & 2D   \\
    2F & 2I_k
\end{array}
\right ), \\
& w_1 = \left (
\begin{array}{cc}
    -I_m & -2D   \\
    0 & I_k
\end{array}
\right ), \hspace{4mm}
w_2 = \left (
\begin{array}{cc}
    I_m & 0   \\
    -2F & -I_k
\end{array}
\right ),
\end{split}
\end{equation}
and where
\begin{equation}
     d_{ij} = \frac{(a_i, b_j)}{(a_i, a_i)},
     \hspace{10mm}
     f_{pq} = \frac{(b_p, a_q)}{(b_p, b_p)}, \vspace{3mm}
\end{equation}
where the $a_i$ and $b_j$ are simple roots in the root systems of
the corresponding to $S_1$ and $S_2$ Kac-Moody Lie algebras
[Kac80]. Let $T = (t_{ij})$ be the diagonal matrix
(\ref{symmetrizable_K}). Then
\begin{equation*}
 (t_{ii}) = \frac{1}{(a_i, a_i)} = \frac{1}{\mathcal{B}(a_i)}.
\end{equation*}
Dividing $T$ and {\bf B} into blocks of size $m\times{m}$ and
$k\times{k}$, we see that
\begin{equation}
  \label{matrix_T}
   T = \left (
  \begin{array}{cc}
    2T_1 & 0 \\
    0   & 2T_2
  \end{array}
    \right ), \hspace{3mm}
   {\bf B} = \left (
  \begin{array}{cc}
    T^{-1}_1 & A \\
    A^t   & T^{-1}_2
  \end{array}
    \right ), \hspace{3mm}
    T_1A = D, \hspace{3mm}
    T_2A^t = F \hspace{1mm}.
\end{equation}

\subsection{Dual graphs and Dual forms}
\index{dual graphs} Every valued graph $\varGamma$ has a dual
graph denoted by $\varGamma^{\vee}$. The dual graph is obtained by
means of transposition $d_{ij} \leftrightarrow d_{ji}$. In other
words, if $K$ is the Cartan matrix for $\varGamma$, then the
Cartan matrix $K^{\vee}$ for $\varGamma^{\vee}$ is $K^t$,
$$
     K^{\vee} =  K^t,
$$
i.e.,
\begin{equation}
  \label{dual_t}
     F^{\vee} = D^t, \hspace{5mm}
     D^{\vee} = F^t\hspace{0.5mm}.
\end{equation}
Therefore,
\begin{equation}
  \label{dual_rel}
     F^{\vee}D^{\vee} = (FD)^t, \hspace{5mm}
     D^{\vee}F^{\vee} = (DF)^t\hspace{0.5mm}.
\end{equation}
For any simply-laced graph, the Cartan matrix is symmetric and $F
= D^t = F^{\vee}$. In this case the graph $\varGamma$ is dual to
itself. Among extended Dynkin diagrams the following pairs of
diagram are dual:
\begin{displaymath}
   (\widetilde{B}_n, \widetilde{C}_n), \hspace{7mm}
   (\widetilde{CD}_n, \widetilde{DC}_n), \hspace{7mm}
   (\widetilde{G}_{21}, \widetilde{G}_{12}), \hspace{7mm}
   (\widetilde{F}_{41}, \widetilde{F}_{42})\hspace{0.5mm}.
\end{displaymath}

Let dual Cartan matrices be factorized by means of the diagonal matrices $T$
and $T^{\vee}$:
\begin{equation*}
   K = T{\bf B}, \hspace{5mm}
   K^{\vee} = T^{\vee}{\bf B}^{\vee}.
\end{equation*}
Then according to \cite[Ch.3, exs.3.1]{Kac93} we have
\begin{equation}
  \label{dual_diag}
      T^{\vee} = T^{-1}\hspace{0.5mm}.
\end{equation}
Since $K^{\vee} = K^t = {\bf B}T = T^{-1}{\bf B}^{\vee}$, we see
that dual Tits forms are related as follows:
\begin{equation}
  \label{dual_tits_forms}
      {\bf B}^{\vee} = T{\bf B}T\hspace{0.5mm}.
\end{equation}

\subsection{Eigenvalues of the Cartan matrix and the Coxeter transformation}
There is a simple relation between eigenvalues of the Cartan
matrix and the Coxeter transformation. Let vector $z \in
\mathcal{E}_\varGamma$ be given in two-component form
\begin{equation}
   z = \left (
   \begin{array}{c}
      x  \\
      y
   \end{array}
   \right ),
\end{equation}
where $\dim\mathcal{E}_\varGamma = k + m = n$, $x \in
\mathcal{E}_{\varGamma_a}$ and $y \in \mathcal{E}_{\varGamma_b}$.
Consider the relation
$$
   {\bf C}z = \lambda{z} \hspace{5mm} \text{ or } \hspace{5mm}
   w_2{z} = \lambda{w}_1{z}.
$$
In the simply-laced case we deduce from (\ref{symmetric_B}) that
\begin{equation}
\label{C_and_B}
   {\bf C}{z} = \lambda{z}
   \hspace{3mm}
   \Longleftrightarrow
   \hspace{3mm}
   \left \{
       \begin{array}{cc}
         \displaystyle\frac{\lambda+1}{2\lambda}x & = Dy     \vspace{3mm} \\
         \displaystyle\frac{\lambda+1}{2}y & = D^{t}x
       \end{array}
   \right.
   \hspace{3mm}
   \Longleftrightarrow
   \hspace{3mm}
   {\bf B}z =
         \frac{\lambda-1}{2}
   \left (
       \begin{array}{c}
         \displaystyle\frac{1}{\lambda}x \vspace{3mm}  \\
         -y
       \end{array}
   \right).
\end{equation}
In the multiply-laced case we deduce that
\index{Coxeter transformation! - multiply-laced case}
\index{Coxeter transformation! - simply-laced case}
\begin{equation}
\label{C_and_B_DF}
   {\bf C}{z} = \lambda{z}
   \hspace{3mm}
   \Longleftrightarrow
   \hspace{3mm}
   \left \{
       \begin{array}{cc}
         \displaystyle\frac{\lambda+1}{2\lambda}x & = Dy     \vspace{3mm} \\
         \displaystyle\frac{\lambda+1}{2}y & = Fx
       \end{array}
   \right.
   \hspace{3mm}
   \Longleftrightarrow
   \hspace{3mm}
   {\bf B}z =
         \frac{\lambda-1}{2}
   \left (
       \begin{array}{c}
         \displaystyle\frac{1}{\lambda}x \vspace{3mm}  \\
         -y
       \end{array}
   \right).
\end{equation}
From (\ref{C_and_B}) and (\ref{C_and_B_DF}) we have
in the simply-laced and multiply-laced cases, respectively:
\index{multiply-laced case}
\index{simply-laced case}
\begin{equation}
  \label{DDt_DtD}
\left \{
  \begin{array}{c}
     DD^{t}x = \displaystyle\frac{(\lambda+1)^2}{4\lambda}x \vspace{3mm} \\
     D^{t}Dy = \displaystyle\frac{(\lambda+1)^2}{4\lambda}y
  \end{array}
 \right.  \hspace{9mm}
 \left \{
  \begin{array}{c}
     DFx = \displaystyle\frac{(\lambda+1)^2}{4\lambda}x \vspace{3mm} \\
     FDy = \displaystyle\frac{(\lambda+1)^2}{4\lambda}y
  \end{array}
 \right.
\end{equation}
Similarly,
\begin{equation}
\label{B_eigenv}
   {\bf B}{z} = \gamma{z}
   \hspace{3mm}
   \Longleftrightarrow
   \hspace{3mm}
   \left \{
       \begin{array}{cc}
         (\gamma-1)x & = Dy  \vspace{3mm} \\
         (\gamma-1)y & = D^{t}x
       \end{array}
   \right.
   \hspace{7mm}
   {\rm resp.} \hspace{3mm}
   \left \{
       \begin{array}{cc}
         (\gamma-1)x & = Dy  \vspace{3mm} \\
         (\gamma-1)y & = Fx;
       \end{array}
   \right.
\end{equation}
and from (\ref{B_eigenv}) we have
\begin{equation}
\left \{
  \begin{array}{c}
     DD^{t}x = (\gamma-1)^2x  \\
     D^{t}Dy = (\gamma-1)^2y
  \end{array}
\right .
 \hspace{7mm}
   {\rm resp.} \hspace{3mm}
\left \{
  \begin{array}{c}
     DFx = (\gamma-1)^2x  \\
     FDy = (\gamma-1)^2y.
  \end{array}
\right .
\end{equation}
By (\ref{symmetrizable_K}) we have
\begin{proposition}[Fixed Points]
\label{fixed_points} 1) The eigenvalues $\lambda$ of the Coxeter
transformation
  and the eigenvalues $\gamma$ of the matrix {\bf B} of the
  Tits form are related as follows
$$
    \frac{(\lambda+1)^2}{4\lambda} = (\gamma-1)^2.
$$

 2) The kernel of the matrix {\bf B} coincides
with the kernel of the Cartan matrix $K$ and coincides with the
space of fixed points of the Coxeter transformation
$$
   \ker{K} = \ker{\bf B} =
     \{ z \mid {\bf C}z = z \}.
$$

3) The space of fixed points of matrix {\bf B} coincides with the
space of anti-fixed points of the Coxeter transformation
$$
   \{ z \mid {\bf B}z = z \} =
   \{ z \mid {\bf C}z = -z \}.
$$
\end{proposition}

For more information about fixed and anti-fixed points of the
powers of the Coxeter transformation,
see Appendix\hspace{1mm}\ref{fixed_anti_fixed}.

\section{An application of the Perron-Frobenius theorem}

\subsection{The pair of matrices $DD^t$ and $D^tD$
(resp. $DF$ and $FD$)}
\begin{remark}
{\rm
The matrices $DD^t$ and $D^tD$ have some nice properties.

(1) They give us all information about the eigenvalues
of Coxeter transformations and Cartan matrices.
 Our results hold for an arbitrary tree $\varGamma$.

(2) The eigenvectors of Coxeter transformations are
combinations of eigenvectors of the matrix $DD^t$ and
eigenvectors of the matrix $D^tD$,
see \S\ref{explicitly}, relations (\ref{case_non_01}),
(\ref{case_1}), (\ref{case_0}).

(3) They satisfy the Perron-Frobenius theorem.
}
\end{remark}

More exactly, properties (1)--(3) will be considered below in
Proposition \ref{gold_pair}, Proposition \ref{pf_pair}, and in
\S\ref{pf_section}, \S\ref{explicitly}.

\begin{proposition}
 \label{gold_pair}
1) The matrices $DD^t$ and $D^tD$ (resp. $DF$ and $FD$)
  have the same non-zero eigenvalues with equal multiplicities.

2) The eigenvalues  $\varphi_i$ of matrices $DD^t$ and $D^tD$
  (resp. $DF$ and $FD$) are non-negative:
  $$
      \varphi_i \geq 0.
  $$

3) The corresponding eigenvalues $\lambda^{\varphi_i}_{1,2}$ of
Coxeter transformations are
  \begin{equation}
    \label{Coxeter_eigenvalues}
        \lambda^{\varphi_i}_{1,2} =
        2\varphi_i - 1 \pm 2\sqrt{\varphi_i(\varphi_i-1)}.
  \end{equation}
The eigenvalues $\lambda^{\varphi_i}_{1,2}$ either lie
  on the unit circle or are real positive numbers. It the latter case
  $\lambda^{\varphi_i}_1$ and $\lambda^{\varphi_i}_2$ are mutually inverse:
 $$
     \lambda^{\varphi_i}_1\lambda^{\varphi_i}_2 = 1.
 $$
\end{proposition}
 \PerfProof
1)   If $DD^tz = \mu{z}$, where $\mu \neq 0$, then $D^tz \neq 0$
and $D^tD(D^tz) = \mu(D^tz)$. We argue similarly for $D^tD$, $DF$,
$FD$. The multiplicities of non-zero eigenvalues coincide since
\begin{equation*}
   x \neq 0,  \hspace{2mm} DD^t{x} = \mu{x} \hspace{3mm}
    \Longrightarrow  \hspace{3mm} D^tD(D^t{x})  =
      \mu(D^t{x}), \text{ and } D^t{x} \neq 0,
\end{equation*}
and
\begin{equation*}
   y \neq 0,  \hspace{2mm} D^tD{y}  = \mu{y} \hspace{3mm}
    \Longrightarrow  \hspace{3mm} DD^t(D{y})  =
       \mu(D{y}), \text{ and } D{y} \neq 0. 
\end{equation*}
\begin{remark}
  {\rm
  The multiplicities of the zero eigenvalue are not equal.
  If $DD^t{x} = 0$ and $x$ is the eigenvector, $x \neq 0$,
  then it is possible that $D^t{x} = 0$ and
  $D^t{x}$ is not an eigenvector of $D^tD$,
  see Remark \ref{remark_star}, 4) below.
  }
\end{remark}

2) The matrices $DD^t$ and $D^tD$ are symmetric and nonnegative
definite. For example,
$$
  <DD^t{x},x>~ = ~<D^t{x}, D^t{x}>~ \geq ~0.
$$
So, if $DD^t{x} = \varphi_i{x}$, then
$$
    <DD^t{x},x>~ = ~\varphi_i<x,x>
$$
and
$$
    \varphi_i  =  \frac{<D^t{x}, D^t{x}>}{<x,x>} \geq  0.
$$
In the multiply-laced case, we deduce from (\ref{matrix_T}) that
the matrix $DF$ is
\begin{equation}
  \label{fd_decomp}
   DF = T_1AT_2A^t.
\end{equation}
Let $\varphi$ be a non-zero eigenvalue for $DF = T_1AT_2A^t$
with eigenvector $x$:
$$
     T_1AT_2A^t{x} ~=~ \varphi{x}.
$$
Since $T$ is a positive diagonal matrix, see \S\ref{cartan}, we have
\begin{equation}
  \label{fd_decomp_2}
   (\sqrt{T_1}AT_2A^t\sqrt{T_1})((\sqrt{T_1})^{-1}{x})
   = \varphi(\sqrt{T_1})^{-1}{x},
\end{equation}
and $\varphi$ is also a non-zero eigenvalue with
eigenvector $(\sqrt{T_1})^{-1}{x}$ for the
matrix $\sqrt{T_1}AT_2A^t\sqrt{T_1}$ which already is symmetric, so
$\varphi \geq 0$.\hspace{3mm} 

\indent 3)
    From (\ref{Coxeter_eigenvalues}) if $0 \leq \varphi \leq 1$
we deduce
that
$$
   {|\lambda_{1,2}^{\varphi_i}|}^2 =
   (2\varphi_i - 1)^2 + 4\varphi_i(1 - \varphi_i) = 1  \hspace{1mm}.
$$
If $\varphi_i > 1$, then
$$
    2\varphi_i - 1 > 2\sqrt{\varphi_i(\varphi_i - 1)}
    \Longrightarrow  \lambda_{1,2}^{\varphi_i} \geq 0  \hspace{1mm}.
$$
Thus,
\begin{equation}
 \label{lambda_phi}
  \begin{array}{cc}
   \lambda_1^{\varphi_i} =
    2\varphi_i - 1 + 2\sqrt{\varphi_i(\varphi_i - 1)} > 1 \hspace{1mm},
    \vspace{3mm} \\
   \lambda_2^{\varphi_i} =
    2\varphi_i - 1 - 2\sqrt{\varphi_i(\varphi_i - 1)} < 1 \hspace{1mm}.
\qed
  \end{array}
\end{equation}

\index{PF-pair} The pair of matrices $(A,B)$ is said to be a {\it
PF-pair} if matrices $A$ and $B$ satisfy the Perron-Frobenius
theorem \footnote{The Perron-Frobenius theorem is well known in
the matrix theory, see \S\ref{sect_pf} and \cite{MM64},
\cite{Ga90}.}.

\begin{proposition} 
  \label{pf_pair}
The matrix pair $(DD^t, D^tD) \hspace{2mm} (resp. \hspace{2mm} (DF, FD))$
is a PF-pair, i.e.,

 $1)~ DD^t$ and $D^tD$
       $(resp.  \hspace{2mm} DF ~and~ FD)$ are non-negative;

 $2)~ DD^t$ and $D^tD$
       $(resp. \hspace{2mm} DF ~and~ FD)$ are indecomposable.
\end{proposition}
 \PerfProof
1) Indeed, in the simply-laced case, the following relation holds
\begin{equation}
  \label{dd_formula}
  \begin{split}
     4(DD^t)_{ij} = &
     4\sum\limits_{p=1}^k{(a_i,b_p)(b_p, a_j)}  = \\
     & \left \{
      \begin{array}{cc}
        s_i  & \text{ the number of edges with a vertex }
        v_i \text{ if } i = j, \\
        1 & \hspace{4mm} \text{ if } \hspace{4mm} |v_i - v_j| = 2, \\
        0 & \hspace{4mm} \text{ if } \hspace{4mm} |v_i - v_j| > 2.
      \end{array}
     \right .
  \end{split}
\end{equation}
In the multiply-laced case, we have
\begin{equation}
    \label{df_formula}
    \begin{split}
     4(DF)_{ij} & =
     4\sum\limits_{p=1}^k
        {\frac{(a_i,b_p)(b_p, a_j)}{(a_i,a_i)(b_p, b_p)}} = \vspace{7mm}\\
    &
       \left  \{
         \begin{array}{cc}
              4\sum\limits_{p=1}^k
        {\displaystyle\frac{(a_i,b_p)^2}{(a_i,a_i)(b_p, b_p)}} =
             4\sum\limits_{p=1}^k {\cos}^2\{a_i, b_p\} & \text{ if }
                 i = j,\vspace{3mm} \\
   4\displaystyle\frac{(a_i,b_p)(b_p, a_j)}{(a_i,a_i)(b_p, b_p)}
                &  \text{ if } |v_i - v_j| = 2, \vspace{3mm} \\
            0  &  \text{ if } |v_i - v_j| > 2.
         \end{array}
       \right .
      \end{split}
\end{equation}

2) Define the distance between two sets of vertices
$A = \{a_i\}_{i \in I}$ and $B = \{b_j\}_{j \in J}$ to be
\begin{equation*}
      \min_{i,j}|a_i - b_j|.
\end{equation*}
If the matrix $DD^t$ is decomposable, then the set of
   vertices \{$v_1,~\ldots ,~v_n$\} can be divided into two
   subsets such that distance between these two subsets is $>2$.
   This contradicts the assumption that $\varGamma$ is connected.
\hspace{3mm}\qedsymbol
\begin{remark}
\label{remark_star}
{\rm
 1) Formula (\ref{dd_formula}) is a particular case of
(\ref{df_formula}), since the angles between the adjacent simple
roots $a_i$ and $b_j$ are $\displaystyle\frac{2\pi}{3}$, so
$\cos\{a_i, b_j\} = \displaystyle\frac{1}{2}$. Of course, the
angles and the lengths of vectors are considered in the sense of
the bilinear form $(\cdot , \cdot)$ from (\ref{bilinear}).

 2) The case $|v_i - v_j| = 2$ from (\ref{df_formula}) can be
expressed in the following form:
\begin{equation}
  \begin{array}{c}
    \displaystyle 4\frac{(a_i,b_p)(b_p, a_j)}{(a_i,a_i)(b_p, b_p)} =
        4\frac{(a_i,b_p)(b_p, a_j)}{|a_i|^2|b_p|^2} = \vspace{3mm} \\
    \displaystyle
4\frac{(a_i,b_p)}{|a_i||b_p|}\frac{(b_p,a_j)}{|a_j||b_p|}\frac{|a_j|}{|a_i|}  =
        4\frac{|a_j|}{|a_i|}\cos\{a_i,b_p\}\cos\{a_j,b_p\}.
  \end{array}
\end{equation}

\begin{figure}[h]
\centering
\includegraphics{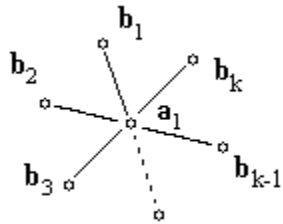}
\caption{\hspace{3mm} The star $*_{k+1}$ with $k$ rays}
\label{small_star}
\end{figure}

3)  One can easily calculate the matrices $DD^t, D^tD, DF, FD$ by means of
  (\ref{dd_formula}) and (\ref{df_formula}).

4) Consider a simple star  $\ast_{k+1}$, Fig.~\ref{small_star}.
  It is bipartite with respect to the following bicolored partition.
  One part of the graph consists of only one vertex $a_1, (m = 1),$
  the other one consists of $k$ vertices \{$b_1,\ldots,b_k$\},
  $n = k+1$. According to (\ref{dd_formula}) the $1\times{1}$ matrix
  $DD^t$ is
  $$
        DD^t = k = n-1,
  $$
and the $k\times{k}$ matrix $D^tD$ is}
  \begin{displaymath}
       D^tD = \left (
        \begin{array}{ccccc}
          1 & 1 & 1 & \dots & 1 \\
          1 & 1 & 1 & \dots & 1 \\
          1 & 1 & 1 & \dots & 1 \\
            &   &   & \dots &   \\
          1 & 1 & 1 & \dots & 1
        \end{array}
        \right ).
  \end{displaymath}
\end{remark}
By Proposition \ref{gold_pair},
heading 2) the matrices $DD^t$ and $D^tD$ have
only one non-zero eigenvalue $\varphi_1 = n-1$. All the other
eigenvalues of $D^tD$ are zeros and the characteristic polynomial of
the $D^tD$ is
$$
         \varphi^{n-1}(\varphi - (n-1)).
$$

\subsection{The Perron-Frobenius theorem applied to $DD^t$ and $D^tD$
  (resp. $DF$ and $FD$)}
  \label{pf_section}
By Proposition \ref{pf_pair} the pairs ($DD^t, D^tD$) (resp. ($DF,
FD$)) are PF-pairs, so we can apply the Perron-Frobenius theorem,
see \S\ref{sect_pf}.
 \index{theorem! - Perron-Frobenius}
 \index{eigenvalues}
\begin{corollary}
  \label{corollary_dominant}
   The matrices  $DD^t$ and $D^tD$ (resp. $DF$ and $FD$)
   have a common simple (i.e., with multiplicity one)
   positive eigenvalue $\varphi_1$.
   This eigenvalue is maximal {\em (called dominant eigenvalue)}:
\begin{equation}
   0 \leq \varphi_i \le \varphi_1,   \hspace{5mm}
   \varphi_1 = \left \{
   \begin{array}{c}
    \displaystyle  \max_{x \geq 0}\min_{1 \leq i \leq m}
      \frac{(DD^tx)_i}{x_i} \hspace{5mm}
      \text{ in the simply-laced case},  \vspace{3mm} \\
    \displaystyle  \max_{x \geq 0}\min_{1 \leq i \leq m}
      \frac{(DFx)_i}{x_i} \hspace{5mm}
      \text{ in the multiply-laced case}. \vspace{3mm}
   \end{array}
   \right .
\end{equation}
There are positive eigenvectors $\mathbb{X}^{\varphi_1}$,
$\mathbb{Y}^{\varphi_1}$ (i.e., non-zero vectors with non-negative
coordinates) corresponding to the eigenvalue $\varphi_1$:
\begin{equation}
  \begin{array}{ccc}
    DD^t\mathbb{X}^{\varphi_1} & = \varphi_1\mathbb{X}^{\varphi_1},
\hspace{5mm}
    D^tD\mathbb{Y}^{\varphi_1} & = \varphi_1\mathbb{Y}^{\varphi_1}
\hspace{1mm}, \\
    DF\mathbb{X}^{\varphi_1} & = \varphi_1\mathbb{X}^{\varphi_1},
\hspace{5mm}
    FD\mathbb{Y}^{\varphi_1} & = \varphi_1\mathbb{Y}^{\varphi_1}
\hspace{1mm}.
   \end{array}
\end{equation}
\end{corollary}

\begin{figure}[h]
\centering
\includegraphics{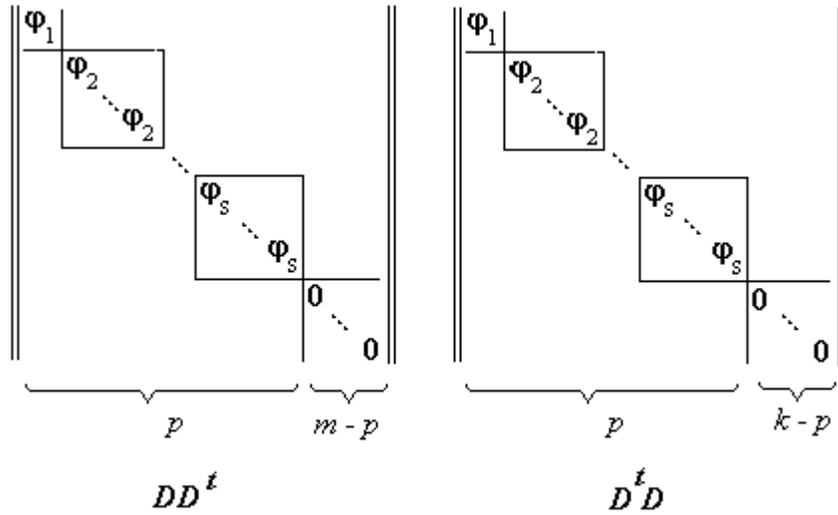}
\caption{\hspace{3mm} The Jordan normal forms of $DD^t$ and $D^tD$}
\label{normal_ddt}
\end{figure}

The matrices $DD^t$ (resp. $D^tD$) are symmetric and can be
diagonalized in the some orthonormal basis of the eigenvectors
from $\mathcal{E}_{\varGamma_a} = \mathbb{R}^h$ (resp.
$\mathcal{E}_{\varGamma_b} = \mathbb{R}^k$). The Jordan normal
forms of these matrices are shown on Fig.~\ref{normal_ddt}.
The normal forms of $DF$ and $FD$ are the same,
since by (\ref{fd_decomp}), (\ref{fd_decomp_2}) the non-zero eigenvalues
of the matrix $DF$ and symmetric matrix $\sqrt{T_1}AT_2A^t\sqrt{T_1}$
coincide. However, the normal bases (i.e., bases which consist of eigenvectors)
for $DF$ and $FD$ are
not necessarily orthonormal, since the eigenvectors of
$DF$ are
obtained from eigenvectors of $\sqrt{T_1}AT_2A^t\sqrt{T_1}$ by
means of the matrix $\sqrt{T_1}$ which does not preserve orthogonality.

\section{The basis of eigenvectors and a theorem on the Jordan form}
\subsection{An explicit construction of the eigenvectors}
 \label{explicitly}
\begin{proposition}
 \label{eigenvectors_creation}
  Let
\begin{equation}
 \label{basis_X}
  \mathbb{X}^{\varphi_1}_1,
  \mathbb{X}^{\varphi_2}_1,...,\mathbb{X}^{\varphi_2}_{t_2},...,
  \mathbb{X}^{\varphi_s}_1,...,\mathbb{X}^{\varphi_s}_{t_s},
  \mathbb{X}^0_1,...,\mathbb{X}^0_{m-p}
\end{equation}
be all the orthonormal eigenvectors for $DD^t$ and
\begin{equation}
  \mathbb{Y}^0_1,...,\mathbb{Y}^0_{k-p}
\end{equation}
be all the orthonormal eigenvectors for $D^tD$ corresponding
to the zero eigenvalue. Then
\begin{equation}
 \label{basis_Y}
  D^t\mathbb{X}^{\varphi_1}_1,
  D^t\mathbb{X}^{\varphi_2}_1,...,D^t\mathbb{X}^{\varphi_2}_{t_2},...,
  D^t\mathbb{X}^{\varphi_s}_1,...,D^t\mathbb{X}^{\varphi_s}_{t_s},
  \mathbb{Y}^0_1,...,\mathbb{Y}^0_{k-p}
\end{equation}
is the set of all orthonormal eigenvectors for $DD^t$.

Bases for $DF, FD$ (not orthonormal) are similarly constructed.
\end{proposition}

 \PerfProof
 Indeed, if $\mathbb{X}^{\varphi_i}$ and $\mathbb{X}^{\varphi_j}$
 are the eigenvectors corresponding to the $\varphi_i$ and $\varphi_j$,
 then the vectors $\mathbb{Y}^{\varphi_i} = D^t\mathbb{X}^{\varphi_i}$
 and $\mathbb{Y}^{\varphi_j} = D^t\mathbb{X}^{\varphi_j}$ are
 eigenvectors of $D^tD$ and
 $$
   <\mathbb{Y}^{\varphi_i}, \mathbb{Y}^{\varphi_j}> =
   <\mathbb{X}^{\varphi_i}, DD^t\mathbb{X}^{\varphi_j}> =
   \varphi_j<\mathbb{X}^{\varphi_i}, \mathbb{X}^{\varphi_j}> = 0 \hspace{1mm}.
   \qed
 $$

Let us construct the eigenvectors for the Coxeter transformation.
We set:

\underline{Case $\varphi_i \neq 0,1$}:
\begin{equation}
 \label{case_non_01}
   z^{\varphi_i}_{r,\nu} =
   \left (
     \begin{array}{c}
        \mathbb{X}^{\varphi_i}_r  \vspace{3mm} \\
       -\displaystyle\frac{2}{\lambda^{\varphi_i}_{\nu}+1}
       D^t\mathbb{X}^{\varphi_i}_r
     \end{array}
   \right ), \hspace{4mm}
   1 \leq i \leq s, \hspace{2mm} 1 \leq r \leq t_i,
   \hspace{2mm} \nu = 1,2 \hspace{1mm}.
\end{equation}
Here $\lambda^{\varphi_i}_{1,2}$ is obtained
by formula (\ref{Coxeter_eigenvalues}).

\underline{Case $\varphi_i = 1$}:
\begin{equation}
 \label{case_1}
   z^1_{r} =
   \left (
     \begin{array}{c}
        \mathbb{X}^1_r  \vspace{3mm} \\
       -D^t\mathbb{X}^1_r
     \end{array}
   \right ), \hspace{4mm}
   \tilde{z}^1_{r} =
   \frac{1}{4}
   \left (
     \begin{array}{c}
        \mathbb{X}^1_r  \vspace{3mm} \\
       D^t\mathbb{X}^1_r
     \end{array}
   \right ), \hspace{4mm}
     1 \leq r \leq t_i.
\end{equation}

\underline{Case $\varphi_i = 0$}:
\begin{equation}
 \label{case_0}
   z^0_{x_\eta} =
   \left (
     \begin{array}{c}
       \mathbb{X}^0_{\eta}  \vspace{3mm} \\
       0
     \end{array}
   \right ), \hspace{2mm}
   1 \leq \eta \leq m-p,
   \hspace{9mm}
   z^0_{y_\xi} =
   \left (
     \begin{array}{c}
       0      \vspace{3mm}     \\
       \mathbb{Y}^0_{\xi}
     \end{array}
   \right ), \hspace{2mm}
   1 \leq \xi \leq k-p.
\end{equation}

\begin{proposition} [\cite{SuSt75, SuSt78}]
 \label{eigenvectors_basis}

1) The vectors (\ref{case_non_01}), (\ref{case_1}), (\ref{case_0})
  constitute a basis in $\mathcal{E}_\varGamma$ over $\mathbb{C}$.

2) The vectors (\ref{case_non_01}) are eigenvectors of the Coxeter
transformation
  corresponding to the eigenvalue $\lambda_{\varphi_i}$:
\begin{equation}
  \label{PK_not_01}
   {\bf C}{z}^{\varphi_i}_{r,\nu} =
      \lambda{z}^{\varphi_i}_{r,\nu}.
\end{equation}

The vectors (\ref{case_1}) are eigenvectors and adjoint vectors
  of the Coxeter transformation corresponding to the eigenvalue $1$:
\begin{equation}
  \label{PK_1}
   {\bf C}z^1_r = z^1_r, \hspace{5mm}
   {\bf C}\tilde{z}^1_r = z^1_r + \tilde{z}^1_r.
\end{equation}

  The vectors (\ref{case_0}) are eigenvectors
  of the Coxeter transformation corresponding to the \\
  eigenvalue $-1$:
\begin{equation}
  \label{PK_0}
   {\bf C}z^0_{x_\eta} = -z^0_{x_\eta}, \hspace{5mm}
   {\bf C}z^0_{y_\xi} = -z^0_{y_\xi}.
\end{equation}
In other words, vectors (\ref{case_non_01}), (\ref{case_1}), (\ref{case_0})
constitute an orthogonal basis which consists of eigenvectors and adjoint
vectors
of the Coxeter transformation. The number
of the adjoint vectors $\tilde{z}^1_r$
is equal to the multiplicity of eigenvalue $\varphi = 1$.
 \end{proposition}

\PerfProof 1) The number of vectors    (\ref{case_non_01}),
(\ref{case_1}), (\ref{case_0}) is $2p + (m-p) + (k-p) = n$ and it
suffices to prove that these vectors are linearly independent. Let
us write down the condition of linear dependence. It splits into 2
conditions: for the $\mathbb{X}$-component and for the
$\mathbb{Y}$-component. The linear independence of vectors
(\ref{case_non_01}), (\ref{case_1}), (\ref{case_0}) follows from
the linear independence of vectors (\ref{basis_X}),
(\ref{basis_Y}). \hspace{3mm} \qedsymbol

2) To prove relation (\ref{PK_not_01}),
the first relation from (\ref{PK_1}) and
relation (\ref{PK_0}), it suffices to check (\ref{C_and_B}).
Let us check that
$$
   {\bf C}\tilde{z}^1_r = z^1_r + \tilde{z}^1_r.
$$
We consider the multiply-laced case.
 Making use of (\ref{matrix_K}) we see that

\begin{displaymath}
   {\bf C} =
     \left (
        \begin{array}{cc}
          4DF - I_m   &  2D \\
          -2F  & I_k
        \end{array}
     \right ).
\end{displaymath}
Then
$$
  {\bf C}\tilde{z}^1_r =
    \frac{1}{4}
    \left (
      \begin{array}{c}
           4DF\mathbb{X}^1_r - \mathbb{X}^1_r  +  2DF\mathbb{X}^1_r \\
           -2F\mathbb{X}^1_r - F\mathbb{X}^1_r
      \end{array}
    \right ) =
    \frac{1}{4}
    \left (
      \begin{array}{c}
            5\mathbb{X}^1_r  \\
           -3\mathbb{X}^1_r
      \end{array}
    \right ) =
   z^1_r + \tilde{z}^1_r. \qed
$$

By heading 2) of Proposition \ref{fixed_points}
and formula (\ref{case_1}) from Proposition \ref{eigenvectors_creation} we have
\begin{corollary}
 \label{diag_form}
The Jordan normal form of the Coxeter transformation is diagonal
if and only if
  the Tits form is nondegenerate \footnote{For the Jordan canonical (normal)
  form,
  see, for example, \cite[Ch.VI]{Ga90} or \cite[Ch.III]{Pr94}.}.  \qedsymbol
\end{corollary}

\subsection{Monotonicity of the dominant eigenvalue}
The following proposition is important for calculation of the
number of $2\times2$ Jordan blocks in the Jordan normal form of
the Coxeter transformation.
\begin{proposition} [\cite{SuSt75, SuSt78}]
  \label{trend_fi}
Let us add an edge to some tree $\varGamma$ and let
$\stackrel{\wedge}{\varGamma}$ be the new graph  (Remark
\ref{operations}). Then:

\index{dominant eigenvalue}
  1) The dominant eigenvalue $\varphi_1$ may only grow:
$$
   \varphi_1(\stackrel{\wedge}{\varGamma}) \geq \varphi_1(\varGamma)\hspace{1mm}.
$$

  2) Let $\varGamma$ be an extended Dynkin diagram, i.e., $\mathcal{B} \geq 0$.
Then the spectra of $DD^t(\stackrel{\wedge}{\varGamma})$ and
$D^tD(\stackrel{\wedge}{\varGamma})$ (resp.
$DF(\stackrel{\wedge}{\varGamma})$ and
 $FD(\stackrel{\wedge}{\varGamma})$) do not contain $1$, i.e.,
$$
   \varphi_i(\stackrel{\wedge}{\varGamma}) \neq 1
   \hspace{3mm} {\rm for \hspace{1mm} all \hspace{1mm}} \varphi_i \in
       Spec(DD^t(\stackrel{\wedge}{\varGamma}))\hspace{0.5mm}.
$$

 3) Let $\mathcal{B}$ be indefinite. Then
$$
      \varphi_1(\stackrel{\wedge}{\varGamma}) > 1\hspace{0.5mm}.
$$
\end{proposition}
 \PerfProof 1) Adding an edge to the vertex $a_i$ we see, according to
(\ref{dd_formula}), that only one element of $DD^t$ changes:
namely, $(DD^t)_{ii}$ changes from $\displaystyle \frac{s_i}{4}$
to $\displaystyle \frac{s_i + 1}{4}$. For the multiply-laced case,
$(DF)_{ii}$ changes by ${\cos}^2(a_i, b_s)$, where $b_s$ is the
new vertex connected with the vertex $a_i$. By Corollary
\ref{corollary_dominant} we have
\begin{equation*}
 \begin{split}
   \varphi_1(\stackrel{\wedge}{\varGamma}) = &
   \max_{x \geq 0}\min_{1 \leq i \leq m}
   \frac{(DF(\stackrel{\wedge}{\varGamma})x)_i}{x_i} \geq \\
   & \min_{1 \leq i \leq m}
   \frac{(DF(\stackrel{\wedge}{\varGamma})\mathbb{X}^{\varphi_1}_i)_i}
            {\mathbb{X}^{\varphi_1}_i}  \geq
   \min_{1 \leq i \leq m}
   \frac{(DF\mathbb{X}^{\varphi_1}_i)_i}
            {\mathbb{X}^{\varphi_1}_i}  =
   \varphi_1\hspace{0.5mm}.
 \end{split}
\end{equation*}

2) The characteristic polynomial of
$DF(\stackrel{\wedge}{\varGamma})$ is
\begin{equation}
  \det|DF(\stackrel{\wedge}{\varGamma}) - \mu{I}| =
    \det|DF - \mu{I}| +
    {\cos}^2\{a_i, b_s\}\hspace{0.5mm}\det|A_i(\mu)|\hspace{1mm},
\end{equation}
where $A_i(\mu)$ is obtained by deleting the $i^{th}$ row and
$i^{th}$ column from the matrix $DF - \mu{I}$. \index{operation
\lq\lq Remove Vertex"} It corresponds to the operation \lq\lq{\it
Remove Vertex}" from the graph $\varGamma$  (Remark
\ref{operations}), the graph obtained by removing vertex is
$\stackrel{\vee}{\varGamma}$, i.e.,
\begin{equation}
  \label{3_char_pol}
  \det|DF(\stackrel{\wedge}{\varGamma}) - \mu{I}| =
    \det|DF - \mu{I}| +
    {\cos}^2\{a_i, b_s\}
     \hspace{0.5mm}\det|DF(\stackrel{\vee}{\varGamma}) - \mu{I}|\hspace{1mm}.
\end{equation}
According to Remark \ref{operations} the quadratic form
$\mathcal{B}(\stackrel{\vee}{\varGamma})$ is positive,
$\stackrel{\vee}{\varGamma}$ is the Dynkin diagram, i.e.,
$\ker{\bf B} = 0$. By Proposition \ref{fixed_points} the Coxeter
transformation for $\stackrel{\vee}{\varGamma}$ does not have
eigenvalue 1. Then by (\ref{Coxeter_eigenvalues}) the
corresponding matrix $DF(\stackrel{\vee}{\varGamma})$ does not
have eigenvalue 1. Thus, in (\ref{3_char_pol}), $\mu = 1$ is a
root of $\det|DF - \mu{I}|$ and is not a root of
$\det|DF(\stackrel{\vee}{\varGamma}) - \mu{I}|$, and therefore
$\mu = 1$ is not a root of $\det|DF(\stackrel{\wedge}{\varGamma})
- \mu{I}|$. The case  ${\cos}^2\{a_i, b_s\} = 0$ is not possible
since $\stackrel{\wedge}{\varGamma}$ is connected.
\hspace{3mm} 

3) By heading 1) adding only one edge to an extended Dynkin
diagram we get $\varphi_1(\stackrel{\wedge}{\varGamma}) \geq
\varphi_1$. By heading 2) $\varphi_1(\stackrel{\wedge}{\varGamma})
> \varphi_1$. But the form $\mathcal{B}$ becomes indefinite after
we add some edges to the extended Dynkin diagram (Remark
\ref{operations}, 4). Thus, $\varphi_1$ is only grows.
\hspace{3mm} \qedsymbol

\begin{proposition}
  The common dominant eigenvalue of $DD^t$ and $D^tD$
  (resp. $DF$ and $FD$) is equal to $1$ if and only if $\varGamma$
  is an extended Dynkin diagram.
\end{proposition}
 \PerfProof 1) Let $\varGamma$ be an extended Dynkin diagram,
then $\mathcal{B} \geq 0$. Since $\ker{\bf B} \neq 0$, we see
by Proposition \ref{fixed_points} that the eigenvalue $\varphi_1$ is
the eigenvalue of the matrices $DF$ and $FD$.
By (\ref{lambda_phi}) we have
$\lambda_1^{\varphi_1} \geq 1$.
Further, since the Weyl group preserves the quadratic form $\mathcal{B}$,
we have
$$
  \mathcal{B}(z^{\varphi_1}) =
  \mathcal{B}({\bf C}z^{\varphi_1}) =
  (\lambda_1^{\varphi_1})^2\mathcal{B}({\bf C}z^{\varphi_1}) .
$$
Therefore, either $\lambda_1^{\varphi_1} = 1$, i.e., $\varphi_1 = 1$,
or $\mathcal{B}(z^{\varphi_1}) = 0$. We will show that in
the latter case $\varphi_1 = 1$, too. By Proposition \ref{eigenvectors_creation}
the vectors $\mathbb{X}^{\varphi_1}$ and $F\mathbb{X}^{\varphi_1}$ have real
coordinates. By (\ref{lambda_phi})
the eigenvalue $\lambda_1^{\varphi_1}$ is also real because
$\varphi_1 \geq 1$. So, the vector $z^{\varphi_1}$ from (\ref{case_1}) is real.
Then from $\mathcal{B}(z^{\varphi_1})$ = 0 we have
$z^{\varphi_1} \in \ker{\bf B}$ and again, $\lambda_1^{\varphi_1} = 1$.

2) Conversely, let  $\lambda_1^{\varphi_1} = 1$. Then, by
Proposition \ref{fixed_points}, $\ker{\bf B} \neq 0$, i.e., the
form $\mathcal{B}$ is degenerate. Let us find whether
$\mathcal{B}$ is definite ($\mathcal{B} \geq 0$) or indefinite. By
heading 3) of Proposition \ref{trend_fi} if $\mathcal{B}$ is
indefinite, then $\lambda_1^{\varphi_1} > 1$. Thus, $\mathcal{B}
\geq 0$ and $\varGamma$ is an extended Dynkin diagram.
\hspace{3mm} \qedsymbol

\begin{corollary}  Let $\mathcal{B} \geq 0$, i.e., let
   $\varGamma$ be an extended Dynkin diagram. Then

1) The kernel of the quadratic form $\mathcal{B}$ is one-dimensional.

2) The vector $z^{\varphi_1}\in \ker{\bf B}$ can be chosen so that
     all its coordinates are positive.
\end{corollary}

\PerfProof
1) ${\varphi_1}$ is a simple eigenvalue of
 $DD^t$ and $D^tD$ (see Corollary \ref{corollary_dominant}).

 2) $\mathbb{X}^{\varphi_1}$ is a positive vector
(Corollary \ref{corollary_dominant}), $D^t\mathbb{X}^{\varphi_1}$ is
a negative vector since $D^t$ is
a nonpositive matrix (\ref{matrix_D}).
So, the vector $z^{\varphi_1}$ from (\ref{case_1}) is positive.
\hspace{3mm} \qedsymbol

The monotonicity of the dominant value $\varphi_1$ and the
corresponding maximal eigenvalue $\lambda^{\varphi_1}$ of the
Coxeter transformation is clearly demonstrated for the diagrams
$T_{2,3,r}$, $T_{3,3,r}$, $T_{2,4,r}$, see Propositions
\ref{polyn_T_1}, \ref{polyn_T_2}, \ref{polyn_T_3} and Tables
\ref{table_char_E_series}, \ref{table_char_E6_series},
\ref{table_char_E7_series} in \S\ref{section_T_pqr}.

 \index{theorem! - on the Jordan form}
\subsection{A theorem on the Jordan form}

Now we can summarize.
\begin{theorem} [\cite{SuSt75, SuSt78, St85}]
\label{th_jordan}
 1) The Jordan form of the Coxeter transformation is diagonal if and only if
  the Tits form is non-degenerate.

 2) If $\mathcal{B} \geq 0$ ($\varGamma$ is an extended Dynkin diagram), then
  the Jordan form of the Coxeter transformation contains only
  one $2\times2$
  Jordan block.  All eigenvalues $\lambda_i$ lie on the unit circle.

 3) If $\mathcal{B}$ is indefinite and degenerate, then
  the number of $2\times2$ Jordan blocks coincides with
  $\dim\ker{\bf B}$. The remaining Jordan blocks are $1\times1$.
  There is a simple maximal eigenvalue $\lambda^{\varphi_1}_1$
  and a simple minimal eigenvalue $\lambda^{\varphi_1}_2$, and
$$
   \lambda^{\varphi_1}_1 > 1, \hspace{5mm} \lambda^{\varphi_1}_2 < 1 .
$$
\end{theorem}

\begin{figure}[h]
\centering
\includegraphics{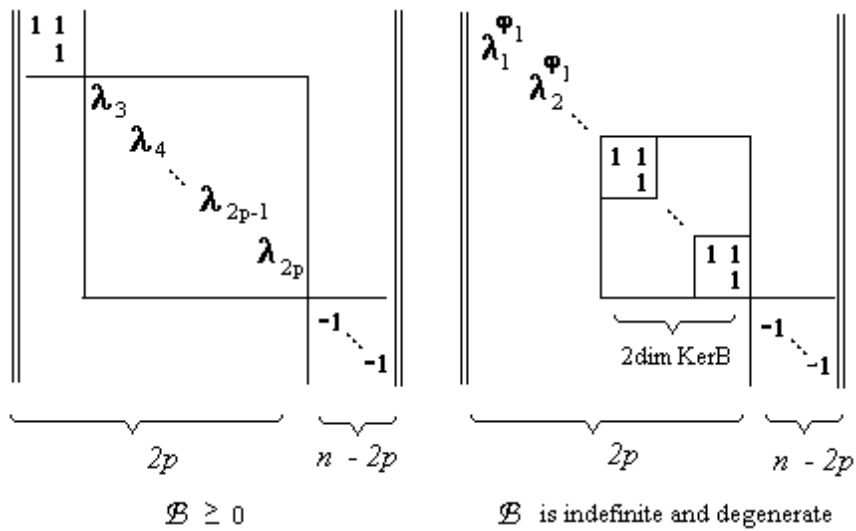}
\caption{\hspace{3mm} The Jordan normal form of the Coxeter transformation}
\label{jordan_form}
\end{figure}

\begin{example}[V.~Kolmykov]
{\rm  The example shows that there is a graph $\varGamma$ with
indefinite and degenerate quadratic form $\mathcal{B}$ such that
$\dim\ker{\bf B}$ is however large (see Fig.
\ref{dim_however_large}) and the Coxeter transformation has
however large number of $2\times2$ Jordan blocks. Consider $n$
instances of the extended Dynkin diagram $\widetilde{D}_4$ with
centers $b_1,...,b_n$ and marked vertices $a_1,...,a_n$. We add
new vertex $b_{n+1}$ and connect it with marked vertices $a_i$ for
$i=1,2...,n$.  The new graph is bipartite: one part consists of
the vertices $b_1,...,b_n, b_{n+1}$;  the matrix $D^tD$ is of size
$(n+1)\times(n+1)$ and has the form

\begin{figure}[h]
\centering
\includegraphics{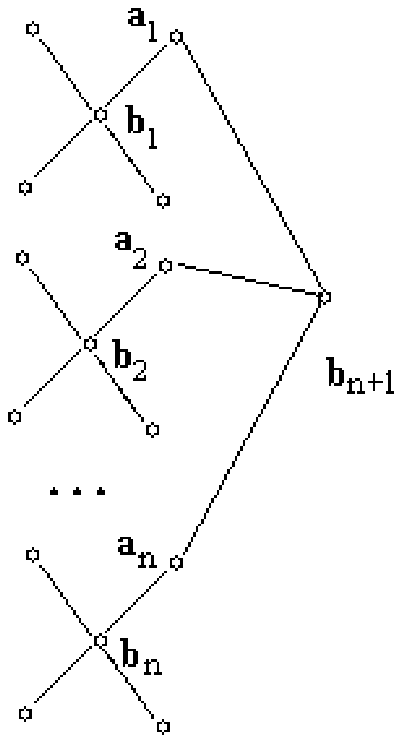}
\caption{\hspace{3mm} A graph $\varGamma$ such that $\dim\ker{\bf
B}$ is however large}
\label{dim_however_large}
\end{figure}

  \begin{displaymath}
       4D^tD = \left (
        \begin{array}{ccccccc}
          n & 1 & 1 & 1 & \dots & 1 & 1 \\
          1 & 4 & 0 & 0 & \dots & 0 & 0 \\
          1 & 0 & 4 & 0 & \dots & 0 & 0 \\
          1 & 0 & 0 & 4 & \dots & 0 & 0 \\
            &   &   &  &  \dots &   &   \\
          1 & 0 & 0 & 0 & \dots & 4 & 0 \\
          1 & 0 & 0 & 0 & \dots & 0 & 4
        \end{array}
        \right ).
  \end{displaymath}
It is easy to show that
 $$
    |4D^tD - \mu{I}| =
    (n - \mu)(4 - \mu)^n - n(4 - \mu)^{n-1} \hspace{0.5mm}.
 $$
 Thus, $\varphi_i = \displaystyle\frac{\mu_i}{4} = 1$ is of multiplicity
 $n - 1$.

  In Proposition \ref{gluing} of \S\ref{shearing} we will see more examples
of graphs $\varGamma(n)$ obtained by gluing $n$ copies of some
graph $\varGamma$. By this proposition every eigenvalue $\lambda$
of $\varGamma$ is also an eigenvalue of $\varGamma(n)$ and the
multiplicity of this eigenvalue for the graph $\varGamma(n)$ is
$(n - 1)\times{m}$, where $m$ is the multiplicity of this
eigenvalue in $\varGamma$. }
\end{example}

%% file: 4shearing.tex

\chapter{\sc\bf Eigenvalues, shearing formulas and diagrams $T_{p,q,r}$}

\setlength{\epigraphwidth}{72mm}

\epigraph{ \ldots most of the fundamental results about simple Lie
algebras, which were discovered by Killing are usually attributed
to \'E.~Cartan. This, despite the meticulousness with which Cartan
noted his indebtedness to Killing. In Cartan's thesis there are 28
references to Lie and 60 to Killing!} {A.~J.~Coleman,
\cite[p.447]{Col89}, 1989.}

\label{chapter_shearing}

\section[The eigenvalues of the affine Coxeter transformation]{The eigenvalues of the affine Coxeter transformation are
roots of unity} \label{roots_of_1}

 The Coxeter transformation corresponding to the
extended Dynkin diagram, i.e., corresponding to affine Kac-Moody
algebra is called the {\it affine Coxeter transformation}.
\begin{theorem} [\cite{SuSt79}, \cite{St82a}, \cite{St85}]
 \label{roots_unity}
 \index{roots of unity}
 The eigenvalues of the affine Coxeter transformation are roots of unity.
 The proper eigenvalues are collected in Table $\ref{table_eigenvalues}$.
\end{theorem}

\PerfProof The eigenvalues for all cases of extended Dynkin
diagrams are easily calculated by means of the generalized
R.~Steinberg's theorem (Theorem \ref{gen_Steinberg}) and Table
\ref{table_char_ext_polynom_Dynkin}. See Remark
\ref{orders_explanations}. \qedsymbol


\begin{remark}
\label{orders_explanations} \rm {According to generalized
R.~Steinberg's theorem
 (Theorem \ref{gen_Steinberg}) orders (not always
different) of eigenvalues (column 3, Table
\ref{table_eigenvalues}) coincide with lengths of branches of the
corresponding Dynkin diagram. If $g$ is the class number
(\ref{class_g}) of the extended Dynkin diagram, the number of
branches of the corresponding Dynkin diagram is $3 - g$. For $g =
0$, we have a simply-laced case of the extended Dynkin diagrams,
i.e., $\widetilde{E}_6$, $\widetilde{E}_7$, $\widetilde{E}_6$,
$\widetilde{D}_n (n \geq 4)$. In this case number of branches is
$3$ . For $g = 1$, we have extended Dynkin diagrams
$\widetilde{F}_{41}$, $\widetilde{F}_{42}$, $\widetilde{CD}_n$,
$\widetilde{CD}_n$, there exists $2$ groups of eigenvalues, see
(\ref{groups_of_chi}). For $g = 2$, we have extended Dynkin
diagrams $\widetilde{G}_{21}$, $\widetilde{G}_{22}$,
$\widetilde{B}_n$, $\widetilde{C}_n$, $\widetilde{BC}_n$ and there
exists only one group of eigenvalues. For $g = 3$, we have
extended Dynkin diagrams $\widetilde{A}_{11}$,
$\widetilde{A}_{12}$, in this case there is only one trivial
eigenvalues $1$}.
\end{remark}

\vfill
\newpage
\begin{table} 
  \centering
  \vspace{2mm}
  \caption{\hspace{3mm}The eigenvalues of affine Coxeter transformations}
  \renewcommand{\arraystretch}{1.5}
  \begin{tabular} {|c|c|c|c|}
  \hline \hline
      Diagram & Eigenvalues $\lambda$
                            & Orders of & Affine Coxeter \cr
       &  & eigenvalues\footnotemark[1] & number $h_a$ \footnotemark[2] \\
  \hline  \hline  
       $\widetilde{E}_6$
     &
       $\begin{array}{c}
        \lambda^1_{1,2} = 1, \\
        \lambda^2_{1,2} = \lambda^3_{1,2} = e^{\pm{2\pi{i}/3}}, \\
        \lambda_7 = -1
       \end{array}
       $
     & $\begin{array}{c}
          \\
          3,3, \\
          2
       \end{array}$
     & 6  \\
  \hline         
        $\widetilde{E}_7$
     &
       $\begin{array}{c}
        \lambda^1_{1,2} = 1, \\
        \lambda^2_{1,2} = e^{\pm{2\pi{i}/3}}, \\
        \lambda^3_{1,2} = e^{\pm{\pi{i}/2}}, \\
        \lambda_{7,8} = -1
       \end{array}
       $
     & $\begin{array}{c}
          \\
          3, \\
          4, \\
          2
       \end{array}$ 
     & 12                  \\
  \hline        
       $\widetilde{E}_8$
     &
       $\begin{array}{c}
        \lambda^1_{1,2} = 1, \\
        \lambda^2_{1,2} = e^{\pm{2\pi{i}/3}}, \\
        \lambda^3_{1,2} = e^{\pm{2\pi{i}/5}}, \hspace{2mm}
        \lambda^3_{1,2} = e^{\pm{4\pi{i}/5}}, \\
        \lambda_9 = -1
       \end{array}
       $
     & $\begin{array}{c}
          \\
          3, \\
          5, \\
          2
       \end{array}$ 
     & 30   \\
  \hline
        $\widetilde{D}_n$
     &
       $\begin{array}{c}
        \lambda^1_{1,2} = 1, \\
        \lambda^2_s = e^{{2s}\pi{i}/(n-2)}, s = 1,2,...,n-3, \\
        \lambda_n = \lambda_{n+1} = -1
       \end{array} $
     & $\begin{array}{c}
          \\
          n-2, \\
          2,2 \\
       \end{array}$ 
     & $\begin{array}{c}
         n-2 \text{ for } n = 2k;   \\
         2(n-2) \text{ for } n = 2k+1
       \end{array}$ \\
  \hline
        $\widetilde{G}_{21},
        \widetilde{G}_{22}$
     &  $\begin{array}{c}
          \lambda^1_{1,2} = 1, \\
          \lambda_3 = -1
        \end{array}$
     & $\begin{array}{c}
          \\
          2 \\
       \end{array}$ 
     & 2 \\
  \hline
        $\widetilde{F}_{41},
        \widetilde{F}_{42} $
     &  $\begin{array}{c}
          \lambda^1_{1,2} = 1, \\
          \lambda^2_{1,2} = e^{{2\pi{i}/3}}, \\
          \lambda_3 = -1
        \end{array}$
     & $\begin{array}{c}
          \\
          3, \\
          2 \\
       \end{array}$
     & 6 \\
  \hline
        $\widetilde{A}_{11},
        \widetilde{A}_{12}$
     & $\lambda^1_{1,2} = 1$
     & 1
     & 1 \\
   \hline
      $\begin{array}{c}
        \widetilde{B}_{n},
        \widetilde{C}_{n}, \\
        \widetilde{BC}_{n}
      \end{array}$
          &
       $\begin{array}{c}
        \lambda^1_{1,2} = 1, \\
        \displaystyle \lambda^2_s = e^{{2s}\pi{i}/n}, s = 1,2,...,n-1
       \end{array} $
     &   $\begin{array}{c}
          \\
          n \\
       \end{array}$ 
     & $n$ \\
   \hline
       $\begin{array}{c}
        \widetilde{CD}_n, \\
        \widetilde{DD}_n
        \end{array} $
     &
      $\begin{array}{c}
        \lambda^1_{1,2} = 1, \vspace{1mm} \\
        \lambda^2_s = e^{{2s}\pi{i}/(n-1)}, s = 1,2,...,n-2 \\
        \lambda_{n+1} = -1
       \end{array} $
     & $\begin{array}{c}
          \\
          n-1, \\
          2 \\
       \end{array}$ 
     & $\begin{array}{c}
         n-1 \text{ for } n = 2k-1; \\
         2(n-1) \text{ for } n = 2k
       \end{array}$ \\
  \hline
\end{tabular}
  \label{table_eigenvalues}
\end{table}

\footnotetext[1]{See Remark \ref{orders_explanations}.}
\footnotetext[2]{See Remark \ref{affine_coxeter_number}.}
 \index{Lie algebra}

\begin{table} 
  \centering
  \vspace{2mm}
  \caption{\hspace{3mm}The Coxeter numbers and affine Coxeter numbers}
  \renewcommand{\arraystretch}{1.5}
  \begin{tabular} {|c|c|c|c|c|}
  \hline \hline
        Extended & Notation in         &  Affine       &             & Dual \footnotemark[2]  \cr
        Dynkin   & context of twisted  & Coxeter       & Coxeter     & Coxeter  \cr
        Diagram  & affine Lie algebra
                   \footnotemark[1]    &  number $h_a$ & number $h$  & number $h^\vee$ \\
  \hline  \hline  
       $\widetilde{E}_6$ & ${E}^{(1)}_6$ & 6 & 12  & 12 \\
  \hline          
       $\widetilde{E}_7$ & ${E}^{(1)}_7$ & 12 & 18 &  18 \\
  \hline          
       $\widetilde{E}_8$ & ${E}^{(1)}_8$ & 30 & 30 & 30 \\
  \hline
        $\widetilde{D}_n$ & ${D}^{(1)}_n$
     & $\begin{array}{c}
         n-2 \text{ for }  n = 2k;    \\
         2(n-2) \text{ for } n = 2k+1
       \end{array}$
     & $2(n-1)$ & $2(n-1)$ \\
  \hline
        $\widetilde{A}_{11}$ & ${A}^{(2)}_2$ & 1 & 3 & 3\\
  \hline
        $\widetilde{A}_{12}$ & ${A}^{(1)}_1$ & 1 & 2 & 2\\
  \hline
        $\widetilde{BC}_{n}$ & ${A}^{(2)}_{2n}$ & $n$ & $2n+1$ & $2n+1$\\
  \hline \hline
        $\widetilde{G}_{22}$ = $\widetilde{G}^{\vee}_{21}$ &
        ${G}^{(1)}_2$ & 2 & 6 & \fbox{4} \\ 
  \hline
        $\widetilde{G}_{21}$ = $\widetilde{G}^{\vee}_{22}$ &
        ${D}^{(3)}_4$ & 2 & \fbox{4} & 6 \\
  \hline \hline
        $\widetilde{F}_{42}$ = $\widetilde{F}^{\vee}_{41}$ &
        ${F}^{(1)}_4$ & 6 & 12 & \fbox{9} \\
  \hline
        $\widetilde{F}_{41}$ = $\widetilde{F}^{\vee}_{42}$ &
        ${E}^{(2)}_6$ & 6 & \fbox{9} & 12 \\
  \hline \hline
        $\widetilde{C}_{n}$ = $\widetilde{B}^{\vee}_{n}$ &
        ${C}^{(1)}_n$ & $n$  & $2n$ & \fbox{n+1} \\
  \hline
        $\widetilde{B}_{n}$ = $\widetilde{C}^{\vee}_{n}$ &
        ${D}^{(2)}_{n+1}$ & $n$  & \fbox{n+1} & $2n$ \\
  \hline \hline
      $\widetilde{CD}_n$ = $\widetilde{DD}^{\vee}_n$
      & ${B}^{(1)}_{n}$
      & $\begin{array}{c}
         n-1 \text{ for }  n = 2k-1;    \\
         2(n-1) \text{ for } n = 2k
       \end{array}$
      & $2n$ & \fbox{2n-1} \\
  \hline
     $\widetilde{DD}_n$ = $\widetilde{CD}^{\vee}_n$ &  ${A}^{(2)}_{2n-1}$
      & $\begin{array}{c}
         n-1 \text{ for }  n = 2k-1;    \\
         2(n-1) \text{ for } n = 2k
       \end{array}$
      & \fbox{2n-1} & $2n$ \\
  \hline
\end{tabular}
  \label{table_Coxeter_numbers}
\end{table}
\index{Coxeter number} \index{affine Coxeter number} \index{root
system}

\begin{remark}
\label{affine_coxeter_number}
 {\rm \indent Let $H$ be the
hyperplane orthogonal to the adjoint vector $\tilde{z}^1$ from
Proposition \ref{eigenvectors_basis}. The vector $\tilde{z}^1$ is
responsible for a $2\times2$ block in the Jordan form of the
affine Coxeter transformation. Due to the presence of a $2\times2$
block, the affine Coxeter transformation is of infinite order in
the Weyl group. The restriction of the Coxeter transformation on
the hyperplane $H$ is, however, of a finite order $h_a$. We call
this number the {\it affine Coxeter number}. The affine Coxeter
number $h_a$ is the least common multiple of orders of eigenvalues
($\neq 1$) of the affine Coxeter transformation (Table
\ref{table_eigenvalues}). We denote by $h$ the {\it Coxeter
number} of the Dynkin diagram. The value $h-1$ is the sum of
coordinates of the highest root $\beta$ of the corresponding root
system, see (\ref{coxeter_result_3}). The imaginary vector $z^1$
depicted on Fig.~\ref{euclidean_diag} coincides with the highest
root $\beta$ extended to the vector with $1$ at the additional
vertex, the one that extends the Dynkin diagram to the extended
Dynkin diagram. Thus, \footnotetext[1]{See,Remark \ref{upper
index_r}.} \footnotetext[2]{For emphasis and to distinguish $h$
from $h^{\vee}$, we frame the value $\min(h, h^{\vee})$ in the
case $h \neq h^{\vee}$.}
\begin{equation}
\label{h_Coxeter_number}
    h = n_1 + \dots + n_l,
\end{equation}
where $n_i$ coordinates of the imaginary vector from
Fig.~\ref{euclidean_diag}. Similarly to (\ref{h_Coxeter_number}),
the {\it dual Coxeter number } is defined as
\begin{equation}
\label{h_Coxeter_number_2}
    h^{\vee} = n^{\vee}_1 + \dots + n^{\vee}_l,
\end{equation}
where $\{n^{\vee}_1,\dots,n^{\vee}_l \}$ are coordinates of the
imaginary vector $h^{\vee}$ of the dual Dynkin diagram,
Fig.~\ref{euclidean_diag}. For values of the dual Coxeter numbers,
see Table \ref{table_Coxeter_numbers}.}
\end{remark}

 \index{twisted affine Lie algebra}
 \index{Lie algebra}
\begin{remark}
  \label{upper index_r}
  {\rm
For the first time, the notation of {\it twisted affine Lie
algebras } from col. 2, Table \ref{table_Coxeter_numbers} appeared
in \cite{Kac69} in the description of finite order automorphisms;
see also \cite{Kac80}, \cite[p.123]{GorOnVi94}, \cite{OnVi90},
\cite{Kac93}. The upper index $r$ in the notation of twisted
affine Lie algebras has an invariant sense: is the order of the
diagram automorphism $\mu$ of $\mathfrak{g}$, where $\mathfrak{g}$
is a complex simple finite dimensional Lie algebra of type $X_N =
A_{l}, D_{l+1}, E_6, D_4$, \cite[Th.8.3]{Kac93}.

The affine Lie algebra associated to a generalized Cartan matrix
of type ${X}^{(1)}_l$  is called a {\it non-twisted affine Lie
algebra}, \cite[Ch.7]{Kac93}.

The affine Lie algebras associated to a generalized Cartan matrix
of type ${X}^{(2)}_l$  and ${X}^{(3)}_l$  are called  {\it twisted
affine Lie algebras}. \cite[Ch.8]{Kac93}.

The corresponding $\mathbb{Z}/r\mathbb{Z}$-gradings are (here
$\bar i\in \mathbb{Z}/r\mathbb{Z}$):
\begin{equation}
  \begin{split}
   & \mathfrak{g} = \mathfrak{g}_{\overline{0}} + \mathfrak{g}_{\overline{1}},
                     \hspace{3mm} \text{for } r=2,
   \\
   & \mathfrak{g} = \mathfrak{g}_{\overline{0}} + \mathfrak{g}_{\overline{1}} +
   \mathfrak{g}_{\overline{2}}, \hspace{3mm}  \text{for } r = 3,
  \end{split}
\end{equation}
see \cite[8.3.1, 8.3.2]{Kac93}, \cite[Ch.3, \S3]{GorOnVi94}
  }
\end{remark}

 \begin{proposition} {\em\bf (\cite[exs.6.3]{Kac93}) }
 \label{prop_lrh}
    Let $A$ be a Cartan matrix of type ${X}^{(r)}_N$ from
    col.2, Table \ref{table_Coxeter_numbers},
    let $l$ = rank of $A$,
    and let $h$ be the Coxeter number. Let $\varDelta$
    be the finite root system of type ${X}_N$. Then
  \begin{equation}
   \label{rlh_relation}
    rlh = |\varDelta|.
  \end{equation}
 \end{proposition}
  \PerfProof
     For $r = 1$, we have H.~S.~M.~Coxeter's proposition (\ref{coxeter_result_1}),
     \S\ref{history_1}. We consider only remaining cases:
 \begin{equation*}
    A^{(2)}_2, \hspace{3mm}
    A^{(2)}_{2n} (l \geq 2), \hspace{3mm}
    A^{(2)}_{2n-1} (l \geq 3), \hspace{3mm}
    D^{(2)}_{n+1} (l \geq 2), \hspace{3mm}
    E^{(2)}_6, \hspace{3mm}
    D^{(3)}_4.
 \end{equation*}
     (see col. 5, Table \ref{table_Coxeter_numbers}).

     (a) $A^{(2)}_2$, $r = 2$;  rank $l = 1$; $h = 3$;
     $X_N = {A}_{2}$,
     $|\varDelta({A}_{2})| = 6$, see \cite[Tab.I]{Bo}.

     (b) ${A}^{(2)}_{2n}$, $r = 2$; rank $l = n$; $h = 2n+1$,
     $X_N = {A}_{2n}$,
     $|\varDelta({A}_{2n})| = 2n(2n+1)$, see \cite[Tab.I]{Bo}.

     (c) ${A}^{(2)}_{2n-1}$, $r = 2$; rank $l = n$; $h = 2n-1$;
     $X_N = {A}_{2n-1}$,
     $|\varDelta({A}_{2n-1})| = 2n(2n-1)$, see \cite[Tab.I]{Bo}.

     (d) ${D}^{(2)}_{n+1}$, $r = 2$;  rank $l = n$; $h = n+1$,
     $X_N = {D}_{n+1}$,
     $|\varDelta({D}_{n+1})| = 2n(n+1)$, see \cite[Tab.IV]{Bo}.

     (e) ${E}^{(2)}_6$, $r = 2$; rank $l = 4$; $h = 9$,
     $X_N = {E}_6$,
     $|\varDelta({E}_6)| = 72$, see \cite[Tab.V]{Bo}.

     (f) ${D}^{(3)}_4$, $r = 3$; rank $l = 2$;  $h = 4$,
     $X_N = {D}_4$,
     $|\varDelta({D}_4)| = 24$, see \cite[Tab.IV]{Bo}.

     Cases (a)-(f) are collected in Table \ref{table_rlh_relation}.
    \qed

\begin{table} 
  \centering
  \vspace{2mm}
  \caption{\hspace{3mm}Proposition \ref{prop_lrh} for $r = 2,3$.}
  \renewcommand{\arraystretch}{1.5}
  \begin{tabular} {|c|c|c|c|c|c|c|}
  \hline \hline
        Extended         & Index & Rank of  & Coxeter     & Root   & $|\varDelta| = |\varDelta(X_N)|$ \cr
        Dynkin           & $r$   &   $A$    & number      & System       &               \cr
        Diagram          &       &          & $h$         & $X_N$ &        \\
  \hline  \hline  
        ${A}^{(2)}_2$ & $2$  & $1$     & $3$        & $A_2$ & $6$ \\
  \hline
        ${A}^{(2)}_{2n}$ & $2$  & $n$   & $2n+1$    & $A_{2n}$ & $2n(2n+1)$ \\
  \hline
        ${A}^{(2)}_{2n-1}$ & $2$  & $n$   & $2n-1$  & $A_{2n-1}$ & $2n(2n-1)$ \\
  \hline
        ${D}^{(2)}_{n+1}$ & $2$ & $n$   & $2n+1$    & $D_{n+1}$ & $2n(2n+1)$ \\
  \hline
        ${E}^{(2)}_6$ & $2$  & $4$     & $9$        & $E_6$ & $72$ \\
  \hline
        ${D}^{(3)}_4$ & $3$  & $2$     & $4$        & $D_4$ & $24$ \\
  \hline
\end{tabular}
  \label{table_rlh_relation}
\end{table}

\section{Bibliographical notes on the spectrum of the Coxeter transformation}
\label{biblio_notes}

The eigenvalues of affine Coxeter transformations were also
calculated by S.~Berman, Y.S.~Lee and R.~Moody in [BLM89]. Theorem
\ref{th_jordan} was also proved by N.~A'Campo \cite{A'C76} and
R.~Howlett \cite{How82}.

Natural difficulties in the study of Cartan matrices and Coxeter transformations
for the graphs containing cycles are connected with the following two facts:

\medskip

1) these graphs have non-symmetrizable Cartan matrices,

2) in general, there are several conjugacy classes of the
Coxeter transformation.

\medskip

We would like to distinguish several works related
to the Coxeter transformation for the graphs with cycles:
C.~M.~Ringel \cite{Rin94}, A.~J.~Coleman \cite{Col89}, Shi Jian-yi \cite{Shi00},
Menshikh and Subbotin \cite{MeSu82}, \cite{Men85}, Boldt and Takane \cite{BT97}.

\index{spectral radius} For generalized Cartan matrices (see
\S\ref{cartan}), i.e., for graphs with cycles, C.~M.~Ringel
\cite{Rin94} showed that the {\it spectral radius} $\rho({\bf C})$
of the Coxeter transformation lies out the unit circle, and is an
eigenvalue of multiplicity one. This result  generalizes Theorem
\ref{th_jordan}, 3), proved in \cite{SuSt75, SuSt78, St85} only
for trees. The spectral radius $\rho({\bf C})$ is used by V.~Dlab
and C.~M.~Ringel to determine the Gelfand-Kirillov dimension of
the preprojective algebras, \cite{DR81}.

\index{Killing polynomials} A.~J.~Coleman \cite{Col89} computed
characteristic polynomials for the Coxeter transformation for all
extended Dynkin diagrams, including the case with cycles
$\widetilde{A}_n$. He baptized these polynomials {\it Killing
polynomials}, (see the epigraph to this chapter). Coleman also
shows that $\widetilde{A}_n$ has $\displaystyle\frac{n}{2}$
spectral conjugacy classes.

V.~V.~Menshikh and V.~F.~Subbotin in \cite{MeSu82}, and
V.~V.~Menshikh in \cite{Men85} established a connection between an
orientation $\varDelta$ of the graph $\varGamma$ and spectral
classes of conjugacy of the Coxeter transformation. For any
orientation ${\varDelta}$ of a given graph $\varGamma$ containing
several cycles, they consider an invariant $R_{\varDelta}$ equal
to the number of arrows directed in a clockwise direction. For any
graph $\varGamma$ containing disjoint cycles, they show that
$R_{\varDelta_1} = R_{\varDelta_2}$ if and only if orientations
$\varDelta_1$ and $\varDelta_2$ can be obtained from each other by
applying a sink-admissible or a source-admissible sequence of
reflections $\sigma_i$, see \S\ref{orientation}, i.e.,
\index{index of the conjugacy class}
$$
   R_{\varDelta_1} =  R_{\varDelta_2}
   \hspace{3mm} \Longleftrightarrow \hspace{3mm}
   \varDelta_1 = \sigma_{i_n}...\sigma_{i_2}\sigma_{i_1}(\varDelta_2)
$$
for any sink-admissible or source-admissible sequence
$i_1,i_2,...,i_n$. Menshikh and Subbotin also showed that two
Coxeter transformations ${\bf C}_{\varDelta_1}$ and ${\bf
C}_{\varDelta_2}$ are conjugate if and only if
\text{$R_{\varDelta_1} = R_{\varDelta_2}$}. The number
$R_{\varDelta}$ is called the {\it index} of the conjugacy class
of the Coxeter transformation. Menshikh and Subbotin  also
calculated the characteristic polynomial of the Coxeter
transformation for every class equivalent to $\varDelta$ for the
extended Dynkin diagram $\widetilde{A}_n$; this polynomial is
\begin{equation}
  \label{char_polyn_An}
     \det|{\bf C} - \lambda{I}| =
     \lambda^{n+1} - \lambda^{n-k+1} - \lambda^{k} + 1,
\end{equation}
where $k = R_{\varDelta}$ is the index of the conjuagacy class of
the Coxeter transformation.

Shi Jian-yi \cite{Shi00} considers conjugacy relation on Coxeter
transformations for the case where $\varGamma$ is just a cycle.
Different equivalence relations on Coxeter transformations are
considered in more difficult cases.  In \cite{Shi00}, Shi also
obtained an explicit formula (\ref{char_polyn_An}).

\index{unicyclic graph} \index{essential cycle} In \cite{Bol96},
\cite{BT97}, Boldt and Takane considered unicyclic graphs.
\begin{definition} [\cite{Bol96}] {\rm
  An {\it essential cycle} of $\varGamma$ is a full subgraph $S$ of
  $\varGamma$ having $m \geq 3$ vertices $x_1,\dots,x_m$
  such that there are edges between $x_i$ and
  $x_{i+1}$ for $i = 1,\dots,m-1$ and also between $x_m$ and
  $x_1$ (thus, $\varGamma$ coincides with $\widetilde{A}_n$). A graph $\varGamma$ is called
  the {\it unicyclic} if it contains precisely one {\it essential cycle}.
  }
 \end{definition}
Boldt and Takane showed how the characteristic polynomial of the
Coxeter transformation for the unicyclic graph $\varGamma$ can be
reduced to the characteristic polynomial for the essential cycle.
They also arrive at the explicit formula (\ref{char_polyn_An}).

\section{Shearing and gluing formulas for the characteristic polynomial}
\label{shearing} The purpose of this section is to prove the
Subbotin-Sumin {\it splitting along the edge} formula for the
characteristic polynomial of the Coxeter transformation
\cite{SuSum82}, to prove its generalization for the multiply-laced
case and to get some of its corollaries that will be used in the
following sections.

\begin{figure}[h]
\centering
\includegraphics{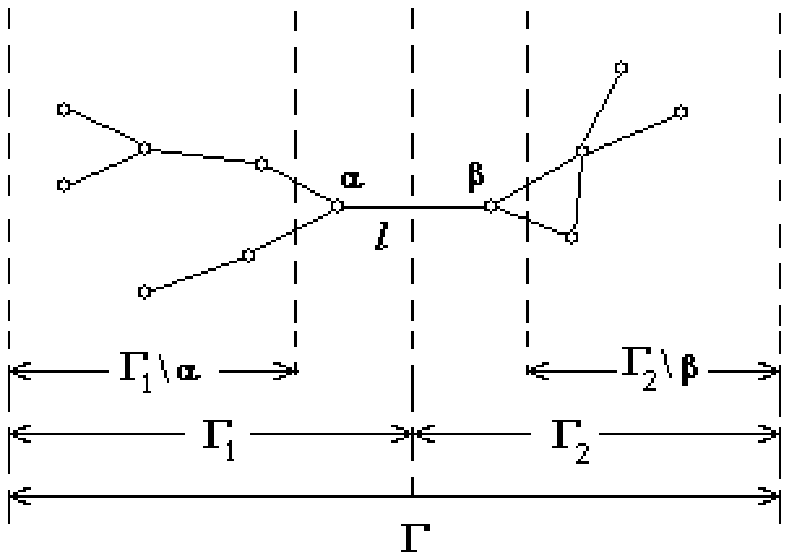}
\caption{\hspace{3mm} A split graph $\varGamma$}
\label{splitted_diagram}
\end{figure}

Let us consider a  characteristic polynomial of the Coxeter transformation
$$
    \mathcal{X}(\varGamma, \lambda) = \det|\varGamma - \lambda{I}|
$$
for the graph $\varGamma$ with a {\it splitting edge}. Recall that
an edge  $l$ is said to be {\it splitting} if deleting it we split
the graph $\varGamma$ into two graphs $\varGamma_1$ and
$\varGamma_2$, see Fig.~\ref{splitted_diagram}.
\index{Subbotin-Sumin formula! - simply-laced case}
\index{formula! - Subbotin-Sumin, simply-laced case}
\begin{proposition} [\cite{SuSum82}]
  \label{split_edge}
  For a given graph $\varGamma$ with a splitting edge $l$, we have
\begin{equation}
  \label{split_edge_eq}
   \mathcal{X}(\varGamma, \lambda) =
       \mathcal{X}(\varGamma_1, \lambda)\mathcal{X}(\varGamma_2, \lambda) -
       \lambda\mathcal{X}(\varGamma_1\backslash\alpha, \lambda)
              \mathcal{X}(\varGamma_2\backslash\beta, \lambda),
\end{equation}
where $\alpha$ and $\beta$ are the endpoints of the deleted edge
$l$.
\end{proposition}
\PerfProof We select a basis in $\mathcal{E}_\varGamma$ such that
\begin{displaymath}
 \begin{split}
   & \mathcal{E}_\varGamma = \mathcal{E}_{\varGamma_1} + \mathcal{E}_{\varGamma_2}, \\
   & \mathcal{E}_{\varGamma_1} = \{a_1,...,a_m\}, \hspace{3mm}
       a_m = \alpha, \hspace{3mm}
       \dim\mathcal{E}_{\varGamma_1} = m,    \\
   & \mathcal{E}_{\varGamma_2} = \{b_1,...,b_k\}, \hspace{3mm}
       b_1 = \beta,  \hspace{3mm}
       \dim\mathcal{E}_{\varGamma_2} = k,
 \end{split}
\end{displaymath}
and a Coxeter transformation {\bf C} such that

\begin{displaymath}
       {\bf C} = {w}_{\varGamma_1}{w}_{\varGamma_2}, \hspace{3mm}
       {w}_{\varGamma_1} = \sigma_{a_m}...\sigma_{a_1}, \hspace{3mm}
       {w}_{\varGamma_2} = \sigma_{b_1}...\sigma_{b_k} \hspace{0.5mm}.
\end{displaymath}
An ordering of reflections $\sigma_{a_i}$ (resp. $\sigma_{b_j}$)
in $w_{\varGamma_1}$ (resp. $w_{\varGamma_2}$) is selected so that
$\sigma_{a_m}$ (resp. $\sigma_{b_1}$) is the last to be executed.
Then the reflections $\sigma_{a_i}$ ($i \neq m$) (resp.
$\sigma_{b_j}$ ($j \neq 1$)) do not act on the coordinates $b_j$
(resp. $a_i$).  Only $\sigma_{a_m}$ acts on the coordinate $b_1$
and $\sigma_{b_1}$ acts on the coordinate $a_m$. We have
\begin{equation}
  \label{w_acts}
\begin{array}{l}
{w}_{\varGamma_1} = \sigma_{a_m}...\sigma_{a_1} = \left (
      \begin{array}{cc}
         {\bf C}_{\varGamma_1} & \delta_{{a_m}{b_1}} \vspace{3mm} \\
          0  &  I
      \end{array}
         \right )
      \begin{array}{c}
        \}m
        \vspace{3mm} \\
        \}k \\
      \end{array}, \vspace{5mm} \\
    {w}_{\varGamma_2} = \sigma_{b_1}...\sigma_{b_k}  = \left (
      \begin{array}{cc}
          I & 0  \vspace{3mm} \\
          \delta_{{b_1}{a_m}} & {\bf C}_{\varGamma_2}
      \end{array}
         \right )
      \begin{array}{c}
        \}m
        \vspace{3mm} \\
        \}k \\
      \end{array},\end{array}
\end{equation}
where ${\bf C}_{\varGamma_1}$ (resp. ${\bf C}_{\varGamma_2}$) is
the Coxeter transformation for the graph $\varGamma_1$ (resp.
$\varGamma_2$) and $\delta_{{a_m}{b_1}}$ (resp.
$\delta_{{b_1}{a_m}}$) means the matrix with zeros everywhere
except 1 in the (${a_m}{b_1}$)th (resp. (${b_1}{a_m}$)th) slot.

\begin{figure}[h]
\centering
\includegraphics{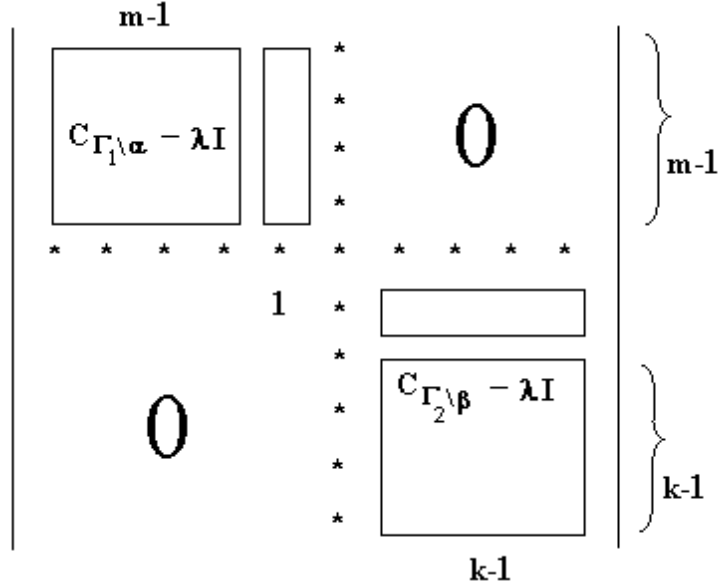}
\caption{\hspace{3mm} The asterisks mark the deleted line $a_m$
and column $b_1$}
\label{expanded_minor}
\end{figure}

Then we have
\begin{equation}
 \label{C_delta_1}
  \begin{array}{c}
   {\bf C} =
     {w}_{\varGamma_1}{w}_{\varGamma_2} = \left (
      \begin{array}{cc}
         {\bf C}_{\varGamma_1} + \delta_{{a_m}{a_m}} &
               \delta_{{a_m}{b_1}}{\bf C}_{\varGamma_2} \vspace{3mm} \\
          \delta_{{b_1}{a_m}}  & {\bf C}_{\varGamma_2}
                \end{array}
         \right ), \vspace{5mm} \\
    \det|{\bf C} - \lambda{I}| =  \left |
      \begin{array}{cc}
         {\bf C}_{\varGamma_1} + \delta_{{a_m}{a_m}} - \lambda{I} &
               \delta_{{a_m}{b_1}}{\bf C}_{\varGamma_2} \vspace{3mm} \\
          \delta_{{b_1}{a_m}}  & {\bf C}_{\varGamma_2} - \lambda{I} \vspace{3mm}
      \end{array}
             \right |.
  \end{array}
\end{equation}
Here, $\delta_{{a_m}{b_1}}{\bf C}_{\varGamma_2}$ is the $b_1$st
line
  (the first line in ${\bf C}_{\varGamma_2}$ for the selected basis)
  of the
  matrix ${\bf C}_{\varGamma_2}$ and we subtract now this line from
  the line $a_m$. Then
\begin{displaymath}
    \det|{\bf C} - \lambda{I}| =  \left |
      \begin{array}{cc}
         {\bf C}_{\varGamma_1} - \lambda{I} &
               \lambda\delta_{{a_m}{b_1}} \vspace{3mm} \\
          \delta_{{b_1}{a_m}}  & {\bf C}_{\varGamma_2} - \lambda{I} \vspace{3mm}
      \end{array}
             \right |
\end{displaymath}
or
\begin{equation}
  \label{C_with_lambda}
    \det|{\bf C} - \lambda{I}| =  \left |
      \begin{array}{cc}
         {\bf C}_{\varGamma_1} - \lambda{I} &
             \begin{array}{cccc}
               0 & 0  & \dots  & 0  \\
               0 & 0  & \dots  & 0  \\
                 &    & \dots  &    \\
               0 & 0  & \dots  & 0  \\
               \lambda & 0  & \dots  & 0
             \end{array} \vspace{3mm} \\
            \begin{array}{cccc}
               0 & \dots & 0 & 1  \\
               0 & \dots & 0 & 0  \\
                 & \dots &   &    \\
               0 & \dots & 0 & 0  \\
               0 & \dots & 0 & 0
             \end{array}
            & {\bf C}_{\varGamma_2} - \lambda{I} \vspace{3mm}
      \end{array}
             \right |.
\end{equation}
Let us expand the determinant $\det|{\bf C} - \lambda{I}|$ in
(\ref{C_with_lambda}) with respect to the minors corresponding to
the line $b_1$ (the line containing $\lambda$ in the column $b_1$,
which is the next after $a_m$). We see that
\begin{equation}
   \det|{\bf C} - \lambda{I}| =
     \det|{\bf C}_{\varGamma_1} - \lambda{I}|\hspace{0.5mm}
     \det|{\bf C}_{\varGamma_2} - \lambda{I}| +
       (-1)^{a_m + (a_m + 1)}\lambda{R},
\end{equation}
where $R$ is the determinant shown on the Fig.~\ref{expanded_minor}.
Expanding the determinant $R$ we get the Subbotin-Sumin
formula (\ref{split_edge_eq}).
\hspace{3mm} \qedsymbol

\index{Subbotin-Sumin formula! - multiply-laced case}
\index{formula! - Subbotin-Sumin, multiply-laced case} Now we will
extend the Subbotin-Sumin formula to the multiply-laced case where
the endpoints $\alpha$ and $\beta$ correspond to simple roots of
different lengths. In this case, formula (\ref{split_edge_eq}) is
modified to (\ref{split_edge_eq_2}).

\begin{proposition}
  \label{split_edge_2}
  For a given graph $\varGamma$ with a splitting weighted edge $l$
  corresponding to roots with different lengths, we have
\begin{equation}
  \label{split_edge_eq_2}
   \mathcal{X}(\varGamma, \lambda) =
       \mathcal{X}(\varGamma_1, \lambda)\mathcal{X}(\varGamma_2, \lambda) -
       \rho\lambda\mathcal{X}(\varGamma_1\backslash\alpha, \lambda)
              \mathcal{X}(\varGamma_2\backslash\beta, \lambda),
\end{equation}
  where $\alpha$ and $\beta$ are the endpoints
  of the deleted edge $l$ and $\rho$ is the ratio of lengths
  of the roots corresponding to these endpoints.
\end{proposition}
\PerfProof
Relation (\ref{C_delta_1}) is modified to the following
\begin{equation*}
   {\bf C} =
      \left (
      \begin{array}{cc}
         {\bf C}_{\varGamma_1} + \rho\delta_{{a_m}{a_m}} &
               \rho\delta_{{a_m}{b_1}}{\bf C}_{\varGamma_2} \vspace{3mm} \\
          \delta_{{b_1}{a_m}}  & {\bf C}_{\varGamma_2}
                \end{array}
         \right ).
\end{equation*}
Multiply the line $b_1$ by $\rho$ and subtract this product from
the line $a_m$, then we obtain
\begin{equation}
 \label{C_delta_2}
    \det|{\bf C} - \lambda{I}| =  \left |
      \begin{array}{cc}
         {\bf C}_{\varGamma_1} + \rho\delta_{{a_m}{a_m}} - \lambda{I} &
               \rho\delta_{{a_m}{b_1}}{\bf C}_{\varGamma_2} \vspace{3mm} \\
          \delta_{{b_1}{a_m}}  & {\bf C}_{\varGamma_2} - \lambda{I} \vspace{3mm}
      \end{array}
             \right | =  \left |
      \begin{array}{cc}
         {\bf C}_{\varGamma_1}  - \lambda{I} &
               \rho\lambda\delta_{{a_m}{b_1}} \vspace{3mm} \\
          \delta_{{b_1}{a_m}}  & {\bf C}_{\varGamma_2} - \lambda{I} \vspace{3mm}
      \end{array}
             \right |
\end{equation}
Further, as in Proposition \ref{split_edge}, we obtain (\ref{split_edge_eq_2}).
\qedsymbol

\begin{corollary}
  \label{split_edge_corol_2}
  Let $\varGamma_2$ in Proposition \ref{split_edge_2}
  be a component containing only one point. Then, the
  following formula holds
\begin{equation}
  \label{split_edge_eq_corol_2}
   \mathcal{X}(\varGamma, \lambda) =
       -(\lambda + 1)\mathcal{X}(\varGamma_1, \lambda) -
       \rho\lambda\mathcal{X}(\varGamma_1\backslash\alpha, \lambda),
\end{equation}
\end{corollary}
\PerfProof
Relation (\ref{w_acts}) is modified to the following
\begin{equation}
\begin{array}{l}
w_{\varGamma_1} = \left (
      \begin{array}{cc}
         {\bf C}_{\varGamma_1} &
             \begin{array}{c}
               0 \\
               0 \\
               \dots \\
               0 \\
               \rho \vspace{3mm} \\
             \end{array} \\
          0  &  1
      \end{array}
         \right )
      \begin{array}{c}
       \left . \begin{array}{c}
                ~\\
                ~\\
                ~\\
                ~\\
                ~\\
               \end{array} \right \}m
        \vspace{3mm} \\
         \hspace{3mm}\}\hspace{1mm}1 \\
      \end{array}, \vspace{5mm} \\
    w_{\varGamma_2} = \left (
      \begin{array}{cc}
          I & 0  \vspace{3mm} \\
          0~0\dots1 & -1
      \end{array}
         \right )
      \begin{array}{c}
        \}m
        \vspace{3mm} \\
        \}1 \\
      \end{array},
\end{array}
\end{equation}
and (\ref{C_delta_2}) is modified to
\begin{equation}
 \label{C_delta_corol_2}
  \begin{array}{c}
   {\bf C} =
      \left (
      \begin{array}{cc}
         {\bf C}_{\varGamma_1} + \rho\delta_{{a_m}{a_m}} &
             \begin{array}{c}
               0 \\
               0 \\
               \dots \\
               0 \\
               -\rho \vspace{3mm} \\
             \end{array} \\
          0~0 \dots 1  & -1
                \end{array}
         \right ), \vspace{5mm} \\
    \det|{\bf C} - \lambda{I}| =  \left |
      \begin{array}{cc}
         {\bf C}_{\varGamma_1} + \rho\delta_{{a_m}{a_m}} - \lambda{I} &
             \begin{array}{c}
               0 \\
               0 \\
               \dots \\
               0 \\
               -\rho \vspace{3mm} \\
             \end{array} \\
          0~0 \dots 1  & -1 - \lambda \vspace{3mm}
      \end{array}
             \right | =  \left |
      \begin{array}{cc}
         {\bf C}_{\varGamma_1} - \lambda{I} &
             \begin{array}{c}
               0 \\
               0 \\
               \dots \\
               0 \\
               \rho\lambda \vspace{3mm} \\
             \end{array} \\
          0~0 \dots 1  & -1 - \lambda \vspace{3mm}
      \end{array}            \right |
  \end{array}
\end{equation}
Further, as in Proposition \ref{split_edge}, we obtain
(\ref{split_edge_eq_corol_2}).
\qedsymbol

The next proposition follows from formula (\ref{split_edge_eq})
and allows one to calculate the spectrum of the graph
$\varGamma(n)$ obtained by gluing $n$ copies of the graph
$\varGamma$.

\begin{figure}[h]
\centering
\includegraphics{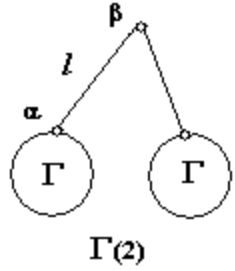}
\caption{\hspace{3mm} Splitting along the edge $l$}
\label{split_along_l}
\end{figure}

\begin{proposition} [\cite{SuSum82}, \cite{KMSS83},
\cite{KMSS83a}]
 \label{gluing}
   Let $*_n$ be a star with $n$ rays coming from a vertex.
   Let $\varGamma(n)$ be the graph obtained from  $*_n$
by gluing to the endpoints of its rays $n$ copies
   of the graph $\varGamma$. Then
\begin{equation}
  \label{n_copies}
   \mathcal{X}(\varGamma(n), \lambda) =
       \mathcal{X}(\varGamma, \lambda)^{n-1}\varphi_{n-1}(\lambda)
\end{equation}
for some polynomial $\varphi_{n-1}(\lambda)$. The polynomials
$\varphi_{n-1}(\lambda)$ are recursively calculated.
\end{proposition}
\PerfProof By splitting along the edge $l$ (\ref{split_edge_eq})
  we get,
  for $n = 2$,  (Fig. \ref{split_along_l})
\begin{displaymath}
  \begin{split}
   & \mathcal{X}(\varGamma(2), \lambda) = \\
   &  \mathcal{X}(\varGamma, \lambda)\mathcal{X}(\varGamma + \beta, \lambda) -
    \lambda\mathcal{X}(\varGamma\backslash\alpha, \lambda)
    \mathcal{X}(\varGamma, \lambda) = \\
   & \mathcal{X}(\varGamma, \lambda)[\mathcal{X}(\varGamma + \beta, \lambda)
       - \lambda\mathcal{X}(\varGamma\backslash\alpha, \lambda)]  = \\
   & \mathcal{X}(\varGamma, \lambda)\varphi_1(\lambda), \hspace{2mm}
     {\rm where} \hspace{2mm} \varphi_1(\lambda) = \\
   & \mathcal{X}(\varGamma + \beta, \lambda)
       - \lambda\mathcal{X}(\varGamma\backslash\alpha, \lambda).
  \end{split}
\end{displaymath}

Let the proposition be true  for $n = r$ and
\begin{displaymath}
  \mathcal{X}(\varGamma(r), \lambda) =
       \mathcal{X}(\varGamma, \lambda)^{r-1}\varphi_{r-1}(\lambda).
\end{displaymath}
Then, for $n = r+1$, we have
\begin{displaymath}
 \begin{split}
  & \mathcal{X}(\varGamma(r+1), \lambda) =
       \mathcal{X}(\varGamma, \lambda)
       \mathcal{X}(\varGamma(r), \lambda)
       - \lambda\mathcal{X}(\varGamma\backslash\alpha, \lambda)
       (\mathcal{X}(\varGamma, \lambda))^r = \\
  &   \mathcal{X}(\varGamma, \lambda)
       \mathcal{X}(\varGamma, \lambda)^{r-1}\varphi_{r-1}(\lambda)
       - \lambda\mathcal{X}(\varGamma\backslash\alpha, \lambda)
     \mathcal{X}(\varGamma, \lambda)^{r}   = \\
  & \mathcal{X}(\varGamma, \lambda))^r
     [\varphi_{r-1}(\lambda) -
       \lambda\mathcal{X}(\varGamma\backslash\alpha, \lambda)].
       \hspace{3mm} \qed
  \end{split}
  \end{displaymath}

\begin{figure}[h]
\centering
\includegraphics{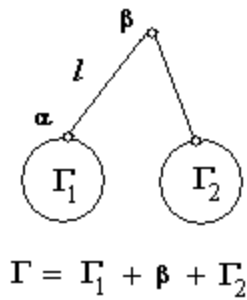}
\caption{\hspace{3mm} Splitting along the edge $l$, $\varGamma_1
\neq \varGamma_2$}
\label{split_along_l_diff_G}
\end{figure}

 The following proposition is due to V.~Kolmykov.
 For the case $\lambda$ = 1, it is formulated in
 \cite[Ch.4, exs.4.12]{Kac93}.
\begin{proposition}
  \label{merge_2G}
  If the spectrum of the Coxeter transformations
  for graphs $\varGamma_1$ and $\varGamma_2$ contains the eigenvalue $\lambda$,
  then this eigenvalue is also the eigenvalue of the graph $\varGamma$
  obtained by gluing as described in Proposition \ref{gluing}.
\end{proposition}
\PerfProof  By splitting along the edge $l$ (\ref{split_edge_eq})
  we get (Fig. \ref{split_along_l_diff_G}
\begin{displaymath}
  \begin{split}
    & \mathcal{X}(\varGamma_1 + \beta + \varGamma_2, \lambda) = \\
    & \mathcal{X}(\varGamma_1, \lambda)\mathcal{X}(\varGamma_2 + \beta, \lambda) -
    \lambda\mathcal{X}(\varGamma\backslash\alpha, \lambda)
    \mathcal{X}(\varGamma_2, \lambda).   \end{split}
\end{displaymath}
If $\mathcal{X}(\varGamma_1, \lambda) |_{\lambda = \lambda_0} = 0$
and $\mathcal{X}(\varGamma_2, \lambda) |_{\lambda = \lambda_0} =
0$, \\
then also $\mathcal{X}(\varGamma, \lambda) |_{\lambda = \lambda_0}
= 0.$ \hspace{3mm} \qed

\index{characteristic polynomial of the Coxeter transformation! -
${A}_n$} Let $\mathcal{X}({A}_n)$ be the characteristic polynomial
of the Coxeter transformation for the Dynkin diagram ${A}_n$.
Then, from the recurrent formula (\ref{split_edge_eq_corol_2}), we
have
\begin{equation}
 \label{recur_chi}
  \begin{split}
    \mathcal{X}({A}_1) = & -(\lambda + 1),  \\
    \mathcal{X}({A}_2) = & \lambda^2 + \lambda + 1,  \\
    \mathcal{X}({A}_3) = & -(\lambda^3 + \lambda^2 + \lambda + 1),  \\
    \mathcal{X}({A}_4) = & \lambda^4 + \lambda^3 + \lambda^2 + \lambda + 1,  \\
    \dots  \\
    \mathcal{X}({A}_n) = & -(\lambda + 1)\mathcal{X}({A}_{n-1}) -
           \lambda\mathcal{X}({A}_{n-2}),
\hspace{5mm} n > 2.
  \end{split}
\end{equation}
\begin{remark}
 \label{up_to_sign}
\rm{
  Since the characteristic polynomial can be chosen up to a factor, we prefer
  to consider the polynomial
\begin{equation}
  \label{def_chi}
   \chi_n = (-1)^n\mathcal{X}({A}_n)
\end{equation}
as the characteristic polynomial of the Coxeter transformation of
${A}_n$, thus the leading coefficient is positive. Pay attention
that the sign $(-1)^n$ for diagrams with $n$ vertices should be
taken into account only in recurrent calculations using formula
(\ref{split_edge_eq_corol_2}) (cf. with calculation for
$\tilde{D}_{n}$ and $\tilde{D}_{4}$ below). }
\end{remark}

\begin{remark}
\label{frame_formula_1}
\index{Frame formula}
  {\rm 1) J.~S.~Frame in \cite[p.784]{Fr51} obtained that
\begin{equation}
   \label{frame}
   \chi({A}_{m+n}) =
     \chi({A}_{m})\chi({A}_{n})
   - \lambda\chi({A}_{m-1})\chi({A}_{n-1}),
\end{equation}
 which easily follows from (\ref{split_edge_eq}).
 Another {\it Frame formula}
 will be proved in Proposition \ref{another_Frame}.

   2) S.~M.~Gussein-Zade in \cite{Gu76} obtained the formula
\begin{equation}
    \mathcal{X}({A}_n) = -(\lambda + 1)\chi_{n-1} - \lambda\chi_{n-2},
\end{equation}
  for the characteristic polynomial of the classical monodromy,
  see (\ref{recur_chi}).
}

\end{remark}

\index{roots of unity} From (\ref{recur_chi}) for $n > 0$, we see
that $\chi_n$ is
  the cyclotomic polynomial whose roots are the primitive $n$th
  roots of unity:
\begin{equation}
 \label{An_chi}
  \begin{split}
    \chi_n = \sum\limits_{i=1}^n\lambda^i =
             \frac{\lambda^{n+1} - 1}{\lambda - 1}.
  \end{split}
\end{equation}

\section{Formulas of the characteristic polynomials for the diagrams $T_{p,q,r}$}
  \label{section_T_pqr}
We give here explicit formulas of characteristic polynomials of
the Coxeter transformations for the three classes of diagrams:
$T_{2,3,r}$, $T_{3,3,r}$, and $T_{2,4,r}$, see \S\ref{hyperbolic},
Fig.~\ref{T_pqr_diagram}.

\subsection{The diagrams $T_{2,3,r}$}

 \index{characteristic polynomial of the Coxeter transformation! -
$T_{2,3,r}$}
 \index{${E}_n$-series}
 The case $T_{2,3,r}$, where $r
\geq 2$, includes diagrams ${D}_5$, ${E}_6$, ${E}_7$, ${E}_8$,
$\tilde{E}_8$, ${E}_{10}$.  Since ${D}_5 = T_{2,3,2} = {E}_5$ and
$\tilde{E}_8 = T_{2,3,6} = {E}_9$, we call these diagrams the
${E}_n$-{\it series}, where $n = r+3$. In Table
\ref{table_char_polynom_Dynkin} and Table
\ref{table_char_ext_polynom_Dynkin} we see that characteristic
polynomials of ${E}_6$, ${E}_7$, ${E}_8$, $\tilde{E}_8$ constitute
a series, and this series can be continued, see (\ref{E_series})
and Table \ref{table_char_E_series}.

\begin{figure}[h]
\centering
\includegraphics{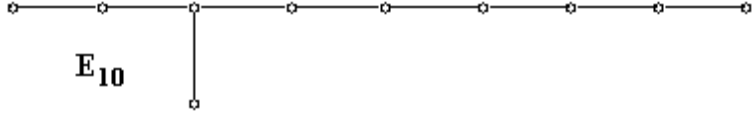}
\caption{\hspace{3mm} The spectral radius attains its minimum on
the graph ${E}_{10}$}
\label{E10_diagram}
\end{figure}

\index{Lehmer's number} \index{theorem! - McMullen}
\begin{remark}
\label{mcmullen} \rm { C.~McMullen observed in \cite{McM02}  that
the spectral radius $\rho({\bf C})$ of the Coxeter transformation
for all graphs with indefinite Tits form attains its minimum when
the diagram is $T_{2,3,7} = {E}_{10}$, see Fig.~\ref{E10_diagram},
\cite{Hir02}, \cite{McM02}, and \S\ref{hyperbolic} about {\it
hyperbolic Dynkin diagrams}. McMullen showed that $\rho({\bf C})$
is the so-called {\it Lehmer's number},
$$
\lambda_{Lehmer} \approx 1.176281...\ , \;\text{see Table
\ref{table_char_E_series}. } $$
For details and definitions, see
\S\ref{sect_numbers}. }
\end{remark}
{\it Lehmer's number} is a root of the polynomial 
$$
   1 + x  - x^3 - x^4 - x^5 - x^6 - x^7 + x^9 + x^{10},
$$
see \cite{Hir02}, \cite{McM02}.
\index{spectral radius}

\begin{proposition}
 \label{polyn_T_1}
The characteristic polynomials of Coxeter transformations for 
the  diagrams $T_{2,3,n}$ are as follows:
\begin{equation}
 \label{E_series}
  \chi(T_{2,3,n-3}) = \chi({E}_n)   =
     \lambda^n + \lambda^{n-1} -
     \sum\limits_{i=3}^{n-3}\lambda^{i} + \lambda + 1,
\end{equation}
where $\chi({E}_n)$ is given up to sign $(-1)^n$, see Remark
\ref{up_to_sign}. In other words,
\begin{equation}
 \label{E_series_1}
  \chi(T_{2,3,n-3}) = \chi({E}_n)   =
\frac{\lambda^{n-2}(\lambda^3 - \lambda - 1) + (\lambda^3 + \lambda^2 -1)}
       {\lambda - 1}.
\end{equation}

The spectral radius $\rho(T_{2,3,n-3})$ converges to the maximal root
$\rho_{max}$
of the equation
\begin{equation}
 \label{zhang_eq}
  \lambda^3 - \lambda - 1 = 0,
\end{equation}
and
\begin{equation}
 \label{zhang_number}
   \rho_{max} =
   \sqrt[3]{\frac{1}{2} + \sqrt{\frac{23}{108}}} +
   \sqrt[3]{\frac{1}{2} - \sqrt{\frac{23}{108}}} \approx 1.324717...\ .
\end{equation}
\end{proposition}
\begin{remark}
\label{rem_Zhang}
\rm{
Note that the spectral radius (\ref{zhang_number}) was obtained
  by Y.~Zhang \cite{Zh89} and used in the study of regular components
  of an Auslander-Reiten quiver. The spectral radius (\ref{zhang_number})
  is also the smallest {\it Pisot number}, see \S\ref{sect_numbers}.
  }
\end{remark}
\index{spectral radius}

\PerfProof Use induction and (\ref{split_edge_eq_corol_2}).  We have

\begin{equation}
 \label{sign_E_series}
  \begin{split}
   & (-1)^{n+2}\chi({E}_{n+2}) =
    -(\lambda + 1)(-1)^{n+1}\chi({E}_{n+1})
   -\lambda(-1)^{n}\chi({E}_n) = \\
   & (-1)^{n+2}[(\lambda + 1)\chi({E}_{n+1}) - \lambda\chi({E}_n)].
  \end{split}
\end{equation}
Thus,
\begin{equation*}
    \chi({E}_{n+2}) =
     (\lambda + 1)\chi({E}_{n+1}) - \lambda\chi({E}_n) =
    \lambda(\chi({E}_{n+1} - \chi({E}_n)) +
     \chi({E}_{n+1}).
\end{equation*}
By the induction hypothesis, (\ref{E_series}) yields
$$
  \lambda(\chi({E}_{n+1} - \chi({E}_n)) =
   \lambda^{n+2} - \lambda^{n} - \lambda^{n-1},
$$
and
\begin{equation*}
  \begin{split}
   & \chi({E}_{n+2}) =
     \lambda^{n+2} - \lambda^{n} - \lambda^{n-1} +
     (\lambda^{n+1} + \lambda^{n} -
     \sum\limits_{i=3}^{n-2}\lambda^{i} + \lambda + 1) = \\
   & \lambda^{n+2} + \lambda^{n+1} -
     \sum\limits_{i=3}^{n-1}\lambda^{i} + \lambda + 1,
  \end{split}
\end{equation*}
and hence (\ref{E_series}) is proved. Further,
\begin{equation*}
\sum\limits_{i=3}^{n-3}\lambda^{i} = \lambda^3\sum\limits_{i=0}^{n-6}\lambda^{i} =
      \lambda^3\frac{\lambda^{n-5} - 1}{\lambda - 1} =
      \frac{\lambda^{n-2} - \lambda^3}{\lambda - 1},
\end{equation*}
and
\begin{equation*}
\begin{split}
\chi({E}_n) =
   &  \lambda^n + \lambda^{n-1} - \frac{\lambda^{n-2} - \lambda^3}{\lambda - 1}
      + \lambda + 1 = \\
   &  \frac{\lambda^{n+1} - \lambda^{n-1} - \lambda^{n-2} + \lambda^3 + \lambda^2 - 1}
      {\lambda - 1} = \vspace{3mm} \\
   & \frac{\lambda^{n-2}(\lambda^3 - \lambda - 1) + (\lambda^3 + \lambda^2 -1)}
       {\lambda - 1}.
\end{split}
\end{equation*}
This proves (\ref{E_series_1}).  For $\lambda > 1$, the
characteristic polynomial $\chi({E}_n)$ converges to
\begin{equation*}
\frac{\lambda^{n-2}(\lambda^3 - \lambda - 1)}
       {\lambda - 1}.
\end{equation*}
By Proposition \ref{convergence_lambda} the spectral radius
$\rho_{max}$ converges to the maximal root
of equation (\ref{zhang_eq}). For the results of calculations 
related to the Proposition \ref{polyn_T_1}, see Table
\ref{table_char_E_series}. \qedsymbol
\begin{table} 
  \centering
  \vspace{2mm}
  \caption{\hspace{3mm}The characteristic polynomials for  $T_{2,3,r}$}
  \renewcommand{\arraystretch}{1.2} 
  \begin{tabular} {|c|c|c|}
  \hline \hline
      $T_{2,3,r}$  &  Characteristic  & Maximal \cr
       and its&     polynomial   & eigenvalue \cr
       well-known  &                  &  outside the \cr
       name        &                  & unit circle        \\
  \hline \hline
     & & \cr
       $(2,3,2)$
     & $\lambda^5 + \lambda^4 + \lambda + 1$
     &  $-$ \cr
     ${D}_5$ & & \\
  \hline
     & & \cr
       $(2,3,3)$
     & $\lambda^6 + \lambda^5 - \lambda^3 + \lambda + 1$
     &  $-$ \cr
     ${E}_6$ & & \\
  \hline
      & & \cr
      $(2,3,4)$
     & $\lambda^7 + \lambda^6 - \lambda^4 - \lambda^3 + \lambda + 1$
     & $-$ \cr
      ${E}_7$ & & \\
  \hline
     & & \cr
       $(2,3,5)$
     & $\lambda^8 + \lambda^7 - \lambda^5 - \lambda^4 -\lambda^3 + \lambda + 1$
     & $-$ \cr
      ${E}_8$ & & \\
  \hline
     & & \cr
       $(2,3,6)$
     & $\lambda^9 + \lambda^8 -
       \lambda^6 - \lambda^5 - \lambda^4 -\lambda^3 + \lambda + 1$
     & $-$ \cr
      $\tilde{E}_8$ & & \\
  \hline
     & & \cr
       $(2,3,7)$
     & $\lambda^{10} + \lambda^9 - \lambda^7
       - \lambda^6 - \lambda^5 - \lambda^4 -\lambda^3 + \lambda + 1$
     & $1.176281...$ \cr
      ${E}_{10}$ & & \\
  \hline
     & & \cr
       $(2,3,8)$
     & $\lambda^{11} + \lambda^{10} - \lambda^8 - \lambda^7
       - \lambda^6 - \lambda^5 - \lambda^4 -\lambda^3 + \lambda + 1$
     & $1.230391...$ \cr
      ${E}_{11}$ & & \\
  \hline
     & & \cr
       $(2,3,9)$
     & $\lambda^{12} + \lambda^{11} - \lambda^9 - \lambda^8 - \lambda^7
       - \lambda^6 - \lambda^5 - \lambda^4 -\lambda^3 + \lambda + 1$
     & $1.261231...$ \cr
      ${E}_{12}$ & & \\
  \hline
     & & \cr
       $(2,3,10)$
     & $ \lambda^{13} + \lambda^{12} - \lambda^{10} -
         \lambda^9 - \lambda^8 - \lambda^7$
          \qquad \qquad \qquad \qquad
     & $1.280638...$ \cr
      ${E}_{13}$
        & \qquad \qquad \qquad \qquad
        $- \lambda^6 - \lambda^5 - \lambda^4 -\lambda^3 + \lambda + 1$ & \\
  \hline
     & & \cr
     \dots & \dots &  \dots \\
     & & \cr
  \hline
     & & \cr
     & & 1.324717... \cr
       $(2,3,n)$
     & $\displaystyle\frac{\lambda^{n-2}(\lambda^3 - \lambda - 1) + (\lambda^3 +
\lambda^2 -1)}
       {\lambda - 1}$
     & \text{as} \cr
      ${E}_{n}$ & & $n \rightarrow \infty$ \\
  \hline \hline
\end{tabular}
  \label{table_char_E_series}
\end{table}

\subsection{The diagrams $T_{3,3,r}$}

 \index{characteristic polynomial of the Coxeter transformation! -
$T_{3,3,r}$}
 \index{${E}_{6,n}$-series}
 The case $T_{3,3,r}$,
where $r \geq 2$, includes diagrams ${E}_6$ and $\tilde{E}_6$, and
so we call these diagrams the ${E}_{6,n}$-{\it series}, where $n =
r-2$. Thus, ${E}_6 =  T_{3,3,2} = {E}_{6,0}$ and $\tilde{E}_6 =
T_{3,3,3} = {E}_{6,1}$. In Table \ref{table_char_polynom_Dynkin}
and Table \ref{table_char_ext_polynom_Dynkin} we see that the
characteristic polynomials of ${E}_6, \tilde{E}_6$ constitute a
series, and this series can be continued, see (\ref{E6_series})
and Table \ref{table_char_E6_series}.

\begin{proposition}
 \label{polyn_T_2}
The characteristic polynomials of Coxeter transformations for 
the  diagrams $T_{3,3,n}$ with $n \geq 3$ are calculated as
follows:
\begin{equation}
 \label{E6_series}
  \chi(T_{3,3,n}) = \chi({E}_{6,n-2})   =
     \lambda^{n+4} + \lambda^{n+3}
     -2\lambda^{n+1}
     -3\sum\limits_{i=4}^{n}\lambda^{i}
     -2\lambda^3 + \lambda + 1,
\end{equation}
where $\chi({E}_{6,n})$ is given up to sign $(-1)^n$, see Remark
\ref{up_to_sign}. In other words,
\begin{equation}
 \begin{split}
 \label{E6_series_1_other}
  \chi(T_{3,3,n}) = & \chi({E}_{6,n-2})   = \\
  &   \frac{\lambda^{n+1}(\lambda^4 - \lambda^2 - 2\lambda - 1)
   + (\lambda^4 + 2\lambda^3 + \lambda^2 - 1)}{\lambda - 1} = \\
  &  \frac{(\lambda^2 + \lambda + 1)
    [\lambda^{n+1}(\lambda^2 - \lambda - 1) + (\lambda^2 + \lambda - 1)]}{\lambda - 1}.
 \end{split}
\end{equation}

The spectral radius $\rho(T_{3,3,n})$ converges to the maximal root
$\rho_{max}$
of the equation
\begin{equation}
 \label{spectr_rad_3}
  \lambda^2 - \lambda - 1 = 0,
\end{equation}
and
\begin{equation}
 \label{spectr_rad_4}
   \rho_{max} =
   \frac{\sqrt{5} + 1}{2} \approx 1.618034...\ .
\end{equation}
\end{proposition}

\begin{remark} \rm{
For $n = 3$, the sum
$\displaystyle\sum\limits_{i=4}^{n}\lambda^{i}$ in
(\ref{E6_series}) disappears and we have
\begin{equation}
 \label{E6_series_1}
 \begin{split}
 & \chi(T_{3,3,3}) = \chi({E}_{6,1}) = \chi(\tilde{E}_6) = \\
 &    \lambda^7 + \lambda^6
     -2\lambda^4
     -2\lambda^3 + \lambda + 1 =
   (\lambda^3 - 1)^2(\lambda + 1).
  \end{split}
\end{equation}
}
\end{remark}

\begin{table} 
  \centering
  \vspace{2mm}
  \caption{\hspace{3mm}The characteristic polynomials for $T_{3,3,r}$}
  \renewcommand{\arraystretch}{1.2} 
  \begin{tabular} {|c|c|c|}
  \hline \hline
      $T_{3,3,r}$  &  Characteristic  & Maximal \cr
     and its &     polynomial   & eigenvalue \cr
       well-known  &                  &  outside the \cr
       name        &                  & unit circle        \\
  \hline \hline
     & & \cr
       $(3,3,2)$
     & $\lambda^6 + \lambda^5 - \lambda^3 + \lambda + 1$
     &  $-$ \cr
     ${E}_6$ & & \\
  \hline
     & & \cr
       $(3,3,3)$
     & $\lambda^7 + \lambda^6 - 2\lambda^4 - 2\lambda^3 + \lambda + 1$
     & $-$ \cr
      $\tilde{E}_6$ & & \\
  \hline
     & & \cr
       $(3,3,4)$
     & $\lambda^8 + \lambda^7 - 2\lambda^5 - 3\lambda^4 - 2\lambda^3 + \lambda + 1$
     & $1.401268...$ \cr
      ${E}_{6,2}$ & & \\
  \hline
     & & \cr
       $(3,3,5)$
     & $\lambda^9 + \lambda^8
       - 2\lambda^6 - 3\lambda^5 - 3\lambda^4 - 2\lambda^3 + \lambda + 1$
     & $1.506136...$ \cr
      ${E}_{6,3}$ & & \\
  \hline
     & & \cr
       $(3,3,6)$
     & $\lambda^{10} + \lambda^9 - 2\lambda^7
       - 3\lambda^6 - 3\lambda^5 - 3\lambda^4 - 2\lambda^3 + \lambda + 1$
     & $1.556030...$ \cr
      ${E}_{6,4}$ & & \\
  \hline
     & & \cr
       $(3,3,7)$
     & $\lambda^{11} + \lambda^{10} - 2\lambda^8 - 3\lambda^7
       - 3\lambda^6 - 3\lambda^5 - 3\lambda^4 - 2\lambda^3 + \lambda + 1$
     & $1.582347...$ \cr
      ${E}_{6,5}$ & & \\
  \hline
     & & \cr
       $(3,3,8)$
     & $\lambda^{12} + \lambda^{11} - 2\lambda^9 - 3\lambda^8 - 3\lambda^7$
       \qquad \qquad \qquad \qquad
     & $1.597005...$ \cr
      ${E}_{6,6}$ &
       \qquad \qquad \qquad \qquad
       $- 3\lambda^6 - 3\lambda^5 - 3\lambda^4 - 2\lambda^3 + \lambda + 1$& \\
  \hline
     & & \cr
       \dots & \dots &  \dots \\
     & & \cr
  \hline
     & & \cr
     & & 1.618034... \cr
       $(3,3,n)$
     & $\displaystyle
      \frac{(\lambda^2 + \lambda + 1)
           [\lambda^{n+1}(\lambda^2 - \lambda - 1) + (\lambda^2 + \lambda - 1)]}
           {\lambda - 1}$
     & \text{as} \cr
      ${E}_{6,n-2}$ & & $n \rightarrow \infty$ \\
  \hline \hline
\end{tabular}
  \label{table_char_E6_series}
\end{table}

 \PerfProof As above in (\ref{E_series}), we
use induction and (\ref{split_edge_eq_corol_2}). So, by
(\ref{split_edge_eq_corol_2}) we have
\begin{equation}
  \label{recurr_E6}
  \begin{split}
     \chi({E}_{6,n+2}) & =
    (\lambda + 1)\chi({E}_{6,n+1}) - \lambda\chi({E}_{6,n}) = \\
    & \lambda(\chi({E}_{6,n+1}) - \chi({E}_{6,n}))
    + \chi({E}_{6,n+1}).
  \end{split}
\end{equation}
By the induction hypothesis and (\ref{E6_series}) we have
\begin{equation}
    \lambda(\chi({E}_{6,n+1}) - \chi({E}_{6,n})) =
    \lambda^{n+8} - \lambda^{n+6} - 2\lambda^{n+5} - \lambda^{n+4},
\end{equation}
and
\begin{equation}
 \begin{split}
   \chi({E}_{6,n+2}) & = \\
  & \lambda^{n+8} - \lambda^{n+6} - 2\lambda^{n+5} - \lambda^{n+4} + \\
    & (\lambda^{n+7} + \lambda^{n+6}
     -2\lambda^{n+4}
     -3\sum\limits_{i=4}^{n+3}\lambda^{i}
     -2\lambda^3 + \lambda + 1) = \\
  & \lambda^{n+8} + \lambda^{n+7} - 2\lambda^{n+5}
     -3\sum\limits_{i=4}^{n+4}\lambda^{i}
     -2\lambda^3 + \lambda + 1;
  \end{split}
\end{equation}
this proves (\ref{E6_series}). Further,

\begin{equation*}
 \begin{split}
  \chi(T_{3,3,n}) =
    & \lambda^{n+4} + \lambda^{n+3} - 2\lambda^{n+1}
      - 3\sum\limits_{i=4}^{n}\lambda^i
      - 2\lambda^3 + \lambda + 1 = \\
    & \lambda^{n+4} + \lambda^{n+3} + \lambda^{n+1}
      - 3\sum\limits_{i=3}^{n+1}\lambda^{i}
      + \lambda^3 + \lambda + 1 = \\ \\
    & \lambda^{n+1}(\lambda^3 + \lambda^2 + 1)
     -3\lambda^3\frac{\lambda^{n-1} -1}{\lambda-1} +\lambda^3 + \lambda + 1,
 \end{split}
\end{equation*}
i.e.,
\begin{equation*}
 \begin{split}
  \chi(T_{3,3,n}) =
    & \frac{\lambda^{n+1}(\lambda^4 - \lambda^2 + \lambda - 1)
     -3\lambda^{n+2} + 3\lambda^3 + (\lambda^4 -\lambda^3 + \lambda^2 - 1)}
     {\lambda - 1} = \\ \\
    & \frac{\lambda^{n+1}(\lambda^4 - \lambda^2 - 2\lambda - 1)
      + (\lambda^4 +2\lambda^3 + \lambda^2 - 1)}
     {\lambda - 1},
 \end{split}
\end{equation*}
and (\ref{E6_series_1_other}) is proved. For $\lambda > 1$, the
characteristic polynomial $\chi(T_{3,3,n})$ converges to
\begin{equation*}
\frac{\lambda^{n+1}(\lambda^4 - \lambda^2 - 2\lambda - 1)}
     {\lambda - 1} =
\frac{\lambda^{n+1}(\lambda^2 - \lambda - 1)(\lambda^2 + \lambda + 1)}
     {\lambda - 1}.
\end{equation*}
By Proposition \ref{convergence_lambda} the spectral radius
$\rho_{max}$ converges to the maximal root of the equation
(\ref{spectr_rad_3}). For the results of calculations related to
the Proposition \ref{polyn_T_2}, see Table
\ref{table_char_E6_series}. \qedsymbol

\subsection{The diagrams $T_{2,4,r}$}

\index{characteristic polynomial of the Coxeter transformation! -
$T_{2,4,r}$} \index{${E}_{7,n}$-series} The case $T_{2,4,r}$,
where $r \geq 2$, includes diagrams ${D}_6, {E}_7, \tilde{E}_7$
and we call these diagrams the ${E}_{7,n}$-{\it series}, where $n
= r-3$. Thus, ${D}_6 = T_{2,4,2} = {E}_{7,-1}$, ${E}_7 = T_{2,4,3}
= {E}_{7,0}$, $\tilde{E}_7 = T_{2,4,4} = {E}_{7,1}$. In Table
\ref{table_char_polynom_Dynkin} and Table
\ref{table_char_ext_polynom_Dynkin} we see that characteristic
polynomials of ${E}_7$, $\tilde{E}_7$ constitute a series, and
this series can be continued, see (\ref{E7_series}) and Table
\ref{table_char_E7_series}.

\begin{proposition}
 \label{polyn_T_3}
The characteristic polynomials of Coxeter transformations for 
diagrams $T_{2,4,n}$, where $n \geq 3$,
  are calculated as follows:
\begin{equation}
 \label{E7_series}
  \chi(T_{2,4,n}) = \chi({E}_{7,n-3})   =
     \lambda^{n+4} + \lambda^{n+3}
     -\lambda^{n+1}
     -2\sum\limits_{i=4}^{n}\lambda^{i}
     -\lambda^3 + \lambda + 1,
\end{equation}
where $\chi({E}_{7,n})$ is given up to sign $(-1)^n$, see Remark
\ref{up_to_sign}.  In other words
\begin{equation}
  \begin{split}
 \label{E7_series_1_other}
  \chi(T_{2,4,n}) = \chi({E}_{7,n-3})   =
 &    \frac{\lambda^{n+1}(\lambda^4 - \lambda^2 - \lambda - 1)
           +(\lambda^4 + \lambda^3 + \lambda - 1)}{\lambda - 1} = \\
 &  \frac{(\lambda+1)(\lambda^{n+1}(\lambda^3 - \lambda^2 -  1)
           +(\lambda^3 + \lambda - 1))}{\lambda - 1}.
 \end{split}
\end{equation}
The spectral radius $\rho(T_{2,4,n})$ converges to the maximal root
$\rho_{max}$
of the equation
\begin{equation}
 \label{spectr_rad_1}
  \lambda^3 - \lambda^2 - 1 = 0,
\end{equation}
and
\begin{equation}
 \label{spectr_rad_2}
   \rho_{max} =
   \frac{1}{3} + \sqrt[3]{\frac{58}{108} + \sqrt{\frac{31}{108}}} +
   \sqrt[3]{\frac{58}{108} - \sqrt{\frac{31}{108}}} \approx 1.465571...\ .
\end{equation}
\end{proposition}
\begin{remark} \rm{
For $n = 3$, the sum
$\displaystyle\sum\limits_{i=4}^{n}\lambda^{i}$ disappears from
(\ref{E7_series}) and we have
\begin{equation}
 \label{E7_series_1}
 \begin{split}
 & \chi(T_{2,4,3}) = \chi({E}_{7,0}) = \chi({E}_7) = \\
 &    \lambda^7 + \lambda^6
     -\lambda^4
     -\lambda^3 + \lambda + 1 =
   (\lambda^4 - 1)(\lambda^3 - 1)(\lambda + 1).
  \end{split}
\end{equation}
}
\end{remark}

\PerfProof
Use induction and (\ref{split_edge_eq_corol_2}).
As above in (\ref{recurr_E6}), we have
\begin{equation}
 \label{recurr_E7}
 \begin{split}
    \chi({E}_{7,n+2}) & = \\
    & (\lambda + 1)\chi({E}_{7,n+1}) - \lambda\chi({E}_{7,n}) = \\
    & \lambda(\chi({E}_{7,n+1}) - \chi({E}_{7,n}))
    + \chi({E}_{7,n+1}).
 \end{split}
\end{equation}
By the induction hypothesis and (\ref{E7_series}) we have
\begin{equation}
    \lambda(\chi({E}_{7,n+1}) - \chi({E}_{7,n})) =
    \lambda^{n+9} - \lambda^{n+7} - \lambda^{n+6} - \lambda^{n+5},
\end{equation}
and we see that
\begin{equation}
\begin{split}
    \chi({E}_{7,n+2}) = &\\
   & \lambda^{n+9} - \lambda^{n+7} - \lambda^{n+6} - \lambda^{n+5} \\
   & + (\lambda^{n+8} + \lambda^{n+7}
     -\lambda^{n+5}
      -2\sum\limits_{i=4}^{n+4}\lambda^{i}
     -\lambda^3 + \lambda + 1) =   \\
   &  \lambda^{n+9} + \lambda^{n+8} - \lambda^{n+6}
     -2\sum\limits_{i=4}^{n+5}\lambda^{i}
     -\lambda^3 + \lambda + 1);
\end{split}
\end{equation}
this proves (\ref{E7_series}). Further,
\begin{equation}
 \begin{split}
 \label{E7_series_lim}
  \chi(T_{2,4,n}) =
    & \lambda^{n+4} + \lambda^{n+3} - \lambda^{n+1}
     -2\sum\limits_{i=4}^{n}\lambda^{i} -\lambda^3 + \lambda + 1 = \\
    & (\lambda^{n+4} + \lambda^{n+3} + \lambda^{n+1})
     -2\sum\limits_{i=3}^{n+1}\lambda^{i} + (\lambda^3 + \lambda + 1) = \\
    & \lambda^{n+1}(\lambda^3 + \lambda^2 + 1)
     -2\frac{\lambda^3(\lambda^{n-1} - 1)}
     {\lambda - 1} + (\lambda^3 + \lambda + 1) = \\ \\
    & \frac{\lambda^{n+1}(\lambda^4 - \lambda^2 + \lambda - 1)
     -2\lambda^{n+2} + 2\lambda^3 + (\lambda^4 - \lambda^3 + \lambda^2 - 1)}
     {\lambda - 1} = \\ \\
    & \frac{\lambda^{n+1}(\lambda^4 - \lambda^2 - \lambda - 1)
      + (\lambda^4 + \lambda^3 + \lambda^2 - 1)}
     {\lambda - 1} = \\ \\
    & \frac{\lambda^{n+1}(\lambda+1)(\lambda^3 - \lambda^2 - 1)
      + (\lambda+1)(\lambda^3 + \lambda - 1)}
     {\lambda - 1},
 \end{split}
\end{equation}
and (\ref{E7_series_1_other}) is proved.
 For $\lambda > 1$, the characteristic
polynomial $\chi(T_{2,4,n})$ converges to
\begin{equation*}
\frac{\lambda^{n+1}(\lambda+1)(\lambda^3 - \lambda^2 - 1)}
     {\lambda - 1}.
\end{equation*}
By Proposition \ref{convergence_lambda}, the spectral radius
$\rho_{max}$ converges to the maximal root of the equation
(\ref{spectr_rad_1}). For the results of calculations related to
the Proposition \ref{polyn_T_3}, see Table
\ref{table_char_E7_series}. \qedsymbol

\index{spectral radius}

\begin{table} 
  \centering
  \vspace{2mm}
  \caption{\hspace{3mm}The characteristic polynomials for  $T_{2,4,r}$}
  \renewcommand{\arraystretch}{1.2} 
  \begin{tabular} {|c|c|c|}
  \hline \hline
      $T_{2,4,r}$  &  Characteristic  & Maximal \cr
      and its&     polynomial   & eigenvalue \cr
       well-known  &                  &  outside the \cr
       name        &                  & unit circle        \\
  \hline \hline
     & & \cr
       $(2,4,2)$
     & $\lambda^6 + \lambda^5 + \lambda + 1$
     &  $-$ \cr
     ${D}_6$ & & \\
  \hline
     & & \cr
       $(2,4,3)$
     & $\lambda^7 + \lambda^6  - \lambda^4  - \lambda^3 + \lambda + 1$
     &  $-$ \cr
     ${E}_7$ & & \\
  \hline
     & & \cr
       $(2,4,4)$
     & $\lambda^8 + \lambda^7 - \lambda^5
        - 2\lambda^4 - \lambda^3 + \lambda + 1$
     & $-$ \cr
      $\tilde{E}_7$ & & \\
  \hline
     & & \cr
       $(2,4,5)$
     & $\lambda^9 + \lambda^8 - \lambda^6
        - 2\lambda^5 -2\lambda^4 - \lambda^3 + \lambda + 1$
     & $1.280638...$ \cr
      ${E}_{7,2}$ & & \\
  \hline
     & & \cr
       $(2,4,6)$
     & $\lambda^{10} + \lambda^9 - \lambda^7
        -2\lambda^6 - 2\lambda^5 -2\lambda^4 - \lambda^3 + \lambda + 1$
     & $1.360000...$ \cr
      ${E}_{7,3}$ & & \\
  \hline
     & & \cr
       $(2,4,7)$
     & $\lambda^{11} + \lambda^{10} - \lambda^8
       -2\lambda^7 -2\lambda^6 - 2\lambda^5
       -2\lambda^4 - \lambda^3 + \lambda + 1$
     & $1.401268...$ \cr
      ${E}_{7,4}$ & & \\
  \hline
     & & \cr
       $(2,4,8)$
     & $\lambda^{12} + \lambda^{11} - \lambda^9
       -2\lambda^8 -2\lambda^7 $
       \qquad \qquad \qquad \qquad
     & $1.425005...$ \cr
      ${E}_{7,5}$ &
       \qquad \qquad \qquad \qquad
       $-2\lambda^6 - 2\lambda^5
       -2\lambda^4 - \lambda^3 + \lambda + 1$ & \\
  \hline
     & & \cr
       \dots & \dots &  \dots \\
     & & \cr
  \hline
     & & \cr
     & & 1.465571... \cr
       $(2,4,n)$
     & $\displaystyle
      \frac{\lambda^{n+1}(\lambda^4 - \lambda^2 - \lambda - 1) +
           (\lambda^4 + \lambda^3 + \lambda -1)}
       {\lambda - 1}$
     & \text{as} \cr
      ${E}_{n}$ & & $n \rightarrow \infty$ \\
  \hline \hline
\end{tabular}
  \label{table_char_E7_series}
\end{table}

The maximal eigenvalues (Tables \ref{table_char_E_series},
\ref{table_char_E6_series}, \ref{table_char_E7_series}) lying out
the unit circle are obtained by means of the online service
``Polynomial Roots Solver'' of ``EngineersToolbox''\footnote{See
  http://www.engineerstoolbox.com.}.

\subsection{Convergence of the sequence of eigenvalues}
\index{spectral radius}

In this section we give a substantiation (Proposition
\ref{convergence_lambda}) of the fact that the maximal eigenvalues
(i.e., spectral radiuses) of the characteristic polynomials for
$T_{p,q,r}$ converges to a real number which is the maximal value
of a known polynomial (see, Propositions \ref{polyn_T_1},
\ref{polyn_T_2}, \ref{polyn_T_3}). We will give this
substantiation in a generic form. Similarly, we have a dual fact
(Proposition \ref{convergence_lambda_2}): the minimal eigenvalues
of the characteristic polynomials for $T_{p,q,r}$ converges to a
real number which is the minimal value of a known polynomial.

\begin{proposition}
  \label{convergence_lambda}
   Let $f(\lambda)$ and  $g(\lambda)$ be some polynomials in
   $\lambda$, and let
  \begin{equation}
      P_n(\lambda) = \lambda^{n}f(\lambda) + g(\lambda)
  \end{equation}
  be a sequence of polynomials in $\lambda$.
  Let $z_n \in \mathbb{R}$ be a root
  of $P_n(\lambda)$,
  lying in the interval $(1,2)$:
  \begin{equation}
    \label{z_interval}
      P_n(z_n) = 0,  \hspace{5mm} z_n \in (1,2), \hspace{3mm} n = 1,2, \dots
  \end{equation}
Suppose that the sequence $\{z_n\}, n = 1,2,\dots$ is
non-decreasing:
    \begin{equation}
      \label{z_inc_seq}
    1 < z_1 \leq z_2 \leq \dots \leq z_n \leq \dots .
    \end{equation}
Then the sequence $z_n$ converges to a real number $z_0$ which is
a root of $f(\lambda)$:
    \begin{equation}
      f(z_0) = 0.
    \end{equation}
\end{proposition}

\PerfProof
  According to (\ref{z_interval}), the sequence (\ref{z_inc_seq})
is non-decreasing and bounded from above. Therefore, there exists
a real number $z_0$ such
  that
  \begin{equation}
    \label{z_0_limit}
      \lim_{k \rightarrow \infty}z_k = z_0.
  \end{equation}
Let us estimate $|f(z_0)|$ as follows:
  \begin{equation}
    \label{estimtion_f_z0}
      |f(z_0)| < |f(z_n)| + |f(z_n) - f(z_0)| .
  \end{equation}
By (\ref{z_0_limit}) we obtain
  \begin{equation}
    \label{estimtion_f_z0_2}
      |f(z_0)| < |f(z_n)| + \varepsilon_n \hspace{1mm}, \hspace{3mm} \text{ where }
      \varepsilon_n = |f(z_n) - f(z_0)| \rightarrow 0 .
  \end{equation}
By (\ref{z_interval}), we have
  \begin{equation}
    \label{f_function_2}
      f(z_n) =  - \frac{g(z_n)}{z_n^n} .
  \end{equation}
Let
  \begin{equation}
    \label{g_function}
      \delta_n = \left |
       \frac{g(z_n)}{z_n^n}
        \right | \hspace{3mm} \text{ for each } n = 1,2,\dots
  \end{equation}
 Since the function $g(z)$ is uniformly bounded on the interval $(1,2)$,
 it follows that $\delta_n  \rightarrow 0$.
By (\ref{estimtion_f_z0_2}) we have
   \begin{equation}
    \label{estimtion_f_z0_4}
      |f(z_0)| <
        \delta_n + \varepsilon_n \rightarrow 0, \hspace{3mm}
        \text{ as } n \rightarrow \infty.
  \end{equation}
Thus, $f(z_0) = 0$. \qedsymbol

\begin{corollary} Proposition \ref{convergence_lambda}
     holds also for the following non-polynomial functions $P_n$:

 1) For
  \begin{equation}
     P_n(\lambda) = \lambda^{n+k}f(\lambda) + g(\lambda),
  \end{equation}
     where $k \in \mathbb{Z}$ is independent of $n$.

  2) For
  \begin{equation}
     P_n(\lambda) = D(\lambda)(\lambda^{n}f(\lambda) + g(\lambda)),
  \end{equation}
     where $D(\lambda)$ is a rational function independent of $n$ and without roots on
     the interval $(1,2)$.
\end{corollary}

\begin{table} 
  \centering
  \vspace{2mm}
  \caption{\hspace{3mm}The diagrams $T_{p,q,r}$ and
     characteristic polynomials $\lambda^{n}f(\lambda) + g(\lambda)$}
  \renewcommand{\arraystretch}{3.9} 
  \begin{tabular} {|c|c|c|c|c|}
  \hline \hline
      Diagram      &  $f(\lambda)$    & $g(\lambda)$ &
      $D(\lambda)$ &  $P_n(\lambda)$  \\
  \hline \hline
       $T(2,3,r)$
     & $\lambda^3 - \lambda - 1$
     & $\lambda^3 + \lambda^2 - 1$
     & $\displaystyle\frac{1}{\lambda-1}$
     & $\lambda^{n-2}f(\lambda) + g(\lambda)$ \\
  \hline
       $T(3,3,r)$
     & $\lambda^2 - \lambda - 1$
     & $\lambda^2 + \lambda - 1$
     & $\displaystyle\frac{\lambda^2 + \lambda + 1}{\lambda-1}$
     & $\lambda^{n+1}f(\lambda) + g(\lambda)$ \\
  \hline
       $T(2,4,r)$
     & $\lambda^3 - \lambda^2 - 1$
     & $\lambda^3 + \lambda - 1$
     & $\displaystyle\frac{\lambda + 1}{\lambda-1}$
     & $\lambda^{n+1}f(\lambda) + g(\lambda)$ \\
  \hline \hline
\end{tabular}
  \label{table_Diagrams T_pqr}
\end{table}

\begin{corollary}
  The maximal values of characteristic polynomials of the Coxeter transformation
  for diagrams $T_{p,q,r}$ satisfy the conditions of Proposition \ref{convergence_lambda},
  see Table \ref{table_Diagrams T_pqr}.
\end{corollary}

\PerfProof Let $\varphi_{1,r} = \varphi_1(T_{p,q,r})$ be the
dominant value of the matrix $DD^t$ (see \S3.3.2) for the diagram
$T_{p,q,r}$.
 According to Proposition \ref{trend_fi} the dominant value
$\varphi_{1,r}$  may only grow, i.e.,
$$
   \varphi_{1,r} \leq \varphi_{1,r+1}.
$$
Therefore, the corresponding eigenvalue $\lambda^{max}_1$ also
only grows. Indeed, by (\ref{Coxeter_eigenvalues}) we get
\begin{equation}
 \label{lambda_max_grows}
  \begin{split}
   \lambda^{max}_{1,r} =
   & 2\varphi_{1,r} - 1 + 2\sqrt{\varphi_{1,r}(\varphi_{1,r} - 1)}
    \leq \\
   & 2\varphi_{1,r+1} - 1 + 2\sqrt{\varphi_{1,r+1}(\varphi_{1,r+1} - 1)}
   = \lambda^{max}_{1,r+1}.
  \end{split}
\end{equation}

In addition, from Proposition \ref{trend_fi} we deduce that
$\varphi_{1,r} > 1$, therefore the corresponding eigenvalue
$\lambda_{1,r}$ of the Coxeter transformation is also $> 1$.

It remains to show that every $\lambda_{1,r} < 2$. This is clear,
because $f(\lambda) > 0$ and $g(\lambda) > 0$ for every $\lambda >
2$, see Table \ref{table_Diagrams T_pqr}. \qedsymbol

\begin{proposition}
  \label{convergence_lambda_2}
   Let $f(\lambda)$ and  $g(\lambda)$ be some polynomials in
   $\lambda$, and let
  \begin{equation}
      P_n(\lambda) = \lambda^{n}f(\lambda) + g(\lambda), \;\text{
      where $n\in\mathbb{Z}_+$}.
  \end{equation}
Let $z_n \in \mathbb{R}$ be a root
  of $P_n(\lambda)$,
  lying in the interval $\left(\displaystyle\frac{1}{2},1\right)$:
  \begin{equation}
    \label{z_interval_05}
      P_n(z_n) = 0,  \hspace{5mm} z_n \in \left(\frac{1}{2},1\right). 
  \end{equation}
Suppose, that the sequence $\{z_n\}_{n\in\mathbb{Z}_+}$ is
non-increasing:
    \begin{equation}
     1 > z_1 \geq z_2 \geq \dots \geq z_n \geq \dots .
    \end{equation}
Then the sequence $z_n$ converges to the some real number $z_0$
which is a root of $g(\lambda)$:
    \begin{equation}
      g(z_0) = 0.
    \end{equation}
\end{proposition}

\PerfProof
  As in Proposition \ref{convergence_lambda}, there is a real number $z_0$ such
  that $z_n \rightarrow  z_0$ as $n \rightarrow \infty$.
Since $g(z_n) \rightarrow g(z_0)$, we obtain
  \begin{equation*}
      |g(z_0)| < |g(z_n)| + \varepsilon_n \hspace{1mm}, \hspace{3mm} \text{ where }
      \varepsilon_n = |g(z_n) - g(z_0)|.
  \end{equation*}
By (\ref{z_interval_05}), we have
  \begin{equation*}
      g(z_n) =  - {z_n^n}{f(z_n)} .
  \end{equation*}
Let
  \begin{equation*}
      \delta_n =  |{z_n^n}{f(z_n)}| \hspace{3mm} \text{ for each } n = 1,2,\dots
  \end{equation*}
 Since the function $f(z_n)$ is uniformly bounded on the interval
$\left(\displaystyle\frac{1}{2},1\right)$, it follows that
$\delta_n \rightarrow 0$. Thus, we have
   \begin{equation*}
      |g(z_0)| <
        \delta_n + \varepsilon_n \rightarrow 0 \hspace{3mm}
        \text{ as } n \rightarrow \infty.
  \end{equation*}
Thus, $g(z_0) = 0$. \qedsymbol

\begin{corollary}
  The minimal values of characteristic polynomials of the Coxeter transformation
  for diagrams $T_{p,q,r}$ satisfy the conditions of
  Proposition \ref{convergence_lambda},
  see Table \ref{table_Diagrams T_pqr}.
\end{corollary}

\PerfProof The minimal value $\lambda^{min}_1$ and the maximal
value $\lambda^{max}_1$ are reciprocal:
$$
    \lambda^{min}_1 = \frac{1}{\lambda^{max}_1}.
$$

 Therefore, by (\ref{lambda_max_grows}), we see that the sequence
of eigenvalues $\lambda_{1, r}^{min}$ is non-increasing:
$$
   \lambda^{min}_{1,r} \geq \lambda^{min}_{1,r+1}.
$$
Since the maximal eigenvalue $\lambda^{max}_{1,r} > 1$, then
$\lambda^{min}_{1,r} < 1$.

It remains to show that $\lambda^{min}_{1,r} >
\displaystyle\frac{1}{2}$ for every $r \in \mathbb{Z}$. But this
is true since $f(\lambda) < 0$ and $g(\lambda) < 0$ for every
$\lambda < \displaystyle\frac{1}{2}$, see Table
\ref{table_Diagrams T_pqr}. \qedsymbol

%% file: 5steinberg.tex

\chapter{\sc\bf R.~Steinberg's theorem, B.~Kostant's construction}
 \label{chapter_steinberg}

\section{R.~Steinberg's theorem and a $(p,q,r)$ mystery}
 \label{Steinberg_simply_laced}
R.~Steinberg in \cite[p.591, $(*)$]{Stb85} observed a property
of affine Coxeter transformations (i.e., transformations corresponding to
extended Dynkin diagrams), which plays the main role in his derivation
of the McKay correspondence.
  Let $(p,q,r)$ be the same as in
  Table \ref{rotation_pol_1}, Table \ref{binary_pol_1} and
  relations (\ref{diophantine_inequality}), (\ref{natural_gen_2}).

\index{theorem! - Steinberg}
  \begin{theorem} [\cite{Stb85}]
    \label{Steinberg}
The affine Coxeter transformation for the extended Dynkin diagram
$\tilde{\varGamma}$ has the same eigenvalues as the product of
three Coxeter transformations of types ${A}_n$, where $n=p-1$,
$q-1$, and $r-1$, corresponding to the branches of the Dynkin
diagram $\varGamma$.
\end{theorem}

Essentially, R.~Steinberg observed that the orders of eigenvalues
of the affine Coxeter transformation corresponding to the {\it
extended} Dynkin diagram $\tilde{\varGamma}$ and given in Table
\ref{table_eigenvalues} coincide with the lengths of branches of
the Dynkin diagram $\varGamma$.

Now, we give the proof of R.~Steinberg's theorem for the simply-laced
  extended Dynkin diagram by using the Subbotin-Sumin splitting
  formula (\ref{split_edge_eq}).

\begin{figure}[h]
\centering
\includegraphics{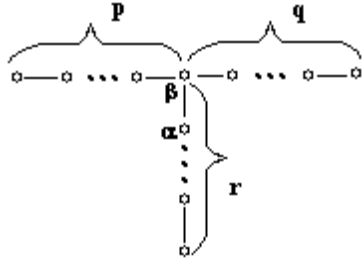}
\caption{\hspace{3mm} Splitting the graph $T_{p,q,r}$ along the
edge $\{ \alpha \beta \}$}
\label{split_Tpqr}
\end{figure}

  In \S\ref{sect_gen_Steinberg},
we generalize this theorem with some modifications
to the multiply-laced case, see Theorem \ref{gen_Steinberg},
Table \ref{table_char_ext_polynom_Dynkin}.

Let $T_{p,q,r}$ be a connected graph with three branches of
lengths $p,q,r$. Let $\chi$ be the characteristic polynomial of
the Coxeter transformation. Split the graph along the edge $\{
\alpha \beta \}$, Fig. \ref{split_Tpqr}. Then, by the
Subbotin-Sumin formula (\ref{split_edge_eq}) we have

\begin{equation}
  \label{sub_sum_1}
  \begin{split}
    \chi =
    & \mathcal{X}({A}_{p+q-1})\mathcal{X}({A}_{r-1}) -
            \lambda\mathcal{X}({A}_{p-1})\mathcal{X}({A}_{q-1})
                 \mathcal{X}({A}_{r-2}) = \\
    & (-1)^{p+q+r}(\chi_{p+q-1}\chi_{r-1} -
             \lambda\chi_{p-1}\chi_{q-1}\chi_{r-2}).
  \end{split}
\end{equation}

By (\ref{An_chi}),  we have, up to a factor $(-1)^{p + q + r}$,
\begin{equation}
 \label{sub_sum_2}
   \chi = \frac{(\lambda^{p+q} -1)(\lambda^{r} -1)}{(\lambda-1)^2} -
          \frac{(\lambda^{p} -1)(\lambda^{q} - 1)(\lambda^{r} -
\lambda)}
               {(\lambda-1)^3}. \vspace{3mm}
\end{equation}

For $p=q$ (e.g., ${E}_6$, $\tilde{E}_6$,
  $\tilde{E}_7$, and ${D}_{n}$ for $n \geq 4$), we have, up to a factor
$(-1)^{2p+r}$,
\begin{equation}
  \label{chi_p_eq_q}
  \begin{split}
   \chi = & \frac{(\lambda^{2p} - 1)(\lambda^{r} -1)}{(\lambda-1)^2} -
          \frac{(\lambda^{p} - 1)^2(\lambda^{r} - \lambda)}
               {(\lambda-1)^3} = \vspace{3mm} \\
   & \frac{(\lambda^{p} -1)}{(\lambda-1)^3}
      ((\lambda^{p} + 1)(\lambda^{r} - 1)(\lambda-1) -
          (\lambda^{p} - 1) (\lambda^{r} - \lambda)).
  \end{split}
\end{equation}

For $q=2p$ (e.g., ${E}_7$, $\tilde{E}_8$), we have, up to a factor
$(-1)^{3p+r}$,
\begin{equation}
  \label{chi_q_eq_2p}
  \begin{split}
   \chi = & \frac{(\lambda^{3p} - 1)(\lambda^{r} -1)}{(\lambda-1)^2} -
          \frac{(\lambda^{2p} - 1)(\lambda^{p} - 1)(\lambda^{r} -
\lambda)}{(\lambda-1)^3} = \vspace{3mm} \\
   & \frac{(\lambda^{p} -1)}{(\lambda-1)^3}
      ((\lambda^{2p} + \lambda^{p} + 1)(\lambda^{r} - 1)(\lambda-1) -
          (\lambda^{2p} - 1) (\lambda^{r} - \lambda)).
  \end{split}
\end{equation}

1) Case $p=q=3$, $r=3$ ($\tilde{E}_6$). From (\ref{chi_p_eq_q}) we
have
\begin{equation}
  \begin{split}
   \chi =
   & \frac{(\lambda^{3} -1)}{(\lambda - 1)^3}
      [(\lambda^{3} + 1)(\lambda^{3} - 1)(\lambda - 1) -
          (\lambda^{3} - 1) (\lambda^{3} - \lambda)] = \\
   &     \frac{(\lambda^{3} - 1)^2}{(\lambda-1)^3}
      [(\lambda^{3} + 1)(\lambda-1) -
          (\lambda^{3} - \lambda)] = \\
   &  \frac{(\lambda^{3} - 1)^2}{(\lambda - 1)^3}
      (\lambda^{4} -2\lambda^3 + 2\lambda -1) = \\
   &     \frac{(\lambda^{3} - 1)^2{(\lambda - 1)^2}(\lambda^2 - 1)}
       {(\lambda-1)^3} =
       (\lambda - 1)^2\chi^2_2\chi_1.
  \end{split}
\end{equation}

Polynomials $\chi_1$ and $\chi_2^2$ have, respectively,
eigenvalues of orders $2,3,3$ which are equal to the lengths of
branches of $E_6$.


2) Case $p=q=4$, $r=2$ ($\tilde{E}_7$). Here, from
(\ref{chi_p_eq_q}) we have
\begin{equation}
  \begin{split}
   \chi =
   & \frac{(\lambda^{4} - 1)}{(\lambda - 1)^3}
      [(\lambda^{4} + 1)(\lambda^{2} - 1)(\lambda - 1) -
          (\lambda^{4} - 1)(\lambda^{2} - \lambda)] = \\
   & \frac{(\lambda^{4} - 1)(\lambda^2-1)(\lambda-1)}{(\lambda-1)^3}
      [(\lambda^{4} + 1) -
          \lambda(\lambda^{2} + 1)] = \\
   & \frac{(\lambda^{4} - 1)(\lambda^3 - 1)
           (\lambda^2-1)(\lambda-1)^2}{(\lambda-1)^3} =
       (\lambda - 1)^2\chi_3\chi_2\chi_1.
   \end{split}
\end{equation}

Polynomials $\chi_1, \chi_2$ and $\chi_3$ have, respectively,
eigenvalues of orders $2,3,4$ which are equal to the lengths of
branches of $E_7$.


3) Case $p=3$, $q=6$, $r=2$ ($\tilde{E}_8$). From
(\ref{chi_q_eq_2p}), we have
\begin{equation}
  \begin{split}
   \chi =   & \frac{(\lambda^{3} -1)}{(\lambda-1)^3}
      [(\lambda^{6} + \lambda^{3} + 1)(\lambda^{2} - 1)(\lambda-1) -
          (\lambda^{6} - 1) (\lambda^{2} - \lambda)] = \\
            & \frac{(\lambda^{3} -1)(\lambda^{2} -1)(\lambda -1)}{(\lambda-1)^3}
      [(\lambda^{6} + \lambda^{3} + 1) -
          \lambda(\lambda^{4} + \lambda^{2} + 1)] = \\
            & \frac{(\lambda^{5} -1)(\lambda^{3} -1)
              (\lambda^{2} - 1)(\lambda - 1)^2}{(\lambda-1)^3} =
            (\lambda - 1)^2\chi_4\chi_2\chi_1.
   \end{split}
\end{equation}

Polynomials $\chi_1, \chi_2$ and $\chi_4$ have, respectively,
eigenvalues of orders $2,3,5$ which are equal to the lengths of
branches of $E_8$.


4) Case $p=q=2$,  $r=2,3,...$ (${D}_{r+2}$). From
(\ref{chi_p_eq_q}) we have
\begin{equation}
  \label{chi_D_r_plus2}
  \begin{split}
     & (-1)^r\mathcal{X}({D}_{r+2}) = \\
     & \frac{(\lambda^{2} - 1)}{(\lambda - 1)^3}
      [(\lambda^{2} + 1)(\lambda^{r} - 1)(\lambda - 1) -
          (\lambda^{2} - 1) (\lambda^{r} - \lambda)] = \\
     & \frac{(\lambda^{2} - 1)(\lambda -1)}{(\lambda-1)^3}
      [(\lambda^{2} + 1)(\lambda^{r} - 1) -
          (\lambda + 1)(\lambda^{r} - \lambda)] =   \\
     & \frac{(\lambda^{2} -1)(\lambda -1)^2
              (\lambda^{r+1} + 1)}{(\lambda-1)^3} =
            (\lambda + 1)(\lambda^{r+1} + 1),
   \end{split}
\end{equation}
and, up to a sign,
\begin{equation}
\chi = (\lambda + 1)(\lambda^{r+1} + 1).
\end{equation}
Pay attention that,  as in (\ref{chi_D_r_plus2}), the sign
$(-1)^r$ should be taken into account only in recurrent
calculations (e.g., for $\tilde{D}_{n}$ and $\tilde{D}_{4}$
below), see Remark \ref{up_to_sign}.

5) Case $\tilde{D}_{r+2}$, contains $r+3$ points.

\begin{figure}[h]
\centering
\includegraphics{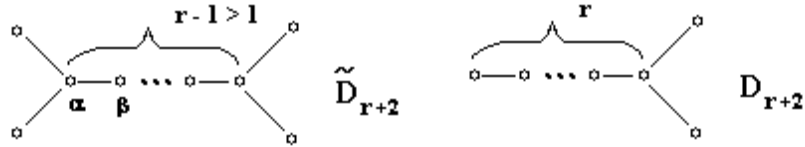}
\caption{\hspace{3mm} Diagrams $\tilde{D}_{r+2}$ and ${D}_{r+2}$}
\label{extended_ordinary_Dn}
\end{figure}

Here, from (\ref{split_edge_eq}) we have
\begin{equation}
  \label{chi_extended_Dr}
   \mathcal{X}(\tilde{D}_{r+2}) = \mathcal{X}(D_{r})\mathcal{X}({A}_3) -
   \lambda\mathcal{X}(D_{r-1})\mathcal{X}({A}_1)^2.
\end{equation}
  From (\ref{chi_extended_Dr}) and (\ref{chi_D_r_plus2})
  we have
\begin{equation*}
  \label{chiextendedDr_1}
  \begin{split}
   & (-1)^{r-1}\mathcal{X}(\tilde{D}_{r+2}) = \\
   & (\lambda + 1)(\lambda^{r-1} + 1)(\lambda^3 + \lambda^2 + \lambda + 1) -
   \lambda(\lambda + 1)(\lambda^{r-2} + 1)(\lambda + 1)^2 = \\
   & (\lambda + 1)^2[(\lambda^{r-1} + 1)(\lambda^2 + 1) -
   (\lambda + 1)(\lambda^{r-1} + \lambda)] = \\
   & (\lambda + 1)^2[(\lambda^{r+1} + \lambda^2 + \lambda^{r-1} + 1) -
   (\lambda^{r} + \lambda^{r-1} + \lambda^2 + \lambda)] = \\
   & (\lambda + 1)^2[(\lambda^{r+1} - \lambda^{r}) -
   (\lambda - 1)] =
   (\lambda + 1)^2(\lambda^{r} - 1)(\lambda - 1) = \\
   & (\lambda - 1)^2\chi^2_1\chi_{r-1},   \end{split}
\end{equation*}
and so
\begin{equation}
  \label{chi_extended_Dr_1}\chi =
  (\lambda - 1)^2\chi^2_1\chi_{r-1}.\end{equation}

Polynomials $\chi_1^2$ and $\chi_{r-1}$ have, respectively,
eigenvalues of orders $2,2,r$ which are equal to the lengths of
branches of $D_{r+2}$.


6) Case $\tilde{D}_4$. Here, by (\ref{chi_D_r_plus2}) we have
\begin{equation*}
\begin{split}
  \mathcal{X}(\tilde{D}_{4}) =
         & \mathcal{X}({D}_4)\mathcal{X}({A}_1) -
          \lambda\mathcal{X}({A}_1)^3 =
          -(\lambda + 1)^2(\lambda^3 + 1) + \lambda(\lambda + 1)^3 = \\
         & -(\lambda + 1)^2(\lambda^3 + 1 - \lambda^2 - \lambda) =
          -(\lambda + 1)^2(\lambda^2 - 1)(\lambda - 1) = \\
         & -(\lambda - 1)^2(\lambda + 1)^3,
         \end{split}
\end{equation*}
and so
\begin{equation}
          \chi = (\lambda - 1)^2\chi^3_1.
\end{equation}

Polynomial $\chi_1^3$ have eigenvalues of orders $2,2,2$ which are
equal to the lengths of branches of $D_4$.


\index{Frame formula} In (\ref{frame}), we saw that the simple
Frame formula is a particular case of the Subbotin-Sumin formula
(\ref{split_edge_eq}). Now we will show another {\it Frame
formula} from \cite[p.785]{Fr51} used by H.~S.~M.~Coxeter in
\cite{Cox51}.
\begin{proposition} [\cite{Fr51}]
\label{another_Frame}
  Let $T_{p,q,r}$ be a connected graph with three branches of lengths
$p, q, r$, where the lengths include the branch point, see
Fig.~\ref{T_pqr_diagram}. Then
\begin{equation}
 \label{frame_formula_2}
   \chi(T_{p,q,r}) =
     \chi({A}_{p+q+r-2})
     - \lambda^2\chi({A}_{p-2})\chi({A}_{q-2})\chi({A}_{r-2}).
\end{equation}
\end{proposition}
\PerfProof
  By (\ref{split_edge_eq}) we have
\begin{equation}
   \chi(T_{p,q,r}) =
     \chi({A}_{p+q-1})\chi({A}_{r-1})
     - \lambda\chi({A}_{p-1})\chi({A}_{q-1})\chi({A}_{r-2}).
\end{equation}
This is equivalent to
\begin{equation}
 \begin{split}
   \chi(T_{p,q,r}) = &
    \chi({A}_{p+q-1})\chi({A}_{r-1})
     - \lambda\chi({A}_{p+q-2})\chi({A}_{r-2}) \\
  &  + \lambda\chi({A}_{p+q-2})\chi({A}_{r-2})
     - \lambda\chi({A}_{p-1})\chi({A}_{q-1})\chi({A}_{r-2}).
 \end{split}
\end{equation}
By the first Frame formula (\ref{frame}) we have
\begin{equation}
   \chi({A}_{p+q+r-2}) =
     \chi({A}_{p+q-1})\chi({A}_{r-1})
     - \lambda\chi({A}_{p+q-2})\chi({A}_{r-2}),
\end{equation}
and
\begin{equation}
 \label{A_pqr_2}
 \begin{split}
  & \chi(T_{p,q,r}) = \\
  &   \chi({A}_{p+q+r-2})
   + \lambda\chi({A}_{p+q-2})\chi({A}_{r-2})
     - \lambda\chi({A}_{p-1})\chi({A}_{q-1})\chi({A}_{r-2}) = \\
  &  \chi({A}_{p+q+r-2}) +
     \lambda\chi({A}_{r-2})
     [\chi({A}_{p+q-2}) - \chi({A}_{p-1})\chi({A}_{q-1})].
 \end{split}
\end{equation}
Again, by  (\ref{frame}) we have
\begin{equation}
  \chi({A}_{p+q-2}) - \chi({A}_{p-1})\chi({A}_{q-1}) =
     - \lambda\chi({A}_{p-2})\chi({A}_{q-2}),
\end{equation}
and from (\ref{A_pqr_2}) we deduce
\begin{equation}
   \chi(T_{p,q,r}) =
     \chi({A}_{p+q+r-2})
     - \lambda^2\chi({A}_{r-2})\chi({A}_{p-2})\chi({A}_{q-2}). \quad \qed
\end{equation}

\section{Characteristic polynomials for the Dynkin diagrams}
 \label{case_Dynkin_diagr}
\index{characteristic polynomial of the Coxeter transformation! -
Dynkin diagrams} In order to calculate characteristic polynomials
of the Coxeter transformations for the Dynkin diagrams, we use the
Subbotin-Sumin formula (\ref{split_edge_eq}), the generalized
Subbotin-Sumin formula (\ref{split_edge_eq_2}), its specialization
(\ref{sub_sum_2}) for arbitrary $(p,q,r)$-trees, and particular
cases of (\ref{sub_sum_2}): formula (\ref{chi_p_eq_q}) for $p=q$,
and (\ref{chi_q_eq_2p}) for $p=2q$.

The results of calculations of this section concerning Dynkin diagrams
are collected in Table \ref{table_char_polynom_Dynkin}.

1) Case $p=q=3$, $r=2$ (${E}_6$). From (\ref{chi_p_eq_q}) we have
\begin{equation}
  \begin{split}
   \chi =
   & \frac{(\lambda^{3} -1)}{(\lambda - 1)^3}
      [(\lambda^{3} + 1)(\lambda^{2} - 1)(\lambda - 1) -
          (\lambda^{3} - 1) (\lambda^{2} - \lambda)] = \\
   &     \frac{(\lambda^{3} - 1)}{(\lambda-1)^3}(\lambda-1)^2
      [(\lambda^{3} + 1)(\lambda+1) -
          \lambda(\lambda^{2} + \lambda + 1)] = \\
   &  \frac{(\lambda^{3} - 1)}{(\lambda - 1)}
      (\lambda^{4} -\lambda^2 + 1) = \\
   &    (\lambda^{2} + \lambda + 1) (\lambda^{4} -\lambda^2 + 1) =
      \lambda^{6} +  \lambda^{5} - \lambda^3 + \lambda + 1.
  \end{split}
\end{equation}
In another form, we have
\begin{equation}
    \chi = \frac{(\lambda^{3} - 1)}{(\lambda - 1)}
      (\lambda^{4} -\lambda^2 + 1) =
      \frac{(\lambda^{6} + 1)}{(\lambda^2 + 1)}
      \frac{(\lambda^{3} - 1)}{(\lambda - 1)}.
\end{equation}

2) Case $p=2$, $q=4$, $r=3$ (${E}_7$).
   From (\ref{chi_q_eq_2p}) we have,  up to a sign,
\begin{equation}
  \begin{split}
   \chi =   & \frac{(\lambda^{2} -1)}{(\lambda-1)^3}
      [(\lambda^{4} + \lambda^{2} + 1)(\lambda^{3} - 1)(\lambda-1) -
          (\lambda^{4} - 1) (\lambda^{3} - \lambda)] = \\
            & \frac{(\lambda^{2} -1)(\lambda - 1)^2}{(\lambda-1)^3}
      [(\lambda^{4} + \lambda^{2} + 1)(\lambda^{2} + \lambda + 1) - \\
      & \qquad \qquad \qquad \qquad
        (\lambda^{3} + \lambda^{2} + \lambda + 1)\lambda(\lambda + 1)] =  \\
            & (\lambda + 1)
      [(\lambda^{4} + \lambda^{2} + 1)(\lambda^{2} + \lambda + 1) -
          (\lambda^{3} + \lambda^{2} + \lambda + 1)(\lambda^2 +
\lambda)] = \\
        & (\lambda + 1)(\lambda^6 - \lambda^3 + 1) =
          \lambda^7 + \lambda^6 - \lambda^4 - \lambda^3 + \lambda + 1.
   \end{split}
\end{equation}
   In another form, we have
\begin{equation}
    \chi = -(\lambda + 1)(\lambda^6 - \lambda^3 + 1) =
           -\frac{(\lambda + 1)(\lambda^{9} + 1)}{(\lambda^3 + 1)}.
\end{equation}

3) Case $p=3$, $q=2$, $r=5$ (${E}_8$). From (\ref{sub_sum_2}) we
have
\begin{equation}
  \begin{split}
   \chi =   & \frac{(\lambda^5 -1)^2}{(\lambda-1)^2}
      -  \frac{(\lambda^3 -1)(\lambda^2 -1)(\lambda^5 - \lambda)}
               {(\lambda-1)^3} = \\
      & (\lambda^4 + \lambda^3 + \lambda^2 + \lambda + 1)^2 - \\
      & \qquad  \qquad ( \lambda^2 + \lambda + 1)(\lambda + 1)
             \lambda(\lambda^3 + \lambda^2 + \lambda + 1) = \\
        &  \lambda^8 + \lambda^7 - \lambda^5 - \lambda^4 - \lambda^3 +
\lambda + 1.
   \end{split}
\end{equation}
Since
\begin{equation*}
    (\lambda^8 + \lambda^7 - \lambda^5 - \lambda^4 - \lambda^3 + \lambda +
1)
    (\lambda^2 - \lambda + 1) = \lambda^{10} - \lambda^5 + 1,
\end{equation*}
we have
\begin{equation}
   \chi =   \frac{\lambda^{10} - \lambda^5 + 1}{\lambda^2 - \lambda + 1} =
            \frac{(\lambda^{15} + 1)(\lambda + 1)}{(\lambda^5 + 1)(\lambda^3
 + 1)}.
\end{equation}

4) For the case ${D}_n$, see (\ref{chi_D_r_plus2}),
\begin{equation}
 \label{D_n}
  \begin{split}
    \mathcal{X}({D}_n) = &  (-1)^n(\lambda + 1)(\lambda^{n-1} + 1),  \\
    \chi = & (\lambda + 1)(\lambda^{n-1} + 1).
  \end{split}
\end{equation}

5) Case ${F}_4$. We use formula (\ref{split_edge_eq_2}), splitting
the diagram ${F}_4$ along the weighted edge into two diagrams
$\varGamma_1 = {A}_2$ and $\varGamma_2 = {A}_2$. Here, $\rho = 2$.
\begin{equation*}
   \chi =
   \chi^2_2 - 2\lambda\chi^2_1 =
   (\lambda^2 + \lambda + 1)^2 - 2\lambda(\lambda + 1)^2 =
   \lambda^4 - \lambda^2 + 1.
\end{equation*}
   In another form, we have
\begin{equation}
    \chi = \lambda^4 - \lambda^2 + 1 =
           \frac{\lambda^{6} + 1}{\lambda^2 + 1}.
\end{equation}

6) Case ${G}_2$. A direct calculation of the Coxeter
transformation gives
\begin{equation*}
   {\bf C} =
     \left (
     \begin{array}{cc}
        -1 & 3   \\
         0 & 1   \\
     \end{array}
     \right )
     \left (
     \begin{array}{cc}
        1 & 0    \\
        1 & -1   \\
     \end{array}
     \right ) =
     \left (
     \begin{array}{cc}
        2 & -3   \\
        1 & -1   \\
     \end{array}
     \right ),
\end{equation*}
and
\begin{equation}
    \chi = \lambda^2 - \lambda + 1 =
           \frac{\lambda^{3} + 1}{\lambda + 1}.
\end{equation}

7) The dual cases ${B}_n$ and ${C}_n = {B}^{\vee}_n$. Since
spectra of the Coxeter transformations of the dual graphs
coincide, we can consider only the case ${B}_n$.

a) Consider ${B}_2$. A direct calculation of the Coxeter
transformation gives
\begin{equation*}
   {\bf C} =
     \left (
     \begin{array}{cc}
        -1 & 2   \\
         0 & 1   \\
     \end{array}
     \right )
     \left (
     \begin{array}{cc}
        1 & 0    \\
        1 & -1   \\
     \end{array}
     \right ) =
     \left (
     \begin{array}{cc}
        1 & -2   \\
        1 & -1   \\
     \end{array}
     \right ),
\end{equation*}
and
\begin{equation}
  \label{B_2}
    \mathcal{X}({B}_2) = \lambda^2  + 1.
\end{equation}

b) Consider ${B}_3$. We use formula (\ref{split_edge_eq_2}),
splitting the diagram ${B}_3$ along the weighted edge
\begin{equation}
  \label{B_3}
  \begin{split}
    \mathcal{X}({B}_3) = & \mathcal{X}({A}_2)\mathcal{X}({A}_1)
            - 2\lambda\mathcal{X}({A}_1) =
            -(\chi_2\chi_1  - 2\lambda\chi_1) = \\
           & -\chi_1(\chi_2 - 2\lambda) =
           -(\lambda + 1)(\lambda^2 + \lambda + 1 - 2\lambda) = \\
           & -(\lambda + 1)(\lambda^2 - \lambda + 1) =
             -(\lambda^3 + 1).
  \end{split}
\end{equation}

c) Case ${B}_n$ for $n \geq 4$.
\begin{equation*}
  \begin{split}
    & \mathcal{X}({B}_n) = \\
    & \chi({A}_{n-3})\chi({B}_3)  -
              \lambda\chi({A}_{n-4})\chi({B}_2) =
             \chi_{n-3}\chi({B}_3)  -
              \lambda\chi_{n-4}\chi({B}_2),
  \end{split}
\end{equation*}
  i.e.,
\begin{equation*}
  \begin{split}
    & (-1)^n\mathcal{X}({B}_n) = \\
           & \frac{\lambda^{n-2} - 1}{\lambda - 1}(\lambda^3 + 1)
    - \lambda\frac{\lambda^{n-3} - 1}{\lambda - 1}(\lambda^2 + 1) =    \\
           & \frac{(\lambda^{n-2} - 1)(\lambda^3 + 1) -
                   (\lambda^{n-3} - 1)(\lambda^3 + \lambda)}
                  {\lambda - 1}   =  \\
           &  \frac{\lambda^{n+1} - \lambda^{n} + \lambda - 1}
              {\lambda - 1}   = \lambda^{n} + 1,
  \end{split}
\end{equation*}
and
\begin{equation}
 \label{B_n}
  \begin{split}
      \mathcal{X}({B}_n) = & (-1)^n(\lambda^{n} + 1).
  \end{split}
\end{equation}
By (\ref{B_2}), (\ref{B_3}) and (\ref{B_n}), we see,
 for the case ${B}_n$, that,  up to a sign,
\begin{equation}
  \chi =  \lambda^{n} + 1.
\end{equation}

\section{A generalization of R.~Steinberg's theorem}
\index{theorem! - Steinberg, multiply-laced case}
 \label{Steinberg_multy}
We will show now that R.~Steinberg's theorem \ref{Steinberg} can be
extended to the multiply-laced case.

\index{branch point}
\index{non-homogeneous point}

\subsection{Folded Dynkin diagrams and branch points}
\label{folded_diag}
\begin{definition} {\rm
A vertex is said to be a {\it branch point} of the Dynkin diagram
in one of the following cases:

(a) if it is the common endpoint of three edges
   (${E}_6, {E}_7, {E}_8, {D}_n$);

(b) if it is the common endpoint of two edges, one of which is
non-weighted and one is weighted (${B}_n, {C}_n, {F}_4$). Such a
vertex is said to be a {\it non-homogeneous point}.

(c) Let, by definition, both points of ${G}_2$ be branch points.
  }
\end{definition}

\begin{figure}[h]
\centering
\includegraphics{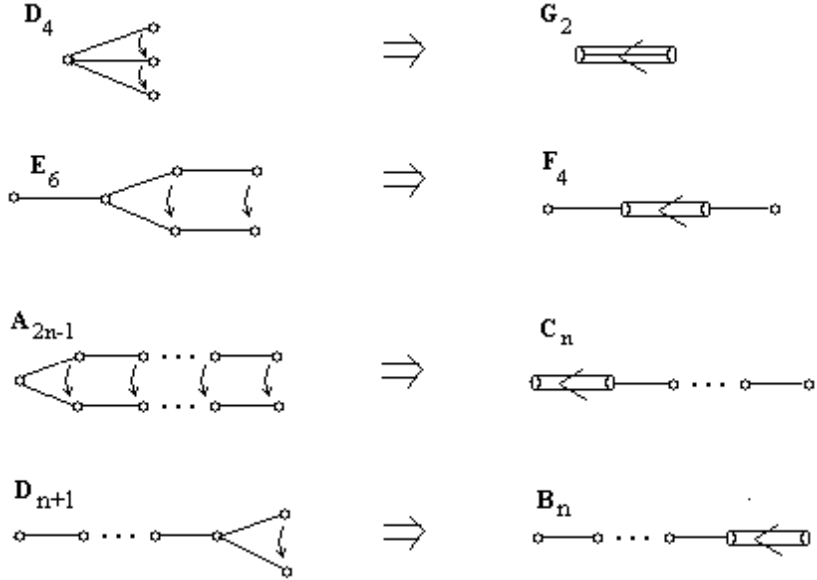}
\caption{\hspace{3mm} Folding operation applied to Dynkin
diagrams}
\label{folding}
\end{figure}

\begin{remark} {\rm
\label{remark_folding} \index{folding operation} \index{folded
Dynkin diagrams} Every multiply-laced Dynkin diagram (and also
every extended Dynkin diagram) can be obtained by a so-called {\it
folding} operation from a simply-laced diagrams, see
Fig.~\ref{folding}, such that the branch point is transformed into
the {\it non-homogeneous point}, see, e.g., I.~Satake,
\cite[p.109]{Sat60}, J.~Tits, \cite[p.39]{Ti66}, P.~Slodowy,
\cite[Appendices I and III]{Sl80}, or more recent works
\cite{FSS96}, \cite{Mohr04}. This fact was our motivation for
considering {\it non-homogeneous points} from (b) as {\it branch
points}.
   }
\end{remark}

Let us divide all extended Dynkin diagrams (trees) into 4 classes 
and let $g$ be the number of the class:

Class $g = 0$ contains all simply-laced diagrams.

Class $g = 1$ contains multiply-laced diagrams which has only one weighted edge
of type (1,2) or (2,1).

Class $g = 2$ contains multiply-laced diagrams with either two weighted edges or
one weighted edge of type (1,3) or (3,1).

Class $g = 3$ contains multiply-laced diagrams with two
vertices.

Thus,
\begin{equation}
\label{class_g}
 \begin{split}
  g = 0 & \text{ for } \tilde{E}_6, \tilde{E}_7, \tilde{E}_8,
                        \tilde{D}_n, \\
  g = 1 & \text{ for }
            \tilde{F}_{41}, \tilde{F}_{42},
            \widetilde{CD}_n, \widetilde{DD}_n,  \\
  g = 2 &  \text{ for }  \tilde{G}_{12}, \tilde{G}_{22},
             \tilde{B}_{n}, \tilde{C}_{n},
             \widetilde{BC}_{n}, \\
  g = 3 &  \text{ for } \tilde{A}_{11}, \tilde{A}_{12}.
  \end{split}
\end{equation}

\subsection{R.~Steinberg's theorem for the multiply-laced case}
  \label{sect_gen_Steinberg}
Now, we can generalize R.~Steinberg's theorem for the
multiply-laced case.

 \index{affine Coxeter transformation}
 \index{characteristic polynomial of the Coxeter transformation! - extended Dynkin diagrams}
 \index{theorem! - Steinberg, multiply-laced case}
\begin{theorem}
 \label{gen_Steinberg}
The affine Coxeter transformation with the extended Dynkin diagram
$\tilde{\varGamma}$ of class $g$ (see (\ref{class_g})) has the
same eigenvalues as the product of the Coxeter transformations of
types ${A}_i$, where $i \in \{p-1, q-1, r-1\}$ and matches to
$(3-g)$ branches of the Dynkin diagram $\varGamma$. In other
words,
\begin{equation}
\label{groups_of_chi}
 \begin{split}
  & \text{ For } g = 0, \text{ the product }
    \chi_{p-1}\chi_{q-1}\chi_{r-1} \text{ is taken}. \\
  & \text{ For } g = 1, \text{ the product }
    \chi_{p-1}\chi_{q-1} \text{ is taken}.  \\
  & \text{ For } g = 2, \text{ the product consists only of }
     \chi_{p-1}. \\
  & \text{ For } g = 3, \text{ the product is trivial (= 1)}. \\
 \end{split}
\end{equation}
The remaining two eigenvalues of the affine Coxeter transformation
are both equal to 1, see Table
\ref{table_char_ext_polynom_Dynkin}.
\end{theorem}

\PerfProof
 For $g = 0$, we have a simply-laced case considered
in R.~Steinberg's theorem (Theorem \ref{Steinberg}).

\underline{Let $g = 1$}.

1) Cases $\tilde{F}_{41}$, $\tilde{F}_{42}$. Since the spectra of
the characteristic polynomials of the Coxeter transformations of
the dual graphs coincide, we consider $\tilde{F}_{41}$. We use
formula (\ref{split_edge_eq_2}), splitting
  the diagram $\tilde{F}_{41}$ along the weighted edge. We have
\begin{equation}
 \begin{split}
  \mathcal{X}(\tilde{F}_{41}) = & -(\chi_3\chi_2 - 2\lambda\chi_2\chi_1)=
          -\chi_2(\chi_3 - 2\lambda\chi_1) =  \\
         & -\chi_2(\lambda^3 + \lambda^2 + \lambda + 1 -
               2\lambda(\lambda + 1)) = \\
         & -\chi_2(\lambda^3 - \lambda^2 - \lambda + 1) = \\
         &  -\chi_2(\lambda^2 - 1)(\lambda - 1) =
           -(\lambda - 1)^2\chi_2\chi_1,
  \end{split}
\end{equation}
and, up to a sign, we have
\begin{equation}
     \chi = (\lambda - 1)^2\chi_2\chi_1.
\end{equation}
Polynomials $\chi_1$ and $\chi_2$ have, respectively, eigenvalues
of orders $2$ and $3$ which are equal to the lengths of the
branches ${A}_1$ and ${A}_2$ of ${F}_{4}$ without non-homogeneous
branch point.

2) Cases $\widetilde{CD}_n$, $\widetilde{DD}_n$. These diagrams
are obtained as extensions of ${B}_n$, see \cite[Tab.II]{Bo}. By
(\ref{split_edge_eq_2}) and by splitting the diagram
$\widetilde{CD}_n$ along the weighted edge, we have
\begin{equation*}
   \mathcal{X}(\widetilde{CD}_n) =
      \mathcal{X}({D}_n)\mathcal{X}({A}_1)
      - 2\lambda\mathcal{X}({D}_{n-1}).
\end{equation*}
By (\ref{D_n})
\begin{equation*}
   \mathcal{X}({D}_n) = (-1)^n(\lambda + 1)(\lambda^{n-1} + 1).
\end{equation*}
Thus,
\begin{equation}
 \begin{split}
  (-1)^{n-1}\mathcal{X}(\widetilde{CD}_n) = & (\lambda + 1)^2(\lambda^{n-1}
+ 1) -
          2\lambda(\lambda + 1)(\lambda^{n-2} + 1) = \\
          & (\lambda + 1)(\lambda^n + \lambda + \lambda^{n-1} + 1
           -2\lambda^{n-1} - 2\lambda) =  \\
          & (\lambda + 1)(\lambda^n - \lambda^{n-1} - \lambda + 1) = \\
          & (\lambda + 1)(\lambda - 1)(\lambda^{n-1} - 1) =
           (\lambda - 1)^2\chi_1\chi_{n-2},
  \end{split}
\end{equation}
and, up to a sign, we have
\begin{equation}
     \chi = (\lambda - 1)^2\chi_1\chi_{n-2}.
\end{equation}
Polynomials $\chi_1$ and $\chi_{n-2}$ have, respectively,
eigenvalues of orders $2$ and $3$ which are equal to the lengths
of the branches ${A}_1$ and ${A}_{n-2}$ of ${B}_{n}$ without
non-homogeneous branch point.

\underline{Let $g = 2$}.

3) Cases $\tilde{G}_{12}$, $\tilde{G}_{22}$. Since the spectra of
the Coxeter transformations of the dual graphs coincide, we
consider $\tilde{G}_{12}$. By (\ref{split_edge_eq_2}), splitting
  the diagram $\tilde{G}_{12}$ along the weighted edge we have
\begin{equation}
 \begin{split}
  \mathcal{X}(\tilde{G}_{12}) = & -(\chi_2\chi_1 - 3\lambda\chi_1) =
         -\chi_1(\chi_2 - 3\lambda) = \\
         & -(\lambda + 1)(\lambda^2 + \lambda + 1 - 3\lambda) = \\
         & -(\lambda - 1)^2(\lambda + 1) =
           -(\lambda - 1)^2\chi_1,
  \end{split}
\end{equation}
and, up to a sign, we have
\begin{equation}
     \chi = (\lambda - 1)^2\chi_1.
\end{equation}

Polynomial $\chi_1$ has the single eigenvalue of order $2$, it
corresponds to the length of the single branch ${A}_1$ of
${G}_{2}$ without non-homogeneous branch point.

4) Cases $\tilde{B}_{n}$, $\tilde{C}_{n}$,
             $\widetilde{BC}_{n}$.
The characteristic polynomials of the Coxeter transformations of
these diagrams coincide. Consider $\tilde{B}_{n}$. By
(\ref{split_edge_eq_2}), splitting
  the diagram $\tilde{B}_{n}$ along the weighted edge we have
\begin{equation}
 \begin{split}
  \mathcal{X}(\tilde{B}_{n}) =
   & \mathcal{X}({B}_n)\mathcal{X}({A}_1)
     - 2\lambda\mathcal{X}({B}_{n-1}),  \text{ i.e., } \\
  (-1)^{n-1}\mathcal{X}(\tilde{B}_{n}) = &
         (\lambda^n + 1)(\lambda + 1) - 2\lambda(\lambda^{n-1} + 1) = \\
         & \lambda^{n+1} + \lambda^n + \lambda + 1
           - 2\lambda^{n-1} -2\lambda = \\
         & \lambda^{n+1} - \lambda^n - \lambda + 1 =
           (\lambda - 1)(\lambda^n - 1),
  \end{split}
\end{equation}
and, up to a sign, we have
\begin{equation}
     \chi = (\lambda - 1)^2\chi_{n-1}.
\end{equation}

\vspace{3mm} Polynomial $\chi_{n-1}$ has the single eigenvalue of
order $n$, it corresponds to the length of the single branch
${A}_{n-1}$ of ${B}_{n}$ without non-homogeneous branch point.

\vspace{3mm}
\underline{Let $g = 3$}.

5) Cases $\tilde{A}_{11}$, $\tilde{A}_{12}$. Direct calculation of
the Coxeter transformation gives \vspace{3mm}
\begin{equation*}
 \begin{split}
 & \text{ for } \tilde{A}_{11} :
   {\bf C} =
     \left (
     \begin{array}{cc}
        -1 & 4   \\
         0 & 1   \\
     \end{array}
     \right )
     \left (
     \begin{array}{cc}
        1 & 0    \\
        1 & -1   \\
     \end{array}
     \right ) =
     \left (
     \begin{array}{cc}
        3 & -4   \\
        1 & -1   \\
     \end{array}
     \right ), \vspace{5mm} \\
 & \text{ for } \tilde{A}_{12}:
   {\bf C} =
     \left (
     \begin{array}{cc}
        -1 & 2   \\
         0 & 1   \\
     \end{array}
     \right )
     \left (
     \begin{array}{cc}
        1 & 0    \\
        2 & -1   \\
     \end{array}
     \right ) =
     \left (
     \begin{array}{cc}
        3 & -2   \\
        2 & -1   \\
     \end{array}
     \right ),
   \end{split}
\end{equation*}
and in both cases we have
\begin{equation}
    \chi = (\lambda - 3)(\lambda + 1) + 4 = (\lambda - 1)^2. \qed \vspace{5mm}
\end{equation}

\section{The Kostant generating function and Poincar\'{e} series}
\label{kostant}
\subsection{The generating function}
\label{generating_fun} \index{Kostant generating function}
\index{Poincar\'{e} series} \index{orbit structure of the Coxeter
transformation}

Let ${\rm Sym}(\mathbb{C}^2)$ be the symmetric algebra over
$\mathbb{C}^2$, in other words, ${\rm Sym}(\mathbb{C}^2) =
\mathbb{C}[x_1, x_2]$, see (\ref{R_algebra}). The symmetric
algebra ${\rm Sym}(\mathbb{C}^2)$ is a graded
$\mathbb{C}$-algebra, see (\ref{R_decomp}):
\begin{equation}
    {\rm Sym}(\mathbb{C}^2) =
\mathop{\oplus}\limits_{n=0}^{\infty}{\rm Sym}^n(\mathbb{C}^2).
\end{equation}

Let $\pi_n$ be the representation of $SU(2)$ in ${\rm
Sym}^n(\mathbb{C}^2)$ induced by its action on $\mathbb{C}^2$. The
set $\{\pi_n\}$, where $n = 0,1,...$ is the set of all irreducible
representations of $SU(2)$, see, e.g., \cite[\S37]{Zhe73}. Let $G$
be any finite subgroup of $SU(2)$, see \S(\ref{finite_subgroups}).
B.~Kostant in \cite{Kos84} considered the question:

\medskip
{\sl how does $\pi_n | G$ decompose for any $n \in \mathbb{N}$?}
\medskip

The answer --- the decomposition $\pi_n | G$ --- is as follows:
\begin{equation}
  \label{decomp_pi_n}
   \pi_n | G = \sum\limits_{i=0}^r{m_i(n)\rho_i},
\end{equation}
where $\rho_i$ are irreducible representations of $G$, considered
in the context of {\it McKay correspondence}, see \S(\ref{McKay}).
Thus, the decomposition (\ref{decomp_pi_n}) reduces the question
to the following one:

\medskip
{\sl what are the multiplicities $m_i(n)$ equal to?}
\medskip

 \index{Lie algebra}
B.~Kostant in \cite{Kos84} obtained the multiplicities $m_i(n)$ by
means the orbit structure of the Coxeter transformation on the
highest root of the corresponding Lie algebra. For further details
concerning this orbit structure and the multiplicities $m_i(n)$,
see \S\ref{orbit_str}.

Note, that multiplicities $m_i(n)$ in (\ref{decomp_pi_n}) are
calculated as follows:
\begin{equation}
  \label{multipl_n}
     m_i(n) = <\pi_n | G , \rho_i>,
\end{equation}
(for the definition of the inner product $<\cdot , \cdot>$, see
(\ref{inner_prod}) ).
\index{induced representation}
\begin{remark}
  \label{remark_decomp}
 {\rm
For further considerations, we extend the relation for
multiplicity (\ref{multipl_n}) to the cases of {\it restricted
representations}
  $\rho_i^\downarrow: = \rho_i\downarrow_H^G$
  and {\it induced representations}
   $\rho_i^\uparrow: = \rho_i\uparrow_H^G$, where $H$ is any subgroup of $G$
  (see \S\ref{section_slodowy}):

\begin{equation}
  \label{decomp_pi_n_2}
   m_i^\downarrow(n) = <\pi_n | H ,\rho_i^\downarrow>,   \hspace{3mm}
   m_i^\uparrow(n) = <\pi_n | G , \rho_i^\uparrow>.
\end{equation}

We do not have any decomposition like (\ref{decomp_pi_n}) neither
for restricted representations $\rho_i^\downarrow$ nor for induced
representations $\rho_i^\uparrow$. Nevertheless, we will sometimes
denote both multiplicities $m_i^\downarrow(n)$ and
$m_i^\uparrow(n)$ in (\ref{decomp_pi_n_2}) by $m_i(n)$ as in
(\ref{decomp_pi_n}). }
\end{remark}

\index{McKay correspondence}
\index{multiplicities $m_i(n)$}

\begin{remark} {\rm
 \label{triv_repr}
1) A representation $\rho: G \longrightarrow GL_k(V)$ defines a
$k$-linear action $G$ on $V$ by \index{trivial representation}
\begin{equation}
  gv = \rho(g)v.
\end{equation}
The pair $(V,\rho)$ is called a $G$-{\it module}. The case where
$\rho(g) = Id_V$ is called the {\it trivial representation} in
$V$. In this case
\begin{equation}
 \label{trivial}
  gv = v \text{ for all } g \in V.
\end{equation}
In (\ref{decomp_pi_n}), the trivial representation $\rho_0$
corresponds to a particular vertex (see \cite{McK80}),
which extends the Dynkin diagram to the extended Dynkin diagram.

2) Let $\rho_0(H)$ (resp. $\rho_0(G)$) be the trivial
representation of any subgroup $H \subset G$ (resp. of group $G$).
The trivial representation $\rho_0(H): H \longrightarrow GL_k(V)$
coincides with the {\it restricted representation}
$\rho_0\downarrow_H^G: G \longrightarrow GL_k(V)$, and the trivial
representation $\rho_0(G): G \longrightarrow GL_k(V)$ coincides
with the {\it induced representation} $\rho_0\uparrow_H^G: H
\longrightarrow GL_k(V)$.}
\end{remark}

Since there is one-to-one correspondence between the $\rho_i$ and
the vertices of the Dynkin diagram, we can define (see
\cite[p.211]{Kos84}) the vectors $v_n$, where $n\in\mathbb{Z}_+$,
as follows:
\begin{equation}
 \label{def_vn}
 v_n  = \sum\limits_{i=0}^r{m_i(n)}\alpha_i,
     \text{ where } \pi_n | G = \sum\limits_{i=0}^r{m_i(n)\rho_i}, 
\end{equation}
where $\alpha_i$ are simple roots of the corresponding extended
Dynkin diagram. Similarly, for  the multiply-laced case, we define
vectors $v_n$ to be:
\begin{equation}
 \label{def_vn_1}
   v_n  = \sum\limits_{i=0}^r{m_i^\uparrow(n)}\alpha_i \hspace{5mm} \text{
   or } \hspace{5mm}
  v_n  = \sum\limits_{i=0}^r{m_i^\downarrow(n)}\alpha_i,
\end{equation}
where the multiplicities $m_i^\uparrow(n)$ and $m_i^\downarrow(n)$
are defined by (\ref{decomp_pi_n_2}). The vector $v_n$ belongs to
the root lattice generated by simple roots. Following B.~Kostant,
we define the generating function $P_G(t)$ for cases
(\ref{def_vn}) and (\ref{def_vn_1}) as follows:
\begin{equation}
 \label{Kostant_gen_func}
    P_G(t) = ([P_G(t)]_0, [P_G(t)]_1, \ldots , [P_G(t)]_r)^t
           := \sum\limits_{n=0}^{\infty}v_n{t^n},
\end{equation}
the components of the vector $P_G(t)$ being the following series
\begin{equation}
 \label{gen_func_i}
    [P_G(t)]_i = \sum\limits_{n=0}^{\infty}\tilde m_i(n){t^n},
\end{equation}
where $i = 0,1,\dots,r$ and $\tilde m_i(n)$ designates $m_i(n),
m_i^\uparrow(n)$ or $m_i^\downarrow(n)$, see
\S\ref{multipl_kostant}.  In particular, for $i=0$, we have
\index{algebra of invariants} \index{symmetric algebra}
\begin{equation}
    [P_G(t)]_0 = \sum\limits_{n=0}^{\infty}m_0(n){t^n},
\end{equation}
where $m_0(n)$ is the multiplicity of the trivial representation
$\rho_0$ (see Remark \ref{triv_repr}) in ${\rm
Sym}^n(\mathbb{C}^2)$. By \S\ref{sect_poincare} the {\it algebra
of invariants} $R^G$ is a subalgebra of the {\it symmetric
algebra} ${\rm Sym}(\mathbb{C}^2)$. Thanks to (\ref{trivial}), we
see that $R^G$ coincides with ${\rm Sym}(\mathbb{C}^2)$, and
$[P_G(t)]_0$ is the Poincar\'{e} series of the algebra of
invariants ${\rm Sym}(\mathbb{C}^2)^G$, i.e.,
\begin{equation}
  \label{poincare_alg_inv}
    [P_G(t)]_0 = P({\rm Sym}(\mathbb{C}^2)^G,t).
\end{equation}
(see \cite[p.221, Rem.3.2]{Kos84}).
\begin{remark}
\label{poincare_m_case} {\rm According to Remark \ref{triv_repr},
heading 2), we have
\begin{equation}
  \label{poincare_alg_inv_2}
    [P_H(t)]_0 = P({\rm Sym}(\mathbb{C}^2)^{\rho_0\downarrow_H^G},t), \hspace{7mm}
    [P_G(t)]_0 = P({\rm Sym}(\mathbb{C}^2)^{\rho_0\uparrow_H^G},t).
\end{equation}
We will need this fact in the proof of W.~Ebeling's theorem for
the multiply-laced case
  in \S\ref{ebeling},
  (see also Remark \ref{remark_decomp} and \S\ref{section_slodowy}).
  }
\end{remark}
The following theorem gives a remarkable formula for calculating
the Poincar\'{e} series for the binary polyhedral groups.
The theorem is known in different forms. B.~Kostant in \cite{Kos84}
shows it in the context of the Coxeter number $h$.

\index{theorem! - Kostant-Kn\"{o}rrer-Gonsales-Verdier}
{\it Theorem K-K-GV (Kostant, Kn\"{o}rrer, Gonsales-Sprinberg, Verdier)}
The Poincar\'{e} series $[P_G(t)]_0$ can be calculated as the following
rational function:
\begin{equation}
  \label{K_K_GV}
      [P_G(t)]_0 = \frac{1 + t^h}{(1 - t^a)(1 - t^b)},
\end{equation}
where
\begin{equation}
  \label{Kostant_numbers_a_b}
      b = h + 2 - a, \text{ and } ab = 2|G|.
\end{equation}
For a proof, see Theorem 1.4 and Theorem 1.8 from \cite{Kos84},
\cite[p.185]{Kn85}, \cite[p.428]{GV83}. We call the numbers $a$
and $b$ the {\it Kostant numbers}. They can be easily calculated,
see Table \ref{Kostant_numbers}, compare also with Table
\ref{binary_pol_1} and Table \ref{klein_rel}. Note, that $a = 2d$,
where $d$ is the maximal coordinate of the {\it nil-root} vector
from the kernel of the Tits form, see \S\ref{roots} and
Fig.~\ref{euclidean_diag}. \index{nil-roots} \index{Kostant
numbers} \index{Coxeter number} \index{binary polyhedral group}
\begin{table} [h]
  \centering
  \vspace{2mm}
  \caption{\hspace{3mm}The binary polyhedral groups (BPG) and the Kostant numbers
$a$, $b$}
   \renewcommand{\arraystretch}{1.5} 
  \begin{tabular} {||c|c|c|c|c|c||}
  \hline \hline
     \quad  Dynkin \quad    &
     \quad Order of \quad   &
     \quad BPG \quad        &
     \quad Coxeter \quad    &
     \quad $a$ \quad        &
     \quad $b$ \quad        \cr
      diagram      &  group   &            &  number        &    &   \\
  \hline \hline
      & & & & & \cr
      ${A}_{n-1}$  & $n$     & $\mathbb{Z}/n\mathbb{Z}$ & $n$ &  2 & $n$ \cr
      & & & & & \\
  \hline
      & & & & & \cr
      ${D}_{n+2}$  & $4n$    & $\mathcal{D}_n$ &   $2n+2$   & 4  & $2n$ \cr
      & & & & & \\
  \hline
      & & & & & \cr
      ${E}_6$  & 24  & $\mathcal{T}$   &   12           & 6  & 8 \cr
      & & & & & \\
  \hline
      & & & & & \cr
      ${E}_7$  & 48  & $\mathcal{O}$   &   18           & 8  & 12 \cr
      & & & & & \\
  \hline
      & & & & & \cr
      ${E}_8$  & 120 & $\mathcal{J}$   &   30           & 12 & 20 \cr
      & & & & & \\
  \hline  \hline
  \end{tabular}
    \label{Kostant_numbers}
\end{table}

\subsection{The characters and the McKay operator}
 \label{sect_char_McKay}
  Let $\chi_1, \chi_2,\dots \chi_r$ be all irreducible $\mathbb{C}$-characters
of a finite group $G$ corresponding to irreducible representations
$\rho_1, \rho_2, \dots, \rho_r$, and let $\chi_1$ correspond to
the trivial representation, i.e., $\chi_1(g) = 1$ for all $g \in
G$.

\index{hermitian inner product} All characters constitute the {\it
character algebra} $C(G)$ of $G$ since $C(G)$ is also a vector
space  over $\mathbb{C}$. An {\it hermitian inner product} $<\cdot
, \cdot>$ on $C(G)$ is defined as follows. For characters $\alpha,
\beta \in C(G)$, let
\begin{equation}
  \label{inner_prod}
   <\alpha, \beta> ~=~ \frac{1}{|G|}
     \sum\limits_{g \in G}\alpha(g)\overline{\beta(g)}
\end{equation}
Sometimes, we will write {\it inner product} $<\rho_i, \rho_j>$ of
the representations meaning actually the inner product of the
corresponding characters $<\chi_{\rho_i}, \chi_{\rho_j}>$.
\index{structure constants} \index{character ring}
\index{characters} Let $z_{ijk} = <\chi_i\chi_j, \chi_k>$, where
 $\chi_i\chi_j$ corresponds to the representation $\rho_i\otimes\rho_j$.
It is known that
$z_{ijk}$ is the multiplicity of the representation $\rho_k$  in
$\rho_i\otimes\rho_j$ and $z_{ijk} = z_{jik}$.
The numbers $z_{ijk}$ are integer and are called the {\it structure constants},
see, e.g., \cite[p.765]{Kar92}.

For every $i \in \{1,...,r\}$, there exists some $\stackrel{\wedge}i \in
\{1,...,r\}$
such that
\begin{equation}
  \chi_{\stackrel{\wedge}{i}}(g) = \overline{\chi_i(g)}
  \text{ for all } g \in G.
\end{equation}
\index{contragredient representation} The character
$\chi_{\stackrel{\wedge}{i}}$ corresponds to the {\it
contragredient representation} $\rho_{\stackrel{\wedge}{i}}$
determined from the relation
\begin{equation}
  \label{def_contragr}
    \rho_{\stackrel{\wedge}{i}}(g) = \rho_i(g)^{'{-1}}.
\end{equation}
We have
\begin{equation}
  \label{char_contragr}
   <\chi_i\chi_j, \chi_k> =
   <\chi_i, \chi_{\stackrel{\wedge}{j}}\chi_k>
\end{equation}
since
\begin{equation*}
 \begin{split}
  <\chi_i\chi_j, \chi_k> = &
   \frac{1}{|G|}\sum\limits_{g \in G}
  \chi_i(g)\chi_j(g)\overline{\chi_k(g)} = \\
  & \frac{1}{|G|}\sum\limits_{g \in G}
  \chi_i(g)(\overline{\overline{\chi_j(g)}\chi_k(g)}) =
   <\chi_i, \chi_{\stackrel{\wedge}{j}}\chi_k>.
 \end{split}
\end{equation*}
\index{unitarity}
\index{unimodularity}
\index{unitary unimodular matrix}
\begin{remark}
\label{group_su2}
{\rm
The group $SU(2)$ is the set of all unitary unimodular $2\times2$ matrices $u$,
  i.e.,
\begin{equation*}
  \label{unit_unim}
  \begin{split}
    &  u^* = u^{-1}  \text{ ({\it unitarity}) }, \\
    &  \det(u) = 1   \text{ ({\it unimodularity}) }.
  \end{split}
\end{equation*}
The matrices $u \in SU(2)$ have the following form:
\begin{equation}
  \label{su2_matr}
  u = \left (
       \begin{array}{cc}
         a     &  b \\
         -b^*  &  a^* \\
       \end{array}
      \right ), \text{ and }
   u^* =      \left (
       \begin{array}{cc}
         a^*  &  -b \\
         b^*  &  a \\
       \end{array}
      \right ),
      \text{ where } aa^* + bb^* = 1,
\end{equation}
see, e.g., \cite[Ch.9, \S6]{Ha89}. The mutually inverse matrices
  $u$ and $u^{-1}$ are
\begin{equation}
  \label{sl2_matr}
  u = \left (
       \begin{array}{cc}
         a  &  b \\
         c  &  d \\
       \end{array}
      \right ), \text{ and }
  u^{-1}=  \left (
       \begin{array}{cc}
         d  &  -b \\
         -c  &  a \\
       \end{array}
      \right ), \text{ where } ad - bc = 1.
\end{equation}

Set
\begin{equation}
  s = \left (
       \begin{array}{cc}
         0  &  -1 \\
         1  &  0 \\
       \end{array}
      \right ), \text{ then }
  s^{-1} = s^3 = \left (
       \begin{array}{cc}
         0  &  1 \\
        -1  &  0 \\
       \end{array}
      \right ).
\end{equation}
For  any $u \in SU(2)$, we have
\begin{equation}
 \label{elem_weyl}
  sus^{-1} = u^{~'{-1}}.
\end{equation}
\index{Weyl element}
The element $s$ is called the {\it Weyl element}.
}
\end{remark}

According to (\ref{def_contragr}) and (\ref{elem_weyl}) we see
that every finite dimensional representation of the group $SL(2,
\mathbb{C})$ (and hence, of $SU(2)$) is equivalent to its
contragredient representation, see \cite[\S37, Rem.3]{Zhe73}. Thus
by (\ref{char_contragr}), for representations $\rho_i$ of any
finite subgroup $G \subset SU(2)$,  we have
\begin{equation}
  \label{char_contragr_1}
   <\chi_i\chi_j, \chi_k> ~=~ <\chi_i, \chi_j\chi_k>.
\end{equation}

Relation (\ref{char_contragr_1}) holds also for characters of
restricted representations $\chi^\downarrow := \chi\downarrow^G_H$
and induced representations
 $\chi^\uparrow := \chi\uparrow^G_H$ (see Remark \ref{remark_decomp} and
\S\ref{section_slodowy}):
\begin{equation}
  \label{char_contragr_2}
\renewcommand{\arraystretch}{1.2}
\begin{array}{l}
<\chi_i\chi^\downarrow_j, \chi_k>_H ~=~ <\chi_i, \chi^\downarrow_j\chi_k>_H,\\
   <\chi_i\chi^\uparrow_j, \chi_k>_G ~=~ <\chi_i, \chi^\uparrow_j\chi_k>_G.
\end{array}
\end{equation}
Indeed, every restricted representation $\chi^\downarrow_j$ (resp.
induced representation $\chi^\uparrow_j$) is decomposed into the
direct sum of irreducible characters $\chi_s \in Irr(H)$ (here,
$Irr(H)$ is the set of irreducible characters of $H$) with some
integer coefficients $a_s$, for example:
\begin{equation}
  \label{char_contragr_3}
   \chi^\downarrow_j = \sum\limits_{\chi_s \in Irr(H)}{a_s}\chi_s,
\end{equation}
and since $\overline{a}_s = a_s$ for all $\chi_s \in Irr(H)$, we
have
\begin{equation}
\begin{split}
  \label{char_contragr_4}
   <\chi_i\chi^\downarrow_j, \chi_k>_H = &
    \sum\limits_{\chi_s \in Irr(H)} <a_s\chi_i\chi_s, \chi_k>_H = \\
   & \sum\limits_{\chi_s \in Irr(H)} a_s<\chi_i, \chi_s\chi_k>_H =  \\
   &  <\chi_i, \sum\limits_{\chi_s \in Irr(H)}a_s\chi_s\chi_k>_H =
    <\chi_i, \chi^\downarrow_j\chi_k>_H.
\end{split}
\end{equation}

\index{McKay-Slodowy correspondence} The matrix of multiplicities
$A := A(G)$ from (\ref{main_McKay}) was introduced by J.~McKay in
\cite{McK80}; it plays the central role in the {\it McKay
correspondence},
 see \S\ref{McKay}. We call this matrix --- or the corresponding operator ---
the {\it McKay matrix} or the {\it McKay operator}. \index{McKay
matrix}\index{McKay operator}

Similarly, let $\tilde{A}$ and $\tilde{A}^\vee$ be matrices
 of multiplicities (\ref{lab_slodowy_1}), (\ref{lab_slodowy_2}).
These matrices were introduced by P.~Slodowy \cite{Sl80} by
analogy with the McKay matrix for the multiply-laced case, see
\S\ref{section_slodowy}. We call these matrices the {\it Slodowy
operators}.

The following result of B.~Kostant \cite{Kos84}, which holds for
the McKay operator holds also for the Slodowy operators.

\begin{proposition}
  \label{kostant_prop}
  If $B$ is either the McKay operator $A$ or the Slodowy operator
  $\tilde{A}$ or $\tilde{A}^\vee$, then
 \begin{equation}
   \label{Kostant_relation}
      Bv_n = v_{n-1} + v_{n+1}.
 \end{equation}

\end{proposition}
\PerfProof
From now on
\begin{equation}
    \rho_i =
     \begin{cases}
        \rho_i &\text{ for } B = A, \vspace{2mm} \\
        \rho_i^\downarrow &\text{ for } B = \tilde{A}, \vspace{2mm} \\
        \rho_i^\uparrow &\text{ for } B = \tilde{A}^\vee, \\
     \end{cases}
 \hspace{7mm}
    m_i(n) =
     \begin{cases}
        m_i(n) &\text{ for } B = A, \vspace{2mm} \\
        m_i^\downarrow(n) &\text{ for } B = \tilde{A}, \vspace{2mm} \\
        m_i^\uparrow(n) &\text{ for } B = \tilde{A}^\vee. \\
     \end{cases}
\end{equation}

By (\ref{def_vn}), (\ref{def_vn_1}),
and by definition of the McKay operator (\ref{main_McKay})
and by definition of the Slodowy operator
(\ref{lab_slodowy_1})  and (\ref{lab_slodowy_2}), we have
\begin{equation}
  \label{McKay_oper}
      Bv_n = B
    \left (
    \begin{array}{c}
       m_0(n) \\
       \dots  \\
       m_r(n)
    \end{array}
    \right ) =
    \left (
    \begin{array}{c}
       \sum a_{0i}m_i(n) \\
       \dots  \\
       \sum a_{ri}m_i(n) \\
    \end{array}
    \right ) =
    \left (
    \begin{array}{c}
       \sum a_{0i}<\rho_i, \pi_n> \\
       \dots  \\
       \sum a_{ri}<\rho_i, \pi_n> \\
    \end{array}
    \right ).
\end{equation}
By (\ref{main_McKay}), (\ref{lab_slodowy_1}) and (\ref{lab_slodowy_2})
we have
\begin{equation*}
  \sum\limits_{i=1}^r a_{0i}<\rho_i, \pi_n> ~=~
  <\sum\limits_{i=1}^r a_{0i}\rho_i, \pi_n> ~=~
  <\rho_{reg}\otimes\rho_i, \pi_n>, \\
\end{equation*}
and from (\ref{McKay_oper}) we obtain
\begin{equation}
  \label{McKay_oper_1}
      Bv_n =
    \left (
    \begin{array}{c}
       <\rho_{reg}\otimes\rho_0, \pi_n> \\
       \dots  \\
       <\rho_{reg}\otimes\rho_r, \pi_n> \\
    \end{array}
    \right ).
\end{equation}
  Here $\rho_{reg}$ is the regular two-dimensional representation which
coincides with the representation $\pi_1$ in ${\rm
Sym}^2(\mathbb{C}^2)$ from \S\ref{generating_fun}. Thus,
\begin{equation}
  \label{McKay_oper_2}
   Bv_n =
    \left (
    \begin{array}{c}
       <\pi_1\otimes\rho_0, \pi_n> \\
       \dots  \\
       <\pi_1\otimes\rho_r, \pi_n> \\
    \end{array}
    \right ).
\end{equation}
From (\ref{char_contragr_1}), (\ref{char_contragr_4}) we obtain
\begin{equation}
  \label{McKay_oper_3}
   Bv_n =
    \left (
    \begin{array}{c}
       <\rho_0, \pi_1\otimes\pi_n> \\
       \dots  \\
       <\rho_r, \pi_1\otimes\pi_n> \\
    \end{array}
    \right ).
\end{equation}
By Clebsch-Gordan formula we have
\begin{equation}
  \label{McKay_oper_4}
   \pi_1\otimes\pi_n = \pi_{n-1} \oplus \pi_{n+1},
\end{equation}
where $\pi_{-1}$ is the zero representation, see
\cite[exs.3.2.4]{Sp77} or \cite[Ch.5, \S6,\S7]{Ha89}. From
(\ref{McKay_oper_3}) and (\ref{McKay_oper_4}) we have
(\ref{Kostant_relation}). \qedsymbol

For the following corollary, see \cite[p.222]{Kos84} and also
\cite[\S4.1]{Sp87}.
\begin{corollary}
Let $x = P_G(t)$ be given by (\ref{Kostant_gen_func}). Then
\begin{equation}
  \label{McKay_oper_5}
  tBx = (1 + t^2)x - v_0, \\
\end{equation}
  where $B$ is either the McKay operator $A$ or the Slodowy operators
  $\tilde{A}$, $\tilde{A}^\vee$.
\end{corollary}
\PerfProof
  From (\ref{Kostant_relation}) we obtain

\begin{equation}
\begin{split}
 & Bx = \sum\limits_{n=0}^{\infty}Bv_n{t^n} =
        \sum\limits_{n=0}^{\infty}(v_{n-1} + v_{n+1}){t^n} =
        \sum\limits_{n=0}^{\infty}v_{n-1}t^n +
        \sum\limits_{n=0}^{\infty}v_{n+1}t^n =  \\
 &  t\sum\limits_{n=1}^{\infty}v_{n-1}t^{n-1} +
    t^{-1}\sum\limits_{n=0}^{\infty}v_{n+1}t^{n+1} =
    tx + t^{-1}(\sum\limits_{n=0}^{\infty}v_n{t^n} - v_0) = \\
 &  tx + t^{-1}x - t^{-1}v_0.  \qed
\end{split}
\end{equation}

\subsection{The Poincar\'{e} series and W.~Ebeling's theorem}
\label{ebeling}
 \index{theorem! - Ebeling}

W.~Ebeling in \cite{Ebl02} makes use of the Kostant relation
(\ref{Kostant_relation}) and deduces a new remarkable fact about
the Poincar\'{e} series, a fact that shows that the Poincar\'{e}
series of a binary polyhedral group (see (\ref{K_K_GV})) is the
quotient of two polynomials: the characteristic polynomial of the
Coxeter transformation and the characteristic polynomial of the
corresponding affine Coxeter transformation, see
\cite[Th.2]{Ebl02}.

We show W.~Ebeling's theorem also for the multiply-laced case, see
Theorem \ref{theorem_ebeling}. The Poincar\'{e} series for the
multiply-laced case is defined by (\ref{poincare_alg_inv_2}).

\begin{table} 
  \centering
  \vspace{2mm}
  \caption{\hspace{3mm}The characteristic polynomials $\chi$, $\tilde{\chi}$ and
the Poincar\'{e} series}
  \renewcommand{\arraystretch}{1.5}  
  \begin{tabular} {|c|c|c|c|}
  \hline \hline
    Dynkin   & Coxeter
             & Affine Coxeter
             & Quotient \cr
    diagram  & transformation $\chi$
             & transformation $\tilde{\chi}$
             & $p(\lambda) = \displaystyle\frac{\chi}{\tilde{\chi}}$ \\
  \hline \hline
     & & & \cr
       ${D}_{4}$
     & $(\lambda + 1)(\lambda^{3} + 1)$
     & $(\lambda - 1)^2(\lambda + 1)^3$
     & $\displaystyle\frac{\lambda^3 + 1}{(\lambda^2 - 1)^2}$ \cr
     & & & \\
  \hline
     & & & \cr
       ${D}_{n+1}$
     & $(\lambda + 1)(\lambda^{n} + 1)$
     & $(\lambda^{n-1} -  1)(\lambda - 1)(\lambda + 1)^2$
     & $\displaystyle\frac{\lambda^{n} + 1}{(\lambda^{n-1} -  1)(\lambda^2 -
1)}$ \cr
     & & & \cr
  \hline
     & & & \cr
       ${E}_6$
     & $\displaystyle\frac{(\lambda^{6} + 1)}{(\lambda^2 + 1)}
       \frac{(\lambda^{3} - 1)}{(\lambda - 1)}$
     & $(\lambda^3 -  1)^2(\lambda + 1)$
     & $\displaystyle\frac{\lambda^6 + 1}{(\lambda^4 -  1)(\lambda^3 - 1)}$ \cr
     & & & \cr
  \hline
     & & & \cr
       ${E}_7$
     & $\displaystyle\frac{(\lambda + 1)(\lambda^{9} + 1)}{(\lambda^3 + 1)}$
     & $(\lambda^4 -  1)(\lambda^3 -  1)(\lambda + 1)$
     & $\displaystyle\frac{\lambda^9 + 1}{(\lambda^4 -  1)(\lambda^6 - 1)}$ \cr
     & & & \cr
  \hline
     & & & \cr
       ${E}_8$
     & $\displaystyle\frac{(\lambda^{15} + 1)(\lambda + 1)}{(\lambda^5 +
1)(\lambda^3  + 1)}$
     & $(\lambda^5 -  1)(\lambda^3 -  1)(\lambda + 1)$
     & $\displaystyle\frac{\lambda^{15} + 1}{(\lambda^{10} -  1)(\lambda^6 -
1)}$ \cr
     & & & \cr
  \hline
     & & & \cr
       ${B}_n$
     & $\lambda^{n} +  1$
     & $(\lambda^{n-1} -  1)(\lambda^2 -1)$
     & $\displaystyle\frac{\lambda^{n} + 1}{(\lambda^{n-1} -  1)(\lambda^2 -
1)}$ \cr
  \hline
     & & & \cr
       ${C}_n$
     & $\lambda^{n} +  1$
     & $(\lambda^n -  1)(\lambda - 1)$
     & $\displaystyle\frac{\lambda^{n} + 1}{(\lambda^n -  1)(\lambda - 1)}$ \cr
  \hline
     & & & \cr
       ${F}_4$
     & $\displaystyle\frac{\lambda^{6} + 1}{\lambda^2 + 1}$
     & $(\lambda^2 -  1)(\lambda^3 - 1)$
     & $\displaystyle\frac{\lambda^6 + 1}{(\lambda^4 -  1)(\lambda^3 - 1)}$ \\
  \hline
     & & & \cr
       ${G}_2$
     & $\displaystyle\frac{\lambda^{3} + 1}{\lambda + 1}$
     & $(\lambda-1)^2(\lambda  + 1)$
     & $\displaystyle\frac{\lambda^3 + 1}{(\lambda^2 - 1)^2}$ \\
  \hline\hline
     & & & \cr
       ${A}_n$
     & $\displaystyle\frac{\lambda^{n+1} - 1}{\lambda - 1}$
     & $(\lambda^{n - k + 1} - 1)(\lambda^k - 1)$
     & $\displaystyle\frac{\lambda^{n+1} - 1}
       {(\lambda - 1)(\lambda^{n - k + 1} - 1)(\lambda^k - 1)}$ \\
  \hline
     & & & \cr
       ${A}_{2n-1}$
     & $\displaystyle\frac{\lambda^{2n} - 1}{\lambda - 1}$
     & $(\lambda^n - 1)^2 \text{ for } k = n$
     & $\displaystyle\frac{\lambda^n + 1}
       {(\lambda^n - 1)(\lambda - 1)}$ \\
  \hline \hline
\end{tabular}
  \label{table_char_polyn_and_Poincare}
\end{table}

\index{theorem! - Ebeling}
\begin{theorem} [generalized W.Ebeling's theorem \cite{Ebl02}]
 \label{theorem_ebeling}
  Let $G$ be a binary polyhedral group and $[P_G(t)]_0$ the Poincar\'{e}
  series (\ref{poincare_alg_inv}) of the algebra of invariants
${\rm Sym}(\mathbb{C}^2)^G$. Then
\begin{equation}
       [P_G(t)]_0 = \frac{\det{M_0}(t)}{\det{M}(t)},
\end{equation}
where
\begin{equation}
   \det{M}(t) = \det|t^{2}I - {\bf C}_a|, \hspace{5mm}
   \det{M_0}(t) = \det|t^{2}I - {\bf C}|,
\end{equation}
{\bf C} is the Coxeter transformation and ${\bf C}_a$ is the
corresponding affine Coxeter transformation.
\end{theorem}

\PerfProof
By (\ref{McKay_oper_5}) we have
\begin{equation}
   [(1 + t^2)I - tB]x = v_0,
\end{equation}
where $x$ is the vector $P_G(t)$ and by Cramer's rule the first coordinate
$P_G(t)$ is
\begin{equation}
       [P_G(t)]_0 = \frac{\det{M_0}(t)}{\det{M}(t)},
\end{equation}
where
\begin{equation}
 \label{M_def_affine}
       \det{M}(t) = \det\left((1 + t^2)I - tB\right),
\end{equation}
and $M_0(t)$ is the matrix obtained by replacing the first column
of $M(t)$ by $v_0 = (1,0,...,0)^t$. The vector $v_0$ corresponds
to the trivial representation $\pi_0$, and by the McKay
correspondence, $v_0$ corresponds to the particular vertex which
extends the Dynkin diagram to the extended Dynkin diagram, see
Remark \ref{triv_repr} and (\ref{def_vn}). Therefore, if
$\det{M}(t)$ corresponds to the affine Coxeter transformation, and
\begin{equation}
 \label{M_affine_cox}
       \det{M}(t) = \det|t^{2}I - {\bf C}_a|,
\end{equation}
then $\det{M_0}(t)$ corresponds to the Coxeter transformation, and
\begin{equation}
 \label{M_cox}
       \det{M_0}(t) = \det|t^{2}I - {\bf C}|.
\end{equation}
So, it suffices to prove (\ref{M_affine_cox}), i.e.,
\begin{equation}
 \label{M_affine_cox_2}
        \det[(1 + t^2)I - tB] = \det|t^{2}I - {\bf C}_a|.
\end{equation}
If $B$ is the McKay operator $A$ given by (\ref{main_McKay}), then
\begin{equation}
 \label{McKay_operator}
       B = 2I - K =
    \left ( \begin{array}{cc}
            2I &  0  \\
            0  &  2I
            \end{array}
    \right ) -
    \left ( \begin{array}{cc}
            2I   & 2D \\
            2D^t &  2I
            \end{array}
    \right ) =
    \left ( \begin{array}{cc}
            0     & -2D \\
            -2D^t &  0
            \end{array}
    \right ),
\end{equation}
where $K$ is a symmetric Cartan matrix (\ref{symmetric_B}).
If $B$ is the Slodowy operator $\tilde{A}$ or $\tilde{A}^\vee$
given by (\ref{lab_slodowy_1}), (\ref{lab_slodowy_2}), then
\begin{equation}
 \label{Slodowy_operator}
       B = 2I - K =
    \left ( \begin{array}{cc}
            2I &  0  \\
            0  &  2I
            \end{array}
    \right ) -
    \left ( \begin{array}{cc}
            2I   & 2D \\
            2F &  2I
            \end{array}
    \right ) =
    \left ( \begin{array}{cc}
            0     & -2D \\
            -2F &  0
            \end{array}
    \right ),
\end{equation}
where $K$ is the symmetrizable but not symmetric Cartan matrix (\ref{matrix_K}).
Thus, in the generic case
\begin{equation}
 \label{M_t}
       M(t) = (1+t^2)I - tB =
    \left ( \begin{array}{cc}
            1+t^2 &  2tD     \\
            2tF  &  1 + t^2
            \end{array}
     \right ).
\end{equation}
Assuming $t \neq 0$ we deduce from (\ref{M_t}) that
\begin{equation}
 \label{M_t_2}
 \begin{array}{cc}
       M(t)
    \left ( \begin{array}{c}
            x    \\
            y
            \end{array}
    \right ) = 0 & \Longleftrightarrow
    \left \{
     \begin{array}{c}
            (1 + t^2)x = -2tDy,  \\
            2tFx = -(1 + t^2)y.
            \end{array}
     \right .  \vspace{5mm} \\
    & \Longleftrightarrow
    \left \{
     \begin{array}{c}
            \displaystyle\frac{(1 + t^2)^2}{4t^2}x = FDy,   \vspace{3mm} \\
            \displaystyle\frac{(1 + t^2)^2}{4t^2}y = DFy.
            \end{array}
     \right .
  \end{array}
\end{equation}
Here we use the Jordan form theory of the Coxeter transformations
constructed in \S\ref{jordan}. According to (\ref{DDt_DtD}),
Proposition \ref{gold_pair} and Proposition \ref{eigenvectors_basis}
we see that $t^2$ is an eigenvalue of the affine Coxeter transformation
${\bf C}_a$, i.e., (\ref{M_affine_cox_2})
together with (\ref{M_affine_cox}) are proved. \qedsymbol

\index{theorem! - Ebeling}
For the results of calculations using W.~Ebeling's theorem, see
Table \ref{table_char_polyn_and_Poincare}.

\begin{remark} {\rm
1) The characteristic polynomials $\chi$ for the Coxeter
transformation and $\tilde{\chi}$ for the affine Coxeter
transformation in Table \ref{table_char_polyn_and_Poincare} are
taken from Tables \ref{table_char_polynom_Dynkin} and
\ref{table_char_ext_polynom_Dynkin}. Pay attention to the fact
that the affine Dynkin diagram for ${B}_n$ is $\widetilde{CD}_n$,
(\cite[Tab.2]{Bo}),  and the affine Dynkin diagram for ${C}_n$ is
$\tilde{C}_n$, (\cite[Tab.3]{Bo}), see Fig.~\ref{euclidean_diag}.

 \index{index of the conjugacy class}
 2) The characteristic
polynomial $\chi$ for the affine Coxeter transformation of ${A}_n$
depends on the {\it index of the conjugacy class} $k$ of the
Coxeter transformation, see (\ref{char_polyn_An}). In the case of
${A}_n$ (for every $k = 1, 2,..., n$) the quotient $p(\lambda) =
\displaystyle\frac{\chi}{\tilde{\chi}}$ contains three factors in
the denominator, and its form is different from (\ref{K_K_GV}),
see Table \ref{table_char_polyn_and_Poincare}.

For the case ${A}_{2n-1}$ and $k = n$, we have
\begin{equation}
 \begin{split}
     p(\lambda) = & \frac{\lambda^{2n} - 1}
                 {(\lambda - 1)(\lambda^{2n-k} - 1)(\lambda^k - 1)} = \vspace{2mm} \\
              & \frac{\lambda^{2n} - 1}
                 {(\lambda - 1)(\lambda^n - 1)(\lambda^n - 1)} =
               \frac{\lambda^n + 1}
                 {(\lambda^n - 1)(\lambda - 1)}
 \end{split}
\end{equation}
and $p(\lambda)$ again is of the form (\ref{K_K_GV}),
see Table \ref{table_char_polyn_and_Poincare}.

\index{folding operation}
3) The quotients $p(\lambda)$ coincide for the following pairs:
\begin{equation}
   \begin{array}{cc}
       {D}_4  \text{ and } {G}_2, \hspace{7mm}
      & {E}_6  \text{ and } {F}_4, \\
       {D}_{n+1} \text{ and } {B}_n  (n \geq 4), \hspace{7mm}
      &  {A}_{2n-1} \text{ and } {C}_n.
   \end{array}
\end{equation}
Note that the second elements of the pairs
are obtained by {\it folding} operation from the first ones,
see Remark \ref{remark_folding}.
}
\end{remark}

%% file: 6regular.tex

\chapter{\sc\bf Regular representations of quivers}
  \label{chap_regular}

\section{The regular representations and the Dlab-Ringel criterion}
\label{criteria_DR} The regular representations are the most
complicated in the category of all representation of the given
quiver. For every Dynkin diagram the category of regular
representations is empty, there are only finite number non-regular
representations (P.Gapriel's theorem \cite{Gab72}). For this
reason, the Dynkin diagrams are called {\it finite type quivers}
in the representation theory of quivers. The regular
representations have been completely described only for the
extended Dynkin diagrams,
 \index{tame quivers} which for this
reason were dubbed {\it tame quivers} in the representation theory
of quivers, (\cite{Naz73}, \cite{DR76}), see \S\ref{tame_quivers}.

\begin{definition}
 \label{def_regular}{\rm
 \index{$\varDelta{'}$-regular vector}
 \index{regular vector in the orientation $\varDelta{'}$}

1) The vector $z \in \mathcal{E}_\varGamma$ is said to be
$\varDelta{'}$-{\it regular} or {\it regular in the orientation}
$\varDelta{'}$ if
\begin{equation}
        {\bf C}^k_{\varDelta{'}}z > 0 \text{ for all } k \in \mathbb{Z}.
\end{equation}

\index{regular representation in the orientation $\varDelta{'}$}
2) The representation $V$ is said to be a {\it regular
representation in the orientation} $\varDelta{'}$ of the quiver
$\varGamma$ if its dimension $\dim{V}$ is a $\varDelta{'}$-regular
vector, i.e.,
\begin{equation}
        {\bf C}^k_{\varDelta{'}}(\dim{V}) > 0 \text{ for all } k \in \mathbb{Z}.
\end{equation}
}
\end{definition}

Let the Coxeter transformation  ${\bf C}_{\varDelta{'}}$
correspond to any orientation $\varDelta{'}$ and let the Coxeter
transformation ${\bf C}_{\varDelta}$ correspond to a bicolored
orientation $\varDelta$. Then ${\bf C}_{\varDelta{'}}$ and ${\bf
C}_{\varDelta}$ are conjugate if the quiver is a tree. Let $T$ be
an element in the Weyl group which relates ${\bf
C}_{\varDelta{'}}$ and ${\bf C}_{\varDelta}$:
\begin{equation}
          {\bf C}_{\varDelta{'}} = T^{-1}{\bf C}_{\varDelta}T.
\end{equation}

\section{Necessary conditions of regularity of the extended Dynkin diagrams}
\index{conditions of regularity! - necessary conditions,
$\mathcal{B} \geq 0$} Let $(\alpha_1, \tilde{\alpha}_1,
\alpha^{\varphi_2}_1, \alpha^{\varphi_2}_2, ...)$ be coordinates
of the vector $Tz$ in the Jordan basis of eigenvectors and adjoint
vectors (\ref{case_non_01}) -- (\ref{case_0}):
\begin{equation}
  \label{Tz}
    Tz = \alpha_1{z}^1 + \tilde{\alpha}_1\tilde{z}^1
+\alpha^{\varphi_2}_1{z}^{\varphi_2}_1 +
\alpha^{\varphi_2}_2{z}^{\varphi_2}_2 + \dots
\end{equation}

According to (\ref{PK_not_01}), (\ref{PK_1}), (\ref{PK_0}), we have
\begin{equation}
    {\bf C}^k_{\varDelta}T{z} = \alpha_1{z}^1 +
\tilde{\alpha}_1(k{z}^1 + \tilde{z}^1) +
\alpha^{\varphi_2}_1(\lambda^{\varphi_2}_1)^k{z}^{\varphi_2}_1 +
\alpha^{\varphi_2}_2(\lambda^{\varphi_2}_2)^k{z}^{\varphi_2}_2 +
 \dots
\end{equation}
and
\begin{equation}
 \label{Tz_delta}
  \begin{split}
    {\bf C}^k_{\varDelta{'}}{z} = &
       T^{-1}{\bf C}^k_{\varDelta{'}}T{z} = \\
                   & \alpha_1{z}^1 +
  \tilde{\alpha}_1(k{z}^1 + {T}^{-1}\tilde{z}^1) + \\
&
\alpha^{\varphi_2}_1(\lambda^{\varphi_2}_1)^k{T}^{-1}{z}^{\varphi_2}_1
+
\alpha^{\varphi_2}_2(\lambda^{\varphi_2}_2)^k{T}^{-1}{z}^{\varphi_2}_2
+ \dots
  \end{split}
\end{equation}

Since $z^1 > 0$,   $k \in \mathbb{Z}$, and
$|\lambda^{\varphi_i}_{1,2}| = 1$, we have the following necessary
condition of $\varDelta{'}$-regularity:
\begin{equation}
   \label{2_coord}
             \tilde{\alpha}_1 = 0,
\end{equation}
where $\tilde{\alpha}_1$ is the coordinate of the vector $Tz$
corresponding to the adjoint basis vector $\tilde{z}^1$.

\underline{a) Simply-laced case}.
The adjoint vector
$\tilde{z}^1$ is orthogonal to vectors
${z}^{\varphi_j}_i$ for $i = 1,2$; the numbers
${\varphi_j}$ are eigenvalues of $DD^t$ since the corresponding
components $\mathbb{X}$ and $\mathbb{Y}$ are orthogonal,
see Proposition $\ref{eigenvectors_creation}$.  Further,
$\tilde{z}^1$ is also orthogonal to ${z}^1$: indeed,
\begin{equation*}
\begin{split}
   & <z^1, \tilde{z}^1> =
    <\mathbb{X}^1, \mathbb{X}^1> - <D^t\mathbb{X}^1, D^t\mathbb{X}^1> = \\
   &    <\mathbb{X}^1, \mathbb{X}^1> - <\mathbb{X}^1, DD^t\mathbb{X}^1> =
    <\mathbb{X}^1, \mathbb{X}^1> - <\mathbb{X}^1, \mathbb{X}^1> = 0.
\end{split}
\end{equation*}

Thus, from (\ref{Tz}) we have
\begin{equation}
  \label{alpha_1_simply}
  <Tz, \tilde{z}^1> =
     \tilde{\alpha}_1<\tilde{z}^1, \tilde{z}^1>
\end{equation}
and (\ref{2_coord}) is equivalent to the following relation:
\begin{equation}
   \label{regular_1}
             <Tz, \tilde{z}^1> = 0.
\end{equation}

\underline{b)  Multiply-laced case}.
   The normal basis of $DF$ (resp. $FD$) is
not orthogonal, however, the dual graph will help us.

\index{conjugate vector}
\begin{definition} {\rm
  \label{conjug_vector}
Let $u$ be an eigenvector of the Coxeter transformation given
by (\ref{case_non_01}), (\ref{case_1}). The vector $\tilde{u}$ is
said to be {\it conjugate} to $u$, if it is obtained
from $u$ by changing the sign of $\mathbb{Y}$-component and
replacing the eigenvalue $\lambda$
by $\displaystyle\frac{1}{\lambda}$.

(It is easily to see that the vectors ${z}^1$ and $4\tilde{z}^1$ are conjugate.)
 }
\end{definition}

\begin{proposition}
  \label{ortogonal_eigenv}
1) Let ${\varphi}$ be an eigenvalue for $DF$ and
${\varphi}^{\vee} \neq \varphi$ an eigenvalue for $D^{\vee}F^{\vee}$.
The eigenvectors of $DF$ and $D^{\vee}F^{\vee}$ corresponding
to these eigenvalues are orthogonal.

2) Let $z_{\varphi}$ be an eigenvector with eigenvalue ${\varphi}
\neq 0$ for the extended Dynkin diagram ${\varGamma}$ and
$z^{\vee}_{\varphi}$ an eigenvector with eigenvalue ${\varphi}$
for the dual diagram ${\varGamma}^{\vee}$. Let
$\tilde{z}^{\vee}_{\varphi}$ be conjugate to
${z}^{\vee}_{\varphi}$. Then ${z}_{\varphi}$ and
$\tilde{z}^{\vee}_{\varphi}$ are orthogonal.
\end{proposition}

\PerfProof 1) Let $DFx = {\varphi}x$ and
${\it D^{\vee}F^{\vee}}x = {\varphi^{\vee}}x^{\vee}$.
Since $\varphi \neq \varphi^{\vee}$,
one of these eigenvalues is $\neq 0$. Let, for example, $\varphi \neq 0$. Then
\begin{equation}
 \label{x_x_vee}
  <x, x^{\vee}> = \frac{1}{\varphi}<DFx, x^{\vee}> =
  \frac{1}{\varphi}<x, (DF)^tx^{\vee}>.
\end{equation}
From (\ref{x_x_vee}) and (\ref{dual_rel}) we have
\begin{equation*}
  <x, x^{\vee}> =
  \frac{1}{\varphi}<x, D^{\vee}F^{\vee}x^{\vee}> =
  \frac{\varphi^{\vee}}{\varphi}<x, x^{\vee}>.
\end{equation*}

Since $\varphi \neq \varphi^{\vee}$, we have $<x, x^{\vee}> = 0$.
\hspace{3mm} \qedsymbol

2) Let us express vectors $z_{\varphi}$ and $\tilde{z}^{\vee}_{\varphi}$ as
follows:
\begin{equation}
  \label{conjugate_z}
    z_{\varphi} =
    \left (
    \begin{array}{c}
          x  \vspace{4mm} \\
          -\displaystyle\frac{2}{\lambda + 1}Fx
    \end{array}
    \right ), \hspace{7mm}
    \tilde{z}^{\vee}_{\varphi} = \left (
    \begin{array}{c}
          x^{\vee}  \vspace{4mm} \\
          \displaystyle\frac{2\lambda}{\lambda + 1}F^{\vee}x^{\vee}
    \end{array}
    \right ).
\end{equation}
According to (\ref{dual_t}) and (\ref{dual_rel}) we have
\begin{equation*}
  \begin{split}
  & <z_{\varphi}, \tilde{z}^{\vee}_{\varphi}> =
  <x, x^{\vee}> -  \frac{4\lambda}{(\lambda + 1)^2}<Fx, F^{\vee}x^{\vee}> = \\
  & <x, x^{\vee}> - \frac{1}{\varphi}<DFx, x^{\vee}> =
    <x, x^{\vee}> - <x, x^{\vee}> = 0.
    \hspace{3mm} \qed
  \end{split}
\end{equation*}
Thus, by Proposition \ref{ortogonal_eigenv} heading 1), the vector
$\tilde{z}^{1\vee}$ is orthogonal to the
${z}^{\varphi}_i$ for $i=1,2$ and by Proposition \ref{ortogonal_eigenv} heading
2), the vector $\tilde{z}^{1\vee}$ is orthogonal to $z^1$.
Therefore from (\ref{Tz}) we deduce that
\begin{equation}
 \label{alpha_1_m}
  <Tz, \tilde{z}^{1\vee}> =
     \tilde{\alpha}_1<\tilde{z}^{1\vee}, \tilde{z}^{1\vee}>
\end{equation}
and (\ref{2_coord}) is equivalent to the following relation:
\begin{equation}
   \label{regular_2}
             <Tz, \tilde{z}^{1\vee}> = 0.
\end{equation}
Relation (\ref{alpha_1_simply}) for the simply-laced case
and relation (\ref{alpha_1_m}) for the multiply-laced case
motivate the following

\begin{definition} {\rm
  \label{rho_defect}
   The linear form $\rho_{\varDelta{'}}(z)$
   given as follows
\begin{equation}
 \label{rho_simply_and_multy}
\rho_{\varDelta{'}}(z) = \begin{cases}<Tz, \tilde{z}^1>&\text{for
a simply-laced
case}, \\
<Tz, \tilde{z}^{1\vee}>& \text{for a multiply-laced case}
 \end{cases}
\end{equation}
is said to be the {\it $\varDelta{'}$-defect} of the vector $z$.
 }
\end{definition}

From (\ref{regular_1}) and (\ref{regular_2})
we derive the following theorem:

 \index{theorem! - on regular representations of quivers}

\begin{theorem} [\cite{SuSt75},\cite{SuSt78}]
 \label{necess_codnd_reg}
If $z$ is a regular vector for the extended Dynkin diagram
$\varGamma$
  in the orientation $\varDelta{'}$, then
\begin{equation}
   \label{regular_3}
      \rho_{\varDelta{'}}(z) = 0.
\end{equation}
\end{theorem}

Below, in Proposition \ref{arbitr_orient}, we will show that the
condition (\ref{regular_3}) is also sufficient, if $z$ is a root
in the root system related to the given extended Dynkin diagram.

\section{The Dlab-Ringel definition of the defect}

The vector $z^1$ from the kernel of the Tits form is the fixed point for
the Weyl group, so $T^{-1}z^1 = z^1$, and from (\ref{Tz}) we deduce
\begin{equation}
  \label{Tz_1}
    z = \alpha_1{z}^1 +\tilde{\alpha}_1{T}^{-1}\tilde{z}^1 +
\alpha^{\varphi_2}_1{T}^{-1}{z}^{\varphi_2}_1 +
\alpha^{\varphi_2}_2{T}^{-1}{z}^{\varphi_2}_2 + \dots
\end{equation}
and from (\ref{Tz_delta}) for $k = h$ being the Coxeter number,
i.e., $(\lambda^{\varphi_2}_{1,2})^k = 1$,  we obtain
\index{Dlab-Ringel formula} \index{formula! - Dlab-Ringel}
\begin{equation}
  \label{Tz_2}
    {\bf C}^k_{\varDelta{'}}z =
                    \alpha_1{z}^1 +
               \tilde{\alpha}_1(k{z}^1 + {T}^{-1}\tilde{z}^1) +
        \alpha^{\varphi_2}_1{T}^{-1}{z}^{\varphi_2}_1 +
       \alpha^{\varphi_2}_2{T}^{-1}{z}^{\varphi_2}_2 + \dots
\end{equation}

From (\ref{Tz_1}) and (\ref{Tz_2}) we get the following formula
due to V.~Dlab and C.~M.~Ringel, see \cite{DR76}: \index{Coxeter
number}
\begin{equation}
  \label{Dlab_Ringel}
    {\bf C}^h_{\varDelta{'}}z =  z + h\tilde{\alpha}_1{z}^1,
\end{equation}
where $h$ is the Coxeter number and $\tilde{\alpha}_1 =
\tilde{\alpha}_1(z)$ is the linear form proportional to
$\rho_{\varDelta{'}}(z)$. By (\ref{alpha_1_simply}) for the
simply-laced case and by (\ref{alpha_1_m}) for the multiply-laced
case the form $\tilde{\alpha}_1(z)$ can be calculated as follows:
\begin{equation}
  \tilde{\alpha}_1(z) =\begin{cases}
\frac{<Tz, \tilde{z}^1>}{<\tilde{z}^1, \tilde{z}^1>}&\text{for a simply-laced
case, }  \vspace{5mm} \\
\frac{<Tz, \tilde{z}^{1\vee}>}
       {<\tilde{z}^{1\vee}, \tilde{z}^{1\vee}>}&\text{for a multiply-laced
case. }
  \end{cases}
\end{equation}

Dlab and Ringel introduced in \cite{DR76} the {\it defect}
$\delta_{\varDelta{'}}$  as a vector from
$\mathcal{E}^*_{\varGamma}$ obtained as a solution of the equation
\begin{equation}
\index{Dlab-Ringel definition of defect}
  \label{DR_defect}
     {\bf C}^*_{\varDelta{'}}\delta_{\varDelta{'}} = \delta_{\varDelta{'}}.
\end{equation}

\begin{proposition} [\cite{St85}]
 \label{defect_coinciding}
  The Dlab-Ringel defect $\delta_{\varDelta{'}}$ given by (\ref{DR_defect})
  coincides (up to coefficient) with the ${\varDelta{'}}$-defect $\rho_{\varDelta{'}}$
given by
Definition $\ref{rho_defect}$.
\end{proposition}
\PerfProof1) Let us show that  $\delta_{\varDelta{'}}$ is obtained
from $\delta_{\varDelta}$ in the same way as $\rho_{\varDelta{'}}$
is obtained from $\rho_{\varDelta}$, so it suffices to prove the
proposition only for bicolored orientations $\varDelta$. Indeed,
$$
  {\bf C}_{\varDelta{'}} = {T}^{-1}{\bf C}_{\varDelta}{T} \hspace{7mm}
\text { implies that }  \hspace{7mm}
  {\bf C}^*_{\varDelta{'}} = {T}^*{\bf C}^*_{\varDelta}{T}^{-1*}.
$$
Since ${\bf C}^*_{\varDelta{'}}\delta_{\varDelta{'}} =
\delta_{\varDelta{'}}$, we obtain
$$
{\bf C}^*_{\varDelta}{T}^{-1*}\delta_{\varDelta{'}} =
        {T}^{-1*}\delta_{\varDelta{'}},
$$
and
$$
        {T}^{-1*}\delta_{\varDelta{'}} = \delta_{\varDelta}, \hspace{7mm}
\text{ i.e., } \hspace{7mm}
      \delta_{\varDelta{'}} = {T}^*\delta_{\varDelta}.
$$
Further,
$$
  <\delta_{\varDelta},{T}z> =
  <{T}^*\delta_{\varDelta},z> =  <\delta_{\varDelta{'}},z>.
$$
Thus,  $\delta_{\varDelta{'}}$ is obtained from
$\delta_{\varDelta}$ in the same way as $\rho_{\varDelta{'}}$ is
obtained from $\rho_{\varDelta}$.

2) Now, let us prove the proposition for the bicolored orientation
$\varDelta$, i.e., let us prove that $\rho_\varDelta$ is
proportional to $\delta_\varDelta$. According to (\ref{C_decomp}),
(\ref{symmetric_B}) and (\ref{matrix_K}) the relation
$$
   {\bf C}^*_{\varDelta}z = z \hspace{7mm} \text { is equivalent to }
   \hspace{7mm}
   w^*_1{z} = w^*_2{z},
$$
which, in turn, is equivalent to the following:
\begin{equation*}
 \left (
  \begin{array}{cc}
    -I     &0 \\
    -2D^t  &I \\
  \end{array}
 \right )
 \left (
  \begin{array}{c}
    x \\
    y \\
  \end{array}
 \right )
 =
 \left (
  \begin{array}{cc}
    I  & -2F \\
    0  &  -I
  \end{array}
 \right )
 \left (
  \begin{array}{c}
    x \\
    y \\
  \end{array}
 \right ).
\end{equation*}
By (\ref{dual_t}) we have
\begin{equation*}
 \left \{
  \begin{array}{c}
    x = F^{t}y = D^{\vee}y, \\
    y = D^{t}x = F^{\vee}x.
  \end{array}
 \right .
\end{equation*}
 Thus,
\begin{equation*}
 z = \left (
  \begin{array}{c}
    x  \\
    F^{\vee}x \\
  \end{array}
 \right ), \text{ where } x = F^{t}D^{t}x = D^{\vee}F^{\vee}x,
\end{equation*}
i.e., $x$ corresponds to $\lambda = 1$, an eigenvalue of {\bf C},
and to $\varphi = 1$, an eigenvalue of $(DF)^{\vee}$. By
(\ref{conjugate_z}) we have $z = \tilde{z}^{1\vee}$, i.e.,
$\delta_\varDelta$ is proportional to $\rho_\varDelta$, see
(\ref{rho_simply_and_multy}). \hspace{3mm} \qedsymbol

\section[Necessary conditions of regularity]{Necessary conditions of regularity
for diagrams with indefinite Tits form} \label{necessary_cond}
\index{conditions of regularity! - necessary conditions,
$\mathcal{B}$ - indefinite} Now, consider the case where
$\mathcal{B}$ is indefinite, i.e., $\varGamma$ is any tree, which
is neither Dynkin diagram nor extended Dynkin diagram.

Let $(\alpha^m_1, {\alpha}^m_2, \alpha^{\varphi_2}_1,
\alpha^{\varphi_2}_2, ...)$ be coordinates of the vector ${T}z$ in
the Jordan basis of eigenvectors and adjoint vectors
(\ref{case_non_01}) -- (\ref{case_0}), where $\alpha^m_1$ and
$\alpha^m_2$ are coordinates corresponding to eigenvectors
 $z^m_1$ and $z^m_2$. The vectors $z^m_1$ and $z^m_2$ correspond to the maximal
eigenvalue
${\varphi}^m = {\varphi}^{max}$ of $DD^t$ and $D^tD$, respectively. Let us
decompose the vector
$Tz$ as follows:
\begin{equation}
  \label{Tz_indef}
    Tz = \alpha^m_1{z}^m_1 +
                     \alpha^m_2{z}^m_2 +
\alpha^{\varphi_2}_1{z}^{\varphi_2}_1 +
\alpha^{\varphi_2}_2{z}^{\varphi_2}_2 + \dots
\end{equation}

According to (\ref{PK_not_01}), (\ref{PK_1}), (\ref{PK_0}) we have
\begin{equation}
 \begin{split}
    & {\bf C}^k_{\varDelta}Tz =  \\
            & (\lambda^{\varphi_m}_1)^k\alpha^m_1{z}^m_1 +
     (\lambda^{\varphi_m}_2)^k\alpha^m_2{z}^m_2 +
(\lambda^{\varphi_2}_1)^k\alpha^{\varphi_2}_1{z}^{\varphi_2}_1 +
(\lambda^{\varphi_2}_1)^k\alpha^{\varphi_2}_2{z}^{\varphi_2}_2 +
\dots
  \end{split}
\end{equation}
and
\begin{equation}
 \label{TCTz_indef}
  \begin{split}
   {\bf C}^k_{\varDelta{'}}z =
       {T}^{-1}{\bf C}^k_{\varDelta{'}}Tz = &
    (\lambda^{\varphi_m}_1)^k\alpha^m_1{T}^{-1}z^m_1 +
    (\lambda^{\varphi_m}_2)^k\alpha^m_2{T}^{-1}z^m_2 + \\
    & (\lambda^{\varphi_2}_1)^k\alpha^{\varphi_2}_1{T}^{-1}z^{\varphi_2}_1 +
    (\lambda^{\varphi_2}_1)^k\alpha^{\varphi_2}_2{T}^{-1}z^{\varphi_2}_2 +
\dots
  \end{split}
\end{equation}

It will be shown in Theorem \ref{th_transf} that the transforming
element $T$ can be modified so that its decomposition does not
contain any given reflection $\sigma_\alpha$. Since the
coordinates of the eigenvectors ${z}^m_{1,2}$ are all positive,
see Corollary \ref{corollary_dominant} and (\ref{PK_not_01}), we
see that each vector ${T}^{-1}{z}^m_{1,2}$ has at least one
positive coordinate. Besides,
\begin{equation}
  \label{lambdas_rel}
   |\lambda^{\varphi_m}_1| > |\lambda^{\varphi_j}_{1,2}| >
|\lambda^{\varphi_m}_2|
\end{equation}
because
\begin{equation*}
  \begin{split}
   & \lambda^{\varphi_m}_1 = \frac{1}{\lambda^{\varphi_m}_2}  \quad \text{ and}\\
   & \lambda^{\varphi_m}_1 = 2\varphi^m - 1 + 2\sqrt{\varphi^m(\varphi^m-1)} >
     2\varphi^i - 1 \pm 2\sqrt{\varphi^i(\varphi^i-1)}.
  \end{split}
\end{equation*}
Thus, since ${T}^{-1}z^m_{1,2}$ have at least one positive
coordinate, we deduce from (\ref{TCTz_indef}) and
(\ref{lambdas_rel}) that
\begin{equation}
  \label{regular_alphas}
   \alpha^{\varphi_m}_1 \ge 0, \hspace{7mm}  \alpha^{\varphi_m}_2 \ge 0.
\end{equation}

As above for the case $\mathcal{B} \geq 0$, let us calculate
$\alpha^{\varphi_m}_{1,2}$. The vector $\tilde{z}^m_1$ conjugate
to the vector $z^m_1$ is orthogonal to the vectors $z^2_{1,2},
z^3_{1,2}, \dots$. Let us show that $\tilde{z}^m_1$ is also
orthogonal to $z^m_1$. Indeed, the vectors $z^m_1$ and
$\tilde{z}^m_1$ can be expressed as follows:
\begin{equation}
    z^m_1 =
    \left (
    \begin{array}{c}
          x^m  \vspace{4mm} \\
          -\displaystyle\frac{2}{\lambda^m_1 + 1}D^tx^m
    \end{array}
    \right ), \hspace{7mm}
    \tilde{z}^m_1 = \left (
    \begin{array}{c}
          x^m  \vspace{4mm} \\
          \displaystyle\frac{2\lambda^m_1}{\lambda^m_1 + 1}D^tx^m
    \end{array}
    \right ),
\end{equation}
where for brevity we designate $\lambda^{\varphi_{max}}_1$ by $\lambda^m_1$,
\;
$z^{\varphi_{max}}_1$ by $z^m_1$, and $x^{\varphi_{max}}$ by $x^m$.
Then,
\begin{equation}
  \begin{split}
    <z^m_1, \tilde{z}^m_1> = & \\
    & <x^m, x^m> - \frac{4\lambda^m_1}{(\lambda^m_1 + 1)^2}<D^tx^m, D^tx^m> = \\
    & <x^m, x^m> - \frac{1}{\varphi_m}<x^m, DD^tx^m> = \\
    &  <x^m, x^m> - <x^m, x^m> = 0.
  \end{split}
\end{equation}
\index{conjugate vector}
The conjugate vector $\tilde{z}^m_1$ is not orthogonal only to
${z}^m_2$ in decomposition (\ref{Tz_indef}),
and
similarly $\tilde{z}^m_2$ is not orthogonal only to ${z}^m_1$.
From (\ref{Tz_indef}) we get
\begin{equation}
  \label{calc_Tz}
   <{T}z, \tilde{z}^m_1> = \alpha^m_2<z^m_2, \tilde{z}^m_1>,
   \hspace{7mm}
   <{T}z, \tilde{z}^m_2> = \alpha^m_1<z^m_1, \tilde{z}^m_2>.
\end{equation}
Let us find $<z^m_2, \tilde{z}^m_1>$ and $<z^m_1, \tilde{z}^m_2>$.
We have:
\begin{equation}
 \label{1_lambda_2}
  \begin{split}
  & <z^m_2, \tilde{z}^m_1> =
   <x^m, x^m> - \frac{4}{(\lambda^m_2 + 1)^2}<x^m, DD^tx^m> = \vspace{3mm} \\
  & <x^m, x^m>(1 - \frac{4}{(\lambda^m_2 + 1)^2}
                   \frac{(\lambda^m_2 + 1)^2}{4(\lambda^m_2)^2}) =
     <x^m, x^m>(1 - (\lambda^m_1)^2).
  \end{split}
\end{equation}
Similarly,
\begin{equation}
 \label{1_lambda_2_1}
  \begin{split}
  & <z^m_1, \tilde{z}^m_2> = <x^m, x^m>(1 - (\lambda^m_2)^2).
  \end{split}
\end{equation}
Thus, from (\ref{calc_Tz}), (\ref{1_lambda_2}) and (\ref{1_lambda_2_1}) we get
\begin{equation}
  \label{calc_Tz_2}
  \begin{split}
 &  <Tz, \tilde{z}^m_1> = \alpha^m_2<x^m, x^m>(1 - (\lambda^m_1)^2), \\
 &  <Tz, \tilde{z}^m_2> = \alpha^m_1<x^m, x^m>(1 - (\lambda^m_2)^2).
  \end{split}
\end{equation}

 \index{theorem! - on regular representations of quivers}

\begin{theorem} [\cite{SuSt75}, \cite{SuSt78}]
  \label{regul_indef_Tits}
  If $z$ is a regular vector for the graph $\varGamma$ with indefinite Tits
  form 
  $\mathcal{B}$ in the orientation $\varDelta{'}$, then
\begin{equation}
   \label{regular_5}
      <Tz, \tilde{z}^m_1> \hspace{3mm} \leq \hspace{3mm} 0,
      \hspace{7mm}
      <Tz, \tilde{z}^m_2> \hspace{3mm} \geq \hspace{3mm} 0.
\end{equation}
\end{theorem}
\PerfProof
 Since $\lambda^m_1 > 1$ and $\lambda^m_2 < 1$,
 the theorem follows from (\ref{calc_Tz_2}) and (\ref{regular_alphas}).
\hspace{3mm} \qedsymbol

We denote the linear form $<Tz, \tilde{z}^m_1>$ (resp. $<Tz,
\tilde{z}^m_2>$) by $\rho^1_{\varDelta{'}}$ (resp.
$\rho^2_{\varDelta{'}}$). Then, conditions (\ref{regular_5}) have
the following form:
\begin{equation}
   \label{regular_6}
      \rho^1_{\varDelta{'}}(z) \hspace{3mm} \leq \hspace{3mm} 0,
      \hspace{7mm}
      \rho^2_{\varDelta{'}}(z) \hspace{3mm} \geq \hspace{3mm} 0.
\end{equation}

A similar results were obtained by Y.~Zhang in
\cite[Prop.1.5]{Zh89}, and by J.~A.~de la Pe\~{n}a, M.~Takane in
\cite[Th.2.3]{PT90}.

For an application of (\ref{regular_5}) to the star graph, see
\S\ref{star_example}.

\section{Transforming elements and sufficient conditions of regularity}
In this section we consider only extended Dynkin diagrams. We will
show in Theorem \ref{necess_eq_suff} that {\it  the necessary
condition of regularity of the vector $z$ (\ref{regular_3})
coincides with the sufficient condition (Proposition
\ref{prop_sufficient}) only if vector $z$ is a positive root in
the corresponding root system.} For an arbitrary vector $z$ this
is not true. Since dimension $\dim{V}$ of the indecomposable
representation $V$ is a positive root (the Gabriel theorem for the
Dynkin diagrams \cite{Gab72}, \cite{BGP73}, the Dlab-Ringel
theorem for the extended Dynkin diagram \cite{DR76}, \cite{Naz73},
the Kac theorem for any diagrams \cite{Kac80}, \cite{Kac82}), then
{\it for the indecomposable representation $V$ the necessary
condition of regularity of the vector $\dim{V}$ coincides with the
sufficient condition}. For an arbitrary decomposable
representation this is not true.

We start from the bicolored orientation $\varDelta$.

\subsection{Sufficient conditions of regularity for the bicolored orientation}

 \index{conditions of regularity! - sufficient conditions,
bicolored orientation}
 \index{root system}

\begin{proposition} [\cite{St82}]
  \label{prop_sufficient}
  Let $\varGamma$ be an extended Dynkin diagram, i.e., $\mathcal{B} \geq 0$.
  Let $z$ be a root in the root system associated with $\varGamma$.
  If the $\varDelta$-defect of the vector $z$ is zero:
\begin{equation}
  \label{suff_cond}
        \rho_\varDelta(z) = 0,
\end{equation}
 then $z$ is regular in the bicolored orientation $\varDelta$.
\end{proposition}
\PerfProof Let
\begin{equation}
    H = \{z \mid z > 0, \; z \text{ is a root},
\;
    \rho_\varDelta(z) = 0 \}.
\end{equation}
It suffices to prove that
\begin{equation}
  \label{suff_inclusion}
    w_1{H} \subset H, \hspace{7mm} w_2{H} \subset H.
\end{equation}
Indeed, if (\ref{suff_inclusion}) holds, then
$$
   {\bf C}_{\varDelta}^kH \subset H \hspace{3mm} \text{for all} \hspace{3mm}
    k \in \mathbb{N},
$$
hence ${\bf C}_{\varDelta}^k{z} > 0$ if $z$ is a positive root
satisfying the condition (\ref{suff_cond}).

So, let us prove, for example, that $w_1{H} \subset H$.
Note that $z$ and $w_1{z}$ are roots simultaneously. Thus,
either $w_1{z} > 0$ or $w_1{z} < 0$.
Suppose, $w_1{z} < 0$.
Together with $z > 0$,  by (\ref{matrix_K}) we have $y = 0$.
Hence
$$
  \rho_\varDelta(z) = <z, \tilde{z}^{1\vee}> =
  <x, \tilde{x}^{1\vee}>.
$$
The coordinates of $\tilde{x}^{1\vee}$ are positive, the
coordinates of $x$ are non-negative. If  $\rho_\varDelta(z) = 0$,
then $x = 0$; this contradicts to the condition $z > 0$.
Therefore, $w_1{z} > 0$. It remains to show that
$$
   \rho_\varDelta(w_1{z}) = 0.
$$
Again, by (\ref{matrix_K}), (\ref{dual_t}) and (\ref{case_1}) we have
\begin{equation}
 \label{wz_z}
  \begin{split}
  & <w_1{z}, \tilde{z}^{1\vee}> =
     <
      \left(
          \begin{array}{c}
            -x - 2Dy \\
            y        \\
          \end{array}
      \right),
      \left(
          \begin{array}{c}
            x^{1\vee}      \\
            F^{\vee}x^{1\vee}    \\
          \end{array}
      \right)
     > = \\
   & -<x, x^{1\vee}> - 2<y, D^t{x}^{1\vee}> +
     <y, F^{\vee}{x}^{1\vee}> = \\
   & -<x, x^{1\vee}> - <y,  F^{\vee}x^{1\vee}> =
      -<z, \tilde{z}^{1\vee}> = 0.
  \end{split}
\end{equation}
Thus, $w_1{H} \subset H$. Similarly, $w_2{H} \subset H$.
\hspace{3mm} \qedsymbol

To prove the sufficient condition of regularity for arbitrary orientation, we
need some properties of transforming elements $T$.

\subsection{A theorem on transforming elements}
\index{transforming elements}
\begin{proposition} [\cite{St82}]
  \label{prop_transf}
Let $\varDelta{'}, \varDelta{''}$ be two arbitrary orientations of
the graph $\varGamma$ that differ by the direction of $k$ edges.
Consider the chain of orientations, in which every two adjacent
orientations differ by the direction of one edge:
\begin{equation}
     \varDelta{'} = \varDelta_0, \varDelta_1, \varDelta_2, \dots,
     \varDelta_{k-1}, \varDelta_k = \varDelta{''}.
\end{equation}
   Then, in the
   Weyl group, there exist elements $P_i$ and $S_i$, where $i = 1,2,...,k$,
such
that
\begin{equation}
 \label{transf_arrow}
  \begin{split}
    & {\bf C}_{\varDelta_0} = P_1{S}_1,  \\
    & {\bf C}_{\varDelta_1} = S_1{P}_1 = P_2{S}_2, \\
    & \dotfill
    \\
    & {\bf C}_{\varDelta_{k-1}} = S_{k-1}P_{k-1} = P_k{S}_k,   \\
    & {\bf C}_{\varDelta_k} = S_k{P}_k.
  \end{split}
\end{equation}
In addition, for each reflection and each $i=1,2,...,k$, this reflection does
not
occur in the decomposition of either $P_i$ or $S^{-1}_i$.
\end{proposition}
\PerfProof
 It suffices to consider the case $k = 1$. Let us consider the graph
$$
     \varGamma_1 \cup \varGamma_2 = \varGamma{\backslash}l.
$$
The graph $\varGamma$ can be depicted as follows:
\begin{equation}
  \begin{array}{cccc}
    \varDelta{'} \hspace{7mm}&  \dots    & \longleftarrow & \dots   \\
                          &  \varGamma_1 &                & \varGamma_2
\vspace{3mm} \\
    \varDelta{''} \hspace{7mm}&  \dots    & \longrightarrow & \dots     \\
                          &  \varGamma_1 &                & \varGamma_2 \\
  \end{array}
\end{equation}
The orientations $\varDelta{''}$ and $\varDelta{'}$ induce the
same orientations on the graphs $\varGamma_1$ and $\varGamma_2$,
and therefore they induce the same Coxeter transformations on
subgraphs $\varGamma_1$ and $\varGamma_2$. Denote the
corresponding Coxeter transformations by ${\bf C}_{\varGamma_1}$
and ${\bf C}_{\varGamma_2}$.  Then
\begin{equation*}
   {\bf C}_{\varDelta{'}} = {\bf C}_{\varGamma_2}{\bf C}_{\varGamma_1}, \hspace{7mm}
   {\bf C}_{\varDelta{''}} = {\bf C}_{\varGamma_1}{\bf C}_{\varGamma_2}.
\end{equation*}
Here,
\begin{equation*}
   P_1 = {\bf C}_{\varGamma_2}, \hspace{7mm}
   S_1 = {\bf C}_{\varGamma_1}. \hspace{6mm} \qed
\end{equation*}

\begin{remark} {\rm
 \label{conj_class}
Observe, that Proposition \ref{prop_transf} gives, in particular,
a simple proof of the fact that all the Coxeter transformations
form one conjugacy class, cf. \cite[Ch.5, \S6]{Bo}.}
\end{remark}

\index{theorem! - on transforming elements}
\begin{theorem} [\cite{St82}]
  \label{th_transf}
  1) Under the condition of Proposition \ref{prop_transf},
$$
  {T}^{-1}{\bf C}_{\varDelta{'}}T = {\bf C}_{\varDelta{''}}
$$
 for the following $k+1$ transforming elements $T:=T_i$:
\begin{equation}
 \label{look_transf}
 \begin{split}
   & T_1 = P_1{P}_2{P}_3...P_{k-2}{P}_{k-1}P_k, \\
   & T_2 = P_1{P}_2{P}_3...P_{k-2}{P}_{k-1}S^{-1}_k, \\
   & T_3 = P_1{P}_2{P}_3...P_{k-2}{S}^{-1}_{k-1}S^{-1}_k, \\
   &    \dots
                               \\
   & T_{k-1} =
            P_1{P}_2{S}^{-1}_3...{S}^{-1}_{k-2}S^{-1}_{k-1}S^{-1}_k,    \\
   & T_k =
            P_1{S}^{-1}_2{S}^{-1}_3...{S}^{-1}_{k-2}{S}^{-1}_{k-1}{S}^{-1}_k, \\
   & T_{k+1} =
      {S}^{-1}_1{S}^{-1}_2{S}^{-1}_3...{S}^{-1}_{k-2}{S}^{-1}_{k-1}{S}^{-1}_k.
 \end{split}
\end{equation}
 In addition, for each reflection $\sigma_\alpha$, there exists a $T_i$ whose
 decomposition does not contain this reflection.

2) The following relation holds:
\begin{equation}
  \label{p_q}
    T_p{T}^{-1}_q = {\bf C}^{q-p}_{\varDelta{'}}.
\end{equation}
\end{theorem}
\PerfProof
1) There are altogether $2^k$ transforming elements of the form
\begin{equation}
  \label{2power_k}
 T = X_1X_2...X_{k-1}X_k, \text{ where } X_i \in \{P_i, S^{-1}_i\}.
\end{equation}
By Proposition \ref{prop_transf}, for each reflection
$\sigma_\alpha$ and for each $i$, we can select $X_i = P_i$ or
$X_i = S^{-1}_i$ such that this reflection does not occur in $T$.
Taking the product of all these elements $X_i$ we obtain the
transforming element $T$ whose decomposition does not contain this
reflection. It remains to show that every transforming element
from the list (\ref{2power_k}) containing $2^k$ elements is of the
form (\ref{look_transf}). By (\ref{transf_arrow}) we have
$$
    S^{-1}_q{P}_{q+1} = P_q{S}^{-1}_{q+1}.
$$
Thus, all symbols $S^{-1}_i$ can be shifted to the right
 and all symbols $P_i$ can be shifted to the left.

2) It suffices to show (\ref{p_q}) for $p < q$. By (\ref{look_transf})
we have
\begin{equation}
 \label{long_prod}
  \begin{split}
   T_p{T}^{-1}_q = &
    (P_1{P}_2...{P}_{k-p+1}{S}^{-1}_{k-p+2}...{S}^{-1}_k)\times
    (S_k...{S}_{k-q+2}{P}_1{P}^{-1}_{k-q+1}...{P}^{-1}_1) = \\
  & P_1{P}_2...{P}_{k-p+1}(S_{k-p+1}{S}_{k-p}...{S}_{k-q+2})
      {P}^{-1}_{k-q+1}...{P}^{-1}_1.
 \end{split}
\end{equation}
Here, $k-p+1 \geq k-q+2$.
In order to simplify (\ref{long_prod}), we use formulas (\ref{transf_arrow}):
\begin{equation}
 \label{simplify_prod}
  \begin{split}
     & P_{k-p}P_{k-p+1}S_{k-p+1}S_{k-p} =
       P_{k-p}(S_{k-p}P_{k-p})S_{k-p} =
      (P_{k-p}S_{k-p})^2, \\
     & \dotfill  \\
     & P_{k-p-1}(P_{k-p}S_{k-p})^2{S}_{k-p-1} = \\
     & P_{k-p-1}(S_{k-p-1}P_{k-p-1})^2{S}_{k-p-1} =
        (P_{k-p-1}S_{k-p-1})^3, \\
     & \dotfill  \\
     & P_{k-q+2}{P}_{k-q+3}...{P}_{k-p+1}
       S_{k-p+1}...{S}_{k-q+2}{S}_{k-q+2} =
       (P_{k-q+2}{S}_{k-q+2})^{q-p} .
  \end{split}
\end{equation}
Thus, from (\ref{long_prod}) and (\ref{simplify_prod}) we deduce:
\begin{equation}
 \label{long_prod_1}
   T_p{T}^{-1}_q =
    P_1{P}_2...{P}_{k-q+1}(P_{k-q+2}{S}_{k-q+2})^{q-p}
    P^{-1}_{k-q+1}...{P}^{-1}_1.
\end{equation}
Again, by (\ref{transf_arrow}) we have
\begin{equation}
 \label{simplify_prod_2}
  \begin{split}
     & {P}_{k-q+1}(P_{k-q+2}{S}_{k-q+2})^{q-p}{P}^{-1}_{k-q+1} = \\
     & \qquad
       P_{k-q+1}(S_{k-q+1}P_{k-q+1})^{q-p}P^{-1}_{k-q+1} =
       (P_{k-q+1}S_{k-q+1})^{q-p}, \vspace{2mm} \\
     & \dotfill \\
     & P_{k-q}(P_{k-q+1}S_{k-q+1})^{q-p}P^{-1}_{k-q} = \\
     & \qquad
       P_{k-q}(S_{k-q}P_{k-q})^{q-p}P^{-1}_{k-q} =
       (P_{k-q}S_{k-q})^{q-p},  \vspace{2mm} \\
     & \dotfill \\
     & P_1(P_2{S}_2)^{q-p}P^{-1}_1 =
       P_1(S_1{P}_1)^{q-p}{P}^{-1}_1 =
       (P_1{S}_1)^{q-p}.
    \end{split}
\end{equation}
Therefore,
\begin{equation}
 \label{simplify_prod_3}
   T_p{T}^{-1}_q = (P_1{S}_1)^{q-p} =
     {\bf C}^{q-p}_{\varDelta_0} = {\bf C}^{q-p}_{\varDelta{'}}.
\qed
\end{equation}

\begin{remark} {\rm
Theorem \ref{th_transf} allows us to select transforming element
$T$ in such a way that its decomposition does not contain any
given refection $\sigma_i$, and therefore $T$ does not change any
given coordinate $i$. This fact was already used once in
\S\ref{necessary_cond} for the proof of the necessary conditions
of regularity for diagrams with indefinite Tits form. Now, we will
use Theorem \ref{th_transf} to carry the sufficient condition of
regularity from a bicolored orientation $\varDelta$ in Proposition
\ref{prop_sufficient} to an arbitrary orientation $\varDelta{'}$,
see Definition \ref{def_regular}. }
\end{remark}

\subsection{Sufficient conditions of regularity for an arbitrary orientation}

\index{conditions of regularity! - sufficient conditions,
arbitrary orientation}
\begin{proposition} [\cite{St82}]
 \label{arbitr_orient}
Let $\varGamma$ be an extended Dynkin diagram. Let $\varDelta =
\varDelta_0$, $\varDelta{'} = \varDelta_k$. If $z$ is a positive
root with zero $\varDelta{'}$-defect, then $z$ is the
$\varDelta{'}$-regular vector.
\end{proposition}

\PerfProof If $z > 0$, then there exists a positive coordinate
$z_\alpha > 0$. Take transforming element $T_i$ whose
decomposition does not contain the reflection $\sigma_\alpha$.
Then $(T_i{z})_\alpha > 0$. Since $T_i{z}$ is the root, we have
$T_i{z} > 0$. The equality $\rho_{\varDelta{'}}({z}) = 0$ means
that $\rho_\varDelta(T_i{z}) = 0$, see Definition
\ref{rho_defect}. Since $T_i{z} > 0$ and $T_i{z}$ is the root, we
see by Proposition \ref{prop_sufficient} that $T_i{z}$ is
$\varDelta$-regular and
\begin{equation}
  \label{t_i_pos}
   {\bf C}^m_\varDelta{T}_i{z} > 0 \text{ for all } m \in \mathbb{N}.
\end{equation}
Suppose that
$$
   u = {T}^{-1}_i{\bf C}^m_\varDelta{T}_i{z} < 0
   \text{ for some } m = m_0.
$$
Then by Theorem \ref{th_transf} there exists $T_j$ such that
$T_j{u} < 0$, i.e.,
\begin{equation}
  \label{t_i_j}
   T_j{T}^{-1}_i{\bf C}^m_\varDelta{T}_i{z} < 0
   \text{ for } m = m_0.
\end{equation}
Again by Theorem \ref{th_transf} we have $T_j{T}^{-1}_i = {\bf
C}^{i-j}_\varDelta$ and by (\ref{t_i_j}) the following relation
holds:
\begin{equation}
  \label{t_i_j_1}
   {\bf C}^{i-j}_\varDelta{\bf C}^m_\varDelta{T}_i{z} < 0
   \text{ for } m = m_0
\end{equation}
that contradicts to (\ref{t_i_pos}). Therefore,
$$
   u = T^{-1}_i{\bf C}^m_\varDelta{T}_i{z} > 0
   \text{ for all } m \in \mathbb{N}.
$$
Thus, $z$ is a $\varDelta{'}$-regular vector. \hspace{3mm}
\qedsymbol

From Theorem \ref{necess_codnd_reg} together with Proposition
\ref{arbitr_orient} we get the following

\index{theorem! - on regular representations of quivers}

\begin{theorem}
\label{necess_eq_suff} The indecomposable representation $V$ of
the graph $\varGamma$ (which is an extended Dynkin diagram) with
orientation $\varDelta{'}$ is regular in the orientation
$\varDelta{'}$ if and only if
$$
    \rho_{\varDelta{'}}(\dim{V}) = 0.
$$
\end{theorem}
\PerfProof
\index{imaginary roots}
Indeed, the dimensions of the indecomposable representations of
the extended Dynkin diagrams are roots, \cite{DR76}, \cite{Kac80}.
However, there are indecomposable representations whose dimensions
are not usual (real) roots but {\it imaginary roots}, see \S\ref{roots}.
They are vectors from the kernel of the Tits form, and are proportional to
vectors $z^{1}$ which are fixed points of the Weyl group. In particular,
$T{z}^{1} = z^{1}$, so
$<z^{1},\tilde{z}^{1\vee}> = 0$ directly implies
$<Tz^{1},\tilde{z}^{1\vee}> = 0$.
\hspace{3mm} \qedsymbol

\subsection{Invariance of the defect}
We will show that the ${\varDelta{'}}$-defect
$\rho_{\varDelta{'}}$ is invariant under the Coxeter
transformation ${\bf C}_{\varDelta{'}}$ and $\rho_{\varDelta{'}}$
does not depend on the choice of the transforming element $T_i$ in
the Weyl group, see Definition \ref{rho_defect} and Theorem
\ref{th_transf}. In other words, the following proposition holds.
\begin{proposition}
 \label{invariance}
1) The Coxeter transformation ${\bf C}_{\varDelta{'}}$ preserves
the linear form $\rho_{\varDelta{'}}$:
\begin{equation}
 \label{invar_rho}
  \rho_{\varDelta{'}}({\bf C}_{\varDelta{'}}z) = \rho_{\varDelta{'}}(z)
  \text{ for any vector } z.
\end{equation}

2) If $T_i$ and $T_j$ are transforming elements defined by (\ref{look_transf}),
i.e.,
\begin{equation}
 \label{condit_transf}
  T_i^{-1}{\bf C}_\varDelta{T_i} = {\bf C}_{\varDelta{'}} \text{ and }
  T_j^{-1}{\bf C}_\varDelta{T_j} = {\bf C}_{\varDelta{'}},
\end{equation}
then we have
\begin{equation}
  \label{invar_T}
  <T_i{z}, \tilde{z}^{1\vee}> ~=~ <T_j{z}, \tilde{z}^{1\vee}>
  \text{ for any vector } z.
\end{equation}
\end{proposition}

\PerfProof First, observe that (\ref{invar_T}) does not hold for
an arbitrary matrix $T$, since the matrix $kT$, where $k\in
\mathbb{R}$, also satisfies (\ref{condit_transf}), but does not
satisfy (\ref{invar_T}).

In (\ref{wz_z}) we showed that
\begin{equation}
 \label{wz_z_1}
  <w_i{z}, \tilde{z}^{1\vee}> ~=~ -<z, \tilde{z}^{1\vee}>.
\end{equation}

By (\ref{wz_z_1}) we have
\begin{equation}
 \label{wz_z_2}
  <z, \tilde{z}^{1\vee}> ~=~ <{\bf C}_{\varDelta}{z}, \tilde{z}^{1\vee}>
  \text{ for any vector } z,
\end{equation}
i.e., (\ref{invar_rho}) holds for any bicolored orientation
$\varDelta$. Since $z$  in (\ref{wz_z_2})is an arbitrary vector,
we have
\begin{equation}
 \label{wz_z_3}
  <Tz, \tilde{z}^{1\vee}> ~=~ <{\bf C}_{\varDelta}{Tz}, \tilde{z}^{1\vee}> ~=~
  <T{\bf C}_{\varDelta{'}}{z}, \tilde{z}^{1\vee}>,
\end{equation}
or
\begin{equation}
 \label{wz_z_4}
  \rho_{\varDelta{'}}(z) = \rho_{\varDelta{'}}({\bf C}_{\varDelta{'}}z),
\end{equation}
i.e., (\ref{invar_rho}) holds for an arbitrary orientation
$\varDelta{'}$.

Further, from (\ref{wz_z_2}) we have
\begin{equation}
 \label{wz_z_5}
  \rho_\varDelta(z) = \rho_\varDelta({\bf C}_\varDelta{z}) =
  \rho_\varDelta({\bf C}_\varDelta^{i-j}{z}).
\end{equation}
By (\ref{p_q})
\begin{equation}
 \label{wz_z_6}
  \rho_\varDelta(z) = \rho_\varDelta(T_j{T_i^{-1}}{z}).
\end{equation}
Substituting $T_i{z}$ instead of $z$ we have
\begin{equation*}
  \rho_\varDelta(T_i{z}) = \rho_\varDelta(T_j{z}),
\end{equation*}
and (\ref{invar_T}) is proved. \qedsymbol


\section{Examples of necessary conditions of regularity}

\index{branch point}
\index{central orientation}
\index{point of
maximum branching}
\begin{definition} {\rm
Let $\varGamma$ be a graph, $x_0 \in \varGamma$ be a point of
maximal branching degree  and let all arrows of $\varGamma$ be
directed to the point $x_0$. The corresponding orientation is said
to be the {\it central orientation} and denoted by $\varDelta_0$.
}
\end{definition}

We will consider the necessary conditions of regularity for some
diagrams in bicolored, central, and other orientations. If an
orientation $\varDelta{'}$ is obtained from another orientation
$\varDelta{''}$ by reversing of all arrows of a graph, then the
corresponding conditions of the regularity coincide:
\begin{equation}
 \rho_{\varDelta{'}} = \rho_{\varDelta{''}} \hspace{3mm}
 \text{if} \hspace{3mm}\mathcal{B} > 0
\end{equation}
or
\begin{equation}
 \rho^1_{\varDelta{'}} = \rho^1_{\varDelta{''}}, \hspace{3mm}
 \rho^2_{\varDelta{'}} = \rho^2_{\varDelta{''}} \hspace{3mm}
 \text{if} \hspace{3mm}\mathcal{B} \text{is indefinite. }\\
\end{equation}

Let us introduce an equivalence relation $\mathcal{R}$ on the set
of orientations of the graph $\varGamma$. Two orientations
$\varDelta{'}$ and $\varDelta{''}$ will be called {\it
equivalent}, if one orientation can be obtained from the other one
by reversing all arrows or by an automorphism of the diagram; for
equivalent orientations, we write:
$$
     \varDelta{'} \equiv \varDelta{''} (mod \mathcal{R}).
$$
Equivalent orientations have the identical conditions of
regularity.

\subsection{The three equivalence classes of orientations of
$\tilde{D}_4$}

\begin{figure}[h]
\centering
\includegraphics{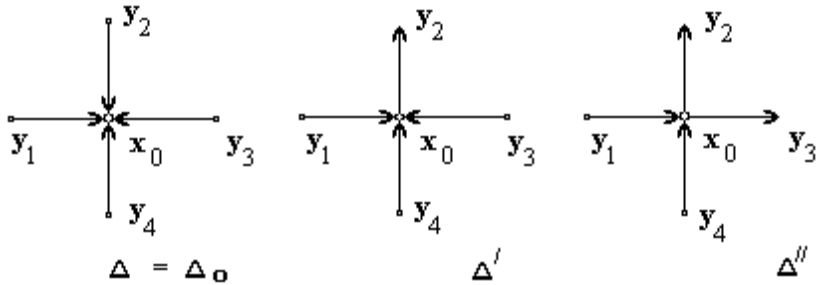}
\caption{\hspace{3mm} For $\tilde{D}_4$, the bicolored
  and central orientations coincide}
\label{bicol_central}
\end{figure}

In the case of $\tilde{D}_4$, the bicolored orientation
$\varDelta$ coincides with the central orientation $\varDelta_0$.
The orientations $\varDelta$, $\varDelta{'}$, and $\varDelta{''}$
cover all possible equivalence classes.

a) The bicolored orientation $\varDelta$, see
Fig.~\ref{bicol_central}.

Here, $T = I$ and by Theorem \ref{necess_codnd_reg} and
Proposition \ref{prop_sufficient} we have
\begin{equation}
 \label{regul_GP}
 z^1 = \left (
 \begin{array}{c}
         2 \\
         1 \\
         1 \\
         1 \\
         1 \\
 \end{array} \right )
 \begin{array}{c}
         x_0 \\
         y_1 \\
         y_2 \\
         y_3 \\
         y_4 \\
 \end{array},
 \hspace{7mm}
 \tilde{z}^1 = \left (
 \begin{array}{c}
         2 \\
         -1 \\
         -1 \\
         -1 \\
         -1 \\
 \end{array} \right ), \hspace{7mm}
 \rho_{\varDelta}(z) = y_1 + y_2 + y_3 + y_4 - 2x_0.
 \hspace{7mm}
\end{equation}

Originally, the linear form $\rho_{\varDelta}(z)$ in
(\ref{regul_GP}) was obtained by I.~M.~Gelfand and V.~A.~Ponomarev
in the work devoted to classifications of quadruples of linear
subspaces of arbitrary dimension \cite{GP72}.

b) The orientation $\varDelta{'}$, see Fig.~\ref{bicol_central}.
Here we have
\begin{equation}
 {\bf C}_{\varDelta} =
 \sigma_{y_4}\sigma_{y_3}\sigma_{y_2}\sigma_{y_1}\sigma_{x_0},
 \hspace{7mm}
 {\bf C}_{\varDelta{'}} =
 \sigma_{y_4}\sigma_{y_3}\sigma_{y_1}\sigma_{x_0}\sigma_{y_2},
 \hspace{7mm}
 T = \sigma_{y_2}.
\end{equation}
Then, by Theorem \ref{necess_codnd_reg} and Proposition
\ref{arbitr_orient} we get the following condition of the
$\varDelta{'}$-regularity:
\begin{equation}
 \label{regul_D4_2}
 Tz = \left (
 \begin{array}{c}
         x_0 \\
         y_1 \\
         x_0 - y_2 \\
         y_3 \\
         y_4 \\
 \end{array} \right ), \hspace{7mm}
 \begin{array}{c}
         y_1 + (x_0 - y_2) + y_3 + y_4 - 2x_0 = 0, \text{ or } \\
         y_1 +y_3 + y_4 = y_2 + x_0.
 \end{array}
\end{equation}

c) The orientation $\varDelta{''}$, see Fig.~\ref{bicol_central}.
Here, $T = \sigma_{y_2}\sigma_{y_3}$ (or
$\sigma_{y_1}\sigma_{y_4}$), We have the following condition of
the $\varDelta{''}$-regularity:
\begin{equation}
 \label{regul_D4_3}
 \begin{array}{c}
         y_1 + (x_0 - y_2) + (x_0 - y_3) + y_4 - 2x_0 = 0, \text{ or }
\\
         y_1 + y_4 = y_2 + y_3.
 \end{array}
\end{equation}

\subsection{The bicolored and central orientations of
$\tilde{E}_6$}

\begin{figure}[h]
\centering
\includegraphics{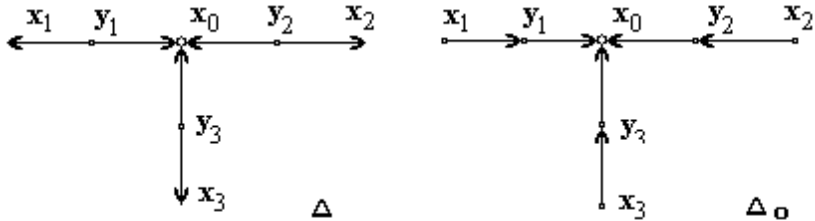}
\caption{\hspace{3mm} For $\tilde{E}_6$, the bicolored and central
orientations}
\label{bicol_central_E6}
\end{figure}

We consider only the bicolored orientation $\varDelta$ and the
central orientation $\varDelta_0$, Fig.~\ref{bicol_central_E6}.
The Coxeter transformations and transforming element $T$ are:
\begin{equation}
 \begin{split}
 & {\bf C}_{\varDelta} =
  \sigma_{y_3}\sigma_{y_2}\sigma_{y_1}
  \sigma_{x_3}\sigma_{x_2}\sigma_{x_1}\sigma_{x_0}, \\
 & {\bf C}_{\varDelta_0} = \sigma_{x_3}\sigma_{x_2}\sigma_{x_1}
  \sigma_{y_3}\sigma_{y_2}\sigma_{y_1}\sigma_{x_0}, \\
 &  T = \sigma_{x_3}\sigma_{x_2}\sigma_{x_1}.
 \end{split}
\end{equation}

We have

\begin{equation}
 \label{regul_GP_3}
 z^1 = \left (
 \begin{array}{c}
         3 \\
         1 \\
         1 \\
         1 \\
         2 \\
         2 \\
         2 \\
 \end{array} \right )
 \begin{array}{c}
         x_0 \\
         x_1 \\
         x_2 \\
         x_3 \\
         y_1 \\
         y_2 \\
         y_3 \\
 \end{array},
 \hspace{7mm}
 \tilde{z}^1 = \left (
 \begin{array}{c}
         3 \\
         1 \\
         1 \\
         1 \\
         -2 \\
         -2 \\
         -2 \\
 \end{array} \right ), \hspace{7mm}
 Tz = \left (
 \begin{array}{c}
         x_0 \\
         y_1 - x_1 \\
         y_2 - x_2 \\
         y_3 - x_3 \\
         y_1 \\
         y_2 \\
         y_3 \\
 \end{array} \right ),
\end{equation}
and by Theorem \ref{necess_codnd_reg} and Proposition
\ref{prop_sufficient} a condition of $\varDelta$-regularity is
\begin{equation}
 \rho_{\varDelta}(z) = 3x_0 + x_1 + x_2 + x_3 + 3x_0 - 2y_1 - 2y_2 - 2y_3.
 \hspace{7mm}
\end{equation}

The choice of element $T$ is ambiguous, but it was shown in
Proposition \ref{invariance} that the regularity condition does
not depend on this choice.

So, a condition of $\varDelta_0$-regularity is
\begin{equation*}
  \begin{split}
     \rho_{\varDelta_0}(z) & = \\
     & 3x_0 + (y_1 - x_1) + (y_2 - x_2) + (y_3 - x_3)
         -2(y_1 + y_2 + y_3) = 0,
   \end{split}
\end{equation*}
 or
\begin{equation}
     \rho_{\varDelta_0}(z) = x_1 + x_2 + x_3 + y_1 + y_2 + y_3 - 3x_0.
\end{equation}

\subsection{The multiply-laced case. The two orientations
 of $\tilde{G}_{21}$ and $\tilde{G}_{22} = {\tilde{G}^{\vee}}_{21}$}

\begin{figure}[h]
\centering
\includegraphics{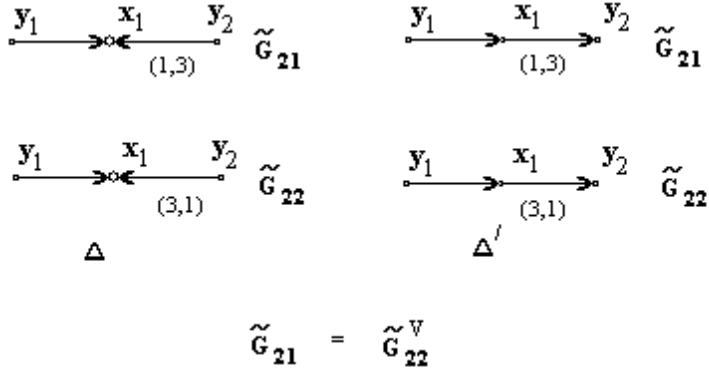}
\caption{\hspace{3mm} The two orientations
 of $\tilde{G}_{21}$ and $\tilde{G}_{22} = {\tilde{G}^{\vee}}_{21}$}
\label{bicol_central_G22}
\end{figure}

We have

\begin{equation}
 \label{regul_GP_4}
 z^1 = \left (
 \begin{array}{c}
         2 \\
         1 \\
         3 \\
 \end{array} \right )
 \begin{array}{c}
         x_1 \\
         y_1 \\
         y_2 \\
 \end{array}, \quad
 \tilde{z}^1 = \left (
 \begin{array}{c}
         2 \\
         -1 \\
         -3 \\
 \end{array} \right ), \quad
 {z}^{1\vee} = \left (
 \begin{array}{c}
         2 \\
         1 \\
         1 \\
 \end{array} \right )
 \begin{array}{c}
         x_1 \\
         y_1 \\
         y_2 \\
 \end{array},   \quad
 {\tilde{z}}^{1\vee} = \left (
 \begin{array}{c}
         2 \\
         -1 \\
         -1 \\
 \end{array} \right ).
\end{equation}

Take vectors $\tilde{z}^1$ and ${\tilde{z}}^{1\vee}$ from
(\ref{regul_GP_4}) and substitute them in (\ref{regular_2}). We
get the following regularity condition $\rho_\varDelta$ (resp.
$\rho^{\vee}_\varDelta$):
\begin{equation}
\begin{split}
& \rho_\varDelta(z) = y_1 + y_2 - 2x_1 \quad \text{ for } \tilde{G}_{21}, \\
& \rho^{\vee}_\varDelta(z) = y_1 + 3y_2 - 2x_1 \quad \text{ for }
\tilde{G}_{22}.
\end{split}
\end{equation}
Now consider the orientation $\varDelta{'}$,
Fig.~\ref{bicol_central_G22}. The Coxeter transformations can be
expressed as follows:
\begin{equation}
 \begin{split}
 & {\bf C}_\varDelta = \sigma_{y_2}\sigma_{y_1}\sigma_{x_1},
 \hspace{7mm}
 {\bf C}^{\vee}_\varDelta =
 \sigma^{\vee}_{y_2}\sigma^{\vee}_{y_1}\sigma_{x_1}^{\vee}, \\
 & {\bf C}_{\varDelta{'}} = \sigma_{y_1}\sigma_{x_1}\sigma_{y_2},
 \hspace{7mm}
 {\bf C}^{\vee}_{\varDelta{'}} =
 \sigma^{\vee}_{y_1}\sigma_{x_1}^{\vee}\sigma^{\vee}_{y_2},
 \end{split}
\end{equation}
Transforming elements $T$ and $T^{\vee}$ are:
\begin{equation}
 T = \sigma_{y_2},
 \hspace{7mm}
 T^{\vee} = \sigma^{\vee}_{y_2}. \\
\end{equation}
and
\begin{equation}
 Tz = \left (
 \begin{array}{c}
         x_1 \\
         y_1 \\
         3x_1 - y_2 \\
 \end{array} \right ), \hspace{7mm}
 T^{\vee}z = \left (
 \begin{array}{c}
         x_1 \\
         y_1 \\
         x_1 - y_2 \\
 \end{array} \right ).
\end{equation}
Finally, we have
\begin{equation}
 \begin{split}
 & \rho_{\varDelta{'}}(z) = y_1 - y_2 + x_1 \quad \text{ for } \tilde{G}_{21}, \\
 & \rho^{\vee}_{\varDelta{'}}(z) = y_1 - 3y_2 + x_1 \quad \text{ for } \tilde{G}_{22}.
 \end{split}
\end{equation}

\subsection{The case of indefinite $\mathcal{B}$. The oriented star
$\ast_{n+1}$}
 \label{star_example}
Consider the oriented star $\ast_{n+1}$ with a bicolored
orientation.

\begin{figure}[h]
\centering
\includegraphics{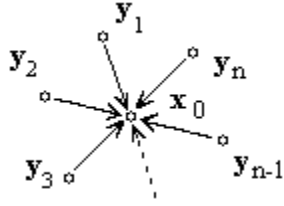}
\caption{\hspace{3mm} A bicolored orientation of the star
$*_{n+1}$}
\label{bicolor_orient}
\end{figure}

According to Remark  \ref{remark_star} the matrix $DD^t$ is a
scalar. By (\ref{dd_formula}) $DD^t = \displaystyle\frac{n}{4}$
and the maximal eigenvalue $\varphi^m = \displaystyle\frac{n}{4}$.
By (\ref{lambda_phi}) we have

\begin{equation}
 \label{lambdas_star}
 \lambda^m_{1,2} = \frac{n - 2 \pm\sqrt{n(n-4)}}{2} \hspace{1mm}.
\end{equation}

\noindent Let $x_m = 1$. Then
$$
 \frac{-2}{\lambda^m_{1,2} + 1}D^tx^m =
 \frac{-2}{\lambda^m_{1,2} + 1}[1,1,...,1]^t.
$$
 Thus, we have the following eigenvectors $z^m_1$, $z^m_2$:

\begin{equation*}
 z^m_1 = \left (
 \begin{array}{c}
    \lambda^m_1 + 1  \\
         1           \\
         1           \\
         ...             \\
         1           \\
 \end{array} \right ), \hspace{3mm}
 z^m_2 = \left (
 \begin{array}{c}
    \lambda^m_2 + 1  \\
         1           \\
         1           \\
         ...             \\
         1           \\
 \end{array} \right ),
\end{equation*}
and their conjugate vectors $\tilde{z}^m_1$, $\tilde{z}^m_2$ (see
Definition \ref{conjug_vector} and Proposition
\ref{ortogonal_eigenv}):
\begin{equation*}
 \tilde{z}^m_1 = \left (
 \begin{array}{c}
    \lambda^m_2 + 1  \\
         -1          \\
         -1          \\
         ...             \\
         -1          \\
 \end{array} \right ), \hspace{3mm}
 \tilde{z}^m_2 = \left (
 \begin{array}{c}
    \lambda^m_1 + 1  \\
         -1          \\
         -1          \\
         ...             \\
         -1          \\
 \end{array} \right ).
\end{equation*}

According to (\ref{regular_5}) we obtain the following condition
of $\varDelta$-regularity:
\begin{equation}
 \label{inequality_lambdas}
 (\lambda^m_2 + 1)x_0 - \sum{y_i} \leq 0, \hspace{5mm}
 (\lambda^m_1 + 1)x_0 - \sum{y_i} \geq 0.
\end{equation}

From (\ref{lambdas_star}) and (\ref{inequality_lambdas}) we deduce
\begin{equation}
 \label{inequality_lambdas_1}
 \begin{split}
 & \frac{n - \sqrt{n(n-4)}}{2}x_0 \leq \sum{y_i} \leq
 \frac{n + \sqrt{n(n-4)}}{2}x_0
     \text{  or} \vspace{7mm} \\
 & \frac{n - \sqrt{n(n-4)}}{2}\sum{y_i} \leq x_0 \leq
 \frac{n + \sqrt{n(n-4)}}{2}\sum{y_i}
  \text { , or} \vspace{7mm} \\
 & |x_0 - \frac{1}{2}\sum{y_i} | \leq \frac{\sqrt{n(n-4)}}{2n}\sum{y_i}
  \text {  or} \vspace{7mm} \\
 & x^2_0 - x_0\sum{y_i} + \frac{1}{4}(\sum{y_i})^2
 \leq \frac{1}{4}(\sum{y_i})^2 - \frac{1}{n}(\sum{y_i})^2
  \text {  or} \vspace{7mm} \\
 & x^2_0 - x_0\sum{y_i}
 \leq - \frac{1}{n}(\sum{y_i})^2
  \text {  or} \vspace{7mm} \\
 & x^2_0 - x_0\sum{y_i} + \sum{y^2_i}
 \leq \sum{y^2_i} - \frac{1}{n}(\sum{y_i})^2
 \end{split}
\end{equation}

Since the left hand side of the latter inequality is the Tits form
$\mathcal{B}$, we obtain the following condition of
$\varDelta$-regularity:
\begin{equation}
 \label{inequality_lambdas_2}
 \mathcal{B}(z)
 \leq \frac{1}{n}\sum\limits_{0 < i < j}(y_i - y_j)^2.
\end{equation}

%% file: 7mckay.tex

\chapter{\sc\bf The McKay correspondence and the Slodowy generalization}

\setlength{\epigraphwidth}{75mm}

\epigraph{We have seen during the past few years a major assault
on the problem of determining all the finite simple groups. ... If
I am right, I foresee new proofs of classification which will owe
little or nothing to the current proofs. They will be much shorter
and will help us to understand the simple groups in a context much
wider than finite group theory.}{J.~McKay, \cite[p.183]{McK80},
1980}

 \label{chapter_slodowy}
\section{Finite subgroups of $SU(2)$ and $SO(3,\mathbb{R})$}
\label{finite_subgroups} Let us consider the special unitary group
$SU(2)$, the subgroup of unitary transformations in
$GL(2,\mathbb{C})$ with determinant 1 and its quotient group
$PSU(2)=SU(2)/\{\pm1 \}$ acting on the complex projective line
$\mathbb{CP}^1$, see \S\ref{sect_projective_line}. The projective
line $\mathbb{CP}^1$ can be identified with the sphere $S^2
\subset \mathbb{R}^3$, see \S\ref{sect_projective_line} and
(\ref{identity_F}).

\index{complex projective line}

The transformations of the sphere $S^2$ induced by elements of
$SO(3,\mathbb{R})$ (orientation preserving rotations of
$\mathbb{R}^3$) correspond under these identifications to
transformations of $\mathbb{CP}^1$ in the group $PSU(2)$:
\begin{equation}
 \label{isom_group}
     PSU(2) \cong SO(3,\mathbb{R}),
\end{equation}
see \cite[Prop.4.4.3]{Sp77}, \cite[\S0.13]{PV94}.

\index{dihedron}
\index{tetrahedron}
\index{octahedron}
\index{dodecahedron}
It is well known that the finite subgroups of $SO(3,\mathbb{R})$
are precisely the rotation groups of the following polyhedra:
the regular $n$-angled pyramid, the $n$-angled dihedron
(a regular plane $n$-gon with two faces), the tetrahedron,
the cube (its the rotation group coincides with the rotation group of the
octahedron),
the icosahedron (its the rotation group coincides with the rotation group of the
dodecahedron).
By (\ref{isom_group}) we have a classification of the
finite subgroups of $PSU(2)$, see Table \ref{rotation_pol_1}.

To get a classification of all finite subgroups of $SU(2)$ we make
use of the double covering
$$
   \pi : SU(2) \longrightarrow SO(3,\mathbb{R}).
$$
If $G$ is a finite subgroup of $SO(3,\mathbb{R})$, we see that the
preimage $\pi^{-1}(G)$ is a finite subgroup of $SU(2)$ and
$|\pi^{-1}(G)| = 2|G|$. \index{binary polyhedral group} The finite
subgroups of $SO(3,\mathbb{R})$ are called {\it polyhedral
groups}, see \text{Table \ref{rotation_pol_1}}. The finite
subgroups of $SU(2)$ are naturally called {\it binary polyhedral
groups}, see Table \ref{binary_pol_1}.
\begin{remark}
 \label{cyclic_group}
{\rm 1) The cyclic group is an exceptional case.
  The preimage of the cyclic group
  $G = \mathbb{Z}/n\mathbb{Z}$ is the even cyclic group
  $\mathbb{Z}/2n\mathbb{Z}$, which is a binary cyclic group.
  The cyclic group $\mathbb{Z}/(2n-1)\mathbb{Z}$ is not a
  preimage with respect to $\pi$, so $\mathbb{Z}/(2n-1)\mathbb{Z}$ is not a
binary polyhedral group. Anyway, the cyclic groups
$\mathbb{Z}/n\mathbb{Z}$ complete the list of finite subgroups of
$SU(2)$.

2) Every finite subgroup of $SL(2,\mathbb{C})$ is conjugate
  to a subgroup in $SU(2)$, so we also have a classification
  of finite subgroups of $SL(2,\mathbb{C})$, \cite{Sp77}.
}
\end{remark}

For every polyhedral group $G$, the axis of rotations under an
element $\gamma \in G$ passes through either the mid-point of a
face, or the mid-point of an edge, or a vertex. We denote the
orders of symmetry of these axes by $p$, $q$, and $r$,
respectively, see \cite{Sl83}. These numbers are listed in Table
\ref{rotation_pol_1}. The triples $p, q, r$ listed in Table
\ref{rotation_pol_1} are exactly the solutions of the diophantine
inequality
\begin{equation}
 \label{diophantine_inequality}
  \frac{1}{p} + \frac{1}{q} + \frac{1}{r} > 1,
\end{equation}
see \cite{Sp77}.

\begin{table} 
  \centering
  \vspace{2mm}
  \caption{\hspace{3mm}The polyhedral groups in $\mathbb{R}^3$}
  \renewcommand{\arraystretch}{1.9}
 \begin{tabular} {||c|c|c|c|c||}
 \hline \hline
   Polyhedron & Orders of symmetries & Rotation group & Group order \\
 \hline \hline
     Pyramid   &  $-$  & cyclic & $n$ \\
 \hline
     Dihedron  &  $n$  2 2  & dihedral & $2n$ \\
 \hline
     Tetrahedron &  $3$  2 3  & $\mathcal{A}_4$ & 12 \\
 \hline
     Cube    &  $4$  2 3  & $\mathcal{S}_4$ & 24 \\
     Octahedron  &  $3$  2 4  & $\mathcal{S}_4$ & 24 \\
 \hline
     Dodecahedron &  $5$  2 3  & $\mathcal{A}_5$ & 60 \\
     Icosahedron &  $3$  2 5  & $\mathcal{A}_5$ & 60  \\
 \hline \hline
 \end{tabular}
\\ \vspace{2mm}
 Here, $\mathcal{S}_m$ (resp. $\mathcal{A}_m$) denotes the
 symmetric, (resp. alternating) \\
 group of all (resp. of all even)
 permutations of $m$ letters.
 \label{rotation_pol_1}
\end{table}

\section{Generators and relations in polyhedral groups}
    \label{gener_rel}
The quaternion group introduced by W.~R.~Hamilton is defined as follows:
It is
generated by three generators $i$, $j$, and $k$ subject to the relations
\begin{equation}
 \label{quaternion}
   i^2 = j^2 = k^2 = ijk = -1;
\end{equation}
for references, see \cite{Cox40}, \cite{CoxM84}. A natural
generalization of the quaternion group is the group generated by
three generators $R$, $S$, and $T$ subject to the relations
\begin{equation}
  \label{natural_gen}
   R^p = S^q = T^r = RST = -1.
\end{equation}
Denote by $\langle p, q, r \rangle$ the group defined by (\ref{natural_gen}).

W.~Threlfall, (see \cite{CoxM84}) has observed that
\begin{equation*}
  \langle 2, 2, n \rangle,
  \langle 2, 3, 3 \rangle,
  \langle 2, 3, 4 \rangle,
  \langle 2, 3, 5 \rangle
\end{equation*}
are the {\it binary polyhedral groups}, of order
\begin{displaymath}
  \displaystyle\frac{4}{\displaystyle\frac{1}{p} + \frac{1}{q} + \frac{1}{r} -
1}
\end{displaymath}
\begin{table} [h]
  \centering
  \vspace{2mm}
  \caption{\hspace{3mm}The finite subgroups of $SU(2)$}
  \index{group! - cyclic}
  \index{group! - binary dihedral}
  \index{group! - binary tetrahedral}
  \index{group! - binary octahedral}
  \index{group! - binary icosahedral}
  \renewcommand{\arraystretch}{1.9}
 \begin{tabular} {||c|c|c|c||}
 \hline \hline
   $\langle l, m, n \rangle$ & Order & Denotation & Well-known name \\
 \hline \hline
      $-$            &  $n$ & $\mathbb{Z}/n\mathbb{Z}$ & cyclic
group \\
 \hline
     $\langle 2, 2, n \rangle$ &    $4n$     & $\mathcal{D}_n$ & binary
dihedral group \\
 \hline
     $\langle 2, 3, 3 \rangle$ &    24  & $\mathcal{T}$ & binary tetrahedral
group \\
 \hline
     $\langle 2, 3, 4 \rangle$ &    48   & $\mathcal{O}$ & binary
octahedral group \\
 \hline
     $\langle 2, 3, 5 \rangle$ &    120  & $\mathcal{J}$ & binary
icosahedral group \\
 \hline \hline
 \end{tabular}
  \label{binary_pol_1}
\end{table}

H.~S.~M.~Coxeter proved in \cite[p.370]{Cox40} that having added one more
generator $Z$ one can replace (\ref{natural_gen})
by the
relations:
\index{group! - $\langle p,q,r \rangle$}
\index{group! - polyhedral}
\begin{equation}
  \label{natural_gen_1}
   R^p = S^q = T^r = RST = Z.
\end{equation}
For the groups $\langle 2, 2, n \rangle$ and $\langle 2, 3, n \rangle$, where $n
= 3,4,5$,
relation (\ref{natural_gen_1}) implies
$$
   Z^2 = 1.
$$
The polyhedral groups from Table \ref{rotation_pol_1} are described by
generations and relations as follows:
\begin{equation}
  \label{natural_gen_2}
   R^p = S^q = T^r = RST = 1,
\end{equation}
see \cite{Cox40}, \cite{CoxM84}.

\section{Kleinian singularities and the Du Val resolution}
 \index{orbit space}
 \index{algebra of invariants}
 \index{symmetric algebra}
Consider the quotient variety $\mathbb{C}^2/G$, where $G$ is a
binary polyhedral group $G \subset SU(2)$ from Table
\ref{binary_pol_1}. According to (\S\ref{orbit_space} and
(\ref{def_orbit_space})) $X = \mathbb{C}^2/G$ is the orbit space
given by the prime spectrum on the algebra of invariants $R^G$:
\begin{equation}
  \label{X_as_Spec}
 X := {\rm Spec}(R^G),
\end{equation}
 where $R = \mathbb{C}[z_1, z_2]$ (see eq.(\ref{R_algebra})) which
coincides with the symmetric algebra ${\rm
Sym}(({\mathbb{C}^2})^{*})$  (see \S\ref{graded_algebras}; for
more details, see \cite{Sp77}, \cite{PV94}, \cite{Sl80}.)

 F.~Klein \cite{Kl1884} observed that the algebra of invariants $\mathbb{C}[z_1, z_2]^G  $
for every binary polyhedral group $G \subset SU(2)$ from Table
\ref{binary_pol_1} can be considered by one approach:

 \index{theorem! - Klein}
\begin{theorem}[F.~Klein, \cite{Kl1884}]
The algebra of invariants $\mathbb{C}[z_1, z_2]^G$ is generated by
3 variables $x, y, z,$ subject to one essential relation
\begin{equation}
 \label{R_curve}
  R(x,y,z) = 0,
\end{equation}
where $R(x,y,z)$ is defined in column 2, Table \ref{klein_rel}. In
other words, the algebra of invariants $\mathbb{C}[z_1, z_2]^G$
coincides with the coordinate algebra (see \S\ref{coord_ring}) of
a curve defined by the eq. (\ref{R_curve}), i.e.,
\begin{equation}
 \mathbb{C}[z_1, z_2]^G \simeq \mathbb{C}[x, y, z]/(R(x,y,z)).
\end{equation}
\end{theorem}

\begin{table} [h]
  \centering
  \vspace{2mm}
  \caption{\hspace{3mm}The relations $R(x, y, z)$ describing the algebra of invariants
      $\mathbb{C}[z_1, z_2]^G$}
  \renewcommand{\arraystretch}{1.9}
 \begin{tabular} {||c|c|c||}
 \hline \hline
    Finite subgroup of $SU(2)$ & Relation $R(x, y, z)$ & Dynkin diagram \\
 \hline \hline
      $\mathbb{Z}/n\mathbb{Z}$ & $x^n + yz$ & $A_{n-1}$ \\
 \hline
      $\mathcal{D}_n$ & $x^{n+1} + xy^2 + z^2$ & $D_{n+2}$ \\
 \hline
      $\mathcal{T}$ & $x^4 + y^3 + z^2$ & $E_6$ \\
 \hline
      $\mathcal{O}$ & $x^3{y} + y^3 + z^2$ & $E_7$ \\
 \hline
      $\mathcal{J}$ & $x^5 + y^3 + z^2$ & $E_8$ \\
 \hline \hline
 \end{tabular}
  \label{klein_rel}
\end{table}

 \index{Kleinian singularity}
 \index{simple surface singularity}
 \index{Du Val singularity}
 \index{rational double point} 

The quotient $X$ from (\ref{X_as_Spec}) has no singularity except
an the origin $O \in \mathbb{C}^3$. The quotient variety $X$ is
called a {\it Kleinian singularity} also known as a {\it Du Val
singularity}, a {\it simple surface singularity} or a {\it
rational double point}.

See Remark \ref{alg_invariants_2}, \cite[\S5]{Sl83}. The quotient
variety $X$ can be embedded as a surface $X \subset \mathbb{C}^3$
with an isolated singularity at the origin.

 \index{theorem! - Sheppard-Todd-Chevalley-Serre}
\begin{remark}
\label{alg_invariants_2} {\rm
    According to Theorem \ref{STCS_theorem}
   (Sheppard-Todd-Chevalley-Serre) the algebra of invariants
   $k[V]^G$ is isomorphic to a polynomial algebra in some number
   of variables if the image of $G$ in $GL(V)$
   is generated by reflections.
   Every binary polyhedral group $G$ is generated by reflections,
   see, e.g., \cite{CoxM84}, therefore the algebra of invariants $k[V]^G$
   for the binary polyhedral group is a polynomial algebra.
   }
\end{remark}

\begin{example}
\label{alg_invariants_3} {\rm
    Consider the cyclic group $G = \mathbb{Z}/r\mathbb{Z}$ of the order $r$.
    The group $G$
    acts on $\mathbb{C}[z_1, z_2]$ as follows:
    \begin{equation}
       (z_1, z_2)  \mapsto (\varepsilon{z_1},
       \varepsilon^{r-1}z_2),
    \end{equation}
   where $\varepsilon = e^{2\pi{i}/r}$,
   and polynomials
    \begin{equation}
        x = z_1z_2, \hspace{5mm} y = -z_1^r, \hspace{5mm} z = z_2^r
    \end{equation}
    are invariant polynomials in $\mathbb{C}[x, y, z]$ which
    satisfy the following relation
    \begin{equation}
       x^r + yz = 0,
    \end{equation}
    see Table \ref{klein_rel}. We have
    \begin{equation*}
       k[V]^G = \mathbb{C}[z_1z_2, z_1^r, z_2^r]
                \simeq \mathbb{C}[x, y, z]/(x^r + yz) .
    \end{equation*}
    see, e.g., \cite[pp.95-97]{Sp77}, \cite[p.143]{PV94}.
   } \vspace{5mm}
\end{example}

\index{exceptional divisor}
 Du Val obtained the following
description of the minimal resolution $\pi: \tilde{X}
\longrightarrow X$ of a Kleinian singularity $X = \mathbb{C}^2/G$,
see \cite{DuVal34}, \cite[\S6.1, \S6.2]{Sl80}
\cite[\S5]{Sl83}\footnote{For more details, see also \cite{Gb02},
\cite{Rie02}, \cite{Hob02}, \cite{Cr01}.}. The {\it exceptional
divisor} (the preimage of the singular point $O$) is a finite
union of complex projective lines: \index{complex projective line}
$$
   \pi^{-1}(O) = L_1 \cup \dots \cup L_n,
   \hspace{3mm} L_i \simeq \mathbb{CP}^1 \text{ for } i =
   1,\dots,n.
$$
The intersection $L_i \cap L_j$ is empty or consists of exactly
one point for $i \neq j$.

 To each complex projective line $L_i$ (which can be identified with the sphere $S^2 \subset
\mathbb{R}^3$, see \S\ref{sect_projective_line}) we assign a
vertex $i$, and two vertices are connected by an edge if the
corresponding projective lines intersect. The corresponding
diagrams are Dynkin diagrams (this phenomenon was observed by Du
Val in \cite{DuVal34} ), see \text{Table \ref{klein_rel}}.

\begin{figure}[h]
\centering
\includegraphics{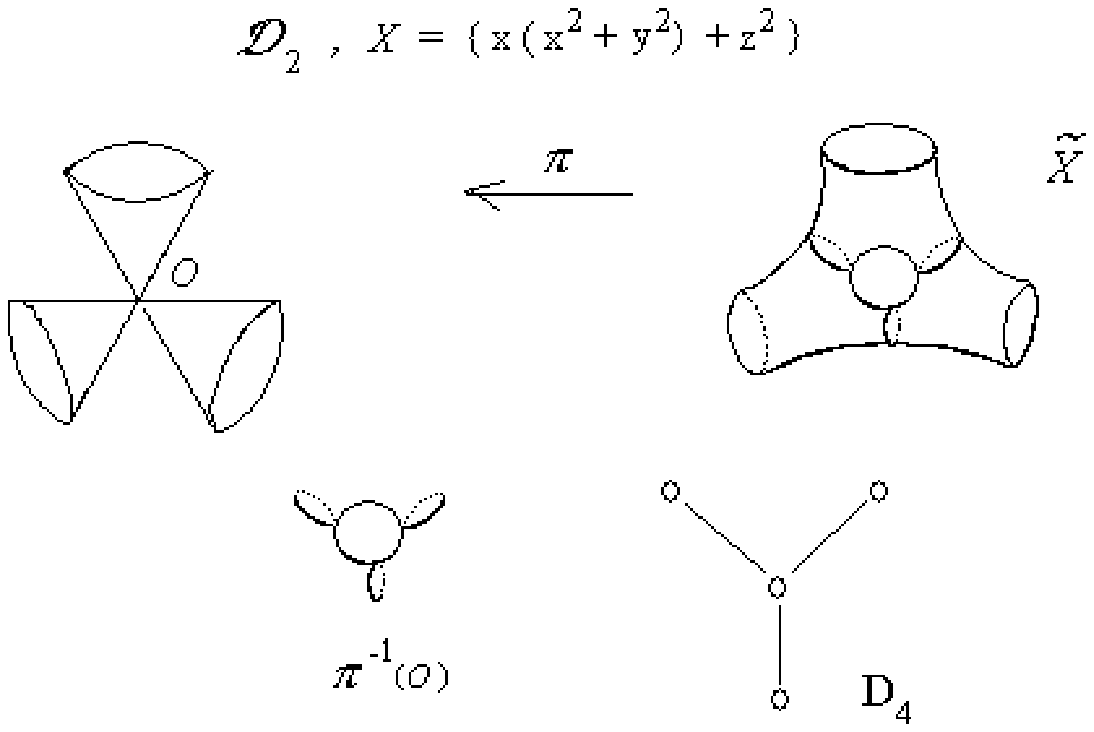}
\caption{\hspace{3mm} The minimal resolution $\pi^{-1}: \tilde{X}
\longrightarrow X$
  for $X = \mathcal{D}_2$}
\label{D2_singularity}
\end{figure}

In the case of the binary dihedral group $\mathcal{D}_2$ the real
resolution of the real variety
$$
   \mathbb{C}^3/R(x,y,z) \cap \mathbb{R}^3
$$
gives a quite faithful picture of the complex situation, the
minimal resolution $\pi^{-1}: \tilde{X} \longrightarrow X$ for $X
= \mathcal{D}_2$ depicted on Fig.~\ref{D2_singularity}. Here
$\pi^{-1}(O)$ consists of four circles, the corresponding diagram
is the Dynkin diagram $D_4$.

\section{The McKay correspondence}
 \label{McKay}
\index{McKay correspondence} Let $G$ be a finite subgroup of
$SU(2)$. Let $\{\rho_0,\rho_1,...\rho_n\}$ be the set of
irreducible finite dimensional complex representations of $G$, of
which $\rho_0$ is the trivial one. Let $\rho: \longrightarrow
SU(2)$ be a faithful
 representation. Then, for each group $G$, we can define
 a matrix $A(G) = (a_{ij})$, by decomposing the tensor products:
\begin{equation}
 \label{main_McKay}
 \rho \otimes \rho_j = \bigoplus\limits_{k=0}^r a_{jk}\rho_k,
  \hspace{7mm} j,k = 0,1,...,r,
\end{equation}
where $a_{jk}$ is the multiplicity of $\rho_k$ in $\rho \otimes
\rho_j$. McKay \cite{McK80} observed that
\begin{center}
{\it The matrix $2I - A(G)$ is the Cartan matrix of the extended
  Dynkin diagram $\tilde\varDelta(G)$ associated to $G$.
  There is one-to-one correspondence between finite subgroups
  of $SU(2)$ and simply-laced extended Dynkin diagrams.
 }\end{center}

This remarkable observation, called the {\it McKay
correspondence}, was based first on an explicit verification
\cite{McK80}.

For the multiply-laced case, the McKay correspondence was extended
by D.~Happel, U.~Preiser, and C.~M.~Ringel in \cite{HPR80}, and by
P.~Slodowy in \cite[App.III]{Sl80}. We consider P.~Slodowy's
approach in \S\ref{slodowy}.

The systematic proof of the McKay correspondence based on the study
of affine Coxeter transformations was given by R.~Steinberg in \cite{Stb85}.

Other explanations of the McKay correspondence have been given by
G.~Gonzalez-Sprinberg and J.-L.~Verdier in \cite{GV83}, by
H.~Kn\"{o}rrer in \cite{Kn85}. A nice review is given by J.~van
Hoboken in \cite{Hob02}.

 \index{McKay operator}
 \index{orbit structure of the Coxeter transformation}
 \index{Lie algebra}
B.~Kostant used the {\it McKay matrix} $A(G)$ (or {\it McKay
operator}) in \cite{Kos84} and showed that the multiplicities
$m_i(n)$ in the decomposition
\begin{equation*}
  \pi_n|G = \sum\limits_{i=0}^r{m_i(n)}\rho_i,
\end{equation*}
(see \S\ref{multipl_kostant}) come in a amazing way from the orbit
structure of the Coxeter transformation on the highest root of the
corresponding Lie algebra $\mathfrak{g}$, see \S\ref{orbit_str}.
More explicitly, the multiplicities $m_i(n)$ come from the orbit
structure of the affine Coxeter transformation on the set of roots
of $\mathfrak{g}$. To calculate these multiplicities, Kostant
employed generating functions and Poincar\'{e} series. We applied
Kostant's technique in \S\ref{kostant} in order to show a relation
between Poincar\'{e} series and a ratio of characteristic
polynomials of the Coxeter transformations.

\section{The Slodowy generalization of the McKay correspondence}
 \label{slodowy}
We consider here the Slodowy generalization \cite{Sl80} of the
McKay correspondence to the multiply-laced case and illustrate
Slodowy's approach with the diagrams $\tilde{F}_{41}$ and
$\tilde{F}_{42}$, see Fig.~\ref{euclidean_diag}.

\index{induced representation}
\index{restricted representation}

Slodowy's approach is based on the consideration of {\it
restricted representations} and {\it induced representations}
instead of an original representation. Let $\rho: G
\longrightarrow GL(V)$ be a representation of a group $G$. We
denote the {\it restricted representation} of $\rho$ to a subgroup
$H \subset G$ by $\rho\downarrow^G_H$, or, briefly,
$\rho^\downarrow$ for fixed $G$ and $H$. Let $\tau: H
\longrightarrow GL(V)$ be a representation of a subgroup $H$. We
denote by $\tau\uparrow^G_H$ the representation {\it induced} by
$\tau$ to a representation of the group $G$ containing $H$; we
briefly write $\tau^\uparrow$ for fixed $G$ and $H$. For a
detailed definition on restricted and induced representations,
see, for example, \cite{Kar92} or \cite{Bak04}.

Let us consider pairs of groups $H \triangleleft G$, where $H$ and $G$ are
binary polyhedral groups from Tables \ref{binary_pol_1}
and \ref{klein_rel}.

\begin{table} [h]
  \centering
  \vspace{2mm}
  \caption{\hspace{3mm}The pairs $H \triangleleft G$ of binary polyhedral groups}
  \renewcommand{\arraystretch}{1.5} 
 \begin{tabular} {||c|c|c|c|c||}
 \hline \hline
   \quad Subgroup \quad &
   \quad Dynkin  \quad  &
   \quad Group \quad    &
   \quad Dynkin \quad   &
   \quad Index \quad    \cr
   \quad $H$  \quad     &
   \quad diagram $\varGamma(H)$ \quad &
   \quad $R$  \quad &
   \quad diagram $\varGamma(G)$ \quad &
   \quad $[G:H]$ \quad \\
 \hline \hline
      & & & & \cr
      $\mathcal{D}_2$ & ${D}_4$ & $\mathcal{T}$ & ${E}_6$ & 3 \cr
      & & & & \\
  \hline
      & & & & \cr
      $\mathcal{T}$  & ${E}_6$ & $\mathcal{O}$ & ${E}_7$ & 2 \cr
      & & & & \\
  \hline
      & & & & \cr
      $\mathcal{D}_{r-1}$ & ${D}_{n+1}$ & $ \mathcal{D}_{2(r-1)}$ &
      ${D}_{2n}$ & 2 \cr
      & & & & \\
 \hline
      & & & & \cr
      $\mathbb{Z}/2r\mathbb{Z}$ & ${A}_{n-1}$ & $\mathcal{D}_r$ &
      ${D}_{r+2}$ & 2 \cr
      & & & & \\
 \hline \hline
 \end{tabular}
  \label{group_pairs}
\end{table}

\subsection{The Slodowy correspondence}
\label{section_slodowy} Let us fix a pair $H \triangleleft G$ from
Table \ref{group_pairs}. We formulate now the essence of the
Slodowy correspondence \cite[App.III]{Sl80}.

 1) Let $\rho_i$, where $i = 1,\dots,n$, be
irreducible representations of $G$; let $\rho^\downarrow_i$ be the
corresponding restricted representations of the subgroup $H$. Let
$\rho$ be a faithful representation of $H$, which may be
considered as the restriction of the fixed faithful representation
$\rho_{f}$ of $G$. Then the following decomposition formula makes
sense
\begin{equation}
 \label{lab_slodowy_1}
     \rho \otimes \rho^\downarrow_i =
          \bigoplus a_{ji} \rho^\downarrow_j
\end{equation}
and uniquely determines an $n\times{n}$ matrix $\tilde{A} =
(a_{ij})$ such that
\begin{equation}
   K = 2I - \tilde{A}
\end{equation}
(see \cite[p.163]{Sl80}), where $K$ is the Cartan matrix of the
corresponding folded extended Dynkin diagram given in Table
\ref{pairs_and_folded_diagr}

2) Let $\tau_i$, where $i = 1,\dots,n$, be irreducible
representations of the subgroup $H$, let $\tau^\uparrow_i$ be the
induced representations of the group $G$. Then the following
decomposition formula makes sense
\begin{equation}
 \label{lab_slodowy_2}
     \rho \otimes \tau^\uparrow_i =
          \bigoplus a_{ij} \tau^\uparrow_j,
\end{equation}
i.e., the decomposition of the induced representation is described
by the matrix $A^\vee = A^t$ which satisfies the relation
\begin{equation}
   K^\vee = 2I - \tilde{A}^\vee
\end{equation}
(see \cite[p.164]{Sl80}), where $K^\vee$ is the Cartan matrix of
the dual folded extended Dynkin diagram given in Table
\ref{pairs_and_folded_diagr}.

\index{Slodowy matrix}
\index{Slodowy correspondence}
\index{McKay-Slodowy correspondence}
We call matrices $\tilde{A}$ and $\tilde{A}^\vee$ the {\it Slodowy matrices},
they are analogs of the McKay matrix. The {\it Slodowy correspondence}
is an analogue to the McKay correspondence for the multiply-laced case, so one
can speak about the {\it McKay-Slodowy correspondence}.

\begin{table} [h]
  \centering
  \renewcommand{\arraystretch}{1.5} 
  \vspace{2mm}
  \caption{\hspace{3mm}The pairs $H \triangleleft G$ and folded extended Dynkin
diagrams}
 \begin{tabular} {||c|c|c||}
 \hline \hline
   \quad Groups \quad  &
   \quad Dynkin diagram \quad  &
   \quad Folded extended  \quad \cr
   \quad $H \triangleleft G$ \quad &
   \quad $\tilde\varGamma(H)$ and $\tilde\varGamma(G)$ \quad &
   \quad Dynkin diagram \quad \\
 \hline \hline
      & & \cr
      $\mathcal{D}_2 \triangleleft \mathcal{T}$ &
      ${D}_4$ and ${E}_6$ &
      $\tilde{G}_{21}$ and $\tilde{G}_{22}$ \cr
      & & \\
 \hline
      & & \cr
      $\mathcal{T} \triangleleft \mathcal{O}$ &
      ${E}_6$ and ${E}_7$ &
      $\tilde{F}_{41}$ and $\tilde{G}_{42}$ \cr
      & & \\
 \hline
      & & \cr
      $\mathcal{D}_{r-1} \triangleleft \mathcal{D}_{2(r-1)}$ &
      ${D}_{n+1}$ and ${D}_{2n}$ &
      $\widetilde{DD}_{n}$ and $\widetilde{CD}_{n}$  \cr
      & & \\
 \hline
      & & \cr
      $\mathbb{Z}/2r\mathbb{Z} \triangleleft \mathcal{D}_r$ &
      ${A}_{n-1}$ and ${D}_{r+2}$ &
      $\tilde{B}_{n}$ and $\tilde{C}_{n}$  \cr
      & & \\
 \hline \hline
 \end{tabular}
  \label{pairs_and_folded_diagr}
\end{table}

\subsection{The binary tetrahedral group and the binary octahedral group}
Now we will illustrate the Slodowy correspondence for the
binary tetrahedral and octahedral groups, i.e., $H$ is
the binary tetrahedral group $\mathcal{T}$ and
$G$ is the binary octahedral group $\mathcal{O}$,
$\mathcal{T} \triangleleft \mathcal{O}$. These groups have orders
$|\mathcal{T}| = 24$ and $|\mathcal{O}| = 48$, see Table \ref{binary_pol_1}.

\index{formula! - Springer}
We will use the Springer formula for elements of the group $\mathcal{O}$
from \cite[\S4.4.11]{Sp77}. Let
\begin{equation}
 a= \left (
  \begin{array}{cc}
     \varepsilon & 0 \\
     0 & \varepsilon^{-1} \\
    \end{array}
 \right ), \hspace{3mm}
 b = \left (
  \begin{array}{cc}
     0 & i \\
     i & 0
    \end{array}
 \right ),  \hspace{3mm}
 c = \frac{1}{\sqrt{2}}\left (
  \begin{array}{cc}
     \varepsilon^{-1} & \varepsilon^{-1} \\
     -\varepsilon & \varepsilon
    \end{array}
 \right ),
\end{equation}
where $\varepsilon = e^{\pi{i}/4}$. Then each of the $48$
different elements $x \in \mathcal{O}$ may be expressed as
follows:
\begin{equation}
 \label{springer_formula}
   x = a^h{b^j}c^l, \hspace{3mm}
   0 \leq h < 8,  \hspace{3mm}
   0 \leq j < 2,  \hspace{3mm}
   0 \leq l < 3.
\end{equation}

\begin{table} [h]
  \centering
  \vspace{2mm}
  \caption{\hspace{3mm}The elements of the binary octahedral group}
  \renewcommand{\arraystretch}{1.5} 
 \begin{tabular} {|c|c|c|}
 \hline \hline
Elements      & Matrix  & Trace    \cr
   $x = a^p{b^j}c^l$ & form   & $tr(x)$   \cr
     $0 \leq p < 8$   &   $u$   &       \\
 \hline \hline
   & & \cr
   $a^p$ &
      $\left (
     \begin{array}{cc}
       \varepsilon^p & 0 \\
       0 & \varepsilon^{-p} \\
       \end{array}
    \right )$ &
     $\displaystyle 2\cos\frac{\pi{p}}{4}$ \cr
     & & \\
  \hline
   & & \cr
     $a^p{b}$ &
      $\left (
    \begin{array}{cc}
       0 & i\varepsilon^p \\
       i\varepsilon^{-p} & 0 \\
      \end{array}
    \right )$ & 0 \cr
    & & \\
  \hline
   & & \cr
     $a^p{c}$ &
      $\displaystyle\frac{1}{\sqrt{2}}\left (
    \begin{array}{cc}
       \varepsilon^{p - 1} & \varepsilon^{p - 1} \\
       -\varepsilon^{-(p-1)} & \varepsilon^{-(p-1)} \\
       \end{array}
    \right )$ =
      $\displaystyle\frac{1}{\sqrt{2}}\left (
    \begin{array}{cc}
       \varepsilon^{p - 1} & \varepsilon^{p - 1} \\
       \varepsilon^{- p + 5} & \varepsilon^{- p + 1} \\
       \end{array}
    \right )$ &
   $\displaystyle \sqrt{2}\cos\frac{\pi{(p-1)}}{4}$ \cr
    & & \\
  \hline
   & & \cr
     $a^p{bc}$ &
      $\displaystyle\frac{1}{\sqrt{2}}\left (
    \begin{array}{cc}
       -i\varepsilon^{p + 1} & i\varepsilon^{p + 1} \\
       i\varepsilon^{-(p + 1)} & i\varepsilon^{-(p + 1)} \\
       \end{array}
    \right )$  =
    $\displaystyle\frac{1}{\sqrt{2}}\left (
    \begin{array}{cc}
       \varepsilon^{p - 1} & \varepsilon^{p + 3} \\
       \varepsilon^{-p + 1} & \varepsilon^{-p + 1} \\
       \end{array}
    \right )$ &
    $\displaystyle \sqrt{2}\cos\frac{\pi{(p - 1)}}{4}$ \cr
    & & \\
  \hline
   & & \cr
     $a^p{c^2}$ &
      $\displaystyle\frac{1}{\sqrt{2}}\left (
    \begin{array}{cc}
       -\varepsilon^{p + 1} & \varepsilon^{p - 1} \\
       -\varepsilon^{-(p - 1)} & -\varepsilon^{-(p + 1)} \\
       \end{array}
    \right )$ =
      $\displaystyle\frac{1}{\sqrt{2}}\left (
    \begin{array}{cc}
       \varepsilon^{p + 5} & \varepsilon^{p - 1} \\
       \varepsilon^{-p + 5} & \varepsilon^{-p + 3} \\
       \end{array}
    \right )$ &
   $\displaystyle -\sqrt{2}\cos\frac{\pi{(p + 1)}}{4}$ \cr
    & & \\
  \hline
   & & \cr
     $a^p{bc^2}$ &
      $\displaystyle\frac{1}{\sqrt{2}}\left (
    \begin{array}{cc}
       -i\varepsilon^{p + 1} & -i\varepsilon^{p - 1} \\
       -i\varepsilon^{-(p - 1)} & i\varepsilon^{-(p + 1)} \\
       \end{array}
    \right )$ =
    $\displaystyle\frac{1}{\sqrt{2}}\left (
    \begin{array}{cc}
       \varepsilon^{p - 1} & \varepsilon^{p - 3} \\
       \varepsilon^{-p - 1} & \varepsilon^{-p + 1} \\
       \end{array}
    \right )$ &
    $\displaystyle \sqrt{2}\cos\frac{\pi{(p - 1)}}{4}$ \cr
     & & \\
 \hline \hline
 \end{tabular}
  \label{elements_bin_octahedral}
\end{table}

The elements $x \in \mathcal{O}$ and their traces are collected
in Table \ref{elements_bin_octahedral}. Observe that
every element $u\in SL(2, \mathbb{C})$ from Table \ref{elements_bin_octahedral}
is of the form
\begin{equation}
     u = \left (
     \begin{array}{cc}
     \alpha & \beta  \\
     -\beta^* & \alpha^*
     \end{array}
     \right ),
\end{equation}
see (\ref{su2_matr}).
We can now
distinguish the elements of $\mathcal{O}$ by their traces and by
means of the 1-dimensional representation $\rho_1$ such that
\begin{equation}
 \rho_1(a) = -1, \hspace{3mm} \rho_1(b) = 1, \hspace{3mm} \rho_1(c) = 1.
\end{equation}
\begin{proposition}
  \label{8_conj_classes}
There are 8 conjugacy classes in the binary octahedral group $\mathcal{O}$.
Rows of Table \ref{separated_by_trace} constitute these conjugacy classes.
\end{proposition}
\PerfProof In order to prove the proposition we need a number of conjugacy
relations.

1) It is easy to verify that
\begin{equation}
 \label{b_swap}
 b \left(
  \begin{array}{cc}
     x & y \\
     u & v \\
    \end{array} \right ) b^{-1} =
  \left(
  \begin{array}{cc}
     v & u \\
     y & x \\
    \end{array} \right )
    \text{ for all } x,y,u,v \in \mathbb{C},
\end{equation}
and
\begin{equation}
 \label{a_conujating}
 a \left(
  \begin{array}{cc}
     x & y \\
     u & v \\
    \end{array} \right ) a^{-1} =
  \left(
  \begin{array}{cc}
     x & \varepsilon^2{y} \\
     \varepsilon^2{u} & v \\
    \end{array} \right )
    \text{ for all } x,y,u,v \in \mathbb{C}.
\end{equation}

\begin{table} [h]
  \centering
  \renewcommand{\arraystretch}{1.3} 
  \vspace{2mm}
  \caption{\hspace{3mm}The conjugacy classes in the binary octahedral group}
 \begin{tabular} {||c|c|c|c|c||}
 \hline \hline
   Trace      & Representation & Conjugacy class & Class & Representative \cr
   $tr(x)$     &  $\rho_1$   & $Cl(g)$ 
     & order & $g$ \\
 \hline \hline
   & & & & \cr
   $2$ & 1 &
      $\left (
     \begin{array}{cc}
       1 & 0 \\
       0 & 1 \\
       \end{array}
    \right )$ & 1 & 1 \cr
   & & & & \\
 \hline
   & & & & \cr
   $-2$ &   1 &
      $b^2 = c^3 = a^4 = \left (
     \begin{array}{cc}
       -1 & 0 \\
       0 & -1 \\
       \end{array}
    \right )$ & 1 & $-1$ \cr
   & & & & \\
 \hline
   & & & & \cr
   $0$ & $-1$ & $ab, a^3b, a^5b, a^7b$, & & \cr
         &  & $a^3bc, a^7bc, ac^2, a^5c^2,$ & 12 & $ab$ \cr
       &  & $a^3c, a^7c, a^3bc^2, a^7bc^2$ & & \cr
   & & & & \\
 \hline
   & & & & \cr
     $0$ & 1 & $a^2, a^6, b, a^2b, a^4b, a^6b$ & 6 & $b$ \cr
   & & & & \\
 \hline
   & & & & \cr
     $-1$ & 1 & $a^4c, a^6c, a^4bc, a^6bc,$ & & \cr
       &  & $c^2, a^6c^2, a^4bc^2, a^6bc^2$ & 8 & $c^2$ \cr
   & & & & \\
 \hline
   & & & & \cr
     $1$ & 1 & $c, a^2c, bc, a^2bc,$ & & \cr
       &  & $a^4c^2, a^2c^2, bc^2, a^2bc^2$ & 8 & $c$ \cr
   & & & & \\
 \hline
   & & & & \cr
     $\sqrt{2}$ & $-1$ & $a, a^7, ac, abc, a^3c^2, abc^2$ & 6 & $a$ \\
 \hline
   & & & & \cr
     $-\sqrt{2}$ & $-1$ & $a^3, a^5, a^5c, a^5bc, a^7c^2, a^5bc^2$ & 6 &
$a^3$ \\
 \hline \hline
 \end{tabular}
  \label{separated_by_trace}
\end{table}

\footnotetext[1]{See Remark \ref{denote_conjugacy_cl}.}

From Table \ref{elements_bin_octahedral} and (\ref{b_swap}) and since
$a^p = a^{p+8}$, we see, for all $p \in \mathbb{Z}$, that
\begin{equation}
 \label{conj_1}
  b(a^p)b^{-1} = a^{-p},
\end{equation}
\begin{equation}
 \label{conj_2}
  b(a^{p}c)b^{-1} = a^{2-p}bc,
\end{equation}
\begin{equation}
 \label{conj_3}
  b(a^{p}bc^2)b^{-1} = a^{2-p}bc^2,
\end{equation}
\begin{equation}
 \label{conj_4}
  b(a^{p}c^2)b^{-1} = a^{6-p}c^2.
\end{equation}
From Table \ref{elements_bin_octahedral} and (\ref{a_conujating}), we see, for
all $p \in \mathbb{Z}$, that
\begin{equation}
 \label{conj_5}
  a(a^{p}c^2)a^{-1} = a^{p-2}bc.
\end{equation}
By (\ref{conj_1})
\begin{equation}
 \label{conj_11}
  a \sim a^7, \hspace{3mm} a^2 \sim a^6, \hspace{3mm} a^3 \sim a^5.
\end{equation}
By (\ref{conj_2})
\begin{equation}
\begin{split}
 \label{conj_21}
  & c \sim a^2bc, \hspace{3mm} ac \sim abc, \hspace{3mm}
  a^2{c} \sim bc, \hspace{3mm} a^3{c} \sim a^{7}bc, \hspace{3mm} \\
  & a^4{c} \sim a^6{bc}, \hspace{3mm} a^5{c} \sim a^5{bc}, \hspace{3mm}
   a^6{c} \sim a^4{bc}, \hspace{3mm} a^7{c} \sim a^3{bc}.
\end{split}
\end{equation}
By (\ref{conj_3})
\begin{equation}
 \label{conj_31}
  bc^2 \sim a^2bc^2, \hspace{3mm} a^{3}bc^2 \sim a^7bc^2, \hspace{3mm}
  a^{4}bc^2 \sim a^{6}bc^2.
\end{equation}
By (\ref{conj_4})
\begin{equation}
 \label{conj_41}
  c^2 \sim a^{6}c^2, \hspace{3mm} a{c}^2 \sim a^{5}c^2, \hspace{3mm}
  a^{2}c^2 \sim a^{4}c^2.
\end{equation}
By (\ref{conj_5})
\begin{equation}
\begin{split}
 \label{conj_51}
  c^2 \sim a^{6}bc, \hspace{3mm} a{c}^2 \sim a^{7}bc, \hspace{3mm}
  a^{2}c^2 \sim bc, \hspace{3mm} a^3{c}^2 \sim abc, \\
  a^4c^2 \sim a^{2}bc, \hspace{3mm} a^5{c}^2 \sim a^{3}bc, \hspace{3mm}
  a^{6}c^2 \sim a^4{bc}, \hspace{3mm} a^7{c}^2 \sim a^5{bc}.
\end{split}
\end{equation}
From (\ref{conj_1}) we have
\begin{equation}
 \label{comm_1}
  ba^{-1} = ab, \hspace{5mm} ba^{-p} = a^{p}b.
\end{equation}
By (\ref{conj_1}) $a^p{b} \sim aa^{p}ba^{-1} = aa^pa{b}$, i.e.,
\begin{equation}
 \label{comm_2}
  a^p{b} \sim a^{p+2}{b},
\end{equation}
i.e.,
\begin{equation}
\begin{split}
 \label{conj_61}
  & b \sim a^2{b} \sim a^4{b} \sim a^6{b}, \\
  & ab \sim a^3{b} \sim a^5{b} \sim a^7{b}.
\end{split}
\end{equation}

2) We need some conjugacy relations performed by $c$.
By (\ref{sl2_matr}) the inverse elements in $SL(2, \mathbb{C})$
are
\begin{equation}
 \label{su2_reverse_matr}
  u^{-1} =   \left (
    \begin{array}{cc}
       d  & -b \\
         -c & a \\
    \end{array}
     \right ) \text { for every }
  u = \left (
    \begin{array}{cc}
       a & b \\
         c & d \\
    \end{array}
     \right ).
\end{equation}
So, we have
\begin{equation}
 \label{c_reverse}
  c = \frac{1}{\sqrt{2}}\left (
      \begin{array}{cc}
       \varepsilon^{-1} & \varepsilon^{-1} \\
         -\varepsilon & \varepsilon
      \end{array} \right ) \text{ and }
  c^{-1} = \frac{1}{\sqrt{2}}\left (
      \begin{array}{cc}
       \varepsilon & -\varepsilon^{-1} \\
         \varepsilon & \varepsilon^{-1}
      \end{array} \right ).
\end{equation}
Thus
\begin{equation}
 \label{c_ap}
  ca^{p}c^{-1} = \frac{1}{2}\left (
      \begin{array}{cc}
       \varepsilon^p + \varepsilon^{-p} & -\varepsilon^{p-2} +
\varepsilon^{-p-2} \\
         -\varepsilon^{p+2} + \varepsilon^{-p+2} & \varepsilon^p +
\varepsilon^{-p}
      \end{array} \right ).
\end{equation}
Since
\begin{equation}
    \varepsilon^2 = i, \hspace{5mm}
    \varepsilon^2 = -i, \hspace{5mm} \varepsilon^4 = \varepsilon^{-4} = -1,
\end{equation}
we obtain from Table \ref{elements_bin_octahedral} and
(\ref{c_ap}):
\begin{equation}
 \label{c_a2}
  ca^{2}c^{-1} = \left (
      \begin{array}{cc}
       0 & -1 \\
         1 & 0
      \end{array} \right ) = a^2b.
\end{equation}
Also, by Table \ref{elements_bin_octahedral} and (\ref{c_ap}) we
have
\begin{equation}
 \label{c_a_1}
  ca^{-1}c^{-1} = \frac{1}{\sqrt{2}}\left (
      \begin{array}{cc}
       1 & 1 \\
         -1 & 1
      \end{array} \right ) = ac, \text{ and } ac^2 = ca^{-1}.
\end{equation}
By (\ref{c_a2}) and (\ref{c_a_1}) we have
\begin{equation}
\begin{split}
 \label{c_a_}
  cac^{-1} = ca^2c^{-1}ca^{-1}c^{-1} & = (a^2b)(ac) = a^2a^{-1}bc = abc, \text{
i.e.,}  \\
  & a \sim abc.
\end{split}
\end{equation}
Further,
\begin{equation}
 \label{c_conj}
  cbc^{-1} = \left (
      \begin{array}{cc}
       \varepsilon^2 & 0    \\
         0 & \varepsilon^{-2}
      \end{array} \right ) = a^2.
\end{equation}
By (\ref{c_a_}) and (\ref{c_conj}) we have
\begin{equation}
\begin{split}
 \label{c_ab}
  cabc^{-1} = cac^{-1}cbc^{-1} & = (abc)(a^2) \sim a^2(abc) = a^3bc,    \\
    & ab \sim a^3bc.
\end{split}
\end{equation}
By (\ref{c_conj})
\begin{equation}
 \label{conj_7}
  cbc^{-1} = a^2, \hspace{3mm} cb^2c^{-1} = a^4,
\end{equation}
so, by (\ref{c_a2}) and (\ref{comm_1}) we have
$ca^6c^{-1} = a^2ba^2ba^2b = a^2a^2b^{-1}ba^2b = a^{6}b$,
i.e.,
\begin{equation}
 \label{conj_71}
  b \sim a^2, \hspace{3mm} b^2 \sim a^4,
\end{equation}
and
\begin{equation}
 \label{conj_72}
  a^2 \sim a^2{b}, \hspace{3mm} a^6 \sim a^{6}b.
\end{equation}
In addition, by (\ref{comm_1}) and (\ref{c_a_1}) we have
\begin{equation}
  a^{p}bc^2 = a^{p+1}a^{-1}bc^2 = a^{p+1}bac^2 = a^{p+1}bca^{-1} \sim
  a^{-1}a^{p+1}bca^{-1}a = a^{p}bc,
\end{equation}
i.e.,
\begin{equation}
  a^{p}bc^2 \sim a^{p}bc,
\end{equation}
and
\begin{equation}
\label{conj_8}
\begin{split}
 & bc \sim bc^2, \hspace{3mm} abc \sim abc^2, \hspace{3mm}
  a^2bc \sim a^2bc^2, \hspace{3mm} a^3bc \sim a^3bc^2,\\
 & a^4bc \sim a^4bc^2, \hspace{3mm} a^5bc \sim a^5bc^2, \hspace{3mm}
  a^6bc \sim a^6bc^2, \hspace{3mm} a^7bc \sim a^7bc^2.
\end{split}
\end{equation}

\begin{remark}
 \label{denote_conjugacy_cl}
 {\rm
 We denote by $Cl(g)$ the conjugacy class containing the element $g \in \mathcal{O}$.
 The union $Cl(b) \cup \{1, -1\}$ constitutes the 8-element subgroup
\begin{equation}
 \label{8_elem_gr}
  \{1, a^2, a^4, a^6, b, a^2b, a^4b, a^6b \} =
  \{1, a^2, -1, -a^2, b, a^2b, -b, -a^2b \}.
\end{equation}
Setting $i = b$, $j = a^2$, $k = a^2b$ we see that
\begin{equation}
\begin{split}
 & i^2 = j^2 = k^2 = -1, \\
 & ji = -ij = k, \hspace{3mm} ik = -ki = j, \hspace{3mm} kj = -jk = i,
\end{split}
\end{equation}
i.e., group (\ref{8_elem_gr}) is the quaternion group $Q_8$, see
(\ref{quaternion}).

According to (\ref{comm_1}), (\ref{c_a2}) and (\ref{c_conj})
we deduce that $Q_8$ is a normal subgroup in
$\mathcal{O}$:
\begin{equation}
 \label{Q_8_normal}
  Q_8 \triangleleft \mathcal{O} .
\end{equation}
 }
\end{remark}

3) Now we can prove conjugacy of elements in every row of Table
\ref{separated_by_trace}.

\underline{Conjugacy class $Cl(ab)$.}
By (\ref{conj_61}) and (\ref{c_ab})
$$
  ab \sim a^3b \sim a^5b \sim a^7b \sim a^3bc.
$$
By (\ref{conj_51}), (\ref{conj_21}) and (\ref{conj_8})
$$
  a^3bc \sim a^5c^2, \hspace{3mm} a^3bc \sim a^7c, \hspace{3mm} a^3bc \sim
a^3bc^2.
$$
By (\ref{conj_31}), (\ref{conj_8}), (\ref{conj_21}) and (\ref{conj_51})
$$
  a^3bc^2 \sim a^7bc^2 \sim a^7bc \sim a^3c, \hspace{3mm} a^7bc \sim ac^2.
$$
Therefore, the elements
\begin{equation*}
  ab, a^3b, a^5b, a^7b, a^3bc, a^7bc, ac^2, a^5c^2,
  a^3c, a^7c, a^3bc^2, a^7bc^2
\end{equation*}
constitute a conjugacy class.

\underline{Conjugacy class $Cl(b)$.}
By (\ref{conj_61}), (\ref{conj_71}) and (\ref{conj_11})
$$
  b \sim a^2b \sim a^4b \sim a^6b,
   \hspace{5mm} b \sim a^2, \hspace{5mm} a^2 \sim a^6,
$$
i.e., the elements
\begin{equation*}
  a^2, a^6, b, a^2b, a^4b, a^6b
\end{equation*}
constitute a conjugacy class.

\underline{Conjugacy class $Cl(c^2)$.}
By (\ref{conj_51}), (\ref{conj_41}), (\ref{conj_21}), (\ref{conj_51}) and
(\ref{conj_8})
$$
  c^2 \sim a^6bc \sim a^4c, \hspace{3mm} c^2 \sim a^6c^2 \sim a^4bc \sim
a^4bc^2,
   \hspace{3mm} a^6bc \sim a^6bc^2, \hspace{3mm} a^6c \sim a^4bc.
$$
Thus, we have a conjugacy class
\begin{equation*}
  c^2, a^6c^2, a^4bc^2, a^6bc^2, a^4c, a^6c, a^4bc, a^6bc.
\end{equation*}

\underline{Conjugacy class $Cl(c)$.}
By (\ref{conj_21}), (\ref{conj_51}), (\ref{conj_41}) and (\ref{conj_31})
$$
  c \sim a^2bc \sim a^4c^2 \sim a^2c^2,
  \hspace{3mm} a^2bc \sim a^2bc^2 \sim bc^2 \sim bc \sim a^2c,
$$
i.e., the elements
\begin{equation*}
  c, a^2c, bc, a^2bc, a^4c^2, a^2c^2, bc^2, a^2bc^2
\end{equation*}
constitute a conjugacy class.

\underline{Conjugacy class $Cl(a)$.}
By (\ref{conj_11}), (\ref{c_a_1}), (\ref{conj_21}), (\ref{conj_51}) and
(\ref{conj_8})
$$
  a \sim a^7 \sim ac \sim abc \sim a^3c^2, \hspace{5mm} abc \sim abc^2,
$$
and we have a conjugacy class
\begin{equation*}
  a, a^7, ac, abc, a^3c^2, abc^2.
\end{equation*}

\underline{Conjugacy class $Cl(a^3)$.} Since $a^7c^2 \sim
c(a^7c^2)c^{-1} = ca^{-1}c$, we have by (\ref{c_a_1}):
$$
a^7c^2 \sim ca^{-1}c = ac^2c = ac^3.
$$
Since $c^3 = a^4 = -1$, we
deduce that
\begin{equation}
 \label{conj_9}
    a^7c^2 \sim a^5.
\end{equation}
By (\ref{conj_11}), (\ref{conj_9}), (\ref{conj_51}), (\ref{conj_21}) and
(\ref{conj_8})
$$
  a^3 \sim a^5 \sim a^7c^2 \sim a^5bc \sim a^5c, \hspace{5mm} a^5bc \sim
a^5bc^2.
$$
The elements
\begin{equation*}
  a^3, a^5, a^5c, a^5bc, a^7c^2, a^5bc^2.
\end{equation*}
constitute a conjugacy class. The proposition \ref{8_conj_classes}
is proved. \qedsymbol

Now consider conjugacy classes in the binary tetrahedral group
$\mathcal{T}$. According to the Springer formula
\cite[\S4.4.10]{Sp77}, the elements of the group $\mathcal{T}$ are
given by (\ref{springer_formula}) with even numbers $h$. In other
words, each of the $24$ different elements $x \in \mathcal{T}$ may
be given as follows:
\begin{equation}
 \label{springer_formula_2}
   x = a^{2h}{b^j}c^l, \hspace{3mm}
   0 \leq h < 4,  \hspace{3mm}
   0 \leq j < 2,  \hspace{3mm}
   0 \leq l < 3.
\end{equation}
There are 24 elements of type (\ref{springer_formula_2}). In the
octahedral group $\mathcal{O}$, the elements
(\ref{springer_formula_2}) constitute 5 conjugacy classes:
$$
   \{1\},\; \{-1\},\; Cl(b),\; Cl(c) \text{ and } Cl(c^2),
$$
see Table \ref{separated_by_trace}. We will see now, that in the
tetrahedral group $\mathcal{T}$, the elements
(\ref{springer_formula_2}) constitute 7 conjugacy classes:
$$
   \{1\},\; \{-1\},\; Cl(b),\; Cl(c),\; Cl(a^4c)\; \text{ and }
   Cl(a^4c^2).
$$
The elements $x \in \mathcal{T}$ and their traces are collected in
Table \ref{separated_by_trace_tetra}. We can now distinguish the
elements of $\mathcal{T}$ by their traces and by means of two
1-dimensional representations $\tau_1$ and $\tau_2$ such that
\begin{equation}
\begin{split}
 & \tau_1(a) = 1, \hspace{3mm} \tau_1(b) = 1, \hspace{3mm} \tau_1(c) = \omega_3,
\\
 & \tau_2(a) = 1, \hspace{3mm} \tau_2(b) = 1, \hspace{3mm} \tau_2(c) =
\omega_3^2,
\end{split}
\end{equation}
where $\omega_3 = e^{2\pi{i}/3}$.
\begin{table} [h]
  \centering
  \vspace{2mm}
  \caption{\hspace{3mm}The conjugacy classes in the binary tetrahedral group}
  \renewcommand{\arraystretch}{1.5} 
 \begin{tabular} {||c|c|c|c|c|c||}
 \hline \hline
   Trace      & Repr.  & Repr.  & Conjugacy class & Class & Representative \cr
   $tr(x)$     & $\tau_1$ & $\tau_2$ &  $Cl(g)$     & order & $g$ \\
 \hline \hline
   & & & & & \cr
   $2$ & 1 & 1 &
      $\left (
     \begin{array}{cc}
       1 & 0 \\
       0 & 1 \\
       \end{array}
    \right )$ & 1 & 1 \cr
   & & & & & \\
 \hline
   & & & & & \cr
   $-2$ &   1 & 1 &
      $\left (
     \begin{array}{cc}
       -1 & 0 \\
       0 & -1 \\
       \end{array}
    \right )$ & 1 & $-1$ \cr
   & & & & & \\
 \hline
   & & & & & \cr
     $0$ & 1 & 1 & $a^2, a^6, b, a^2b, a^4b, a^6b$ & 6 & $b$ \cr
   & & & & & \\
 \hline
   & & & & & \cr
     $-1$ & $\omega_3$ & $\omega_3^2$ & $a^4c, a^6c, a^4bc, a^6bc,$   & 4 &
$a^4c = -c$ \cr
   & & & & & \\
 \hline
   & & & & & \cr
     $-1$ & $\omega_3^2$ & $\omega_3$ & $c^2, a^6c^2, a^4bc^2, a^6bc^2$ & 4
& $c^2$ \cr
   & & & & & \\
 \hline
   & & & & & \cr
$1$ & $\omega_3$ & $\omega_3^2$ & $c, a^2c, bc, a^2bc$ & 4 & $c$
 \cr
   & & & & & \\
 \hline
   & & & & & \cr
$1$ & $\omega_3^2$ & $\omega_3$ & $a^4c^2, a^2c^2, bc^2, a^2bc^2$
& 4 & $a^4c^2 = -c^2$ \cr
   & & & & & \\
 \hline \hline
 \end{tabular}
  \label{separated_by_trace_tetra}
\end{table}

\begin{proposition}
 \label{7_conj_classes}
There are 7 conjugacy classes in the binary tetrahedral group $\mathcal{T}$.
The rows of Table \ref{separated_by_trace_tetra} constitute these conjugacy
classes.
\end{proposition}
\PerfProof We will prove that the classes
$Cl(b)$, $Cl(a^4c)$, $Cl(c^2)$, $Cl(c)$, and $Cl(a^4c^2)$
constitute all conjugacy classes.

\underline{Conjugacy class $Cl(b)$.}
Observe that relation (\ref{conj_61}) uses conjugating element $a$ which
does not lie in $\mathcal{T}$. We conjugate by $a^{2p}$:
\begin{equation}
 \label{conj_a2}
  b \sim a^2ba^{-2} = a^4b, \hspace{3mm}
  a^2b \sim a^6a^2ba^{-6} = a^6b.
\end{equation}
Then by (\ref{conj_a2}), (\ref{conj_11}), (\ref{c_a2}) and (\ref{conj_71}) we
have
\begin{equation}
 b \sim a^4b, \hspace{3mm} a^2b \sim a^6b, \hspace{3mm}
 a^2 \sim a^6, \hspace{3mm} a^2 \sim a^2b, \hspace{3mm} b \sim a^2,
\end{equation}
i.e., the elements $b$, $a^2$, $a^6$, $a^2b$, $a^4b$, and $a^6b$
constitute a conjugacy class.

\underline{Conjugacy class $Cl(a^4c)$.}
By (\ref{conj_21})
conjugating by $b$ we get
\begin{equation}
 \label{conj_a2_2}
  a^4bc \sim a^6c, \hspace{3mm}
  a^6bc \sim a^4c.
\end{equation}
Further,
\begin{equation}
 \label{conj_a2_c}
 a^2ca^{-2} = \left (
   \begin{array}{cc}
      \varepsilon^{-1} & \varepsilon^3 \\
      -\varepsilon^{-3} & \varepsilon
     \end{array} \right ) =
     \left ( \begin{array}{cc}
      \varepsilon^{-1} & \varepsilon^3 \\
      \varepsilon & \varepsilon
     \end{array} \right ),
\end{equation}
and by Table \ref{elements_bin_octahedral} we have
\begin{equation}
 \label{conj_a2_c_1}
  a^2ca^{-2} = bc, \hspace{3mm} c \sim bc,
\end{equation}
where conjugating element is $a^2$. From (\ref{conj_a2_c_1}) we obtain
\begin{equation}
 \label{conj_a2_c_2}
  a^6bc = a^6(a^2ca^{-2}) = ca^{-2} \sim a^6(ca^{-2})a^{-6} = a^6c.
\end{equation}
 So,
\begin{equation}
 \label{conj_a2_c_3}
  a^6bc \sim a^6c,
\end{equation}
and together with (\ref{conj_a2_2}) we see that the elements
$a^4c$, $a^4bc$, $a^6c$, and $a^6bc$ are conjugate.

\underline{Conjugacy class $Cl(c^2)$.}
By (\ref{conj_41}) and (\ref{conj_31})
\begin{equation}
 \label{conj_c2_class_1}
  c^2 \sim a^6c^2, \hspace{3mm} a^4bc^2 \sim a^6bc^2.
\end{equation}
Further, by (\ref{conj_a2_c_1}), (\ref{c_conj}) and (\ref{comm_1})
\begin{equation}
 \label{conj_c2_class_2}
  a^2c^2a^{-2} = (a^2ca^{-2})^2 = (bc)^2 =
  b(cb)c = b(a^2c)c = a^{-2}bc^2 = a^6bc^2,
\end{equation}
i.e.,
\begin{equation}
 \label{conj_c2_class_3}
  c^2 \sim a^6bc^2,
\end{equation}
and, with (\ref{conj_c2_class_1}), we deduce that
the
elements $c^2$, $a^6c^2$, $a^4bc^2$, and $a^6bc^2$ are conjugate.

\underline{Conjugacy class $Cl(c)$.}
By (\ref{conj_a2_c_1}), (\ref{conj_21})
\begin{equation}
 \label{conj_c_class_1}
  c \sim bc, \hspace{3mm} c \sim a^2bc, \hspace{3mm}
    bc \sim a^2c,
\end{equation}
i.e., $c$, $bc$, $a^2bc$, and $a^2c$ are conjugate.

\underline{Conjugacy class $Cl(a^4c^2)$.}
By (\ref{conj_21}), (\ref{conj_31})
\begin{equation}
 \label{conj_a4_c2_class}
  a^4c^2 \sim a^2c^2, \hspace{3mm} bc^2 \sim a^2bc^2.
\end{equation}
Further, by (\ref{conj_7}) we have
\begin{equation}
 \label{conj_a4_c2_class_2}
  bc^2 \sim c(bc^2)c^{-1} = (cb)c = (a^2c)c = a^2c.
\end{equation}
Thus, the elements $a^4c^2$, $a^2c^2$, $bc^2$, $a^2c$ are
conjugate. The proposition \ref{7_conj_classes} is proved.
\qedsymbol

\subsection{Representations of the binary octahedral and
tetrahedral groups}

\begin{table} [h]
  \centering
  \vspace{2mm}
  \caption{\hspace{3mm}The characters of the binary octahedral group}
  \renewcommand{\arraystretch}{1.5}
 \begin{tabular} {||c|c|c|c|c|c|c|c|c|c||}
 \hline \hline
   Character & \multicolumn{8}{c|}
         {Conjugacy class $Cl(g)$ and its order $|Cl(g)|$ under it} & Note on \\
 \cline{2-9}
   $\psi_i$ & $Cl(1)$ & $Cl(-1)$ & $Cl(ab)$ &
         $Cl(b)$ & $Cl(c^2)$ & $Cl(c)$ & $Cl(a)$ & $Cl(a^3)$ & represent. \cr
        &  1  & 1   & 12   & 6   & 8   & 8   & 6  &  6 & $\rho_i$ \\
 \hline \hline
     $\psi_0$ &  1  & 1   & 1   & 1   & 1   & 1   & 1  &  1 & trivial \\
 \hline
     $\psi_1$ &  1  & 1   & $-1$   & 1   & 1   & 1   & $-1$  & $-1$ &
$\rho_1(a) = -a$ \\
 \hline
     $\psi_2$ &  2  & 2   & 0   & 2   & $-1$   & $-1$  & 0  & 0 &
$\gamma_2\pi_2$ \\
 \hline
     $\psi_3$ &  2  & $-2$  & 0   & 0   & $-1$   & 1
          & $\sqrt{2}$ & $-\sqrt{2}$ & faithful \\
 \hline
     $\psi_4$ &  2  & $-2$  & 0   & 0   & $-1$   & 1
          & $-\sqrt{2}$ & $\sqrt{2}$ & faithful \\
 \hline
     $\psi_5$ &  3  & 3  & $-1$   & $-1$   & 0   & 0   & 1 & 1 &
$\gamma_5\pi_{56}$ \\
 \hline
     $\psi_6$ &  3  & 3  &  1   & $-1$   & 0   & 0   & $-1$ & $-1$ &
$\gamma_6\pi_{56}$ \\
 \hline
     $\psi_7$ &  4  & $-4$  &  0   & 0   & 1   & $-1$   & 0 & 0 &
$\rho_2\otimes\rho_3$ \\
 \hline \hline
 \end{tabular}
  \label{characters_O}
\end{table}

\begin{proposition}
\label{charachters_O_prop}
 The group $\mathcal{O}$ has the following $8$ irreducible
representations.

1) Two 1-dimensional representations:
\begin{equation}
 \label{constr_rho_01}
\begin{split}
 & \rho_0(a) = \rho_0(b) = \rho_0(c) = 1; \\
 & \rho_1(a) = -1, \hspace{3mm} \rho_1(b) = \rho_1(c) = 1,
\end{split}
\end{equation}

2) Two faithful 2-dimensional representations:
\begin{equation}
 \label{constr_rho_34}
\begin{split}
 & \rho_3(a) = a, \hspace{3mm} \rho_3(b) = b, \hspace{3mm} \rho_3(c) = c; \\
 & \rho_4(a) = -a, \hspace{3mm} \rho_4(b) = b, \hspace{3mm} \rho_4(c) = c.
\end{split}
\end{equation}

3) The 2-dimensional representation $\rho_2$ constructed by means
of an epimorphism to the symmetric group $S_3$
\begin{equation}
  \mathcal{O} \longrightarrow \mathcal{O}/Q_8 \simeq S_3,
\end{equation}
where $Q_8$ is the quaternion group.

4) Two 3-dimensional representations $\rho_5$ and $\rho_6$
constructed by means of an epimorphism to the symmetric group
$S_4$
\begin{equation}
   \mathcal{O} \longrightarrow \mathcal{O}/\{1,-1\} \simeq S_4.
\end{equation}
The representations $\rho_5$ and $\rho_6$ are related as follows
\begin{equation}
 \rho_6(a) = -\rho_5(a), \hspace{3mm} \rho_6(b) = \rho_5(a),
       \hspace{3mm} \rho_6(c) = \rho_5(c).
\end{equation}

5) The 4-dimensional representation $\rho_7$ constructed as the
tensor product $\rho_2\otimes\rho_3$ (it coincides with
$\rho_2\otimes\rho_4$).

The characters of representations $\rho_i$ for $i=0,\dots,7$ are
collected in Table \ref{characters_O}.
\end{proposition}

\PerfProof 1) and 2) are clear from constructions (\ref{constr_rho_01}) and
(\ref{constr_rho_34}).

3) We construct the third 2-dimensional representation $\rho_2$ by using
the homomorphism
\begin{equation}
 \label{factor_Q8_1}
 \pi_2: \mathcal{O} \longrightarrow \mathcal{O}/Q_8,
\end{equation}
see (\ref{Q_8_normal}). The quotient group $\mathcal{O}/Q_8$ is isomorphic
to the symmetric group $S_3$ consisting of $6$ elements. The cosets of
$\mathcal{O}/Q_8$ are
\begin{equation}
 \begin{split}
  & \{ Q_8, cQ_8, c^2Q_8, abQ_8, acQ_8, ac^2Q_8 \} = \\
  & \{ \{1\}, \{c\}, \{c^2\}, \{ab\}, \{ac\}, \{ac^2\} \}.
 \end{split}
\end{equation}

By (\ref{c_a_}), (\ref{c_a_1}) and (\ref{conj_4})
\begin{equation}
\begin{split}
 & (ab)(ab) = aa^{-1}b^2 = -1 \in Q_8,
       \text{ i.e., } \{ab\}^2 = \{1\}, \\
 & (ac)(ac) = a(ca)c = a(abc^2)c = -a^2b \in Q_8,
       \text{ i.e., } \{ac\}^2 = \{1\}, \\
 & (ac^2)(ac^2) = ca^{-1}(ac^2) = c^3 = -1 \in Q_8,
       \text{ i.e., } \{ac^2\}^2 = \{1\}, \\
 & (ab)(ac) = bc \in cQ_8
       \text{ i.e., } \{ab\}\{ac\} = \{c\}, \\
 & (ac)(ab) = a(ca)b = a(abc^2)b
       = -a^2b(c^2b^{-1}) = \\
     & -a^2b(b^{-1}a^6c^2) = -c^{2} \in c^2Q_8,
       \text{ i.e., } \{ab\}\{ac\} = \{c\}^2, \\
 & (c)(ab) = (ca)b = (abc^2)b
       = -ab(c^2b^{-1}) = \\
     & -ab(b^{-1}a^6c^2) = a^3c^{2} \in Q_8ac^2 = ac^2Q_8,
       \text{ i.e., } \{c\}\{ab\} = \{ac^2\}, \\
\end{split}
\end{equation}
and so on. Thus,
\begin{equation}
 \label{isom_S3}
\begin{split}
 & \mathcal{O}/Q_8 \simeq S_3, \\
 & \{ab\} \simeq (12), \hspace{3mm}
    \{ac\} \simeq (13), \hspace{3mm}
    \{ac^2\} \simeq (23), \\
 &  \{c\} \simeq (123), \hspace{3mm}
    \{c^2\} \simeq (132).
\end{split}
\end{equation}
The symmetric group $S_3$ has a 2-dimensional representation
$\gamma_2$ such that
\begin{equation*}
\begin{split}
 & tr~\gamma_2(12) = tr~\gamma_2(13) = tr~\gamma_2(23) = 0, \hspace{3mm}
   \text{ i.e., } \\
 & tr~\gamma_2\{ab\} = tr~\gamma_2\{ac\} = tr~\gamma_2\{ac^2\} = 0, \\
 & tr~\gamma_2(123) = tr~\gamma_2(132) = -1,
   \text{ i.e., } tr~\gamma_2\{c\} = tr~\gamma_2\{c^2\} = -1, \\
 &  tr~\gamma_2(1) = 2,
\end{split}
\end{equation*}
see, e.g., \cite[\S32]{CR62}. We consider now the representation
$\rho_2$ as the composition of epimorphism $\pi_2$ and the
representation $\gamma_2$, i.e.,
\begin{equation}
  \rho_2 = \gamma_2 \pi_2.
\end{equation}
So, for all $u \in Cl(ab) \cup Cl(a) \cup Cl(a^3)$, we see   that
$$
\pi_2(u) \in \{ab\} \cup \{ac\} \cup \{ac^2\}\;\text{ and
$tr~\rho_2(u) = 0$.}
$$
For all $ u \in Cl(c) \cup Cl(c^2)$, we see that
$$
\pi_2(u) \in \{c\} \cup \{c^2\}\;\text{ and $tr~\rho_2(u) = -1$.}
$$
Finally, for all $u \in Cl(1) \cup Cl(-1) \cup Cl(b)$ we see that
$$
\pi_2(u) \in \{1\}\;\text{  and $tr~\rho_2(u) = 2$.}
$$
Thus we obtain the row of characters $\psi_3$.

\begin{table} [h]
  \centering
  \vspace{2mm}
  \caption{\hspace{3mm}The 3-dimensional characters of $S_4$}
  \renewcommand{\arraystretch}{1.5}
 \begin{tabular} {||c|c|c|c|c|c||}
 \hline \hline
   Character  & \multicolumn{5}{c||}
         {Conjugacy class $C_i \subset S_4$ and} \cr
          & \multicolumn{5}{c||}{its order $|C_i|$ under it} \\
 \cline{2-6}
 & ~$C_1~$ & ~$C_2$~ & ~$C_3$~ & ~$C_4$~ & ~$C_5$~  \cr
              &  1  & 6   & 8   & 6   & 3    \\
 \hline \hline
     $\gamma_5$ &  3  & $-1$   & 0   & 1   & $-1$      \\
 \hline
     $\gamma_6$ &  3  &  1   & 0   & $-1$  & $-1$     \\
 \hline \hline
 \end{tabular}
  \label{characters_S4}
\end{table}

4) We construct representations $\rho_5$ and $\rho_6$ by means of
the epimorphism
\begin{equation}
   \pi_{56} : \mathcal{O} \longrightarrow \mathcal{O}/\{1,-1\}.
\end{equation}

The epimorphism $\pi_{56}$ is well-defined because the subgroup
$\{1,-1\}$ is normal:
\begin{equation}
     \{1,-1\} \triangleleft \mathcal{O}.
\end{equation}
The quotient group $\mathcal{O}/\{1,-1\}$ is the $24$-element
octahedral group coinciding with the symmetric group $S_4$. By
\cite[\S32]{CR62}  $S_4$ has two 3-dimensional representations,
$\gamma_5$ and $\gamma_6$, with characters as in Table
\ref{characters_S4}.

In Table \ref{characters_S4} we give conjugacy classes $C_i$ of
the group $S_4$ together with the number of elements of these
classes:
\begin{equation}
\begin{split}
  & C_1 = \{1\}, \\
  & C_2 = \{ab, a^3b, a^3c, a^3bc, ac^2, a^3bc^2\}, \\
  & C_3 = \{c, c^2, bc, a^2c, a^2bc, bc^2, a^2c^2, a^2bc^2 \}, \\
  & C_4 = \{a, a^3, ac, abc, a^3c^2, abc^2 \}, \\
  & C_5 = \{a^2, b, a^2b\}.
\end{split}
\end{equation}
We have
\begin{equation}
  \rho_5 = \gamma_5 \pi_{56}, \hspace{3mm}
  \rho_6 = \gamma_6 \pi_{56}.
\end{equation}
For $ u \in Cl(1) \cup Cl(-1)$, we see that
$$
\pi_{56}(u) = 1\;\text{ and $tr~\rho_5(u) = 3$.}
$$
For all $u \in Cl(ab)$, we see  that
$$
\pi_{56}(u) \in C_2\;\text{ and $tr~\rho_5(u) = -1$.}
$$
For $ u \in Cl(b)$, we see that
$$
\pi_{56}(u) \in C_5\;\text{ and $tr~\rho_5(u) = -1$.}
$$
For $u \in Cl(c) \cup Cl(c^2)$ , we see that
$$
\pi_{56}(u) \in C_3\;\text{ and $tr~\rho_5(u) = 0$.}
$$
Finally, for every $u \in Cl(a) \cup Cl(a^2)$ we see that
$$
\pi_{56}(u) \in C_4\;\text{ and $tr~\rho_5(u) = 1$.}
$$
Thus we obtain the row of characters $\psi_5$.

Note, that $\rho_6$ can be obtained from $\rho_5$ by the following
relations:
\begin{equation}
  \rho_6(a) = -\rho_5(a), \hspace{3mm}
  \rho_6(b) = \rho_5(b), \hspace{3mm}
  \rho_6(c) = \rho_5(c).
\end{equation}

5) Finally, the 4-dimensional representation $\rho_7$ is
constructed as either of the tensor products $\rho_2\otimes\rho_3$
or $\rho_2\otimes\rho_4$. Observe that $\psi_2\psi_3$ and
$\psi_2\psi_4$ have the same characters,
 see Table \ref{characters_rho_7}.

\begin{table} [h]
  \centering
  \vspace{2mm}
  \caption{\hspace{3mm}The character of the 4-dimensional representation $\rho_7$}
  \renewcommand{\arraystretch}{1.5}
 \begin{tabular} {||c|c|c|c|c|c|c|c|c||}
 \hline \hline
       & $Cl(1)$ & $Cl(-1)$ & $Cl(ab)$ &
         $Cl(b)$ & $Cl(c^2)$ & $Cl(c)$ & $Cl(a)$ & $Cl(a^3)$  \\
 \hline \hline
     $\psi_7 = \psi_2\psi_3 = \psi_2\psi_4$ &  4  & -4   & 0   & 0   & 1   & -1   & 0  &  0  \\
 \hline \hline
 \end{tabular}
  \label{characters_rho_7}
\end{table}

The irreducibility of $\rho_2\otimes\rho_3$ follows from the fact
that
$$
    < \psi_2\psi_3, \psi_2\psi_3 > = \frac{16 + 16 + 8 + 8}{48} = 1.
\qed
$$

Select $\rho_3$ as a faithful representation of $\mathcal{O}$ from
the McKay correspondence. All irreducible representations $\rho_i$
of $\mathcal{O}$ (see Proposition \ref{charachters_O_prop} and
Table \ref{characters_O}) can be placed in vertices of the
extended Dynkin diagram $\tilde{E}_7$, see (\ref{mckay_graph_E7}):
\begin{equation}
 \label{mckay_graph_E7}
 \begin{array}{cccc cccc cccc c}
 \rho_0 & \line(1,0){15} & \rho_3 & \line(1,0){15} & \rho_5 & \line(1,0){15} &
 \rho_7 & \line(1,0){15} &
 \rho_6 & \line(1,0){15} & \rho_4 & \line(1,0){15} & \rho_1 \\
 & & & & & & \line(0,1){20} & & & & &  \\
 & & & & & & \rho_2 & & & & & & \\
 \end{array}
\end{equation}
Then, according to the McKay correspondence we have the following
decompositions of the tensor products $\rho_3\otimes\rho_i$:
\begin{equation}
\label{rho_tensor_1}
\begin{split}
 & \rho_3\otimes\rho_0 = \rho_3, \\
 & \rho_3\otimes\rho_1 = \rho_4, \\
 & \rho_3\otimes\rho_2 = \rho_7, \\
 & \rho_3\otimes\rho_3 = \rho_0 + \rho_5, \\
 & \rho_3\otimes\rho_4 = \rho_1 + \rho_6, \\
 & \rho_3\otimes\rho_5 = \rho_3 + \rho_7, \\
 & \rho_3\otimes\rho_6 = \rho_4 + \rho_7, \\
 & \rho_3\otimes\rho_7 = \rho_2 + \rho_5 + \rho_6. \\
\end{split}
\end{equation}

\begin{table} [h]
  \centering
  \vspace{2mm}
  \caption{\hspace{3mm}The characters of the binary tetrahedral group}
  \renewcommand{\arraystretch}{1.5}
 \begin{tabular} {||c|c|c|c|c|c|c|c|c||}
 \hline \hline
   Character & \multicolumn{7}{c|}
         {Conjugacy class $Cl(g)$ and its order $|Cl(g)|$ under it}  & Note on \\
 \cline{2-8}
   $\chi_i$ & $Cl(1)$ & $Cl(-1)$ & $Cl(b)$ &
         $Cl(c)$ & $Cl(c^2)$ & $Cl(-c)$ & $Cl(-c^2)$ & $\tau_i$ \cr
        &  1  & 1   & 6   & 4  & 4   & 4   & 4  &  \\
 \hline \hline
     $\chi_0$ &  1  & 1   & 1   & 1  & 1   & 1   & 1  & trivial \\
 \hline
     $\chi_1$ &  1  & 1   & 1   & $\omega_3$  & $\omega_3^2$
                       & $\omega_3$  & $\omega_3^2$  & $\tau_1(c) =
\omega_3$   \\
 \hline
     $\chi_2$ &  1  & 1   & 1   & $\omega_3^2$ & $\omega_3$
                       & $\omega_3^2$ & $\omega_3$   & $\tau_2(c) =
\omega_3^2$   \\
 \hline
     $\chi_3$ &  2  & $-2$  & 0   & 1   & $-1$   & $-1$  & 1  & faithful \\
 \hline
     $\chi_4$ &  2  & $-2$  & 0   & $\omega_3$  & $-\omega_3^2$
                       & $-\omega_3$ & $\omega_3^2$ & $\tau_3\otimes\tau_1$
\\
 \hline
     $\chi_5$ &  2  & $-2$  & 0   & $\omega_3^2$ & $-\omega_3$
                       & $-\omega_3^2$ & $\omega_3$ & $\tau_3\otimes\tau_2$
\\
 \hline
     $\chi_6$ &  3  & 3  & $-1$   & 0   & 0   & 0  & 0  & $\gamma_6\pi_6$ \\
 \hline \hline
 \end{tabular}
  \label{characters_T}
\end{table}

In Table \ref{characters_T} we assume that $\omega_3 =
e^{2\pi{i}/3}$.

\begin{proposition}
\label{charachters_T_prop}
The group $\mathcal{T}$ has the
following 7 irreducible
 representations:

1) Three 1-dimensional representations
\begin{equation}
\begin{split}
 & \tau_0(a) = \tau_0(b) = \tau_0(c) = 1, \\
 & \tau_1(a) = \tau_1(b) = 1, \hspace{3mm} \tau_1(c) = \omega_3, \\
 & \tau_2(a) = \tau_2(b) = 1, \hspace{3mm} \tau_2(c) = \omega_3^2, \\
\end{split}
\end{equation}
 Representations $\tau_1$ and $\tau_2$ can be constructed by using an
epimorphism
onto the alternating group $A_4$:
\begin{equation}
 \label{epim_A4}
  \mathcal{T} \longrightarrow \mathcal{T}/\{1,-1\} = A_4.
\end{equation}

2) The faithful 2-dimensional representation
\begin{equation}
  \tau_3(a) = a, \hspace{3mm} \tau_3(b) = b, \hspace{3mm}
  \tau_0(c) = c, \\
\end{equation}

3) Two 2-dimensional representation $\tau_4$ and $\tau_5$
constructed as tensor products
\begin{equation}
  \tau_4 = \tau_3\otimes\tau_1, \hspace{3mm}
  \tau_5 = \tau_3\otimes\tau_2.
\end{equation}

4) The 3-dimensional representation $\tau_6$ constructed by using
an epimorphism (\ref{epim_A4}) onto the alternating group $A_4$.

The characters of representations $\tau_i$ for $i = 0,\dots,6$ are
collected in Table \ref{characters_T}.

\end{proposition}
\PerfProof 1), 2) and 3) are easily checked.

4) Consider the 3-dimensional representation of the alternating
group $A_4$ with the character given in Table
\ref{3_character_A4}, see, e.g., \cite[\S32]{CR62}:

\begin{table} [h]
  \centering
  \vspace{2mm}
  \caption{\hspace{3mm}The 3-dimensional character of $A_4$}
  \renewcommand{\arraystretch}{1.5}
 \begin{tabular} {||c|c|c|c|c||}
 \hline \hline
   Character  & \multicolumn{4}{c||}
      {Conjugacy class $C_i \subset A_4$ and} \cr
      & \multicolumn{4}{c||}{its order $|C_i|$ under it} \\
     \cline{2-5}
    & \hspace{2mm} $C_1$ \hspace{2mm}
    & \hspace{2mm} $C_2$ \hspace{2mm}
    & \hspace{2mm} $C_3$ \hspace{2mm} & $C_4$  \cr
         &  1  & 3   & 4   & 4   \\
 \hline \hline
     $\gamma_6$ &  3  & $-1$   & 0   & 0   \\
 \hline \hline
 \end{tabular}
  \label{3_character_A4}
\end{table}

In Table \ref{3_character_A4} we have $C_1 = \{1\}$, $C_2$
contains only elements of order $2$, and $C_3$, $C_4$ contain only
elements of order $3$, see (\ref{4_conj_classes}).
\begin{equation}
 \label{4_conj_classes}
\begin{split}
 & C_1 = \{1\},  \\
 & C_2 = \{a^2, b, a^2b\}, \\
 & C_3 = \{c, a^2c, bc, a^2bc \}, \\
 & C_4 = \{c^2, a^2c^2, bc^2, a^2bc^2 \}.
\end{split}
\end{equation}

If $\pi_6$ is an epimorphism (\ref{epim_A4}), then
\begin{equation}
  \tau_6 = \gamma_6\pi_6,
\end{equation}
and
\begin{equation}
\renewcommand{\arraystretch}{1.2}
\begin{array}{ll}
\tau_6: Cl(1) \cup Cl(-1) \longrightarrow C_1, & tr~\tau_6(u) = 3, \\
\tau_6: Cl(b) \longrightarrow C_2, & tr~\tau_6(u) = -1, \\
\tau_6: Cl(c) \cup Cl(-c) \longrightarrow C_3,& tr~\tau_6(u) = 0, \\
\tau_6: Cl(c^2) \cup Cl(-c^2) \longrightarrow C_4, & tr~\tau_6(u) = 0. \\
\end{array}
\end{equation}
Thus we get the last row in Table \ref{characters_T}.
\qedsymbol

Select $\tau_3$ as a faithful representation of $\mathcal{T}$ from
the McKay correspondence. All irreducible representations $\tau_i$
of $\mathcal{T}$  (see Proposition \ref{charachters_T_prop} and
Table \ref{characters_T}) can be placed in vertices of the
extended Dynkin diagram $\tilde{E}_6$, see (\ref{mckay_graph_E6}):
\begin{equation}
 \label{mckay_graph_E6}
 \begin{array}{cccc cccc c}
 \tau_1 & \line(1,0){15} & \tau_4 & \line(1,0){15} &
 \tau_6 & \line(1,0){15} & \tau_5 & \line(1,0){15} & \tau_2  \\
 & & & & \line(0,1){20} & & &  \\
 & & & & \tau_3 & & &  \\
 & & & & \line(0,1){20} & & &  \\
 & & & & \tau_0 & & &  \\
 \end{array}
\end{equation}
Then, according to the McKay correspondence, we have the following
decompositions of the tensor products $\tau_3\otimes\tau_i$:
\begin{equation}
\begin{split}
 & \tau_3\otimes\tau_0 = \tau_3, \\
 & \tau_3\otimes\tau_1 = \tau_4, \\
 & \tau_3\otimes\tau_2 = \tau_5, \\
 & \tau_3\otimes\tau_3 = \tau_0 + \tau_6, \\
 & \tau_3\otimes\tau_4 = \tau_1 + \tau_6, \\
 & \tau_3\otimes\tau_5 = \tau_2 + \tau_6, \\
 & \tau_3\otimes\tau_6 = \tau_3 + \tau_4 + \tau_5. \\
\end{split}
\end{equation}

\subsection{Induced and restricted representations}
\index{induced representation}
\index{restricted representation}

Let us denote the characters of induced and restricted representations
of $\mathcal{T}$ and $\mathcal{O}$ as follows:
\begin{equation}
\begin{split}
 \text{irreducible representations } \tau \text{ of } & \mathcal{T}
     \Longleftrightarrow \text{ characters are denoted by } \chi,    \\
 \text{irreducible representations } \rho \text{ of } & \mathcal{O}
     \Longleftrightarrow \text{ characters are denoted by } \psi,    \\
 \text{induced representations } \tau\uparrow^\mathcal{O}_\mathcal{T}
     \text{ of } & \mathcal{O}
     \Longleftrightarrow \text{ characters are denoted by } \chi^\uparrow,
\\
 \text{restricted representations } \rho\downarrow^\mathcal{O}_\mathcal{T}
     \text{ of } & \mathcal{T}
     \Longleftrightarrow \text{ characters are denoted by } \psi^\downarrow.
\\
\end{split}
\end{equation}

\begin{table} [h]
  \centering
  \vspace{2mm}
  \caption{\hspace{3mm}The restricted characters of the binary octahedral group}
  \renewcommand{\arraystretch}{1.9}
 \begin{tabular} {||c|c|c|c|c|c|c|c||}
 \hline \hline
   Character & \multicolumn{7}{c||}
         {Conjugacy class $Cl(g)$ and order of class $|Cl(g)|$}  \\
 \cline{2-8}
   $\psi_i^\downarrow$
        & $Cl_{\mathcal{O}}(1)$ & $Cl_{\mathcal{O}}(-1)$ & $Cl_{\mathcal{O}}(b)$
          & \multicolumn{2}{c|}{$Cl_{\mathcal{O}}(c^2)$}
              & \multicolumn{2}{c||}{$Cl_{\mathcal{O}}(c)$}  \cr
        & 1 & 1 & 6
          & \multicolumn{2}{c|}{8}
              & \multicolumn{2}{c||}{8}  \cr
 \cline{2-8}
              & $Cl_{\mathcal{T}}(1)$   & $Cl_{\mathcal{T}}(-1)$
              & $Cl_{\mathcal{T}}(b)$   & $Cl_{\mathcal{T}}(-c)$
              & $Cl_{\mathcal{T}}(c^2)$  & $Cl_{\mathcal{T}}(c)$
              & $Cl_{\mathcal{T}}(-c^2)$  \cr
        &  1  & 1   & 6   & 4  & 4   & 4   & 4  \\
 \hline \hline
     $\psi_1^\downarrow = \psi_0^\downarrow$
          &  1  & 1  & 1  & 1  & 1  & 1  & 1   \\
 \hline
     $\psi_2^\downarrow$
          &  2  & 2  & 2  & $-1$ & $-1$  & $-1$ & $-1$   \\
 \hline
     $\psi_4^\downarrow = \psi_3^\downarrow$
          &  2  & $-2$ & 0  & $-1$ & $-1$  & 1 & 1   \\
 \hline
     $\psi_6^\downarrow = \psi_5^\downarrow$
          &  3  & 3  & $-1$ & 0 & 0  & 0 & 0  \\
 \hline
     $\psi_7^\downarrow$
          &  4  & $-4$ & 0  & 1 & 1  & $-1$ & $-1$    \\
 \hline \hline
 \end{tabular}
  \label{restr_char_O}
\end{table}
Consider the restriction of the binary octahedral group $\mathcal{O}$ onto the
binary tetrahedral subgroup $\mathcal{T}$. Then the conjugacy
classes $Cl(a)$, $Cl(a^3)$ and $Cl(ab)$ disappear, and the remaining
5 classes split into 7 conjugacy classes, see Table \ref{restr_char_O}.
We denote the conjugacy classes of $\mathcal{O}$ by $Cl_\mathcal{O}$ and
conjugacy classes of $\mathcal{T}$ by $Cl_\mathcal{T}$.

Now consider the restricted representations $\psi_i^\downarrow$ from
$\mathcal{O}$ onto $\mathcal{T}$. By Table \ref{characters_O}
$\psi_0$ (resp. $\psi_3$ or $\psi_5$)
differs from $\psi_1$ (resp. $\psi_4$ or $\psi_6$)
only on $Cl(a)$, $Cl(a^3)$, and $Cl(ab)$, so we have the following coinciding
pairs
of restricted representations:
\begin{equation}
\psi_0^\downarrow = \psi_1^\downarrow, \quad
\psi_3^\downarrow = \psi_4^\downarrow, \quad
\psi_5^\downarrow = \psi_6^\downarrow.
\end{equation}

The values of characters $\psi_i^\downarrow$ for $i = 0, 2, 3, 5, 7$    are
easily obtained
from the corresponding characters $\psi_i$.

\begin{figure}[h]
\centering
\includegraphics{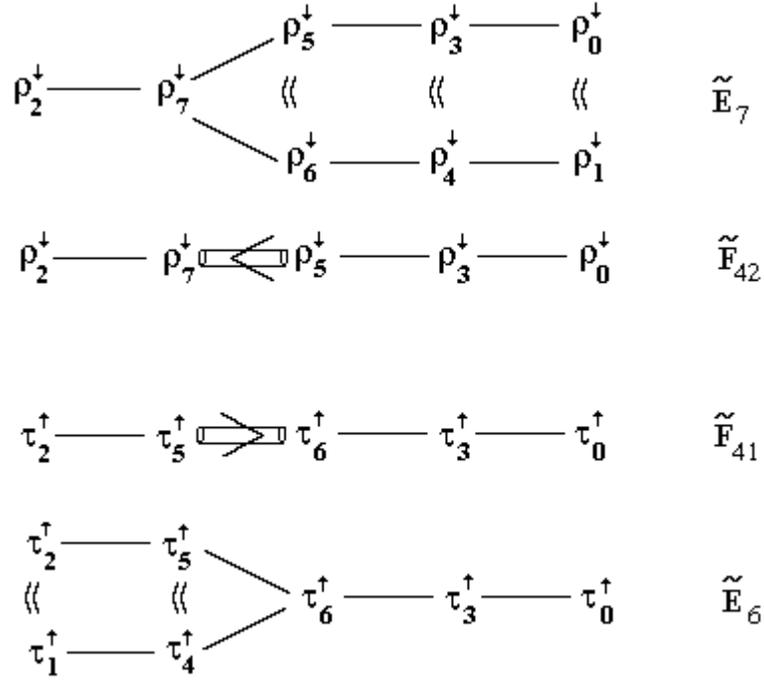}
\caption{\hspace{3mm} The induced and restricted representations
     of $\mathcal{T} \triangleleft \mathcal{O}$ }
\label{induced_restr_repr}
\end{figure}

Observe that $\rho_3^\downarrow = \tau_3$ is a faithful
representation of $\mathcal{T}$ with character $\psi_3^\downarrow
= \chi_3$. All irreducible representations $\rho_i^\downarrow$ for
$i = 0, 2, 3, 5, 7$ can be placed in vertices of the extended
Dynkin diagram $\tilde{F}_{42}$, see
Fig.~\ref{induced_restr_repr}. From Table \ref{restr_char_O} we
have
\begin{equation}
\label{tensor_prod_Slodowy_1}
\begin{split}
 & \tau_3\otimes\rho_0^\downarrow = \rho_3^\downarrow\otimes\rho_0^\downarrow =
\rho_3^\downarrow, \\
 & \tau_3\otimes\rho_2^\downarrow = \rho_3^\downarrow\otimes\rho_2^\downarrow =
\rho_7^\downarrow, \\
 & \tau_3\otimes\rho_3^\downarrow = \rho_3^\downarrow\otimes\rho_3^\downarrow =
\rho_0^\downarrow + \rho_5^\downarrow, \\
 & \tau_3\otimes\rho_5^\downarrow = \rho_3^\downarrow\otimes\rho_5^\downarrow =
\rho_3^\downarrow + \rho_7^\downarrow, \\
 & \tau_3\otimes\rho_7^\downarrow = \rho_3^\downarrow\otimes\rho_7^\downarrow =
\rho_2^\downarrow + 2\rho_5^\downarrow. \\
\end{split}
\end{equation}
Thus, the decompositions (\ref{tensor_prod_Slodowy_1}) constitute the following
matrix
\begin{equation}
\label{matrix_Slodowy_1}
  \tilde{A} = \left ( \begin{array}{ccccc}
        0 & 1 & 0 & 0 & 0 \vspace{2mm} \\
                1 & 0 & 2 & 0 & 0 \vspace{2mm} \\
                0 & 1 & 0 & 1 & 0 \vspace{2mm} \\
                0 & 0 & 1 & 0 & 1 \vspace{2mm} \\
                0 & 0 & 0 & 1 & 0
       \end{array}
  \right )
  \begin{array}{c}
    \rho_2^\downarrow \vspace{2mm} \\
    \rho_7^\downarrow \vspace{2mm} \\
    \rho_5^\downarrow \vspace{2mm} \\
    \rho_3^\downarrow \vspace{2mm} \\
    \rho_0^\downarrow
  \end{array}
  \hspace{7mm}
\end{equation}
where the row associated with $\rho_i^\downarrow$ for $i =
2,7,5,3,0$ consists of the decomposition coefficients
(\ref{tensor_prod_Slodowy_1}) of $\tau_3\otimes\rho_i^\downarrow$.
The matrix $\tilde{A}$ satisfies the relation \index{Slodowy
matrix} \index{Slodowy correspondence}
\begin{equation}
 \label{slodowy_1}
   \tilde{A} = 2I - K,
\end{equation}
where $K$ is the Cartan matrix for the extended Dynkin diagram
$\tilde{F}_{42}$. The matrix $K$ in (\ref{slodowy_1}) does not
correspond to the bicolored order of vertices (see
\S\ref{bicolored_C}), but can be obtained from the bicolored one
by permutations of rows and columns. We call the matrix
$\tilde{A}$ the {\it Slodowy matrix}; it is an analog of the McKay
matrix for the multiply-laced case.

\index{formula! - Frobenius reciprocity} \index{Frobenius
reciprocity formula} Now we move on to the dual case and consider
induced representations. To obtain induced representations
$\chi\uparrow^\mathcal{O}_\mathcal{T}$, we use the {\it Frobenius
reciprocity formula} connecting restricted and induced
representations, see, e.g., \cite[21.16]{JL2001}
\begin{equation}
 \label{reciprocity}
 <\psi, \chi\uparrow^G_H>_G ~=~ <\psi\downarrow^G_H, \chi>_H.
\end{equation}
By (\ref{reciprocity}) we have the following expression for the
characters of induced representations
\begin{equation}
 \label{reciprocity_2}
 \chi\uparrow^G_H ~=~
 \sum\limits_{\psi_i \in Irr(G)}
 <\psi_i, \chi\uparrow^G_H>_G\psi_i ~=~
 \sum\limits_{\psi_i \in Irr(G)}<\psi\downarrow^G_H, \chi>_H\psi_i.
\end{equation}
Thus, to calculate the characters $\chi^\uparrow =
\chi\uparrow^\mathcal{O}_\mathcal{T}$, we only need to calculate
the inner products
\begin{equation}
 \label{inner_prod_1}
 <\psi^\downarrow, \chi>_\mathcal{T} ~=~
 <\psi\downarrow^\mathcal{O}_\mathcal{T}, \chi>_\mathcal{T}.
\end{equation}

\begin{table} [h]
  \centering
  \vspace{2mm}
  \caption{\hspace{3mm}The inner products $<\psi^\downarrow, \chi>$}
  \renewcommand{\arraystretch}{1.5}
 \begin{tabular} {||c|c|c|c|c|c||}
 \hline \hline
        & $\psi^\downarrow_1 = \psi^\downarrow_0$
        & $\psi^\downarrow_2$
        & $\psi^\downarrow_3 = \psi^\downarrow_4$
        & $\psi^\downarrow_5 = \psi^\downarrow_6$
        & $\psi^\downarrow_7$ \\
 \hline \hline
     $\chi_0$ &  1  & 0   & 0   & 0  & 0  \\
 \hline
     $\chi_1$ &  0  & 1   & 0   & 0  & 0  \\
 \hline
     $\chi_2$ &  0  & 1   & 0   & 0  & 0  \\
 \hline
     $\chi_3$ &  0  & 0   & 1   & 0  & 0  \\
 \hline
     $\chi_4$ &  0  & 0   & 0  & 0  & 1  \\
 \hline
     $\chi_5$ &  0  & 0   & 0  & 0  & 1  \\
 \hline
     $\chi_6$ &  0  & 0   & 0  & 1  & 0  \\
 \hline \hline
 \end{tabular}
  \label{table_inner_prod}
\end{table}

One can obtain the inner products (\ref{inner_prod_1}) from
Tables \ref{characters_T} and \ref{restr_char_O}.
The results are given in Table \ref{table_inner_prod}.
Further, from Table \ref{table_inner_prod} and (\ref{reciprocity_2}) we deduce
\begin{equation}
\label{chi_induced_1}
\begin{split}
     \chi_0^\uparrow = & \psi_0 + \psi_1, \\
     \chi_1^\uparrow = \chi_2^\uparrow = & \psi_2, \\
     \chi_3^\uparrow = & \psi_3 + \psi_4, \\
     \chi_4^\uparrow = \chi_5^\uparrow = & \psi_7, \\
     \chi_6^\uparrow = & \psi_5 + \psi_6, \\
\end{split}
\hspace{20mm}
\begin{split}
     \tau_0^\uparrow = & \rho_0 + \rho_1, \\
     \tau_1^\uparrow = \tau_2^\uparrow = & \rho_2, \\
     \tau_3^\uparrow = & \rho_3 + \rho_4, \\
     \tau_4^\uparrow = \tau_5^\uparrow = & \rho_7, \\
     \tau_6^\uparrow = & \rho_5 + \rho_6. \\
\end{split}
\end{equation}

Let us find the tensor products $\rho_{f}\otimes\tau_i^\uparrow =
\rho_3\otimes\tau_i^\uparrow$, where $\rho_{f}$ means the faithful
representation of $\mathcal{O}$.
 By (\ref{rho_tensor_1}) and (\ref{chi_induced_1})
we have
\begin{equation}
\label{tensor_prod_Slodowy_2}
\begin{split}
  \rho_3 \otimes \tau_0^\uparrow = &
  \rho_3 \otimes (\rho_0 + \rho_1) =
    \rho_3 \otimes \rho_0 + \rho_3 \otimes \rho_1 =
    \rho_3 + \rho_4 = \tau_3^\uparrow, \\
  \rho_3 \otimes \tau_1^\uparrow = & \rho_3 \otimes \tau_2^\uparrow =
  \rho_3 \otimes \rho_2 = \rho_7 = \tau_4^\uparrow = \tau_5^\uparrow, \\
  \rho_3 \otimes \tau_3^\uparrow = &
  \rho_3 \otimes (\rho_3 + \rho_4) = \\
 & \rho_3 \otimes \rho_3 + \rho_3 \otimes \rho_4 =
  (\rho_0 + \rho_5) + (\rho_1 + \rho_6) =
  \tau_0^\uparrow + \tau_6^\uparrow, \\
  \rho_3 \otimes \tau_4^\uparrow = & \rho_3 \otimes \tau_5^\uparrow =
  \rho_3 \otimes \rho_7 = \rho_2 + \rho_5 + \rho_6 =
  \tau_1^\uparrow + \tau_6^\uparrow =
  \tau_2^\uparrow + \tau_6^\uparrow, \\
  \rho_3 \otimes \tau_6^\uparrow = &
  \rho_3 \otimes (\rho_5 + \rho_6) =
  \rho_3 \otimes \rho_5 + \rho_3 \otimes \rho_6 = \\
 & (\rho_3 + \rho_7) + (\rho_4 + \rho_7) =
  \tau_3^\uparrow + 2\tau_4^\uparrow =
  \tau_3^\uparrow + 2\tau_5^\uparrow.
\end{split}
\end{equation}

Here $\rho_3$ is the faithful representation of $\mathcal{O}$ with
character $\psi_3$. All irreducible representations
$\tau_i^\uparrow$ for $i = 0, 2, 3, 5, 6$ can be placed in
vertices of the extended Dynkin diagram $\tilde{F}_{41}$, see
Fig.~\ref{induced_restr_repr}. For other details, see
\cite[App.III, p. 164]{Sl80}.

The decompositions (\ref{tensor_prod_Slodowy_2}) constitute the following matrix
\begin{equation}
\label{matrix_Slodowy_2}
  \tilde{A}^\vee = \left ( \begin{array}{ccccc}
        0 & 1 & 0 & 0 & 0 \vspace{2mm} \\
                1 & 0 & 1 & 0 & 0 \vspace{2mm} \\
                0 & 2 & 0 & 1 & 0 \vspace{2mm} \\
                0 & 0 & 1 & 0 & 1 \vspace{2mm} \\
                0 & 0 & 0 & 1 & 0
       \end{array}
  \right )
  \begin{array}{c}
    \tau_2^\uparrow \vspace{2mm} \\
    \tau_5^\uparrow \vspace{2mm} \\
    \tau_6^\uparrow \vspace{2mm} \\
    \tau_3^\uparrow \vspace{2mm} \\
    \tau_0^\uparrow
  \end{array}
  \hspace{7mm}
\end{equation}
where the row associated with $\tau_i^\uparrow$ for $i =
2,5,6,3,0$ consists of the decomposition coefficients
(\ref{tensor_prod_Slodowy_2}) of $\rho_3 \otimes \tau_i^\uparrow$.
As in (\ref{slodowy_1}), the matrix $\tilde{A}^\vee$ satisfies the
relation
\begin{equation}
 \label{slodowy_2}
   \tilde{A}^\vee = 2I - K^\vee,
\end{equation}
where $K^\vee$ is the Cartan matrix for the extended Dynkin
diagram $\tilde{F}_{41}$. We see that the matrices $\tilde{A}$ and
$\tilde{A}^\vee$ are mutually transposed:
\begin{equation}
  \tilde{A}^t = \tilde{A}^\vee.
\end{equation}

\index{Slodowy matrix}
\index{Slodowy correspondence}
As in (\ref{slodowy_1}), we call the matrix $\tilde{A}^\vee$ the {\it Slodowy
matrix}.

%% file: 8hurwitz.tex
\chapter{\sc\bf Miscellaneous}

 \section{The triangle groups and Hurwitz groups}
  \label{hurwitz_gr}
 \index{triangle group}
 \index{group! - triangle}
 \index{group! - $(p,q,r)$}
 \index{group! - Hurwitz}

The group generated by $X, Y, Z$ that satisfy the relations
 (\ref{natural_gen_2})
 \begin{equation*}
 X^p = Y^q = Z^r = XYZ = 1
 \end{equation*}
 is said to be a {\it triangle group}.
As it was mentioned in \S\ref{gener_rel},
 the finite polyhedral groups from Table \ref{rotation_pol_1} are
 triangle groups. Set
 \begin{equation*}
\mu(p,q,r):=\frac{1}{p} + \frac{1}{q} + \frac{1}{r}.
 \end{equation*}
\index{group! - soluble}
 \index{group! - insoluble}
 \index{group! - $(2,3,7)$}
 The triangle group is finite if and only if $\mu(p,q,r) > 1$.
 There are only three triangle groups
 \begin{equation*}
 (2,4,4), \hspace{3mm} (2,3,6), \hspace{3mm} (3,3,3),
 \end{equation*}
for which $\mu(p,q,r) = 1$.
 These groups are infinite and soluble. For references, see
 \cite{Con90}, \cite{Mu01}.
 For all other triangle groups, we have
\begin{equation*}
 \mu(p,q,r) < 1.
\end{equation*}
These groups are infinite and insoluble. The value $1 -
\mu(p,q,r)$ attains the minimum value $\displaystyle\frac{1}{42}$
at $(2,3,7)$. Thus, $(2,3,7)$ is, in a sense, the minimal infinite
insoluble triangle group. The importance of the triangle group
$(2,3,7)$ is revealed by the following theorem due to Hurwitz.
\index{theorem! - Hurwitz}
\begin{theorem}[Hurwitz, \cite{Hur1893}]
If $X$ is a compact Riemann surface of genus $g > 1$, then $|AutX|
\leq 84(g-1)$, and moreover, the upper bound of this order is
attained if an only if $|AutX|$ is a homomorphic image of the
triangle group $(2,3,7)$.
\end{theorem}
For further references, see \cite{Con90}, \cite{Con03}.

A {\it Hurwitz group} is any finite nontrivial quotient of the
triangle group $(2,3,7)$. In other words, the finite group $G$ is
the Hurwitz group if it has generators $X, Y, Z \in G$ such that
\begin{equation*}
 X^2 = Y^3 = Z^7 = XYZ = 1.
\end{equation*}

M.~Conder writes that the significance of the Hurwitz groups
\lq\lq \ldots is perhaps best explained by referring to some
aspects of the theory of Fuchsin groups, hyperbolic geometry,
Riemann surfaces...", see \cite[p.359]{Con90} and a bibliography
cited there.

\index{theorem! - Conder}
 M.~Conder \cite{Con80} using the method of coset graphs developed
by G.~Higman have shown that
the alternating group $A_n$ is a Hurwitz group for all $n \geq 168$.

\index{theorem! - Lucchini-Tamburini-Wilson}
 Recently, A.~Lucchini, M.~C.~Tamburini and J.~S.~Wilson
showed that most finite simple classical groups of sufficiently
large rank are Hurwitz groups, see \cite{LuT99}, \cite{LuTW00}.
For example, the groups $SL_n(q)$ are Hurwitz, for all $n > 286$
\cite{LuTW00}, and the groups $Sp_{2n}(q),  SU_{2n}(q)$ are
Hurwitz, for all $n > 371$ \cite{LuT99}. (Mentioned groups act in
the n-dimensional vector space over the field $F_q$ of the prime
characteristic $q$).

\index{theorem! - Wilson on Monster}
\index{group! - Monster}
 The sporadic groups have been treated in a series of papers
by Woldar and others, see \cite{Wi01} for a survey and references.
R.~A.~Wilson shows in \cite{Wi01}, that the Monster is also a
Hurwitz group.

\section{Algebraic integers}
  \label{sect_numbers}
\index{algebraic number} \index{algebraic integer}
\index{conjugates of an algebraic number} If $\lambda$ is a root
of the polynomial equation
\begin{equation}
   \label{monic}
     a_n{x}^n + a_n{x}^{n-1} + \dots + a_1{x} + a_0 = 0,
\end{equation}
where $a_i$ for $i = 0,1,\dots,n$ are integers and $\lambda$
satisfies no similar equation of degree $<n$, then $\lambda$ is
said to be an {\it algebraic number} of degree $n$. If $\lambda$
is an algebraic number and $a_n = 1$, then $\lambda$ is called an
{\it algebraic integer}.

 \index{polynomial monic}
 \index{polynomial monic integer}

 A polynomial $p(x)$ in which the
coefficient of the highest order term is equal to 1 is called the
{\it monic polynomial}. The polynomial (\ref{monic}) with integer
coefficients and $a_n = 1$ is the {\it monic integer polynomial}.

The algebraic integers of degree 1 are the ordinary integers
(elements of $\mathbb{Z}$). If $\alpha$ is an algebraic number of
degree $n$ satisfying the polynomial equation
$$
   (x - \alpha)(x - \beta)(x - \gamma) \dots = 0,
$$
then there are $n-1$ other algebraic numbers $\beta$, $\gamma$,
... called the {\it conjugates} of $\alpha$. Furthermore, if
$\alpha$ satisfies any other algebraic equation, then its
conjugates also satisfy the same equation.

\begin{definition} \rm{
An algebraic integer $\lambda > 1$ is said to be a {\it Pisot
number} if its conjugates (other then $\lambda$ itself) satisfy
$|\lambda^{'}| < 1$. }
\end{definition}

\index{Pisot number}
\index{Zhang number} The smallest Pisot
number,
\begin{equation}
  \label{pisot_zhang}
          \lambda_{Pisot} \approx 1.324717...,
\end{equation}
is a root of $\lambda^3 - \lambda - 1 = 0$ (for details and
references, see \cite{McM02}). This number appears in Proposition
\ref{polyn_T_1} as a limit of the spectral radius
$\rho(T_{2,3,n})$ as $n \rightarrow \infty$. This number was also
obtained by Y.~Zhang [Zh89] and used in the study of regular
components of an Auslander-Reiten quiver, so we call the number
(\ref{pisot_zhang}) the {\it Zhang number}.

\begin{definition} {\rm
Let $p(x)$ be a monic integer polynomial, and define its {\it
Mahler measure} to be
\begin{equation}
\Vert p(x)\Vert  = \prod\limits_{\beta} |\beta| ,
\end{equation}
where $\beta$ runs over all (complex) roots of $p(x)$ outside the
unit circle. }
\end{definition}
It is well known that $\Vert p(x)\Vert  = 1$ if and only if all
roots of $p(x)$ are roots of unity. In 1933, Lehmer \cite{Leh33}
asks whether for each $\varDelta \ge 1$, there exists an algebraic
integer such that
\begin{equation}
 \label{numb_lehmer}
    1 < \Vert \alpha \Vert < 1 + \varDelta
\end{equation}

\index{Mahler measure}
\index{Lehmer's number}

Lehmer found polynomials with smallest Mahler measure for small
degrees and stated in \cite[p.18]{Leh33} that the polynomial with
minimal root $\alpha$ (in the sense of \ref{numb_lehmer}) he could
find is the polynomial of degree $10$:
\begin{equation}
  \label{lehmer_polynom}
   1 + x  - x^3 - x^4 - x^5 - x^6 - x^7 + x^9 + x^{10},
\end{equation}
see \cite{Hir02}, \cite{McM02}; cf. Remark \ref{mcmullen}. Outside
the unit circle, the polynomial (\ref{lehmer_polynom}) has only
one root
\begin{equation}
  \label{lehmer_number}
  \lambda_{Lehmer} \approx 1.176281...
\end{equation}
The number (\ref{lehmer_number}) is called {\it Lehmer's number};
see Proposition \ref{polyn_T_1}, Remark \ref{mcmullen} and Table
\ref{table_char_E_series}.

\begin{definition} \rm{
An algebraic integer $\lambda > 1$ is said to be a {\it Salem
number} if its conjugates $\lambda'$ (other then $\lambda$ itself)
satisfy $|\lambda'| \leq 1$ and include
$\displaystyle\frac{1}{\lambda}$.
  }
\end{definition}
It is known that every Pisot number is a limit of Salem numbers;
for details and references, see \cite[p.177]{McM02}.
Conjecturally, Lehmer's number (\ref{lehmer_number}) is the
smallest Salem number, \cite{Leh33}, \cite{GH01}.

The positive root of the quadratic equation $\lambda^2 - \lambda -
1 = 0$ is a well-known constant
\begin{equation}
   \lambda_{Golden} \approx 1.618034...,
\end{equation}
called the {\it Golden mean} or {\it Divine proportion}. This
number appears in Proposition \ref{polyn_T_2} as a limit of the
spectral radius $\rho(T_{3,3,n})$ as $n \rightarrow \infty$.

The smallest Mahler measure among polynomials of degree at most
$6$ is
\begin{equation}
    M_6 = ||x^6 - x^4 + x^3 - x^2 + 1|| \approx 1.401268...,
\end{equation}
see \cite[p.1700]{Mos98}. This number appears in Proposition
\ref{polyn_T_2} as a root of 
\begin{equation*}
    \chi(T_{3,3,4}) = x^8 + x^7 - 2x^5 - 3x^4 - 2x^3 + x + 1.
\end{equation*}
This is the polynomial of minimal degree among polynomials
$\chi(T_{3,3,n})$, where $n = 4,5,6,\dots $ with indefinite Tits
form, see Proposition \ref{polyn_T_2} and Table
\ref{table_char_E6_series}.

\section{The Perron-Frobenius Theorem}
  \label{sect_pf}
 \index{matrix! - positive} \index{matrix! - nonnegative}
 \index{matrix! - reducible} \index{matrix! - irreducible}
 \index{theorem! - Perron-Frobenius}

We say that a matrix is {\it positive} (resp. {\it nonnegative})
if all its entries are positive (resp. nonnegative). We use the
notation $A > 0$ (resp. $A \geq 0$) for positive (resp.
nonnegative) matrix. A square $n\times{n}$ matrix $A$ is called
{\it reducible} if the indices 1, 2, ..., $n$ can be divided into
the disjoint union of two nonempty sets $\{i_1, i_2,\dots,i_p\}$
and $\{j_1,j_2,\dots,j_q\}$ (with $p+q=n$) such that
\begin{equation*}
  a_{{i_\alpha}{j_\beta}} = 0, \text{ for }
  \alpha = 1,\dots,p \text{ and } \beta = 1,\dots,q.
\end{equation*}
In other words, $A$ is {\it reducible} if there exists a
permutation matrix $P$, such that
\begin{equation*}
  PAP^t = \left (
    \begin{array}{cc}
      B & 0 \\
      C & D
    \end{array}
      \right ),
\end{equation*}
where $B$ and $D$ are square matrices.
 A square matrix which is not reducible is said to be {\it irreducible}.

 \index{theorem! - on fixed and anti-fixed points}

\begin{theorem} [Perron-Frobenius]
   \label{pf_theorem}
 Let $A$ be $n\times{n}$ nonnegative irreducible matrix. Then the following holds:

 1) One of its eigenvalues $\lambda$ is positive and greater than or equal to
 (in absolute value) all other eigenvalues $\lambda_i$:
\begin{equation*}
   |\lambda_i| \leq \lambda, \hspace{3mm} i = 1,2,\dots,n.
\end{equation*}

 2) There is a positive eigenvector $z$ corresponding to that eigenvalue:
\begin{equation*}
   Az = \lambda{z}, \text{ where } z =(z_1,\dots,z_n)^t \text{ and }  z_i > 0
   \text { for } i = 1,2,\dots,n.
\end{equation*}
Such an eigenvalue $\lambda$ is called the {\em dominant
eigenvalue} of $A$.

 3) The eigenvalue $\lambda$ is a simple root of the characteristic equation of $A$.
\end{theorem}

The following important corollary from the Perron-Frobenius
theorem holds for the eigenvalue $\lambda$:
\begin{equation*}
  \begin{split}
  & \displaystyle \lambda = \max_{z \geq 0} \min_i \frac{(Az)_i}{z_i}
     \hspace{7mm} (z_i \neq 0),       \\
  & \displaystyle \lambda = \min_{z \geq 0} \max_i \frac{(Az)_i}{z_i}
     \hspace{7mm} (z_i \neq 0).
  \end{split}
\end{equation*}

For details, see \cite{MM64}, \cite{Ga90}.

\section{The complex projective line and stereographic projection}
  \label{sect_projective_line}
\index{complex projective line}
   The $n$-dimensional {\it complex projective space} $\mathbb{C}P^n$
   is the set of all complex lines in $\mathbb{C}^{n+1}$ passing
   through the origin. A point $z \in \mathbb{C}P^n$ is defined
   up to the {\it complex} factor $w$. In other words, two
   points
$$
   z_1 = (z^0_1, z^1_1, z^2_1,\dots, z^n_1 ),\hspace{3mm}
   z_2 = (z^0_2, z^1_2, z^2_2,\dots, z^n_2 ) \in \mathbb{C}^{n+1}
$$
   lie on the same line, if
\begin{equation}
  \label{complex_mult}
  \begin{split}
   & z_2 = wz_1 \text{ for some complex
    multiplier } w \in \mathbb{C}, \text{ i.e. } \\
   & z^i_2 = wz^i_1\;\text{ for $i = 0,1,\dots, n$}.
  \end{split}
\end{equation}
The points (\ref{complex_mult}) constitute an equivalence class
 denoted by $[ z^0 : z^1 : \dots : z^n ]$.

 Clearly, the {\it complex projective line}
 $\mathbb{C}P^1$ is the set of all lines in $\mathbb{C}^2$.
 By (\ref{complex_mult}), the points of $\mathbb{C}P^1$ are classes of complex
 pairs $z = [z^0, z^1]$ up to a factor $w \in \mathbb{C}$.

The correspondence
\begin{equation}
  \label{complex_2}
     [z^0 : z^1] \mapsto z^0/z^1
\end{equation}
sets the following bijection maps:
\begin{equation}
 \begin{split}
  & \mathbb{C}P^1 \backslash \{[1:0]\} \Longleftrightarrow \mathbb{C},
    \text{ and } \\
  & \mathbb{C}P^1 \Longleftrightarrow \mathbb{C} \cup \infty.
 \end{split}
\end{equation}

Let $z_1, z_2$ be two vectors from $\mathbb{C}^2$. Define a map
\begin{equation}
   F : \mathbb{C}^2 \longrightarrow \mathbb{R}^3
\end{equation}
by setting
\begin{equation}
  \label{rule_F}
   F(z_1, z_2) :=
   \left (
     \frac{z_1\overline{z}_2 + \overline{z}_1{z}_2}
          {z_1\overline{z}_1 + \overline{z}_2{z}_2},
     \frac{z_1\overline{z}_2 - \overline{z}_1{z}_2}
          {i(z_1\overline{z}_1 + \overline{z}_2{z}_2)},
     \frac{z_1\overline{z}_1 - \overline{z}_2{z}_2}
          {z_1\overline{z}_1 + \overline{z}_2{z}_2}
   \right ),
\end{equation}
(see, for example, \cite{Alv02}).

Let $z_2 = wz_1$,  then
\begin{equation}
  \label{map_F}
   F(z_1, z_2) =
   \left (
     \frac{\overline{w} + w}{1 + |w|^2},
     \frac{\overline{w} - w}{1 + |w|^2},
     \frac{1 - |w|^2}{1 + |w|^2}
   \right ).
\end{equation}
If $w = u + iv$, then
\begin{equation}
 \label{map_to_sphere}
   F(z_1, z_2) =
   \left (
     \frac{2u}{1 + u^2 + v^2},
     \frac{2v}{1 + u^2 + v^2},
     \frac{1 - u^2 - v^2}{1 + u^2 + v^2}
   \right ).
\end{equation}

It is easily to see that
\begin{equation*}
  (2u)^2 + (2v)^2 + (1 - u^2 - v^2)^2 = (1 + u^2 + v^2)^2,
\end{equation*}
so $F(z_1, z_2)$ maps every vector $[z_1, z_2]$ to a point
(\ref{map_to_sphere}) on the unit sphere $S^2$ in $\mathbb{R}^3$.
Another vector $[z_3, z_4]$ defines the same line in the
$\mathbb{C}P^1$ if and only if
\begin{equation*}
    z_4/z_3 = z_2/z_1 = w,
\end{equation*}
i.e., the map F (\ref{map_F}) defines a bijection from the complex
projective line $\mathbb{C}P^1$ and the unit sphere in
$\mathbb{R}^3$:
\begin{equation}
  \label{identity_F}
     F : \mathbb{C}P^1 \Longleftrightarrow S^2.
\end{equation}

\begin{figure}[h]
\centering
\includegraphics{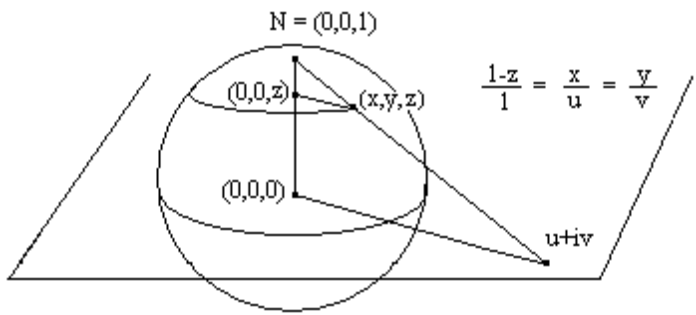}
\caption{\hspace{3mm} The stereographic projection}
\label{stereogr}
\end{figure}

Let $(x, y, z)$ be a point on $S^2$ distinct from the {\it north
pole} $N = \{0,0,1\}$. The {\it stereographic projection} is the
map
\begin{equation}
     S : S^2 \backslash N \longrightarrow \mathbb{C}
\end{equation}
defined by
\begin{equation}
     S(x, y, z) :=  \frac{x}{1-z} + {i}\frac{y}{1-z}. \vspace{3mm}
\end{equation}
see Fig.~\ref{stereogr}. For details, see, for example,
\cite[\S2.9]{Jen94}.

Consider the composition $S(F(z_1,z_2))$. By (\ref{rule_F}), we
have
\begin{equation}
     1 - z = \frac{2\overline{z}_2{z_2}}{z_1\overline{z}_1 +
     \overline{z}_2{z}_2},
\end{equation}
and
\begin{equation}
\begin{split}
    & \frac{x}{1 - z} =
      \frac{z_1\overline{z}_2 + \overline{z}_1{z}_2}
      {2\overline{z}_2{z_2}} =
      \frac{1}{2}
      \left (
      \frac{z_1}{z_2} + \frac{\overline{z_1}}{\overline{z_2}}
      \right ), \\ \\
    & \frac{y}{1 - z} =
      \frac{z_1\overline{z}_2 - \overline{z}_1{z}_2}
      {2\overline{z}_2{z_2}} =
            \frac{1}{2}
      \left (
      \frac{z_1}{z_2} - \frac{\overline{z_1}}{\overline{z_2}}
      \right ).
\end{split}
\end{equation}
Let
$$
  \frac{z_1}{z_2} = w = u + iv.
$$
Then
\begin{equation}
     \frac{x}{1 - z} =
      \frac{1}{2}\left ( w + \overline{w} \right ) = u,
      \hspace{5mm}
     \frac{y}{1 - z} =
      \frac{1}{2}\left ( w - \overline{w} \right ) = v,
\end{equation}
and
\begin{equation}
   S(F(z_1,z_2)) = u + iv = w = \frac{z_1}{z_2}.
\end{equation}

\section{Prime spectrum, coordinate ring, orbit space}
  \label{Hilbert_th}

\subsection{Hilbert's Nullstellensatz (Theorem of zeros)}
 \label{Hilbert_zeros}

From now on, we assume that $R$ is a commutative ring with
identity.

 \index{maximal ideal}
 \index{zero set}
A proper ideal
$\mathfrak{m}$ of $R$ is said to be {\it maximal} if
$\mathfrak{m}$ is not a proper subset of any other proper ideal of
$R$. An ideal $\mathfrak{m} \subset R$ is maximal if and only if
the quotient ring $R/\mathfrak{m}$ is a field. For example, every
ideal $p\mathbb{Z}$ is maximal in the ring of integers
$\mathbb{Z}$, if $p$ is a prime number, and, in this case, the
quotient ring $\mathbb{Z}/p\mathbb{Z}$ is a field.

 \index{prime ideal}
 \index{integral domain}
A proper ideal
$\mathfrak{p} $ of a commutative ring $R$ is called a {\it prime
ideal} if the following condition holds:
\begin{equation}
   \text{for any } a,b \in R, \text{ if } a\cdot{b} \in \mathfrak{p},
   \text{ then either } a \in \mathfrak{p} \text{ or } b \in \mathfrak{p}.
\end{equation}
An ideal $\mathfrak{p} \subset R$ is prime if and only if the
quotient ring $R/\mathfrak{p}$ is an {\it integral domain} (i.e.,
the commutative ring which has no divisors of $0$). Examples of
prime ideals:
\begin{equation}
  \begin{split}
    & \text{ ideal } (5,\sqrt{6}) \text{ in the ring } \mathbb{Z}[\sqrt{6}],\\
    & \text{ ideal } (x) \text { in the ring }\mathbb{Z}[x],\\
    & \text{ ideal } (y + x + 1) \text{ in the ring } \mathbb{C}[x,y].
  \end{split}
\end{equation}
The maximal ideals are prime since the fields are integral
domains, but not conversely:
\begin{equation}
   \fbox{maximal ideals} \subset \fbox{prime ideals} ~~.
\end{equation}
 The ideal $(x) \in \mathbb{Z}[x]$ is prime, but not maximal, since, for
 example:
\begin{equation}
    (x) \subset (2,x) \subset \mathbb{Z}[x].
\end{equation}

 Let $k$ be an algebraically closed field (e.g., the complex field
$\mathbb{C}$), and let $I$ be an ideal in $k[x_1,\dots,x_n]$.
Define $\mathbb{V}_k(I)$, the {\it zero set} of $I$ , by
\begin{equation}
  \label{zero_set}
    \mathbb{V}_k(I) = \{(a_1,\dots,a_n) \in k  ~|~ f(a_1,\dots.a_n) = 0
      \text{ for all }  f \in I \}.
\end{equation}

 \index{affine space}
 \index{Zariski's topology}
Denote by $\mathbb{A}^n_{k}$ (or just $\mathbb{A}^n$) the
$n$-dimensional {\it affine space} over the field $k$. The {\it
Zariski topology} on $\mathbb{A}^n_k$  is defined to be the
topology whose closed sets are the zero sets $\mathbb{V}_k(I)$.

For any ideal $I$ of the commutative ring $R$, the radical
 $\sqrt{I}$ of $I$, is the set
\begin{equation}
    \{a \in R ~|~ a^n \in I \text { for some integer } n > 0 \}.
\end{equation}

 \index{radical}
 \index{radical ideal}
The radical of an ideal $I$ is always an ideal of $R$. If $I =
\sqrt{I}$, then $I$ is called a {\it radical ideal}. The prime
ideals are radical:
\begin{equation}
   \fbox{prime ideals} \subset \fbox{radical ideals} ~~.
\end{equation}

 \index{theorem! - Hilbert's Nullstellensatz}
 \index{Nullstellensatz}

\begin{theorem}
  \label{hilbert_null}
{\em (Hilbert's Nullstellensatz, \cite[Ch.2]{Re88}, or
\cite[pp.9-11]{Mum88}) }

  Let $k$ be an algebraically closed field.

 (1) The maximal ideals in the ring $A = k[x_1,\dots,x_n]$
   are the ideals
\begin{equation}
  \label{closed_points}
    \mathfrak{m}_{P} = (x_1 - a_1,\dots,x_n - a_n)
\end{equation}
 for some point $P = (a_1,\dots,a_n)$. The ideal
 $\mathfrak{m}_{P}$
 coincides with the ideal $I(P)$ of all functions which are zero in
 $P$.

 (2) If $\mathbb{V}_k(I) = \emptyset$, then $I = k[x_1,\dots,x_n].$

 (3) For any ideal $J \subset A$
\begin{equation}
     I(\mathbb{V}_k(J)) = \sqrt{J},
\end{equation}
\end{theorem}

 \index{affine variety}

The set $X \subset \mathbb{A}^n_{k}$ is called the {\it affine
variety} if $X = \mathbb{V}_k(I)$ for some ideal $I \subset A$,
see (\ref{zero_set}). An affine variety $X \subset \mathbb{A}^n_k$
is called {\it irreducible affine variety} if there does not exist
a decomposition
\begin{equation*}
  X = X_1 \cup X_2,
\end{equation*}
where $X_1$, $X_2$ are two proper subsets of $X$. For example, the
affine variety
$$
  X = \{(x,y) \subset \mathbb{A}^2_{\mathbb{C}} ~|~ xy = 0 \}
$$
is decomposed to the sum of
$$
  X_1 = \{(x,y) \subset \mathbb{A}^2_{\mathbb{C}} ~|~ x = 0 \}
  \text{ and }
  X_2 = \{(x,y) \subset \mathbb{A}^2_{\mathbb{C}} ~|~ y = 0 \}.
$$
\begin{proposition}
  \label{irreducible_set}
   Let $X \subset \mathbb{A}^n_k$ be an affine variety and $I = I(X)$
   the corresponding ideal, i.e., $X = \mathbb{V}_k(I)$, see
   (\ref{zero_set}). Then
\begin{equation*}
    X \text{ is irreducible } \Longleftrightarrow I(X)
    \text{ is the prime ideal}.
\end{equation*}
\end{proposition}
For a proof of this proposition, see, e.g., \cite[\S3.7]{Re88}.

One of the important corollaries of Hilbert's Nullstellensatz is
the following correspondence \cite[\S3.10]{Re88} between
subvarieties
 $X \subset\mathbb{A}^n_{k}$ and ideals $I \subset A$:
\begin{equation}
\begin{array}{ccc}
    \{ \text{ radical ideals }  \} & \Longleftrightarrow & \{ \text{ affine variety } \} \\
       \bigcup &  & \bigcup \\
    \{ \text{ prime ideals } \} & \Longleftrightarrow & \{ \text{ irreducible affine variety }  \} \\
       \bigcup &  & \bigcup \\
    \{ \text{ maximal ideals } \} & \Longleftrightarrow & \{ \text{ points }  \} \\
\end{array}
\end{equation}

\subsection{Prime spectrum}
 \label{prime_spec}

 \index{prime spectrum ${\rm Spec}(R)$}
 The {\it prime spectrum} ${\rm Spec}(R)$ of the  commutative
ring $R$ is defined to be the set of proper prime ideals of $R$:
\begin{equation}
     \{ \mathfrak{p} \subset R ~|~ \mathfrak{p}
        \text{ is a prime ideal of } R \}.
\end{equation}
The ring $R$ itself is not counted as a prime ideal, but ($0$), if
prime, is counted.

 \index{closed sets of ${\rm Spec}(R)$}
 \index{closed points}
 \index{closure}
 \index{generic point}
 A topology is imposed on ${\rm Spec}(R)$ by defining
the sets of {\it closed sets}. The closed set $V(J)$ is the set of
the form
\begin{equation}
  \label{closed_set}
     \{ \mathfrak{p}  ~|~ \mathfrak{p} \supseteq J
        \text{ for some ideal } J \subset R \}.
\end{equation}
This set is denoted by $V(J)$. The {\it closure}
$\overline{\mathcal{P}}$ of the subset $\mathcal{P}$ is the
intersection of all closed sets containing $\mathcal{P}$:
\begin{equation}
  \label{closure}
    \overline{\mathcal{P}} =
     \bigcap\limits_{V(J) \supset \mathcal{P}}V(J).
\end{equation}

Consider the closure of points of the topological space ${\rm
Spec}(R)$. For a point of ${\rm Spec}(R)$ which is the prime ideal
$\mathfrak{p}$, we have the closure
\begin{equation}
  \label{closure_of_point}
    \overline{\mathfrak{p}} =
     \bigcap\limits_{V(J) \supset \mathfrak{p}}V(J) =
     V(\mathfrak{p}),
\end{equation}
which consists of all prime ideals $\mathfrak{p}^{'} \supset
\mathfrak{p}$. In particular, the closure
$\overline{\mathfrak{p}}$ consists of one point if and only if the
ideal $\mathfrak{p}$ is maximal. The point which coincides with
its closure is called the {\it closed point}. Thus, closed points
in ${\rm Spec}(R)$ correspond one-to-one with maximal ideals.

In the conditions of Hilbert's Nullstellensatz every maximal ideal
$J$ is $\mathfrak{m}_P$ (\ref{closed_points}) defined by some
point $P = (a_1,\dots,a_n) \in k^n$. Thus, in the finitely
generated ring $A = k[x_1,\dots,x_n]$ over an algebraically closed
field $k$ the closed points in ${\rm Spec}(A)$ correspond to the
maximal ideals and correspond to the points $P = (a_1,\dots,a_n)
\in k^n$.

There exist non-closed points in the topological space ${\rm
Spec}(R)$. Let $R$ has no divisors of $0$. Then the ideal (0) is
prime and is contained in all other prime ideals. Thus, the
closure $\overline{(0)}$ in (\ref{closure_of_point}) consists of
all prime ideals, i.e. coincides with the space ${\rm Spec}(R)$.
The point $(0)$ is everywhere dense point in ${\rm Spec}(R)$. The
everywhere dense is called the {\it generic point}.

\begin{definition} {\rm
Let $\mathcal{P}$ is an irreducible closed subset of
 ${\rm Spec}(R)$. A point $a \in \mathcal{P}$ is called
the {\it generic point of $\mathcal{P}$} if the closure
$\overline{a}$ coincides with $\mathcal{P}$. }
\end{definition}

\begin{proposition}{\em (\cite[p.126]{Mum88})}
   If $x \in {\rm Spec}(R)$, then the closure $\overline{\{x\}}$ of $\{x\}$
is irreducible and $x$ is a generic point of $\overline{\{x\}}$.
Conversely, every irreducible closed subset $\mathcal{P} \subset
{\rm Spec}(R)$ equals $V(J)$ for some prime ideal $J \subset R$
and $J$ is its unique generic point.
\end{proposition}

 \index{distinguished open set of ${\rm Spec}(R)$}
 A {\it distinguished open set} of ${\rm Spec}(R)$ is defined
 to be an open set of the form
\begin{equation}
     {\rm Spec}(R) := \{ \mathfrak{p} \in {\rm Spec}(R) ~|~ f
     \notin \mathfrak{p} \}.
\end{equation}
for any element $f \in R$.
\begin{example} {\rm
 1) If R is a field, then ${\rm Spec}(R)$ has just one point $(0)$.

 2) Let $R = k[X]$ be a polynomial ring in one variable $x$.
 Then ${\rm Spec}(R)$ is the affine line $A^1_k$ over $k$. There exist
 two types of prime ideals: $(0)$ and $(f(x))$, where $f$ is
 an irreducible polynomial, corresponding to the prime ideal $(f(x))$.
 For algebraically closed $k$, the closed points are all of the
 form $(X - a)$. The point $(0)$ is generic.

 3) ${\rm Spec}(\mathbb{Z})$ consists of closed points for every prime
 ideal (p), plus ideal (0).
}
\end{example}

\subsection{Coordinate ring}
 \label{coord_ring}
 \index{coordinate ring}

Let $V \subset \mathbb{A}^n_k$ be an affine variety and $I(V)$ be
an ideal of $V$. The {\it coordinate ring} of the affine variety
$V$ is the quotient ring denoted by $k[V]$ and defined as follows:
\begin{equation}
   k[V] := k[x_1,\dots,x_n]/I(V).
\end{equation}

A {\it regular  function} on the affine variety $V$ is the
restriction to $V$ of a polynomial in $x_1,\dots,x_n$, modulo
$I(V)$ (i.e., modulo functions vanishing on $V$). Thus, the
regular functions on the affine variety $V$ are elements of the
coordinate ring $k[V]$.
\begin{example} {\rm
1) Let $f(x,y)$ be a complex polynomial function from
$\mathbb{C}^2$. The coordinate ring of a plane curve defined by
the equation $f(x,y) = 0$ in $\mathbb{A}^2_{\mathbb{C}}$ is
\begin{equation}
   \mathbb{C}[x, y]/(f(x,y)),
\end{equation}
where $(f(x,y))$ is the ideal generated by the polynomial
$f(x,y)$.

2) Consider two algebraic curves: $y = x^r (r \in \mathbb{N})$,
and $y = 0$. The coordinate rings of these curves are isomorphic:
\begin{equation*}
  \begin{array}{cccc}
     & y = x^r  &  \text{ corresponds to }&  \mathbb{C}[x,y]/(y - x^r)
                \simeq \mathbb{C}[x,x^r] \simeq \mathbb{C}[x],
     \\
     & y = 0 &  \text{ corresponds to }&  \mathbb{C}[x,y]/(y)
                \simeq \mathbb{C}[x,x] \simeq \mathbb{C}[x].
  \end{array}
\end{equation*}
Thus, in the sense of algebraic geometry,
 curves $y = x^r$ and $y = 0$ are the same.

 3) On the other hand
\begin{equation*}
     \mathbb{C}[x,y]/(y^2 - x^3) \ncong \mathbb{C}[x],
\end{equation*}
and the curves $y^2 = x^3$ and $y = 0$ are not the same (in the
sense of algebraic geometry). Really, there exists the isomorphism
$$
  x \longmapsto T^2, \hspace{5mm} y \longmapsto  T^3,
$$
and
\begin{equation*}
     \mathbb{C}[x,y]/(y^2 - x^3) \simeq \mathbb{C}[T^2, T^3]
         \subset \mathbb{C}[T].
\end{equation*}

 \index{Neile's parabola}
 \index{semicubical parabola}
The affine variety $X = \{(x,y) ~|~ y^2 = x^3 \}$ is called {\it
Neile's parabola}, or the {\it semicubical parabola}. }
\end{example}

\begin{remark}{\rm
Consider the coordinate ring $R = k[V]$ of the affine variety $V$
over algebraically closed $k$. Then,}
\begin{center}{\it
 the prime spectrum ${\rm Spec}(R)$ contains exactly \\
 the same information as the variety V}.
\end{center}
{\rm  Really, according to Hilbert's Nullstellensatz, the maximal
ideals of $k[V]$ correspond one-to-one with points of $V$:
\begin{equation*}
          v \in V \Longleftrightarrow  \mathfrak{m}_v \subset k[V].
\end{equation*}
Besides, according to Proposition \ref{irreducible_set} every
other prime ideal $\mathfrak{p} \subset k[V]$ is the intersection
of maximal ideals corresponding to the points of some irreducible
subvariety $Y \subset V$:
\begin{equation*}
          \mathfrak{p}_{Y} = \bigcap\limits_{v \in Y}\mathfrak{m}_v.
\end{equation*}
}
\end{remark}

\subsection{Orbit space}
  \label{orbit_space}
 \index{orbit space}
 \index{quotient space}
 \index{Zariski's topology}
 \index{geometric invariant theory(GIT)}

Let $G$ be a finite group acting on affine variety $V$ and $k[V]$
be the coordinate ring of $V$. The problem here is that the
quotient space $V/G$ may not exist, even for very trivial group
actions, see Example \ref{non_closed_orbit}. To find how this
problem is solved in the geometric invariant theory (GIT), see
\cite{Dol03} or \cite{Kr85}.

\begin{example}
  \label{non_closed_orbit}
{\rm Consider the multiplicative group $G = k^*$ on the affine
line $A^1_k$. An orbit space for the group $G$ consists of two
orbits: $\{0\}$ and $A^1_k -\{0\}$. The second orbit $A^1_k
-\{0\}$ is not closed subset in the Zariski's topology, and the
first orbit $\{0\}$ is contained in the closure of the orbit
$A^1_k -\{0\}$. Thus, the point $\{0\}$ is the generic point
(\S\ref{prime_spec}). }
\end{example}

 \index{categorical quotient}

The following two definitions of the {\it categorical quotient}
and the {\it orbit space} can be found, e.g., in the
P.~E.~Newstead textbook \cite[p.39]{Newst78}.

\begin{definition}{\rm
   Let $G$ be an algebraic group acting on a variety $V$.
   A {\it categorical quotient} of $V$ by $G$ is a pair
   $(Y, \varphi)$, where $Y$ ia  variety and $\varphi$ is a
   morphism $\varphi : V \longrightarrow Y$ such that

   (i) $\varphi$ is constant on the orbits of the action;

   (ii) for any variety $Y^{'}$ and morphism
   $\varphi^{'} : V \longrightarrow Y^{'}$ which is constant on
   orbits, there is a unique morphism $\psi : Y \longrightarrow Y^{'}$
   such that $\psi \circ \varphi = \varphi^{'}$.
   }
\end{definition}
\begin{definition}{\rm
  A categorical quotient of $V$ by $G$ is called an {\it orbit space}
  if $\varphi^{-1}(y)$ consists of a single orbit for all $y \in
  Y$. The {\it orbit space} is denoted by $V/G$.
   }
\end{definition}
    The {\it orbit space} $X = V/G$ is an affine variety
whose points correspond one-to-one with orbits of the group
action.
\begin{proposition}
  Let $G$ be a finite group of order $n$ acting on affine variety $V$.
  Assume that the characteristic of $k$ does not divide $n$.  Set
\begin{equation}
 \label{def_orbit_space}
     X := {\rm Spec}~k[V]^G.
\end{equation}
Then $X$ is an orbit space for the action $G$ on $V$.
\end{proposition}
For a proof of this theorem see, e.g., \cite[Ch.1 \S2,
Ex.11]{Shaf88} or \cite[Th.1.6]{Rom00}.

\section{The orbit structure of the Coxeter transformation}
 \label{orbit_str}

 \index{highest root}
 \index{orbit structure of the Coxeter transformation}
 \index{Lie algebra}
 \index{Cartan subalgebra}
 Let $\mathfrak{g}$ be the simple complex Lie algebra
 of type $A, D$ or $E$, and $\mathfrak{h}$ be a Cartan
 subalgebra of $\mathfrak{g}$. Let $\mathfrak{h}^{\vee}$ be
 the dual space to $\mathfrak{h}$ and $\alpha_i \in
 \mathfrak{h}^{\vee}, i = 1,\dots,l$ be an ordered set of simple
 positive roots.
 Here, we follow B.~Kostant's description \cite{Kos84}
 of the orbit structure of the Coxeter transformation ${\bf C}$
 on the highest root in the root system of $\mathfrak{g}$.
 We consider a bipartite graph and a bicolored Coxeter
 transformation from \S\ref{bicolor_part}, \S\ref{conjugacy}.
 Let $\beta$ be the highest root of ($\mathfrak{h}, \mathfrak{g}$),
 see \S\ref{roots_of_1}. Between two bicolored Coxeter
 transformations (\ref{C_decomp}) we select $w_1, w_2$ such that
 \begin{equation}
         w_2\beta = \beta.
 \end{equation}
Consider, for example, the Dynkin diagram $E_6$. Here,
\begin{equation}
    w_1 =  \left (
      \begin{array}{cccccc}
        1 & & &  &  &  \\
        & 1 & &  &  &  \\
        & & 1 &  &  &  \\
        1 & 1 & 0 & -1 & &  \\
        1 & 0 & 1 & & -1 & \\
        1 & 0 & 0 & & & -1 \\
      \end{array}
       \right )
      \begin{array}{c}
         x_0 \\
         x_1 \\
         x_2 \\
         y_1 \\
         y_2 \\
         y_3 \\
      \end{array},
       \hspace{5mm}
   w_2 = \left (
      \begin{array}{cccccc}
        -1 & & & 1 & 1 & 1 \\
        & -1 & & 1 & 0 & 0 \\
        &  & -1 & 0 & 1 & 0 \\
        &  &    & 1 & & \\
        &  &    & & 1 & \\
        &  &    & & & 1 \\
      \end{array}
       \right )
      \begin{array}{c}
         x_0 \\
         x_1 \\
         x_2 \\
         y_1 \\
         y_2 \\
         y_3 \\
      \end{array},
\end{equation}
where any vector $z \in \mathfrak{h}^{\vee}$ and the highest
vector $\beta$ are:
\begin{equation}
  \label{def_E6_vectors}
      \begin{array}{ccccc}
        x_1 & y_1 & x_0 & y_2 & x_2 \\
            &     & y_3 &     &
      \end{array},
      \hspace{5mm}
      \beta =
      \begin{array}{ccccc}
        1 & 2 & 3 & 2 & 1 \\
        &     & 2 &     &
      \end{array} ,
      \hspace{5mm}
      \text{ or } \hspace{5mm}
      \left (
      \begin{array}{c}
         3 \\
         1 \\
         1 \\
         2 \\
         2 \\
         2 \\
      \end{array}
      \right )
      \begin{array}{c}
         x_0 \\
         x_1 \\
         x_2 \\
         y_1 \\
         y_2 \\
         y_3 \\
      \end{array},
\end{equation}

Further, following B.Kostant \cite[Th.1.5]{Kos84} consider the
alternating products $\tau^{(n)}$:
\begin{equation}
  \begin{split}
    & \tau^{(1)} = w_1, \\
    & \tau^{(2)} = {\bf C} = w_2w_1, \\
    & \tau^{(3)} = w_1{\bf C} = w_1w_2w_1, \\
    & \dots, \\
    & \tau^{(n)} =
       \begin{cases}
         {\bf C}^k = w_2w_1\dots{w_2}w_1 \text{ for } n = 2k, \\
         w_1{\bf C}^k = w_1w_2w_1\dots{w_2}w_1 \text{ for } n = 2k+1, \\
       \end{cases}
  \end{split}
\end{equation}
and the orbit of $\tau^{(n)}$ on the highest vector $\beta$:
\begin{equation}
   \beta_n = \tau^{(n)}\beta,
\end{equation}
where $n = 1,\dots,h-1$ ($h$ is the Coxeter number,
see(\ref{roots_of_1})).

\begin{theorem}{\em\bf (B.~Kostant, \cite[Theorems 1.3, 1.4, 1.5, 1.8]{Kos84})}

  1)  There exist $z_j \in \mathfrak{h}^{\vee}, j = 1,\dots,h-1$
   and even integers $2 \le a \le b \le h$
   (see (\ref{Kostant_numbers_a_b}) and Table
   \ref{Kostant_numbers})such that the generating functions
   $P_{G}(t)$
   (see (\ref{def_generating_funct}), (\ref{def_P_G_2}), or
         (\ref{Kostant_gen_func}), (\ref{gen_func_i}) )
   are obtained as follows:
  \begin{equation}
    \label{series_P_G_i}
      [P_{G}(t)]_i =
       \begin{cases}
         \displaystyle\frac{1 + t^h}
            {(1 - t^a)(1 - t^b)} \text{ for } i = 0, \vspace{5mm} \\
         \displaystyle\frac{\sum\limits_{j=1}^{h-1}z_j{t}^j}
            {(1 - t^a)(1 - t^b)} \text{ for } i = 1,\dots,r. \\
       \end{cases}
  \end{equation}
  Indices $i = 1,\dots,r$ enumerate vertices of the
  Dynkin diagram and the coordinates of the vectors $z_n$; $i = 0$ corresponds to the additional vertex,
  the one that extends the Dynkin diagram to the extended Dynkin
  diagram.

  2) Vectors $z_n$ (we call these vectors the assembling vectors)
  are obtained by means of the orbit of $\tau^{(n)}$ on the highest vector $\beta$
  as follows:
\begin{equation}
  z_n = \tau^{(n-1)}\beta - \tau^{(n)}\beta.
\end{equation}

  3) On has
\begin{equation}
  z_g = 2\alpha_{*}, \text{ where } g = \frac{h}{2},
\end{equation}
and $\alpha_{*}$ is the simple root corresponding to the branch
point for diagrams $D_n$, $E_n$ and to the midpoint for the
diagram $A_{2m-1}$. In all these cases $h$ is even, and $g$ is
integer. The diagram $A_{2m}$ has been excluded.

  4) The series of assembling vectors $z_n$ is symmetry:
\begin{equation}
  z_{g+k} = z_{g-k} \text{ for } k = 1,\dots,g.
\end{equation}
\end{theorem}

In the case of the Dynkin diagram $E_6$ vectors $\tau^{(n)}\beta$
are given in Table \ref{table_orbit_max_root}, and the assembling
vectors $z_n$ are given in Table \ref{table_vectors_z}. Vector
$z_6$ coincides with $2\alpha_{x_0}$, where $\alpha_{x_0}$ is the
simple root corresponding to the vertex $x_0$, see
(\ref{def_E6_vectors}). From Table \ref{table_vectors_z} we see,
that
\begin{equation}
    z_1 = z_{11}, \hspace{5mm}
    z_2 = z_{10}, \hspace{5mm}
    z_3 = z_9 , \hspace{5mm}
    z_4 = z_8, \hspace{5mm}
    z_5 = z_7.
\end{equation}

\begin{table} 
 \centering
 \vspace{2mm}
 \caption{\hspace{3mm}Orbit of the Coxeter transformation on the highest root}
 \renewcommand{\arraystretch}{1.3}
  \begin{tabular} {||c|c|c||}
  \hline \hline
      $\beta$   &
      $\tau^{(1)}\beta = w_1\beta$ &
      $\tau^{(2)}\beta = {\bf C}\beta$ \\
  \hline
      $\begin{array}{ccccc}
        1 & 2 & 3 & 2 & 1 \\
        &     & 2 &     &
      \end{array}$
&     $\begin{array}{ccccc}
        1 & 2 & 3 & 2 & 1 \\
        &     & 1 &     &
      \end{array}$
&     $\begin{array}{ccccc}
        1 & 2 & 2 & 2 & 1 \\
        &     & 1 &     &
      \end{array}$
      \\
  \hline \hline
     $\tau^{(3)}\beta = w_1{\bf C}\beta$ &
     $\tau^{(4)}\beta = {\bf C}^2\beta$ &
     $\tau^{(5)}\beta = w_1{\bf C}^2\beta$
     \\
  \hline
     $\begin{array}{ccccc}
        1 & 1 & 2 & 1 & 1 \\
        &     & 1 &     &
      \end{array}$
&     $\begin{array}{ccccc}
        0 & 1 & 1 & 1 & 0 \\
        &     & 1 &     &
      \end{array}$
&     $\begin{array}{ccccc}
        0 & 0 & 1 & 0 & 0 \\
        &     & 0 &     &
      \end{array}$
      \\
  \hline \hline
      $\tau^{(6)}\beta = {\bf C}^3\beta$ &
      $\tau^{(7)}\beta = w_1{\bf C}^3\beta$ &
      $\tau^{(8)}\beta = {\bf C}^4\beta$ \\
  \hline
      $\begin{array}{ccccc}
        0 & 0 & -1 & 0 & 0 \\
        &     & 0 &     &
      \end{array}$
&     $\begin{array}{ccccc}
        0 & -1 & -1 & -1 & 0 \\
        &     & -1 &     &
      \end{array}$
&      $\begin{array}{ccccc}
        -1 & -1 & -2 & -1 & -1 \\
        &     & -1 &     &
      \end{array}$
      \\
  \hline \hline
      $\tau^{(9)}\beta = w_1{\bf C}^4\beta$  &
      $\tau^{(10)}\beta = {\bf C}^5\beta$ &
      $\tau^{(11)}\beta = w_1{\bf C}^5\beta$ \\
  \hline
      $\begin{array}{ccccc}
        -1 & -2 & -2 & -2 & -1 \\
        &     & -1 &     &
      \end{array}$
&     $\begin{array}{ccccc}
        -1 & -2 & -3 & -2 & -1 \\
        &     & -1 &     &
      \end{array}$
&     $\begin{array}{ccccc}
        -1 & -2 & -3 & -2 & -1 \\
        &     & -2&     &
      \end{array}$ \\
  \hline \hline
\end{tabular}
  \label{table_orbit_max_root}
\end{table}

\begin{table} 
 \centering
 \vspace{2mm}
 \caption{\hspace{3mm}Assembling vectors $z_n = \tau^{(n-1)}\beta -
            \tau^{(n)}\beta$}
 \renewcommand{\arraystretch}{1.3}
  \begin{tabular} {||c|c|c||}
  \hline \hline
      $z_1 = \beta - w_1\beta$   &
      $z_2 = w_1\beta - {\bf C}\beta$ &
      $z_3 = {\bf C}\beta - w_1{\bf C}\beta$ \\
  \hline
      $\begin{array}{ccccc}
        0 & 0 & 0 & 0 & 0 \\
        &     & 1 &     &
      \end{array}$
&     $\begin{array}{ccccc}
        0 & 0 & 1 & 0 & 0 \\
        &     & 0 &     &
      \end{array}$
&     $\begin{array}{ccccc}
        0 & 1 & 0 & 1 & 0 \\
        &     & 0 &     &
      \end{array}$
      \\
  \hline \hline
      $z_4 = w_1{\bf C}\beta - {\bf C}^2\beta$   &
      $z_5 = {\bf C}^2\beta - w_1{\bf C}^2\beta$ &
      $z_6 = w_1{\bf C}^2\beta - {\bf C}^3\beta$ \\
  \hline
      $\begin{array}{ccccc}
        1 & 0 & 1 & 0 & 1 \\
        &     & 0 &     &
      \end{array}$
&     $\begin{array}{ccccc}
        0 & 1 & 0 & 1 & 0 \\
        &     & 1 &     &
      \end{array}$
&     $\begin{array}{ccccc}
        0 & 0 & 2 & 0 & 0 \\
        &     & 0 &     &
      \end{array}$
      \\
  \hline \hline
      $z_7 = {\bf C}^3\beta - w_1{\bf C}^3\beta$   &
      $z_8 = w_1{\bf C}^3\beta - {\bf C}^4\beta$ &
      $z_9 = {\bf C}^4\beta - w_1{\bf C}^4\beta$ \\
  \hline
      $\begin{array}{ccccc}
        0 & 1 & 0 & 1 & 0 \\
        &     & 1 &     &
      \end{array}$
&     $\begin{array}{ccccc}
        1 & 0 & 1 & 0 & 1 \\
        &     & 0 &     &
      \end{array}$
&     $\begin{array}{ccccc}
        0 & 1 & 0 & 1 & 0 \\
        &     & 0 &     &
      \end{array}$
      \\
  \hline \hline
      $z_{10} = w_1{\bf C}^4\beta - {\bf C}^5\beta$   &
      $z_{11} = {\bf C}^5\beta - w_1{\bf C}^5\beta$ & \\
  \hline
      $\begin{array}{ccccc}
        0 & 0 & 1 & 0 & 0 \\
        &     & 0 &     &
      \end{array}$
&     $\begin{array}{ccccc}
        0 & 0 & 0 & 0 & 0 \\
        &     & 1 &     &
      \end{array}$
&     \\
  \hline \hline
\end{tabular}
  \label{table_vectors_z}
\end{table}

Denote by $z(t)_i$ the polynomial $\sum\limits_{j=1}^{h-1}z_jt^j$
from (\ref{series_P_G_i}). In the case of $E_6$ we have:
\begin{equation}
 \label{def_zt_2}
  \begin{split}
    & z(t)_{x_1} = t^4 + t^8,  \\
    & z(t)_{y_1} = t^3 + t^5 + t^7 + t^9, \\
    & z(t)_{y_2} = t^3 + t^5 + t^7 + t^9, \\
    & z(t)_{x_2} = t^4 + t^8, \\
    & z(t)_{y_3} = t + t^5 + t^7 + t^{11}.
  \end{split}
\end{equation}
The Kostant numbers $a,b$ (see Table \ref{Kostant_numbers}) for
$E_6$ are $a = 6$, $b = 8$. From (\ref{series_P_G_i}) and
(\ref{def_zt_2}), we have

\begin{equation}
 \begin{split}
 & [P_G(t)]_{x_1} = [P_G(t)]_{x_2} = \frac{t^4 + t^8}{(1 - t^6)(1 - t^8)},
    \\ \\
 & [P_G(t)]_{y_1} = [P_G(t)]_{y_2} = \frac{t^3 + t^5 + t^7 + t^9}{(1 - t^6)(1 - t^8)},
    \\ \\
 & [P_G(t)]_{y_3} \hspace{21mm} = \frac{t + t^5 + t^7 + t^{11}}{(1 - t^6)(1 - t^8)}.
 \end{split}
\end{equation}
Since
\begin{equation}
     1 - t^6 = \sum\limits_{n=0}^{\infty}t^{6n}, \hspace{5mm}
     1 - t^8 = \sum\limits_{n=0}^{\infty}t^{8n},
\end{equation}
we have
\begin{equation}
 \label{calc_poincare}
 \begin{split}
 & [P_G(t)]_{x_1} = [P_G(t)]_{x_2} =
           \sum\limits_{i=0, j= 0}^{\infty}(t^{6i + 8j + 4} + t^{6i + 8j + 8}),
    \\ \\
 & [P_G(t)]_{y_1} = [P_G(t)]_{y_2} =
          \sum\limits_{i=0, j= 0}^{\infty}(t^{6i + 8j + 3} + t^{6i + 8j + 5} +
                                           t^{6i + 8j + 7} + t^{6i + 8j + 9}),
    \\ \\
 & [P_G(t)]_{y_3} \hspace{21mm} =
          \sum\limits_{i=0, j= 0}^{\infty}(t^{6i + 8j + 1} + t^{6i + 8j + 5}+
                                           t^{6i + 8j + 7} + t^{6i + 8j + 11}),
 \end{split}
\end{equation}

Recall, that $m_{\alpha}(n), \alpha = x_1, x_2, y_1, y_2, y_3$ are
the multiplicities of the indecomposable representations
$\rho_{\alpha}$ of $G$ (considered in the context of McKay
correspondence, \S(\ref{McKay})) in the decomposition of $\pi_n|G$
(\ref{decomp_pi_n}). These multiplicities are the coefficients of
the Poincar\'{e} series (\ref{calc_poincare}), see (\ref{def_vn}),
(\ref{Kostant_gen_func}), (\ref{gen_func_i}). For example,
\begin{equation}
  [P_G(t)]_{x_1} = [P_G(t)]_{x_2} =
    t^4 + t^8 + t^{10} + t^{12} + t^{14} + 2t^{16} + t^{18} + 2t^{20} + \dots
\end{equation}
\begin{equation}
  m_{x_1} = m_{x_2} =
  \begin{cases}
     0 ~\text{ for } n = 1,2,3,5,6 \text{ and } n = 2k+1, k \geq 3, \\
     1 ~\text{ for } n = 4,8,10,12,14,18, \dots \\
     2 ~\text{ for } n = 16,20, \dots \\
     \dots
  \end{cases}
\end{equation}
In particular representations $\rho_{x_1}(n)$ and $\rho_{x_2}(n)$
do not enter in the decomposition of $\pi_n$ of $SU(2)$ (see
\S\ref{generating_fun}) for all odd $n$.

\section{Fixed and anti-fixed points of the Coxeter transformation}
\subsection{Chebyshev polynomials and the McKay-Slodowy matrix}
\label{fixed_anti_fixed}
 \index{fixed points}
 \index{anti-fixed points}

 We continue to use block arithmetic for $(m+k)\times(m+k)$ matrices
 as in \S\ref{bicolored_C}, namely:
\begin{equation}
 \label{matrix_K_blocks}
K = \left (
\begin{array}{cc}
 2I & 2D \\
 2F & 2I
\end{array}
\right ), \hspace{5mm} w_1 = \left (
\begin{array}{cc}
 -I & -2D \\
 0 & I
\end{array}
\right ), \hspace{5mm} w_2 = \left (
\begin{array}{cc}
 I & 0 \\
 -2F & -I
\end{array}
\right ).
\end{equation}
 Let
\begin{equation}
 \beta = \frac{1}{2}K - I = \left (
 \begin{array}{cc}
 0 & D \\
     F & 0
 \end{array}
 \right ).
\end{equation}
The matrix $\beta$ coincides (up to a factor
$\displaystyle-\frac{1}{2}$) with the McKay matrix from
\S\ref{McKay} or with the Slodowy matrix in the multiply-laced
case, (\ref{slodowy_1}), (\ref{slodowy_2}).

Introduce two recursive series $f_p(\beta)$ and $g_p(\beta)$ of
polynomials in $\beta$:
\begin{equation}
 \label{def_f}
 \begin{array}{l}
    f_0(\beta) = 0,\\
f_1(\beta) = K = 2\beta + 2I, \vspace{2mm} \\
    f_{p+2}(\beta) = 2\beta{f}_{p+1}(\beta) - f_p(\beta)\; \text{ for }
    p \geq 0
    \end{array}
\end{equation}
and
\begin{equation}
 \label{def_g}
 \begin{array}{l}
    g_0(\beta) = 2I,\\
g_1(\beta) = K -2I = 2\beta, \vspace{2mm} \\
    g_{p+2}(\beta) = 2\beta{g}_{p+1}(\beta) - g_p(\beta)\; \text{ for }
    p \geq 0.
    \end{array}
\end{equation}
For example,
\begin{equation}
 \label{f2_f3}
 \begin{array}{l}
 f_2(\beta) = 2\beta(2\beta + 2I)
 = (K - 2I)K = 4\left (
 \begin{array}{cc}
        DF & D \\
        F & FD
 \end{array}
 \right ),  \vspace{2mm} \\
     f_3(\beta) = (K - 2I)(K - 2I)K - K =
     2\left (
 \begin{array}{cc}
     4DF - I & 4DFD - D \\
            4FDF - F & 4FD - I
 \end{array}
 \right ), \vspace{2mm}
 \end{array}
\end{equation}
\begin{equation}
 \label{g2_g3}
 \begin{array}{l}
 g_2(\beta) = 2(2\beta^2 - I) = 2\left (
 \begin{array}{cc}
        2DF - I & 0 \\
        0 & 2FD - I
 \end{array}
 \right ),  \vspace{3mm} \\
     g_3(\beta) = 4\beta(2{\beta}^2 - I) - 2\beta =
     2\left (
 \begin{array}{cc}
     0 & 8DFD - 6D \\
            8FDF - 6F & 0
 \end{array}
 \right ). \vspace{2mm}
 \end{array}
\end{equation}
\begin{remark}  {\rm Let us formally substitute $\beta = \cos{t}$. Then

1) The polynomials $g_p(\beta)$ are Chebyshev polynomials of the
first kind (up to a factor $\displaystyle\frac{1}{2}$):
$$
 g_p(\cos{t}) = 2\cos{pt} \hspace{0.5mm}.
$$
Indeed,
$$
 g_1(\beta) = 2\beta = 2\cos{t}, \hspace{2mm}
 g_2(\beta) = 2(2\beta^2-1) = 2\cos{2t}, \hspace{2mm}
$$
and
$$
 g_{p+2}(\beta) =
 2(2\cos{t}\cos(p+1)t - \cos{pt}) =
 2\cos(p+2)t . \vspace{2mm}
$$

2) The polynomials $f_p(\beta)$ are Chebyshev polynomials of the
second kind (up to a factor $2\beta + 2$):
$$
 f_p(\cos{t}) = 2(\cos{t} + 1)\frac{\sin{pt}}{\sin{t}} .
$$
Indeed,
$$
 f_1(\beta) = 2\beta + 2I = 2\cos{t} + 2, \hspace{2mm}
 f_2(\beta) = 2\beta(2\beta+2I) =
 4\cos{t}(\cos{t} + 1),
$$
and \vspace{2mm}
$$
 f_{p+2}(\beta) =
 2\frac{\cos{t} + 1}{\sin{t}}
 (2\cos{t}\sin{(p+1)t} - \sin{pt}) =
 2(\cos{t} + 1)\frac{\sin{(p+2)t}}{\sin{t}} .
 \vspace{2mm}
$$
}
\end{remark}
\subsection{A theorem on fixed and anti-fixed points} The purpose
of this section is to prove the following
\begin{theorem} [\cite{SuSt82}]
 \label{fixed_main_proposition}
 1) The fixed points of the powers of the Coxeter transformation
 satisfy the following relations
 \begin{displaymath}
 {\bf C}^p{z} = z \hspace{2mm}
     \Longleftrightarrow \hspace{2mm}
     f_p(\beta){z} = 0 .
 \end{displaymath}

 2) The anti-fixed points of the powers of the Coxeter transformation
 satisfy the following relations
 \begin{displaymath}
 {\bf C}^p{z} = -{z} \hspace{2mm}
     \Longleftrightarrow \hspace{2mm}
     g_p(\beta){z} = 0 .
 \end{displaymath}
\end{theorem}

In the proof of Theorem \ref{fixed_main_proposition} we will need
some properties of block matrices. Define the involution $\sigma$
by setting
\begin{equation}
 \sigma : A = \left (
     \begin{array}{cc}
        P & Q \\
        S & T
     \end{array}
     \right ) \hspace{2mm}\mapsto
 \hspace{2mm} \stackrel{\line(1,0){7}}{A} \hspace{1mm} =
 \left (
     \begin{array}{cc}
        P & Q \\
        -S & -T
     \end{array}
     \right ).
 \end{equation}
The map $\sigma$ is a linear operator because $
\stackrel{\line(1,0){67}}{\alpha_1A_1 + \alpha_2A_2} =
 \alpha_1\stackrel{\line(1,0){7}}{A_1} +
 \alpha_2\stackrel{\line(1,0){7}}{A_2}
$.

The involution $\sigma$ does not preserve products:
$$
 \stackrel{\line(1,0){18}}{AR} \hspace{3mm} \neq \hspace{3mm}
 \stackrel{\line(1,0){7}}{A}\hspace{0.5mm}\stackrel{\line(1,0){7}}{R},
$$
but
\begin{equation}
 \stackrel{\line(1,0){18}}{AR} \hspace{2mm} = \hspace{2mm}
 \stackrel{\line(1,0){7}}{A}\hspace{0.5mm}R \hspace{5mm}
\text{ for any } R = \left (
 \begin{array}{cc}
                     X & Y \\
                     U & V
                     \end{array}
 \right ),
\end{equation}
and
\begin{equation}
 \label{inv_RA}
 R\stackrel{\line(1,0){7}}{A} \hspace{2mm} =
 -\stackrel{\line(1,0){18}}{RA} \text{ for }
 R = \left (
 \begin{array}{cc}
     0 & Y \\
            U & 0
 \end{array}
     \right ).
\end{equation}
The following relations are easy to check:

\begin{equation}
 \label{C_w1w2}
 \begin{split}
& {\bf C} = w_1{w_2} =
 \left (
 \begin{array}{cc}
     4DF - I & 2D \\
            -2F & -I
 \end{array}
     \right ), \\
& \dotfill \\
& {\bf C}^{-1} = w_2{w_1} =
 \left (
 \begin{array}{cc}
     -I & -2D \\
            2F & 4FD - I
 \end{array}
     \right ), \vspace{3mm}
 \end{split}
\end{equation}

\begin{equation}
 \begin{split}
 \label{w2c_and_w1C}
 & w_2{\bf C} =
 \left (
 \begin{array}{cc}
     4DF - I & 2D \\
            -8FDF+4F & -4FD+I
 \end{array}
     \right ),   \\  \\
 & w_1{\bf C}^{-1} =
 \left (
 \begin{array}{cc}
     -4DF + I & 4D - 8DFD\\
            2F & 4FD - I
 \end{array}
     \right ),
 \end{split}
\end{equation}
         \vspace{7mm}
\begin{equation}
 \label{w1_w2}
\begin{array}{c}
 w_1 - w_2 =
 4\left (
 \begin{array}{cc}
     -2I & -2D \\
            2F & 2I
 \end{array}
     \right )   = -\stackrel{\line(1,0){10}}{K}, \hspace{5mm}
 w_1 + w_2 = -(K -2I),
 \end{array}    \vspace{3mm}
\end{equation}

\begin{equation}
 \begin{split}
 \label{C_C1}
 & {\bf C} - {\bf C}^{-1} =
 4\left (
 \begin{array}{cc}
     DF & D \\
            -F & -FD
 \end{array}
     \right )   =
 4\stackrel{\line(1,0){67}}{\left (
 \begin{array}{cc}
     DF & D \\
            F & FD
 \end{array}
     \right )}, \\ \\ 
 & {\bf C} + {\bf C}^{-1} =
 2\left (
 \begin{array}{cc}
     2DF-I & 0 \\
            0 & 2FD-I
 \end{array}
     \right ),
     \end{split}
\end{equation}
     \vspace{7mm}
\begin{equation}
 \label{w2c_w1C}
 \begin{split}
 w_2{\bf C} - w_1{\bf C}^{-1} =
 & 2\left (
     \begin{array}{cc}
     4DF - I & 4DFD - D \\
        -(4FDF - F) & -(4FD - I)
     \end{array}
    \right ) =  \vspace{3mm} \\ \\
 & 2\stackrel{\line(1,0){137}}{
 \left (
     \begin{array}{cc}
     4DF - I & 4DFD - D \\
        4FDF - F & 4FD - I
     \end{array}
    \right )},
 \end{split}
\end{equation}
 \vspace{7mm}
\begin{equation}
 \label{w2c_plus_w1C}
 w_2{\bf C} + w_1{\bf C}^{-1} =
 \left (
     \begin{array}{cc}
     0 & -8DFD + 6D \\
        -8FDF + 6F & 0
     \end{array}
    \right ),
\end{equation}
     \vspace{7mm}
\begin{equation}
 \label{beta_w}
 2\beta{w}_2 + I = -{\bf C}, \hspace{5mm}
 2\beta{w}_1 + I = -{\bf C}^{-1}. \vspace{3mm}
\end{equation}

\begin{proposition}
 \label{prop_pm}
 The following relations hold

\begin{equation}
 \label{fixed_prop_1}
 w_2 - w_1 =
 \stackrel{\line(1,0){25}}{f_1(\beta)},
\end{equation}

\begin{equation}
 \label{fixed_prop_2}
 {\bf C} - {\bf C}^{-1} =
 \stackrel{\line(1,0){25}}{f_2(\beta)},
\end{equation}

\begin{equation}
 \label{fixed_prop_3}
 w_2{\bf C} - w_1{\bf C}^{-1} =
 \stackrel{\line(1,0){25}}{f_3(\beta)},
\end{equation}

\begin{equation}
 \label{fixed_prop_4}
 w_2 + w_1 = -{g}_1(\beta),
\end{equation}

\begin{equation}
 \label{fixed_prop_5}
 {\bf C} + {\bf C}^{-1} =
 g_2(\beta),
\end{equation}

\begin{equation}
 \label{fixed_prop_6}
 w_2{\bf C} + w{\bf C}^{-1} =
 -{g}_3(\beta).
\end{equation}
\end{proposition}

 \PerfProof The proof of this theorem follows
 from comparing relations (\ref{def_f})-(\ref{def_g})
 and relations (\ref{C_C1})--(\ref{w2c_plus_w1C}).

 (\ref{fixed_prop_1}) follows from (\ref{def_f}) and (\ref{w1_w2}).

 (\ref{fixed_prop_2}) follows from (\ref{f2_f3}) and (\ref{C_C1}).

 (\ref{fixed_prop_3}) follows from (\ref{f2_f3}) and (\ref{w2c_w1C}).

 (\ref{fixed_prop_4}) follows from (\ref{def_g}) and (\ref{w1_w2}).

 (\ref{fixed_prop_5}) follows from (\ref{g2_g3}) and (\ref{C_C1}).

 (\ref{fixed_prop_6}) follows from (\ref{g2_g3}) and (\ref{w2c_plus_w1C}).
 \qedsymbol
\begin{proposition}
\label{fm_gm} The following relations hold:

\begin{equation}
\stackrel{\line(1,0){35}}{f_{2m}(\beta)} = {\bf C}^m - {\bf
C}^{-m}, \hspace{5mm} \stackrel{\line(1,0){40}}{f_{2m+1}(\beta)} =
w_2{\bf C}^m - w_1{\bf C}^{-m},
\end{equation}

\begin{equation}
g_{2m}(\beta) = {\bf C}^m + {\bf C}^{-m}, \hspace{5mm}
g_{2m+1}(\beta) = -(w_2{\bf C}^m + w_1{\bf C}^{-m}).
\end{equation}
\end{proposition}

\PerfProof 1) Relations (\ref{fixed_prop_1})--(\ref{fixed_prop_3})
realize the basis of induction.\\
By (\ref{inv_RA}), by the induction hypothesis, and by
(\ref{beta_w}) we have
\begin{displaymath}
 \begin{split}
 & \stackrel{\line(1,0){40}}{f_{2m+2}(\beta)} ~=~
 \stackrel{\line(1,0){60}}{2\beta{f}_{2m+1}(\beta)} -
 \stackrel{\line(1,0){35}}{{f}_{2m}(\beta)} ~=~
 -2\beta\stackrel{\line(1,0){40}}{{f}_{2m+1}(\beta)} -
 \stackrel{\line(1,0){30}}{{f}_{2m}(\beta)} = \\
 & -2\beta(w_2{\bf C}^m - w_1{\bf C}^{-m}) -
        ({\bf C}^m - {\bf C}^{-m}) = \\
 & -(2\beta{w}_2 + I){\bf C}^m + (2\beta{w}_1 + I){\bf C}^{-m} ~=~
 {\bf C}^{m+1} - {\bf C}^{-(m+1)}.
 \end{split}
\end{displaymath}
In the same way we have
\begin{displaymath}
 \begin{split}
 & \stackrel{\line(1,0){40}}{f_{2m+3}(\beta)} ~=~
 -2\beta\stackrel{\line(1,0){40}}{f_{2m+2}(\beta)} -
 \stackrel{\line(1,0){40}}{f_{2m+1}(\beta)} = \\
 & -2\beta({\bf C}^{m+1} - {\bf C}^{-(m+1)}) -
        w_2{\bf C}^m - w_1{\bf C}^{-m} ~= \\
 & -(2\beta{w}_1 + I){w}_2{\bf C}^m +
 (2\beta{w}_2 + I){w}_1{\bf C}^{-m} = \\
 & {\bf C}^{-1}{w}_2{\bf C}^m - {\bf C}{w}_1{\bf C}^{-m} ~=~
 w_1{\bf C}^{m+1} - w_2{\bf C}^{-(m+1)}.
 \end{split}
\end{displaymath}

 2) Relations (\ref{fixed_prop_4})-(\ref{fixed_prop_6})
realize the basis of induction. By the induction hypothesis and by
(\ref{beta_w}) we have
\begin{displaymath}
 \begin{split}
 & g_{2m+2}(\beta) ~=~
 2\beta{g}_{2m+1}(\beta) -
 g_{2m}(\beta) =     \\
 & -2\beta(w_2{\bf C}^m + w_1{\bf C}^{-m}) -
        ({\bf C}^m + {\bf C}^{-m}) = \\
 & -(2\beta{w}_2 + I){\bf C}^m - (2\beta{w}_1 + I){\bf C}^{-m} =
 {\bf C}^{m+1} + {\bf C}^{-(m+1)}.
 \end{split}
\end{displaymath}
By analogy, we have
\begin{displaymath}
 \begin{split}
 & g_{2m+3}(\beta) ~=~
 2\beta{g}_{2m+2}(\beta) -
 {g}_{2m+1}(\beta) = \\
 & -2\beta({\bf C}^{m+1} + {\bf C}^{-(m+1)}) -
        (w_2{\bf C}^m + w_1{\bf C}^{-m}) = \\
& -(2\beta{w}_1 + I){w}_2{\bf C}^m -(2\beta{w}_2 + I){w}_1{\bf
C}^{-m} = {w}_1{\bf C}^{m+1} + {w}_2{\bf C}^{-(m+1)}. \qed
\end{split}
\end{displaymath}
Theorem \ref{fixed_main_proposition} now follows from Proposition
\ref{fm_gm}.